\documentclass{memo-l}
\usepackage{amsmath,amssymb,amsthm,stmaryrd}
\usepackage{lscape}
\usepackage[all]{xy}
\numberwithin{section}{chapter}
\numberwithin{equation}{section}
\newtheorem{thm}[equation]{Theorem}
\newtheorem{cor}[equation]{Corollary}
\newtheorem{prop}[equation]{Proposition}
\newtheorem{lem}[equation]{Lemma}
\newtheorem*{thm*}{Theorem}
\newtheorem{ques}[equation]{Question}
\theoremstyle{definition}
\newtheorem{defn}[equation]{Definition}
\theoremstyle{remark}
\newtheorem{rmk}[equation]{Remark}
\newtheorem{exam}[equation]{Example}

\newtheorem{const}[equation]{Construction}

\makeindex

\def\co{\colon\thinspace}
\newcommand{\mb}[1]{\mathbb{#1}}
\newcommand{\mf}[1]{\mathfrak{#1}}
\newcommand{\mc}[1]{\mathcal{#1}}
\newcommand{\mbf}[1]{\mathbf{#1}}

\newcommand{\Ext}{{\rm Ext}}

\newcommand{\Hom}{{\rm Hom}}

\newcommand{\colim}{\varinjlim}
\DeclareMathOperator*{\hocolim}{hocolim}
\DeclareMathOperator*{\holim}{holim}

\newcommand{\overfrom}{\mathop\leftarrow}
\newcommand{\into}{\hookrightarrow}

\newcommand{\longoverto}{\mathop{\longrightarrow}}
\newcommand{\Map}{{\rm Map}}
\newcommand{\End}{{\rm End}}
\DeclareMathOperator{\Spec}{Spec}

\DeclareMathOperator{\Lie}{Lie}
\DeclareMathOperator{\Ker}{Ker}
\newcommand{\Cl}{{\rm Cl}}
\newcommand{\mmod}{{\sslash}}

\newcommand{\xym}[1]{
\vskip 0.7pc
\centerline{\xymatrix{#1}}
\vskip 0.7pc
}

\newcommand{\abs}[1]{\lvert #1 \rvert}
\newcommand{\norm}[1]{\lVert #1 \rVert}
\newcommand{\bra}[1]{\langle #1 \rangle}
\newcommand{\br}[1]{\overline{#1}}

\newcommand{\ul}[1]{\underline{#1}}
\newcommand{\td}[1]{\widetilde{#1}}

\newcommand{\ZZ}{\mathbb{Z}}
\newcommand{\RR}{\mathbb{R}}
\newcommand{\CC}{\mathbb{C}}
\newcommand{\QQ}{\mathbb{Q}}

\newcommand{\FF}{\mathbb{F}}
\newcommand{\GG}{\mathbb{G}}
\newcommand{\MS}{\mathbb{S}}
\newcommand{\AF}{\mathbb{A}}

\newcommand{\calO}{\mathcal{O}}
\newcommand{\Sh}{\mathit{Sh}}
\newcommand{\TAF}{\mathit{TAF}}
\newcommand{\TMF}{\mathit{TMF}}

\DeclareMathOperator{\Aut}{Aut}
\DeclareMathOperator*{\Tot}{Tot}

\DeclareMathOperator{\AbVar}{AbVar}
\DeclareMathOperator{\Top}{Top}
\usepackage{graphicx}

\DeclareMathOperator{\inv}{inv}
\DeclareMathOperator{\Tr}{Tr}

\DeclareMathOperator{\VSub}{VSub}
\DeclareMathOperator{\Pic}{Pic}
\DeclareMathOperator{\Def}{Def}
\DeclareMathOperator{\Spf}{Spf}
\DeclareMathOperator{\Alg}{Alg}
\DeclareMathOperator{\im}{Im}
\DeclareMathOperator{\Set}{Set}
\DeclareMathOperator{\Sp}{Sp}
\DeclareMathOperator{\Pre}{Pre}
\DeclareMathOperator{\Shv}{Shv}
\DeclareMathOperator{\CoInd}{CoInd}
\DeclareMathOperator{\Ind}{Ind}
\DeclareMathOperator{\Res}{Res}
\DeclareMathOperator{\Orb}{Orb}
\DeclareMathOperator{\Fr}{Fr}
\DeclareMathOperator{\Cor}{Cor}

\begin{document}

\frontmatter

\title{Topological Automorphic Forms}

\author{Mark Behrens}
\address{Department of Mathematics, Massachusetts Institute of Technology, 
Cambridge, MA 02140}
\email{mbehrens@math.mit.edu}
\thanks{The first author was partially supported by NSF Grant 
\#0605100, a grant from the Sloan
foundation, and DARPA}

\author{Tyler Lawson}
\address{Department of Mathematics, University of Minnesota, 
Minneapolis, MN 55455}
\email{tlawson@math.umn.edu}
\thanks{The second author was partially supported by NSF Grant \#0402950}

\date{November 7, 2007}

\begin{abstract}
We apply a theorem of J.~Lurie to produce cohomology theories associated to 
certain Shimura varieties of type $U(1,n-1)$. 
These cohomology theories of topological automorphic forms ($\TAF$) are 
related 
to Shimura varieties in the same way that $\TMF$ is related to the moduli 
space of elliptic curves. We study the cohomology operations on these 
theories, and relate them to certain Hecke algebras. We compute the 
$K(n)$-local homotopy types of these cohomology theories, and determine that 
$K(n)$-locally these spectra are given by finite products of homotopy fixed 
point spectra of the Morava E-theory $E_n$ by finite subgroups of the Morava 
stabilizer group. We construct spectra $Q_U(K)$ for compact open subgroups
$K$ of certain adele groups, generalizing the spectra $Q(\ell)$ 
studied by the first author in the modular case. We show that the spectra 
$Q_U(K)$ admit finite resolutions by the spectra $\TAF$, arising from the 
theory of buildings. We prove that the $K(n)$-localizations of the spectra 
$Q_U(K)$ are finite products of homotopy fixed point spectra of $E_n$ with 
respect to certain arithmetic subgroups of the Morava stabilizer groups, 
which N.~Naumann has shown (in certain cases) to be dense. Thus the spectra 
$Q_U(K)$ approximate the $K(n)$-local sphere to the same degree that the 
spectra $Q(\ell)$ approximate the $K(2)$-local sphere.
\end{abstract}

\maketitle

%    Dedication.  If the dedication is longer than a line or two,
%    remove the centering instructions and the line break.
%\cleardoublepage
%\thispagestyle{empty}
%\vspace*{13.5pc}
%\begin{center}
%  Dedication text (use \\[2pt] for line break if necessary)
%\end{center}
%\cleardoublepage

%    Change page number to 7 if a dedication is present.
\setcounter{page}{4}

\tableofcontents

\chapter*{Introduction}
\label{chap:intro}

\section{Background and motivation}\label{sec:motivation}

\subsection*{The chromatic filtration}

Let $X$ be a finite spectrum, and let $p$ be a prime number.  
The chromatic tower of $X$ is the
tower of Bousfield localizations
$$ \cdots \rightarrow X_{E(3)} \rightarrow X_{E(2)} \rightarrow X_{E(1)}
\rightarrow X_{E(0)} $$
with respect to the $p$-primary 
Johnson-Wilson spectra $E(n)$ (where $E(0)$ is the
Eilenberg-MacLane spectrum $H\QQ$).  The chromatic convergence theorem of
Hopkins and Ravenel \cite{ravenelorange} states that the $p$-localization
$X_{(p)}$ is recovered by taking the homotopy
inverse limit of this tower, and that the homotopy groups of $X_{(p)}$ 
are given
by
$$ \pi_* X_{(p)} = \varprojlim_n \pi_* X_{E(n)}. $$

The salient feature of this approach is that the homotopy groups of the
monochromatic layers \index{monochromatic layer} $M_n X$ \index{2MNX@$M_n X$}, given by the homotopy fibers
$$ M_n X \rightarrow X_{E(n)} \rightarrow X_{E(n-1)}, $$
fit into periodic families.  

We pause to explain how this works.
The periodicity theorem of Hopkins 
and Smith \cite{hopkinssmith} implies that for a cofinal collection of
indices $I = (i_0, \ldots, i_{n})$ there exist finite
complexes $M(I)^0 = M(p^{i_0}, \ldots, v_n^{i_n})^0$ \index{2MI0@$M(I)^0$}
inductively given by fiber
sequences:
$$ 
M(p^{i_0}, \ldots, v_{n}^{i_{n}})^0 \rightarrow 
M(p^{i_0}, \ldots, v_{n-1}^{i_{n-1}})^0 \xrightarrow{v_n^{i_n}} 
\Sigma^{-i_n\abs{v_n}} M(p^{i_0}, \ldots, v_{n-1}^{i_{n-1}})^0
$$
Here, $v_n^{i_n}$ is a $v_n$-self map --- it induces an isomorphism on
Morava $K(n)$-homology.  (The $0$ superscript is used 
to indicate that we have
arranged for the top cell of these finite complexes to be in dimension $0$.) 
The $n$th monochromatic layer admits the following
alternate description: there is an equivalence
\begin{equation}\label{eq:monochromatic}
M_n X \simeq \hocolim_{I = (i_0, \ldots, i_{n-1})} 
X_{E(n)} \wedge M(I)^0. 
\end{equation}
The homotopy colimit is taken with respect to an appropriate cofinal
collection of indices $(i_0, \ldots, i_{n-1})$.
Suppose that
$$ \alpha : S^k \rightarrow M_n X $$
represents a non-trivial 
element of $\pi_k M_n X$.  Then, using (\ref{eq:monochromatic}), there
exists a sequence $I = (i_0, \ldots, i_{n-1})$ such that $\alpha$ factors as a
composite
$$ S^k \xrightarrow{\tilde{\alpha}} X_{E(n)} \wedge M(I)^0 
\rightarrow M_n X. $$
Let
$$ v_n^{i_n} : \Sigma^{i_n \abs{v_n}} M(I)^0
\rightarrow M(I)^0 $$
be a $v_n$-self map.  The map $v_n^{i_n}$ is an equivalence after smashing with
$X_{E(n)}$ --- this is because $v_n^{i_n}$ is a $K(n)$-equivalence, and the
smash product of an $E(n)$-local spectrum with a type $n$ complex is
$K(n)$-local.  For each integer $s$, the composite
$$ \alpha_s : 
S^{k + i_ns\abs{v_n}} \xrightarrow{\tilde{\alpha}} 
\Sigma^{i_ns\abs{v_n}} X_{E(n)} \wedge M(I)^0 \xrightarrow{(v_n^{i_n})^s} 
X_{E(n)}  \wedge M(I)^0 \rightarrow M_n X $$
gives rise to a family of elements $\alpha_s \in \pi_* M_n X$ which
specializes to $\alpha$ for $s = 0$.  

There are homotopy pullback squares:
$$
\xymatrix{
X_{E(n)} \ar[r] \ar[d] &
X_{K(n)} \ar[d]
\\
X_{E(n-1)} \ar[r] &
(X_{K(n)})_{E(n-1)}
}
$$
In particular, taking the fibers of the vertical arrows, 
there is an equivalence $M_n X \simeq M_n (X_{K(n)})$, and
understanding the $n$th monochromatic layer of $X$ 
is equivalent to understanding
the $K(n)$-localization of $X$.

The idea of the chromatic tower originates with Jack Morava, who also
developed computational techniques for understanding the homotopy groups of
$X_{K(n)}$ and $M_n X$.
For simplicity, let us specialize our discussion to the case where $X$ is
the sphere spectrum $S$.
Let $E_n$ \index{2En@$E_n$} be the Morava $E$-theory spectrum associated to the 
Honda height
$n$ formal group $H_n$ \index{2Hn@$H_n$} over $\br{\FF}_p$. 
It is the Landweber exact cohomology theory whose formal group is
the Lubin-Tate universal deformation of $H_n$, with coefficient ring
$$ (E_n)_* = W(\br{\FF}_p)[[u_1, \ldots, u_{n-1}]][u^{\pm 1}] $$
(here $W(\br{\FF}_p)$ \index{2Wk@$W(k)$} is the Witt ring of $\br{\FF}_p$).
Let 
$$ \GG_n = \MS_n \rtimes Gal $$
denote \index{2Gn@$\GG_n$} the extended Morava stabilizer group.  
Here $Gal$ \index{2Gal@$Gal$} 
denotes the
absolute Galois group of $\FF_p$, and $\MS_n = \Aut(H_n)$ 
\index{2Sn@$\MS_n$} 
is the 
$p$-adic analytic
group of automorphisms of the formal group $H_n$.
Morava's change of rings theorems state that the Adams-Novikov spectral
sequences for $M_n S$ and $S_{K(n)}$ take the form:
\begin{align}
\label{eq:Mnss}
H^*_c(\GG_n ; (E_n)_*/(p^\infty, \ldots, u_{n-1}^\infty)) &  
\Rightarrow \pi_* M_n S \\
\label{eq:Knss}
H^*_c(\GG_n ; (E_n)_*) & \Rightarrow \pi_* S_{K(n)}
\end{align}
Thus the homotopy groups of $M_n S$ and $S_{K(n)}$  
are intimately related to the formal 
moduli of commutative $1$-dimensional formal
groups of height $n$.  

The action of $\GG_n$ on the ring $(E_n)_*$ is understood
\cite{devinatzhopkinsaction}, but is very complicated.  The computations of
the $E_2$-terms of the spectral sequences (\ref{eq:Mnss}) and (\ref{eq:Knss}) 
are only known for small
$n$, and even in the case of $n = 2$ the 
computations of Shimomura, Wang, and Yabe
\cite{shimomurawang}, \cite{shimomurayabe}, 
while quite impressive, are difficult to comprehend fully. 

Goerss and Hopkins
\cite{goersshopkins} extended work of Hopkins and Miller \cite{rezk} to
show that $E_n$ is an $E_\infty$-ring spectrum, and that the group $\GG_n$ acts
on $E_n$ by $E_\infty$-ring maps.  Devinatz and Hopkins \cite{devinatzhopkins} 
gave a construction
of continuous homotopy fixed points of $E_n$ with respect to closed
subgroups of $\GG_n$, and refined Morava's change of rings theorem
(\ref{eq:Knss}) to prove
that the natural map
\begin{equation}\label{eq:DHequiv}
S_{K(n)} \xrightarrow{\simeq} E_n^{h\GG_n}
\end{equation}
is an equivalence.   This result represents a different sort of
computation: the \emph{homotopy type} of the $K(n)$-local sphere is given
as the hypercohomology of the group $\GG_n$ with coefficients in the
\emph{spectrum} $E_n$.

\subsection*{Finite resolutions of $S_{K(n)}$}

In \cite{ghmr}, \cite{henn}, Goerss, Henn, Mahowald, and Rezk proposed a
conceptual framework in which to understand the computations of Shimomura,
Wang, and Yabe.  Namely,
they produced finite 
decompositions of the $K(n)$-local sphere into homotopy
fixed point spectra of the form $(E_n^{hF})^{hGal}$ for various
\emph{finite} subgroups $F$ of the group $\MS_n$ (\cite{ghmr} treats the
case $n = 2$ and $p=3$, and \cite{henn} generalizes this to
$n=p-1$ for $p > 2$).  In their work, explicit
finite resolutions of the trivial $\MS_n$-module $\ZZ_p$ by permutation
modules of the form $\ZZ_p[[\MS_n/F]]$ are produced.  Obstruction theory is
then used to realize these resolutions in the category of spectra.  
These approaches often yield very efficient resolutions.
However, the algebraic resolutions are 
non-canonical.  Moreover, given an algebraic resolution, there can be 
computational
difficulties in showing that the obstructions to a topological realization
of the resolution vanish.

When the $1$-dimensional formal group $H_n$ occurs as the formal 
completion of a
$1$-dimensional commutative group scheme $A$ in characteristic $p$, 
then global methods may be used to study $K(n)$-local homotopy theory.  The
formal group is naturally contained in the $p$-torsion of the group
scheme $A$.  The geometry of the torsion at primes $\ell \ne p$ can be
used to produce \emph{canonical} 
finite resolutions analogous to those discussed above, as
we now describe.

\subsubsection*{The $K(1)$-local sphere and the $J$ spectrum}

In the case of $n = 1$, the formal completion of the
multiplicative group $\GG_m/\br{\FF}_p$ is isomorphic to the height $1$
formal group $H_1$. 
The $J$ 
\index{2J@$J$}
\index{J spectrum@$J$ spectrum} spectrum is given by the fiber sequence
\begin{equation}\label{eq:Jsequence}
J \rightarrow \mathit{KO}_p \xrightarrow{\psi^\ell - 1} \mathit{KO}_p
\end{equation}
where
$\ell$ is a topological generator of $\ZZ_p^\times$
(respectively $\ZZ_2^\times/\{\pm 1\}$ if $p = 2$),
and $\psi^\ell$ is the
$\ell$th Adams operation.
The natural map
\begin{equation}\label{eq:Adams-Baird}
S_{K(1)} \xrightarrow{\simeq} J
\end{equation}
was shown to be an equivalence by Adams and Baird, and also by Ravenel
\cite{bousfield}, \cite{ravenel84}.
The group of endomorphisms
$$ \End(\GG_m)[1/\ell]^\times = (\ZZ[1/\ell])^\times = \pm \ell^\ZZ $$
is a dense subgroup of the Morava stabilizer group $\MS_1 = \ZZ_p^\times$.
The $p$-adic $K$-theory spectrum is given by
$$ KO_p \simeq (E_1^{h\{\pm 1\}})^{hGal} $$
and under this equivalence, the fiber sequence (\ref{eq:Jsequence}) is
equivalent to the fiber sequence of homotopy fixed point spectra:
\begin{equation}\label{eq:Esequence}
(E_1^{\pm \ell^\ZZ})^{hGal} \rightarrow (E_1^{h\{\pm 1\}})^{hGal}
\xrightarrow{[\ell]-1} (E_1^{h\{\pm 1\}})^{hGal}.
\end{equation}
We see that there is an equivalence
\begin{equation}\label{eq:Jequiv}
J \simeq (E_1^{h\pm \ell^\ZZ})^{hGal}.
\end{equation}
Equivalence~(\ref{eq:Adams-Baird}) is then a statement concerning
the equivalence
$$ (E_1^{h\MS_1})^{hGal} \xrightarrow{\simeq} (E_1^{h\pm \ell^\ZZ})^{hGal} $$
of two homotopy fixed point spectra.

\subsubsection*{The $K(2)$-local sphere and the spectrum $Q(\ell)$}

A similar treatment of the case $n = 2$ was initiated by Mahowald and Rezk
\cite{mahowaldrezk}, and carried out further in previous work of the authors.
The multiplicative group is replaced with a supersingular
elliptic curve $C/\br{\FF}_p$.
The formal completion of $C$ is isomorphic to the height $2$
formal group $H_2$. 
The first author \cite{mark}
defined a spectrum $Q(\ell)$ \index{2Ql@$Q(\ell)$} to be the totalization of a
semi-cosimplicial spectrum arising from analogs of the Adams operations
occurring in the theory of topological modular
forms:
\begin{equation}\label{eq:Qelldef}
Q(\ell) = \holim \left(
\TMF \Rightarrow \TMF_0(\ell) \times \TMF \Rrightarrow
\TMF_0(\ell) \right).
\end{equation}
For simplicity, assume that $p$ is odd.
The authors \cite{marktyler} 
showed that if $\ell$ is a 
topological generator
of $\ZZ_p^\times$, the group
\begin{equation}\label{eq:gammadef}
\Gamma = (\End(C)[1/\ell])^\times
\end{equation}
is dense in the Morava stabilizer group $\MS_2$.  The first author
\cite{markbldg} used the
action of $\Gamma$ on the building for $GL_2(\QQ_\ell)$ to show there is an
equivalence
\begin{equation}\label{eq:Qellequiv}
Q(\ell)_{K(2)} \xrightarrow{\simeq} (E_2^{h\Gamma})^{hGal}.
\end{equation}
However, the spectrum $Q(\ell)_{K(2)}$ is not equivalent to $S_{K(2)}
\simeq
(E_2^{h\MS})^{hGal}$.  Rather, the conjecture is that there is a fiber
sequence
\begin{equation}\label{eq:cofiberconj}
D_{K(2)}Q(\ell)_{K(2)} \rightarrow S_{K(2)} \rightarrow Q(\ell)_{K(2)}.
\end{equation}
This conjecture was verified in \cite{mark} in the case of $p = 3$ and 
$\ell = 2$.

\subsection*{Periodic families in the stable stems and arithmetic
congruences}

Miller, Ravenel, and Wilson \cite{millerravenelwilson}, generalizing the
work of Adams, Smith, and Toda, 
suggested a prominent source of $v_n$-periodic families of
elements in the stable homotopy groups of spheres.  Namely, for $I =
(p^{i_0}, \ldots, v_{n-1}^{i_{n-1}})$, suppose that the associated complex
$M(I)^0$ admits a $v_n$-self map
$$ v_n^{i_n} : \Sigma^{i_n \abs{v_n}} M(I)^0 \rightarrow M(I)^0. $$
The composite
$$ S^{i_n s\abs{v_n}- \norm{I}} \xrightarrow{i} \Sigma^{i_n s \abs{v_n}}M(I)^0 
\xrightarrow{(v_n^{i_n})^s} M(I)^0
\xrightarrow{j} S, $$
where $i$ is the inclusion of the bottom cell, $\norm{I}$ is the quantity
$$ \norm{I} = i_1\abs{v_1} + i_2\abs{v_2} + \cdots + i_{n-1}\abs{v_{n-1}} + n,$$
and 
\index{2I@$\norm{I}$}
$j$ is the projection
onto the top cell,
gives the \emph{Greek letter element}  \index{Greek letter element}: 
$$ \alpha^{(n)}_{si_n/i_{n-1}, \ldots, i_0} \in \pi_* S. $$
(This \index{1alphan@$\alpha^{(n)}_{si_n/i_{n-1}, \ldots, i_0}$}
element is not uniquely
defined --- it depends inductively on the choices of self maps $(v_1^{i_1},
\ldots, v_n^{i_n})$.)  The image of this element in $\pi_* S_{E(n)}$ 
factors through
$M_nS$ as the composite
$$ S^{i_ns\abs{v_n}-\norm{I}} \rightarrow 
\Sigma^{i_ns\abs{v_n}} M(I)^0_{E(n)} 
\xrightarrow{(v_n^{i_n})^s} M(I)^0_{E(n)}
\hookrightarrow M_n S.
$$
This composite is necessarily 
detected in the $0$-line of the spectral sequence~\ref{eq:Mnss} by a
\emph{non-zero} element:
\begin{equation}\label{eq:invariants}
x_{i_n/i_{n-1}, \ldots, i_0} \in 
H^0(\GG_n, (E_n)_*/(p^\infty, \ldots, u_{n-1}^\infty))
\end{equation}
The \index{2Xin@$x_{i_n/i_{n-1}, \ldots, i_0}$} 
existence of these invariants determines which
sequences $(i_0, \ldots, i_n)$ \emph{could} give rise to a Greek letter
element $\alpha^{(n)}_{i_n/i_{n-1}, \ldots, i_0}$.  

There are
three important questions one can ask:

\begin{ques}\label{ques:Greek1}
For which sequences $(i_0, \ldots, i_n)$ do there exist invariants
$$ x_{i_n/i_{n-1}, \ldots, i_0} \in H^0(\GG_n, (E_n)_*/(p^\infty, \ldots,
u_{n-1}^\infty))? $$
\end{ques}

\begin{ques}\label{ques:Greek2}
Does there exist a corresponding complex $M(I)^0$ and $v_n$-self map
$v_n^{i_n}$?
\end{ques}

\begin{ques}\label{ques:Greek3}
If so, is the corresponding Greek letter element
$\alpha^{(n)}_{i_n/i_{n-1}, \ldots, i_0}$ non-trivial?
\end{ques}

Question~\ref{ques:Greek1} is a computation, which, as we will describe
below, is known for $n = 1,2$.  We will explain how the $J$-spectrum, in
the case $n = 1$, and the spectrum $Q(\ell)$ for $n = 2$, relate these
computations to certain arithmetic congruences.
Question~\ref{ques:Greek2} 
is a very difficult question in computational homotopy
theory: it was completely answered for $n = 1$ by Adams \cite{adamsJXIV} and 
Mahowald \cite{mahowald}, while even for $p \ge 5$, there are only
\emph{partial}
results for $n = 2$, due to Smith, Zahler, and Oka (see
\cite[Sec.~5.5]{ravenelgreen}).
Question~\ref{ques:Greek3} is more tractable.  
For instance, it is completely solved
in \cite{millerravenelwilson} and \cite{shimomura} for $n = 2$  
using the chromatic spectral sequence.

\subsubsection*{The elements $\alpha_{i/j}$ and Bernoulli numbers}

Assume that $p > 2$.  In what follows, let $\nu_p$ denote $p$-adic
valuation.  The homotopy elements $\alpha_{i/j}$ generate the image of the
classical $J$ homomorphism.  The corresponding invariants $x_{i/j}$ are
classified by the following theorem \cite[Thm.~4.2]{millerravenelwilson}.

\begin{thm}
There exists an invariant
$$ x_{i/j} \in H^0_c(\GG_1, (E_1)_{2t}/(p^\infty)) $$
of order $p^j$ if and only if $t = (p-1)i$ and $j \le \nu_p(i)+1$.
\end{thm}

The orders of the generators of these groups of invariants are related to
the $p$-adic valuations of denominators of Bernoulli numbers.  Let $B_i$
denote the $i$th Bernoulli number.

\begin{lem}[Lipshitz-Sylvester {\cite[Lem.~B.2]{milnorstasheff}}]
For every pair of integers $k$ and $n$, 
the quantity $k^{n}(k^{n} - 1)B_n/n$ is an integer.
\end{lem}

In particular, if $k$ is chosen to be prime to $p$, then we deduce that
$$ (k^n - 1)B_n/n $$
is $p$-integral.  If $\ell$ is chosen such that the subgroup $\pm \ell^\ZZ$
is dense in $\ZZ_p^\times$, there is an isomorphism
$$ H_c^0(\ZZ_p^\times, (E_1)_{2t}/(p^\infty))^{Gal} \xrightarrow{\cong} 
H^0(\pm \ell^\ZZ, (E_1)_{2t}/(p^\infty))^{Gal}. $$
The sequence (\ref{eq:Esequence}) implies that, for $t$ even, there is an exact
sequence
$$
0 \rightarrow 
H^0(\pm \ell^\ZZ, (E_1)_{2t}/(p^\infty))^{Gal} \rightarrow \ZZ/p^\infty
\xrightarrow{\ell^t-1} \ZZ/p^\infty. 
$$
In particular, the image of $B_t/t$ in $\ZZ/p^\infty$ is annihilated by
$\ell^t-1$, and $B_t/t$ gives an element $y_t$ of $H^0(\GG_1,
(E_1)_{2t}/(p^\infty))$.  In fact, these invariants $y_t$ are
generators.  In summary, there is a correspondence
$$  \alpha_{i/j}  \leftrightarrow  B_t/t \in \QQ/\ZZ_{(p)}  $$
for $t = (p-1)i$ and $j = \nu_p(i)+1$.

\subsubsection*{The elements $\alpha_{i/j}$ and Eisenstein series}

We now assume, for simplicity, that $p > 3$.
We now explain how the Eisenstein series give an alternative approach
to the invariants $x_{i/j}$.  By explicit calculation, 
there is an abstract isomorphism
$$ H^0(\GG_1, (E_1)_{2t}/(p^\infty)) \cong H^0(\GG_2,
(E_2[v_1^{-1}])_{2t}/(p^\infty) ). $$
(Presumably this is related to Hopkins' chromatic splitting conjecture.) 
Choosing $\ell$ as before, the group $\Gamma$ of (\ref{eq:gammadef}) is
a dense subgroup of $\MS_2$, and so there is an isomorphism
$$ H^0(\MS_2, (E_2[v_1^{-1}])_{2t}/(p^\infty))^{Gal} \cong 
H^0(\Gamma, (E_2[v_1^{-1}])_{2t}/(p^\infty))^{Gal}. $$
The equivalence (\ref{eq:Qellequiv}) relates this computation to the
computation of the equalizers
$$ \mc{A}_{t/j} \rightarrow \TMF_{2t}/p^j
\underset{d_1}{\overset{d_0}{\rightrightarrows}}
\TMF_0(\ell)_{2t}/p^j \oplus \TMF_{2t}/p^j  $$
where the maps $d_0$ and $d_1$ are the maps induced from the
semi-cosimplicial coface maps of (\ref{eq:Qelldef}).

For a $\ZZ[1/N]$-algebra $R$, let $M_t(\Gamma_0(N))_R$
\index{2MtG0NR@$M_t(\Gamma_0(N))_R$}
denote the weight
$t$ modular forms for the congruence subgroup $\Gamma_0(N)$ defined over
the ring $R$.
For $p > 3$ there are isomorphisms
\begin{align*}
(\TMF_{(p)})_{2t} & \cong M_t(\Gamma_0(1))_{\ZZ_{(p)}}[\Delta^{-1}], \\
(\TMF_0(\ell)_{(p)})_{2t} & \cong M_t(\Gamma_0(\ell))_{\ZZ_{(p)}}[\Delta^{-1}].
\end{align*}
The coface map $d_1$ acts on a modular form $f \in
M_t(\Gamma_0(1))_{\ZZ_{(p)}}[\Delta^{-1}]$ by
$$ d_1(f) = (f,f) \in M_t(\Gamma_0(\ell))_{\ZZ_{(p)}}[\Delta^{-1}] \oplus
M_t(\Gamma_0(1))_{\ZZ_{(p)}}[\Delta^{-1}]. $$
The coface map $d_0$ acts by
$$ d_0(f) = (\ell^t V_\ell(f),  \ell^t f) $$
where, on the level of $q$-expansions, the operator $V_\ell$ is given by
$$ (V_\ell (f))(q) = f(q^\ell). $$
Thus, \index{2Vl@$V_\ell$}
we have
$$
\mc{A}_{(t;j)} = \left\{
f \in M_t(\Gamma_0(1))[\Delta^{-1}]/p^j \: : \: 
\begin{array}{l}
\mathrm{(i)} \: (\ell^t - 1) f \equiv 0 \mod p^j \\
\mathrm{(ii)} \: \ell^t V_\ell(f) - f \equiv 0 \mod p^j
\end{array}
\right\}.
$$
Condition (i) implies that $\mc{A}_{(t;j)}$ cannot have elements of order
$p^j$ 
unless $t =
(p-1)p^{j-1}s$.  If this is satisfied, Condition (ii) on $f$ 
is equivalent to the assertion that its $q$-expansion satisfies
$$ f(q) \equiv \text{constant} \mod p^j. $$
The Eisenstein series $E_t \in M_t(\Gamma_0(1))_{\ZZ}$ have $q$-expansions
$$ E_t(q) = 1 - \frac{2t}{B_t} \left( \sum_{i \ge 1} \sigma_{t-1}(i) q^i
\right) $$
where
$$ \sigma_k(i) := \sum_{d \vert i}d^k. $$
Our remarks concerning the $p$-adic valuations of the 
denominators of $B_t/t$ imply that for $j \le \nu_p(i) + 1$, the Eisenstein
series $E_{(p-1)i}$ are generators of $\mc{A}_{(t;j)}$ of order $p^j$.  
Thus, for $t = (p-1)i$ and $j = \nu_p(i)+1$, there is a
correspondence
$$ \alpha_{i/j} \leftrightarrow E_t \in \mc{A}_{(t;j)}. $$

\begin{rmk}
This technique of computing the $1$-line of the Adams-Novikov spectral
sequence using Eisenstein series is first seen in the work of A.~Baker
\cite{bakerANSS}.  Baker uses Hecke operators, and his results may be
derived from the discussion above by exploiting the relationship between
the Hecke operator $T_\ell$ and the verschiebung $V_\ell$.
\end{rmk}

\subsubsection*{The elements $\beta_{i/j,k}$ and congruences of modular
forms}

We continue to assume that $p > 3$.  The sequences $(j,k,i)$ giving invariants
$x_{i/j,k}$ corresponding to homotopy elements $\beta_{i/j,k} \in S_{E(2)}$ are
completely classified:

\begin{thm}[Miller-Ravenel-Wilson
\cite{millerravenelwilson}]\label{thm:millerravenelwilson}
The exists an invariant
$$ x_{i/j,k} \in H^0_c(\GG_2, (E_2)_{2t}/(p^\infty, v_1^{\infty})) $$
of order $p^k$ if and only if
\begin{enumerate}
\item $t = (p^2-1)i - (p-1)j$, 
\item $1 \le j \le p^{\nu_p(i)} + p^{\nu_p(i)-1} - 1$ ($j = 1$ if $\nu_p(i)
= 0$),
\item If $m$ is the unique number satisfying
$$ p^{\nu_p(i)-m-1}+ p^{\nu_p(i)-m-2} - 1 < j \le 
p^{\nu_p(i)-m} + p^{\nu_p(i)-m-1} - 1
$$
then $k \le \min \{ \nu_p(j)+1, m+1 \}$.
\end{enumerate}
\end{thm}

We will explain how these elaborate conditions on $(j,k,i)$ given in
Theorem~\ref{thm:millerravenelwilson} reflect a
congruence phenomenon occurring amongst $q$-expansions of modular forms.
The discussion of congruence
properties of Eisenstein series in the last section specializes to give:
$$ E_{p-1}(q) \equiv 1 \mod p $$
and thus
$$ E_{p-1}^{p^{k-1}}(q) \equiv 1 \mod p^k. $$
By the $q$-expansion principle, multiplication by $E_{p-1}^{p^{k-1}}$ gives
an \emph{injection}:
$$ \cdot E_{p-1}^{p^{k-1}}: M_{t}(\Gamma_0(N))_{\ZZ/p^k} \hookrightarrow 
M_{t + (p-1)p^{k-1}}(\Gamma_0(N))_{\ZZ/p^k} $$
for any $N$ coprime to $p$.  In fact, there is a converse:

\begin{thm}[Serre {\cite[4.4.2]{katzpadic}}]
Let $N$ be coprime to $p$.
Suppose that for $t_1 < t_2$ 
we are given modular forms $f_i \in M_{t_i}(\Gamma_0(N))_{\ZZ_p}$  whose
$q$-expansions satisfy
$$ f_1(q) \equiv f_2(q) \mod p^k. $$
Then $t_1$ and $t_2$ are congruent modulo $(p-1)p^{k-1}$, and
$$ f_2 \equiv E_{p-1}^{\frac{t_2-t_1}{p-1}} \cdot f_1 \mod p^k. $$
\end{thm}

The Hasse invariant $v_1 \in M_{p-1}(\Gamma_0(N))_{\FF_p}$ lifts to
$E_{p-1}$ (see \cite[2.1]{katzpadic} --- it is necessary that $p > 3$).  
This allows us to reinterpret the
groups $H^0(\GG_2, (E_2)_{t}/(p^\infty, u_1^\infty))$ in terms of
$p$-adic congruences of modular forms.  If $\ell$ is chosen such that 
the group $\Gamma$
(\ref{eq:gammadef}) is dense in the Morava stabilizer group $\MS_2$, there
is an isomorphism
$$ H^0(\MS_2, (E_2)_{2t}/(p^\infty, u_1^\infty)) \cong H^0(\Gamma,
(E_2)_{2t}/(p^\infty, u_1^\infty)). $$
The equivalence (\ref{eq:Qellequiv}), together with the semi-cosimplicial
resolution (\ref{eq:Qelldef}), allow us to deduce that there is an
isomorphism
$$ H^0(\Gamma,
(E_2)_{2t}/(p^\infty, u_1^\infty)) = \colim_{j,k} \mc{B}_{(t;j,k)} $$
where the colimit is taken over pairs $(j,k)$ where $j \equiv 0 \mod
(p-1)p^{k-1}$ and the groups $\mc{B}_{(t;j,k)}$ for such $(j,k)$ are given by
the equalizer diagram:
$$
\mc{B}_{(t;j,k)} \rightarrow
\begin{array}{c}
\frac{M_t(\Gamma_0(1))_{\ZZ/p^k}[\Delta^{-1}]}{M_{t-j}
(\Gamma_0(1))_{\ZZ/p^k}[\Delta^{-1}]}
\end{array}
\underset{d_1}{\overset{d_0}{\rightrightarrows}}
\begin{array}{c}
\frac{M_t(\Gamma_0(\ell))_{\ZZ/p^k}[\Delta^{-1}]}{M_{t-j}
(\Gamma_0(\ell))_{\ZZ/p^k}[\Delta^{-1}]} \\
\oplus \\
\frac{M_t(\Gamma_0(1))_{\ZZ/p^k}[\Delta^{-1}]}{M_{t-j}
(\Gamma_0(1))_{\ZZ/p^k}[\Delta^{-1}]}
\end{array}
$$
Using our explicit description of the cosimplicial coface maps $d_0$ and
$d_1$, we see that
$$
\mc{B}_{(t;j,k)} \cong 
\left\{f \in 
\begin{array}{c}
\frac{M_t(\Gamma_0(1))_{\ZZ/p^k}[\Delta^{-1}]}{M_{t-j}
(\Gamma_0(1))_{\ZZ/p^k}[\Delta^{-1}]}
\end{array}
\: : \: 
\begin{array}{l}
\mathrm{(i)} \: 
(\ell^t-1)f(q) \equiv h(q) \mod p^k, \\
\: \qquad h \in M_{t-j}(\Gamma_0(1))[\Delta^{-1}]. \\
\mathrm{(ii)} \:
(\ell^t f(q^\ell) - f(q) \equiv g(q) \mod p^k, \\
\: \qquad g \in M_{t-j}(\Gamma_0(\ell))[\Delta^{-1}].
\end{array}
\right\}
$$
This translates into the following conclusion: there exist level $1$ 
modular forms
$f_t \in M_t(\Gamma_0(1))[\Delta^{-1}]$ such that $f_t(q)$ is not
congruent to any form of lower weight, and such that the existence of
$\beta_{i/j,k} \in \pi_*S_{E(2)}$ is \emph{equivalent} to the existence of a
congruence
$$ \ell^t f(q^\ell) - f(q) \equiv g(q) \mod p^k \quad \text{for} \quad 
g \in M_{t-j}(\Gamma_0(\ell))[\Delta^{-1}]. $$
for $t = (p^2-1)i$.

\begin{rmk}
In \cite{laures} G.~Laures defined the $f$ invariant, which gives an
embedding of the $2$-line of the Adams-Novikov spectral sequence into a
ring of divided congruences of modular forms.  The authors are unclear at
this time how Laures' congruences are related to the congruences described
here.
\end{rmk}

\section{Subject matter of this book}

The purpose of this book is to propose entries to the third column of the
following table:
\vspace{10pt}

\centerline{
\begin{tabular}{c|ccc}
& $n = 1$ & $n = 2$ & $n \ge 3$ \\
\hline
\\
$E$ & $KO$ & $\TMF$ & ? 
\\ \\
$Q$ & $J$ & $Q(\ell)$ & ?
\end{tabular}
}
(In fact, our constructions will also specialize to $n = 1,2$.)  Each of
the columns of the table gives spectra $E$ and $Q$, such that:
\begin{enumerate}
\item $E$ is an $E_\infty$-ring spectrum 
associated to a $1$ dimensional formal group of height $n$
residing in the formal completion of a commutative group scheme.
\item The localization $E_{K(n)}$ is a product of homotopy fixed point
spectra of $E_n$ for finite subgroups of $\MS_n$.
\item $Q$ admits a finite semi-cosimplicial resolution by
$E_\infty$-spectra related to $E$.
\item The localization $Q_{K(n)}$ is a homotopy fixed point spectrum of
$E_n$ for an infinite discrete group which is a dense subgroup of $\MS_n$.
\end{enumerate}

The moduli space of elliptic curves will be replaced by
certain \emph{PEL Shimura varieties} of type $U(1,n-1)$.  On first encounter,
these moduli spaces seem complicated and unmotivated (at least from the
point of view of homotopy theory).  There
are many excellent references available (for the integral models we are using in
this book the reader is encouraged to consult \cite{kottwitz},
\cite{harristaylor}, and \cite{hida}), but these can be initially 
inaccessible to homotopy theorists.  We therefore devote the first few
chapters of this book to an expository summary of the Shimura varieties we
are considering, singling out only those aspects relevant for our
applications to $K(n)$-local homotopy theory.

We briefly summarize the contents of this book.

\subsection*{Topological automorphic forms}

To obtain formal groups of height greater
than $2$, we must study abelian varieties $A/\br{\FF}_p$ \index{2A@$A$} of 
dimension greater than $1$.
Since only the moduli of $1$-dimensional formal groups seems to be relevant
to homotopy theory, we must introduce a device which canonically splits
off a $1$-dimensional formal summand from the formal completion
$\widehat{A}$ \index{2A@$\widehat{A}$}.  
This will be accomplished by fixing an imaginary quadratic
extension $F$ \index{2F@$F$} of $\QQ$, with ring of integers $\mathcal{O}_F$ 
\index{2OF@$\mc{O}_F$}
in which $p$ splits as $uu^c$ \index{2U@$u$} \index{2Uc@$u^c$}, 
and insisting that
$A$ admit complex multiplication by $F$, given by a ring homomorphism
$$i\co \mathcal{O}_F \hookrightarrow \End(A). $$
The \index{2I@$i$}
splitting $F_p \cong F_u \times F_{u^c}$ induces a splitting
$$ \widehat{A} \cong \widehat{A}_u \times \widehat{A}_{u^c}. $$
We \index{2Au@$\widehat{A}_u$}
shall insist that $\widehat{A}_u$ have dimension $1$.  The other summand
$\widehat{A}_{u^c}$ will be controlled by means of a compatible
polarization --- which is given by a certain kind of prime to $p$ isogeny
$$ \lambda \co A \rightarrow A^\vee $$
which \index{1lambda@$\lambda$}
induces an isomorphism
$$ \lambda_* \co A(u) \rightarrow A(u^c)^\vee $$
between the \emph{$p$-divisible groups} 
which extend the formal groups $\widehat{A}_u$,
$\widehat{A}_{u^c}^\vee$.

Moduli stacks of triples $(A,i,\lambda)$ (over $\ZZ_p$) where $A$ has dimension $n$ 
are given by Shimura stacks \index{Shimura stack}
$\Sh$ \index{2Sh@$\Sh$} 
associated to rational forms of the unitary group $U(1,n-1)$. 
The
formal groups $\widehat{A}_u$ associated to points of these Shimura
stacks attain heights up to and including height $n$.  There are two
caveats: 
\begin{enumerate}
\item In order to allow for more general forms of $U(1,n-1)$, we will
actually 
fix a maximal order $\mathcal{O}_B$ \index{2OB@$\mc{O}_B$} in   
a central simple algebra $B$ \index{2B@$B$} over $F$ of dimension $n^2$, which is
split at $u$.  The algebra 
$B$ shall be endowed with a positive involution $\ast$ which
restricts to conjugation on $F$.  We shall actually
consider moduli of triples $(A,i,\lambda)$ where $A$ is an abelian variety
of dimension $n^2$ with $\calO_B$-action given by an inclusion of rings
$$ i \co \calO_B \hookrightarrow \End(A) $$
and $\lambda$ is a compatible polarization.  The case of $n$-dimensional
polarized abelian varieties with complex multiplication by $F$ is recovered
by choosing $B$ to be the matrix algebra $M_n(F)$ by Morita equivalence.

\item Unlike the case of supersingular elliptic curves, there are
infinitely many isogeny classes of triples $(A,i,\lambda)$ for which
$\widehat{A}_u$ has height $n$.  Fixing an isogeny class is tantamount to
fixing a compatible pairing $\bra{-,-}$ \index{0bra@$\bra{-,-}$} 
on the rank $1$ left $B$-module $V = B$. \index{2V@$V$}
The points in the
isogeny class are those points $(A,i,\lambda)$ whose 
$\lambda$-Weil pairing on the $\ell$-adic Tate modules at
each of the places $\ell$ is $B$-linearly similar to $(V_\ell, \bra{-,-})$.
\end{enumerate}

For these more general $n^2$-dimensional 
$B$-linear polarized abelian varieties $(A, i,
\lambda)$, the formal group $\widehat{A}$ inherits an action of the
$p$-complete algebra $B_p \cong
M_n(F_p)$.  Fixing a rank $1$ projection $\epsilon \in M_n(F_p)$ 
\index{1epsilon@$\epsilon$}
gives a
splitting 
$$ \widehat{A} \cong (\epsilon \widehat{A}_u \times \epsilon
\widehat{A}_{u^c})^n,
$$
and we insist that the summand $\epsilon \widehat{A}_u$ is
$1$-dimensional.

If the algebra 
$B$ is a division algebra, then $\Sh$ possesses an \'etale projective
cover (this is the case considered in \cite{harristaylor}).  
The Shimura varieties $\Sh$ 
associated to other $B$ need not be compact,
but we do not pursue compactifications here.

Lurie has announced an extension of the Goerss-Hopkins-Miller theorem
to $p$-divisible groups.
Lurie's theorem gives rise to a homotopy sheaf of $E_\infty$-ring
spectra $\mc{E}$ \index{2E@$\mc{E}$} on $\Sh$ in the \'etale topology.  
The global sections 
of $\mc{E}$ give rise to a $p$-complete spectrum $\TAF$ \index{2TAF@$\TAF$} 
of topological
automorphic forms, \index{topological
automorphic forms} and a descent spectral sequence
$$ H^s(\Sh, \omega^{\otimes t}) \Rightarrow \pi_{2t-s}(\TAF). $$
Here, $\omega$ \index{1omega@$\omega$} 
is a line bundle on $\Sh$ such that over $\CC$, the
holomorphic sections of $\omega^{\otimes t}$ give rise to certain
determinant weight holomorphic automorphic forms.

\subsection*{Hypercohomology of adele groups}

Let $\AF$ \index{2A@$\AF$} denote the rational adeles.  If $S$ is a set of
places of $\QQ$, we shall use the convention that places in the superscript
are omitted and places in the subscript are included:
\begin{align*}
\AF_S & = {\prod_{v \in S}}' \QQ_v, \\
\AF^S & = {\prod_{v \not\in S}}' \QQ_v.
\end{align*}
If \index{2AS@$\AF^S$} \index{2AS@$\AF_S$}
$S$ is a set of finite rational primes, 
then this convention naturally extends to
the profinite integers:
\begin{align*}
\ZZ_S & = {\prod_{p \in S}} \ZZ_p, \\
\ZZ^S & = {\prod_{p \not\in S}} \ZZ_p.
\end{align*}
This \index{2ZS@$\ZZ^S$} \index{2ZS@$\ZZ_S$}
notational philosophy will be extended to other contexts as we see
fit.

Let $GU/\QQ$ \index{2GU@$GU$} 
be the algebraic group of similitudes of our fixed pairing
$\bra{-,-}$
on $V$:
$$ GU(R) = \{ g \in \Aut_B(V \otimes_\QQ R) \: : \: (gx, gy) = \nu(g)(x,y),
\: \nu(g) \in R^\times \}. $$
Let \index{1nu@$\nu$}
$U$ be the subgroup of isometries.
\index{2U@$U$}
Let $L$ \index{2L@$L$} be the $\mathcal{O}_B$-lattice in $V = B$ given by the 
maximal order $\calO_B$.
Let $K^p_0$ 
\index{2Kp0@$K^p_0$}
be the compact open subgroup of the group $GU(\AF^{p,\infty})$  
given by
$$ K^p_0 = \{ g \in GU(\AF^{p,\infty}) \: : \: g(L \otimes \ZZ^p) = L \otimes
\ZZ^p \}. $$
(Here, $\AF^{p,\infty} = \prod'_{\ell \ne p} \QQ_\ell$ is the ring of finite
adeles away from $p$, and $\ZZ^p = \prod_{\ell \ne p} \ZZ_\ell$ 
is the subring of integers.)
To each open subgroup $K^p$ 
\index{2Kp@$K^p$}
of $K^p_0$ one associates an \'etale cover
$$ \Sh(K^p) \rightarrow \Sh. $$
The \index{2ShKp@$\Sh(K^p)$}
stack $\Sh(K^p)$ is the moduli stack
of prime-to-$p$
isogeny classes of tuples of the form
$(A,i,\lambda,[\eta]_{K^p})$ where $i$ is a
compatible inclusion
$$ i \co \calO_{B,(p)} \hookrightarrow \End(A)_{(p)} $$
and
$$ \eta \co (V \otimes_\QQ \AF^{p,\infty}, \bra{-,-}) \xrightarrow{\cong} 
(V^p(A), \lambda\bra{-,-}) $$
represents \index{1eta@$\eta$} \index{etaKp@$[\eta]_{K^p}$}
a $K^p$ orbit $[\eta]_{K^p}$ of $B$-linear similitudes.
Associated to the \'etale cover $\Sh(K^p)$ is a spectrum $\TAF(K^p)$.
\index{2TAFKp@$\TAF(K^p)$}
Choosing $K^p = K^p_0$ recovers the spectrum $\TAF$.  

Associated to the tower of covers $\{\Sh(K^p)\}$ (as $K^p$ varies over the
open subgroups of $K^p_0$) is a filtered system of $E_\infty$-ring 
spectra $\{\TAF(K^p) \} $.  The colimit
$$ V_{GU} = \varinjlim_{K^p} \TAF(K^p) $$
admits \index{2VGU@$V_{GU}$}
an action of the group $GU(\AF^{p,\infty})$ by $E_\infty$-ring
maps, giving it the structure of a \emph{smooth $G$-spectrum}.
\index{smooth!$G$-spectrum}
The spectra $\TAF(K^p)$ are recovered by taking \emph{smooth homotopy fixed
points}
\index{smooth!homotopy fixed points}
$$ \TAF(K^p) \simeq V_{GU}^{hK^p}. $$

More generally, for any set of primes $S$ not containing $p$, and any
compact open subgroup $K^{p,S}$ contained in $GU(\AF^{p,S,\infty})$, we may
consider the homotopy fixed point spectrum
$$ Q_{GU}(K^{p,S}) := V_{GU}^{hK^{p,S}_+} $$
where \index{2QGU@$Q_{GU}$}
$$ K^{p,S}_+ := K^{p,S}GU(\AF_S) \subseteq GU(\AF^{p,\infty}). $$
Letting \index{2KpS+@$K^{p,S}_+$}
$GU^1(\AF_S)$ \index{2GU1@$GU^1$} be the subgroup
$$ GU^1(\AF_S) = \ker \left( GU(\AF_S) \xrightarrow{\nu} {\prod_{\ell \in
S}}'
\QQ_\ell^\times \xrightarrow{\nu_\ell} \prod_{\ell \in S} \ZZ \right) $$
of similitudes whose similitude norm has valuation $0$ at every place in
$S$, it is also convenient to consider the fixed point spectrum
$$ Q_{U}(K^{p,S}) := V_{GU}^{hK^{p,S}_{1,+}} $$
where \index{2QU@$Q_U$}
$$ K^{p,S}_{1,+} := K^{p,S}GU^1(\AF_S) \subseteq GU(\AF^{p,\infty}). $$
These \index{2KpS1+@$K^{p,S}_{1,+}$}
spectra should be regarded as generalizations of the spectra
$Q(\ell)$ to our setting.
As the results described in the next section demonstrate, the spectrum
$Q_U$ seems to be more closely related to the $K(n)$-local sphere than the
spectrum $Q_{GU}$. 

\subsection*{$K(n)$-local theory}

Let $\Sh(K^p)^{[n]}$ \index{2ShKpn@$\Sh(K^p)^{[n]}$} 
be the (finite) locus of $\Sh(K^p) \otimes_{\ZZ_p} \FF_p$ 
where the formal group
$\epsilon \widehat{A}_u$ is of height $n$.  
The $K(n)$-localization of the spectrum 
$\Sh(K^p)$ is given by the following theorem.

\begin{thm*}[Corollary~\ref{cor:K(n)localTAFQ}]
There is an equivalence
\begin{equation}\label{eq:K(n)localTAF}
\TAF(K^p)_{K(n)} \simeq \left( \prod_{(A,i, \lambda, [\eta]_{K^p}) 
\in \Sh(K^p)^{[n]}(\br{\FF}_p)}
E_n^{h\Aut(A,i,\lambda, [\eta]_{K^p})} \right)^{hGal}.
\end{equation}
\end{thm*}

Fix a $\br{\FF}_p$-point $(A,i,\lambda, [\eta]_{K^p})$ of $\Sh^{[n]}$, and
choose a representative $\eta$ of the $K^p$-orbit $[\eta]_{K^p}$.  
Let $\Gamma$ \index{1gamma@$\Gamma$} be the group of prime-to-$p$ 
quasi-isometries of $(A,i,\lambda)$.
$$ \Gamma = \{ \gamma \in \End_B(A)_{(p)}^\times \: : \: \gamma^\vee 
\lambda \gamma =
\lambda \} $$
The fixed representative $\eta$ induces an embedding
$$ \Gamma \hookrightarrow GU^1(\AF^{p,S,\infty}). $$
Let $\Gamma(K^{p,S})$ \index{1gammaK@$\Gamma(K)$} denote the subgroup
$$ \Gamma(K^{p,S}) = \Gamma \cap K^{p,S}. $$
The group $\Gamma(K^{p,S})$ acts naturally on the formal group
$$ \epsilon \widehat{A}_u \cong H_n $$
by isomorphisms, giving an embedding
$$ \Gamma(K^{p,S}) \hookrightarrow \MS_n. $$
We define a spectrum $E(\Gamma(K^{p,S}))$ 
\index{2EGKpS@$E(\Gamma(K^{p,S}))$}
to be the 
homotopy fixed point
spectrum
$$ E(\Gamma(K^{p,S})) = E_n^{h\Gamma(K^{p,S})}. $$
There is a natural map
\begin{equation}\label{eq:approximation}
E_n^{h\MS_n} \rightarrow E(\Gamma(K^{p,S})).
\end{equation}
The spectrum $E(\Gamma(K^{p,S}))$ 
should be regarded as some sort of approximation to the
$K(n)$-local sphere. 
At this time we are unable to even formulate a 
conjecture (analogous to (\ref{eq:Adams-Baird}) in the case of $n=1$ and 
(\ref{eq:cofiberconj}) in the case of $n=2$) quantifying the
difference between the spectrum $E(\Gamma(K^{p,S}))$ 
and the $K(n)$-local sphere.  Naumann \cite{naumann} has proven
that in certain circumstances, there exist small sets of primes $S$ 
such that
$\Gamma(K^{p,S})$ is of small index in the Morava stabilizer group $\MS_n$ 
(see Theorem~\ref{thm:naumann}).  The group $\Gamma$ is always dense in $\MS_n$
(Proposition~\ref{prop:plocaldense}).

The $K(n)$-localization of the spectrum 
$Q_U(K^{p,S})$ is related to Morava
$E$-theory by the
following theorem.

\begin{thm*}[Corollary~\ref{cor:K(n)localTAFQ}]
There is an equivalence
\begin{equation}\label{eq:K(n)localQU}
Q_U(K^{p,S})_{K(n)} \simeq \left( \prod_{[g] \in \Gamma \backslash 
GU^1(\AF^{p,S,\infty})/K^{p,S}}
E(\Gamma(gK^{p,S}g^{-1})) \right)^{hGal}.
\end{equation}
\end{thm*}

This theorem is an analog of (\ref{eq:Jequiv}) in the case $n=1$, and
(\ref{eq:Qellequiv}) in the case of $n=2$.

\subsection*{Building decompositions}

For simplicity, we now consider the case where the set $S$ consists of a single prime $\ell
\ne p$, where $B$ is split over all of the places that divide $\ell$.
The group $\Gamma(K^{p,\ell})$ sits naturally inside the group $U(\QQ_\ell)$ 
through
its action on the Tate module $V_\ell(A)$ through $B$-linear isometries.
Let $\mathcal{B}(U)$ be the Bruhat-Tits 
building for $U(\QQ_\ell)$.  The finite dimensional 
complex $\mc{B}(U)$ \index{2BU@$\mc{B}(U)$} is contractible, and
the group $U(\QQ_\ell)$ acts on it with compact open stabilizers.  This
action extends naturally to an action of the group $GU^1(\QQ_\ell)$.

\begin{thm*}[Proposition~\ref{prop:aut}]
The group $\Gamma(K^{p,\ell})$ acts on $\mc{B}(U)$ with finite stabilizers,
and these stabilizers are the automorphism groups of points of
$\Sh(K^p)^{[n]}(\br{\FF}_p)$ for various compact open subgroups $K^p$ of
$GU^1(\AF^{p,\infty})$.
\end{thm*}

Therefore, the virtual cohomological dimension of $\Gamma(K^{p,\ell})$ 
is equal to the
dimension of $\mc{B}(U)$.  We have
$$
\dim \mc{B}(U) = 
\begin{cases}
n & \text{$\ell$ splits in $F$}, \\
\text{$n/2$ or $(n-2)/2$} & \text{$\ell$ does not split in $F$, $n$ even}, \\
(n-1)/2 & \text{$\ell$ does not split in $F$, $n$ odd}. \\
\end{cases}
$$

We see that the virtual cohomological dimension of $\Gamma(K^{p,\ell})$ is 
maximized when $\ell$ splits in $F$.  In this case, there is an isomorphism
$$ U(\QQ_\ell) \cong GL_n(\QQ_\ell) $$
and the building $\mc{B}(U)$ is especially easy to understand.  
Even in the split case, unless $n = 1$, the map (\ref{eq:approximation})
cannot be an equivalence, because $\Gamma(K^{p,\ell})$ has virtual cohomological
dimension $n$, whereas $\MS_n$ has virtual cohomological dimension $n^2$. 

The action of $GU^1(\QQ_\ell)$ on $\mc{B}(U)$ gives rise to the following
theorem.

\begin{thm*}[Theorem~\ref{thm:mainthm}]
There is a
semi-cosimplicial spectrum $Q_{U}(K^{p,\ell})^\bullet$ 
\index{2Q@$Q^\bullet$}
of length $d = \dim \mc{B}(U)$
whose $s$th term ($0 \le s \le d$)
is given by
$$ Q_U(K^{p,\ell})^s = \prod_{[\sigma]} \TAF(K(\sigma)). $$
The \index{2Ks@$K(\sigma)$} \index{2Kls@$K_\ell(\sigma)$}
product ranges over $U(\QQ_\ell)$ orbits of 
$s$-simplices $[\sigma]$ in the building $\mathcal{B}(U)$.  
The groups $K(\sigma)$ are given by
$$ K(\sigma) = K^p K_\ell(\sigma) $$
where $K_\ell(\sigma)$ is the subgroup of $GU^1(\QQ_\ell)$ which stabilizes
$\sigma$.
There is an equivalence
\begin{equation}\label{eq:QUequiv}
Q_U(K^{p,\ell}) \simeq \Tot Q_U(K^{p,\ell})^\bullet.
\end{equation}
\end{thm*} 

\begin{rmk}
Theorem~\ref{thm:mainthm} gives a similar semi-cosimplicial resolution of
the spectrum $Q_{GU}(K^{p,\ell})$.
\end{rmk}

The spectra $Q_U(K^{p,\ell})$ 
should be regarded as a height $n$ analog of the $J$-theory
spectrum.  
The approach taken in this paper towards defining $Q_U(K^{p,\ell})$ is different
from that taken in defining the spectra $J$ and $Q(\ell)$.  The spectra $J$ and
$Q(\ell)$ were \emph{defined} as the totalization of an appropriate
semi-cosimplicial spectrum (\ref{eq:Jsequence}), (\ref{eq:Qelldef}).  Here,
we define $Q_U(K^{p,\ell})$ to be a hypercohomology spectrum, and then show
that it \emph{admits} a semi-cosimplicial decomposition.

\subsection*{Potential applications to Greek letter elements}

We outline some potential applications the spectra $\TAF$ and $Q_U$ could
have in the study of Greek letter elements in chromatic
filtration greater than $2$.  Specifically, we will discuss the
applicability of our methods towards resolving Questions~\ref{ques:Greek1}
and \ref{ques:Greek3}.  We do not anticipate our constructions to be
relevant towards the resolution of Question~\ref{ques:Greek2}.

\subsubsection*{Existence of $x_{i_n/i_{n-1}, \ldots, i_0}$}

Suppose that the set of primes $S$ and the compact open subgroup 
$K^{p,S}$ are chosen such that the corresponding
subgroup
$$ \Gamma(K^{p,S}) \subset \MS_n $$
is dense.  For instance, $S$ may consist of a single prime $\ell$
(see Theorem~\ref{thm:naumann}) or $S$ may always be taken to 
consist of all primes different
from $p$ (see Proposition~\ref{prop:plocaldense}).
Then the natural map
$$ H^0(\MS_n, (E_n)_{2t}/(p^\infty, \ldots, u_{n-1}^{\infty})) 
\rightarrow 
H^0(\Gamma(K^{p,S}), (E_n)_{2t}/(p^\infty, \ldots, u_{n-1}^{\infty}))$$
is an isomorphism.

Suppose that $S$ consists of a single prime $\ell$ which splits in $F$.
Then we have $U(\QQ_\ell) \cong GL_n(\QQ_\ell)$.
The equivalences (\ref{eq:K(n)localTAF}), (\ref{eq:K(n)localQU}), and
(\ref{eq:QUequiv}) may be combined to show that the existence of 
$x_{i_n/(i_{n-1}, \ldots, i_0)}$ is, under suitable hypotheses, 
detected in the equalizer of
the maps
$$ \pi_{2t} M_n \TAF(K^{p,\ell}K_\ell)
\underset{d_1}{\overset{d_0}{\rightrightarrows}}
\prod_{[\sigma]} \pi_{2t} M_n \TAF(K^{p,\ell}K_{\ell}(\sigma)). $$
Here $M_n$ denotes the $n$th monochromatic fiber, $K_{\ell} \subset
GL_n(\QQ_\ell)$ is a 
stabilizer of the unique orbit of $0$-simplices in the building
$\mc{B}(GL_n(\QQ_\ell))$, and the product ranges over the orbits $[\sigma]$ of
$1$-simplices in $\mc{B}(GL_n(\QQ_{\ell}))$, with stabilizers
$K_\ell(\sigma) \subset GL_n(\QQ_\ell)$.  It is our hope that this
relationship will lead to a description of the invariants
$x_{i_n/i_{n-1}, \ldots, i_0}$ (Question~\ref{ques:Greek1}) 
related to arithmetic properties of
automorphic forms analogous to the descriptions for $n = 1,2$ given in
Section~\ref{sec:motivation}.

\subsubsection*{Non-triviality of Greek letter elements}

Miller, Ravenel, and Wilson \cite{millerravenelwilson} construct Greek
letter elements in the $n$-line of the Adams-Novikov spectral sequence in
the following manner.  
There is a composition of connecting homomorphisms
$$
\partial_n: H^0(\GG_n, (E_n)_*/(p^\infty, \ldots, u_{n-1}^\infty)) \rightarrow
\Ext^{n,*}_{BP_*BP}(BP_*, BP_*).
$$
Given
$x_{i_n/i_{n-1}, \ldots, i_0} \in H^0(\GG_n, (E_n)_*/(p^\infty, \ldots,
u_{n-1}^\infty))$, its image under $\partial_n$ is the element in the
Adams-Novikov spectral sequence which will detect the corresponding Greek
letter element $\alpha^{(n)}_{i_n/i_{n-1}, \ldots, i_0}$ if it exists.
Question~\ref{ques:Greek3} is answered in the affirmative if 
$$ \partial_n(x_{i_n/i_{n-1}, \ldots, i_0}) \ne 0. $$
Because the length of the semi-cosimplicial resolution of $Q(\ell)$ is $n$, 
there is a corresponding sequence of boundary homomorphisms for the
spectrum $Q_U(K^{p,\ell})$, where Question~\ref{ques:Greek3} may have a
more tractable solution.

\section{Organization of this book}

Many chapters of this book are expository in nature, and are meant to
serve as a convenient place for the reader who is a homotopy theorist to
assimilate the necessary background information in a motivated and direct
manner.

In Chapter~\ref{chap:pdiv} we briefly describe the theory of $p$-divisible
groups, and recall their classification up to isogeny.

In Chapter~\ref{chap:tatehonda} we review the Honda-Tate classification of
abelian varieties up to isogeny over $\br{\FF}_p$.
Our aim is to establish that there is an obvious choice of 
isogeny class which supplies
height $n$ formal groups of dimension $1$.  We explain the generalization
to $B$-linear abelian varieties.

In Chapter~\ref{chap:tate} we describe the notion of a level structure, and
explain how homomorphisms of abelian varieties and abelian schemes may be
understood through the Tate representation.

In Chapter~\ref{chap:pol} we introduce the notion of a polarization of an
abelian variety, and its associated Weil pairing on the Tate module.  
We review the classification of polarizations up to
isogeny, which is essentially given by the classification of certain
alternating forms on the Tate module.

In Chapter~\ref{chap:galois} we explain some basic facts about alternating
and hermitian forms, and their groups of isometries and similitudes.  We
summarize the classification of alternating forms.  This classification
completes the classification of polarizations, and is subsequently applied
to the analysis of the height $n$ locus $\Sh^{[n]}$ of $\Sh$.

The classification of polarized abelian varieties up to isogeny suggests
a natural moduli problem to study, which we explain in
Chapter~\ref{chap:shimura}.  There are two equivalent formulations of the
moduli problem.  The moduli problem is representable by the Shimura stack
$\Sh$.

In Chapter~\ref{chap:deformations}, we summarize the Grothendieck-Messing
theory of deformations of $p$-divisible groups.  We then explain how this
is related to the deformation theory of abelian varieties via Serre-Tate
theory.  This machinery is then 
applied to understand the deformation theory of
mod $p$ points of $\Sh$.

In Chapter~\ref{chap:taf}, we summarize Lurie's generalization of the
Hopkins-Miller theorem.  We then use the deformation theory to apply this
theorem to our Shimura stacks to define spectra of 
topological automorphic forms.

Chapter~\ref{chap:autform} summarizes the classical holomorphic automorphic
forms associated to $\Sh$, and explains the relationship to topological
automorphic forms.

In Chapter~\ref{chap:smGspt}, we pause to rapidly develop the notion of a
smooth $G$-spectrum for $G$ a locally compact totally disconnected group.
We define smooth homotopy fixed points, and study the behavior of this
construction under restriction and coinduction. 
We prove that a smooth $G$-spectrum is determined up to equivalence by its
homotopy fixed points for compact open subgroups of $G$.  We define
transfer maps and observe that the homotopy fixed points of a smooth
$G$-spectrum define a Mackey functor from a variant of the Burnside
category to the stable homotopy category.

In Chapter~\ref{chap:ops} we describe $E_\infty$-operations on $\TAF(K^p)$ 
given by
elements of the group 
$GU(\AF^{p,\infty})$.  This gives the data to produce a smooth
$GU(\AF^{p,\infty})$-spectrum $V_{GU}$.  
We then describe how this structure extends to make the functor $\TAF(-)$ a
Mackey functor with values in the stable homotopy category.
We discuss the obstructions to deducing that the Hecke
algebra for the pair $(GU(\AF^{p,\infty}), K^p)$ acts on $\TAF(K^p)$
through cohomology operations, and give necessary conditions for these
obstructions to vanish.

In Chapter~\ref{chap:buildings}, we recall explicit lattice chain models for
the buildings for the groups $GU(\QQ_\ell)$ and $U(\QQ_\ell)$.

In Chapter~\ref{chap:Qspectra}, the spectra $Q_{GU}$ and $Q_U$ are defined.
We then describe the semi-cosimplicial resolution of $Q_{U}$.

In Chapter~\ref{chap:K(n)local}, we describe the height $n$ locus
$\Sh^{[n]}$.  We then relate the $K(n)$-localization of $\TAF$ and $Q_U$ to
fixed point spectra of Morava $E$-theory.

In Chapter~\ref{chap:examples}, 
we study the example of $n = 1$.  In this case the spectra $\TAF$ are
easily described as products of fixed points of the $p$ completion of
the complex $K$-theory spectrum.  The spectrum
$Q_U(K^{p,\ell})$ is observed to be equivalent to a product of 
$J$-spectra for $F$
and $\ell$ chosen suitably well.  The reader might find it helpful to 
look ahead to this section to see how the theory plays out in this
particular context.

\section{Acknowledgments}

The cohomology theory $\TAF$ studied in this book is not the first
generalization of $\TMF$ to higher chromatic levels.  Homotopy fixed points
$EO_n$ of $E_n$ by finite subgroups of the Morava stabilizer group 
were introduced by Hopkins and Miller, and have been studied by Hill, Mahowald,
Nave, and others.
In  
certain instances, the relationship of these spectra to moduli spaces
arising in algebraic geometry have also been considered by Hopkins, Gorbounov, 
Mahowald and Lurie in the case of $n = p-1$, and Ravenel in more general
cases.

The idea of associating cohomology theories to Shimura varieties also does
not originate with the authors --- Goerss, Lurie, and Morava have
considered such matters.  The authors are especially grateful to Jacob Lurie, 
who generously shared his work and results with the authors.  

The authors benefited greatly from conversations with Frank Calegari,
Daniel Davis, Johan
de Jong, Paul Goerss, Mike Hopkins, Haynes Miller, Catherine O'Neil, and 
Akshay Venkatesh.  The authors would like to thank Niko Naumann for his
comments, his observant corrections, 
and for keeping us informed of his timely work.
The first author also expresses his gratitude to Robert Kottwitz, 
for both
introducing him to the subject of Shimura varieties, and sharing his
time and knowledge.

\aufm{Mark Behrens \\ Tyler Lawson}

\mainmatter

\chapter{$p$-divisible groups}
\label{chap:pdiv}

If $A$ is an $n$-dimensional abelian variety over a field $k$, then
it has a group $A[p^\infty]$ \index{2Ap@$A[p^\infty]$} of $p$-torsion points defined over an
algebraic closure $\bar k$.  If $k$ has characteristic distinct from
$p$, this is a group abstractly isomorphic to $(\mb
Z/p^\infty)^{2n}$ with an action of the Galois group of $k$.

However, if $p$ divides the characteristic of $k$, this breaks down.
For example, an elliptic curve over a finite field can either be
ordinary (having $p$-torsion points isomorphic to the group $\mb
Z/p^\infty$) or supersingular (having no non-trivial 
$p$-torsion points whatsoever).
In these cases, however, these missing torsion points appear in the
formal group of the abelian variety.  
If the elliptic curve is ordinary, the formal group has
height $1$, whereas if the elliptic curve
is supersingular, the formal group has height $2$.  The missing rank in the
$p$-torsion arises as height in the formal group.

The correct thing to do in this context is to instead consider the
inductive system of group schemes given by the kernels of the homomorphism
$$ [p^i]\co A \rightarrow A. $$
The \index{2Pi@$[p^i]$}
ind-finite group schemes that arise are called $p$-divisible groups.

\section{Definitions}

\begin{defn}
A {\em $p$-divisible group\/} 
\index{p-divisible group@$p$-divisible group}
of height \index{p-divisible group@$p$-divisible group!height} 
$h$ over a scheme $S$ is a sequence of
commutative group schemes over $S$
\[
\{1\} = \mb G_0 \into \mb G_1 \into \mb G_2 \into \cdots
\]
such that each $\mb G_i$ is locally free of rank $p^{ih}$ over $S$, 
and such that for each $i \ge 0$ the sequence
$$ 0 \rightarrow \GG_i \hookrightarrow \GG_{i+1} \xrightarrow{[p^i]}
\GG_{i+1} $$
is exact.
\end{defn}

If $S$ is affine, this corresponds to a pro-system $\{R_i\}$
of bicommutative Hopf algebras over a ring $R$, with $R_i$ finitely generated
projective of rank $p^{ih}$.

\begin{exam}
Any finite group $G$ gives rise to a finite group scheme $G_S$ over $S$
$$ G_S = \coprod_{G} S. $$
The inductive system $\{ ((\ZZ/p^i)^h)_S \}$ gives a $p$-divisible group
$((\ZZ/p^\infty)^h)_S$ of height $h$.
\end{exam}

\begin{exam}\label{exam:formal}
Let $R$ be a $p$-complete ring, and let $\GG$ be a formal group over $R$.
Let $\br{\GG}$ denote the reduction of $\GG$ to $k = R/(p)$, and assume
that $\br{\GG}$ has constant and finite height $h$.  The the $p$-series of
$\br{\GG}$ takes the form
$$ [p]_{\br{\GG}}(x) = u x^{p^h} + \cdots \mod p$$
with $u \in k^\times$.
The $p^i$-series of $\br{\GG}$ is then given by
$$ [p^i]_{\br{\GG}}(x) = u' x^{p^{ih}} + \cdots \mod p $$
for $u' \in k^\times$.
We deduce, using Weierstrass preparation, that
the ring of functions of the formal group $\GG$ is given by
$$ \calO_\GG = R[[x]] = \lim_i R[[x]]/([p^i]_{\GG}(x)). $$
Therefore, $\GG$ may be regarded as a $p$-divisible group of height $h$
over $\Spec (R)$.
\end{exam}

\begin{exam}
If $A/S$ is an abelian scheme of dimension $n$, the inductive system of 
kernels $\{ A[p^i] \}$
\index{2Api@$A[p^i]$}
of the
$p^i$th power maps forms a $p$-divisible group 
\index{p-divisible group@$p$-divisible group!of an abelian scheme}
of height
$2n$, which we denote by $A(p)$.
\index{2Ap@$A(p)$}
\end{exam}

\begin{exam}
Suppose $E$ is an even periodic $p$-complete ring
spectrum.  The spectrum $E$ is complex orientable, and so 
associated to $E$ we have a formal group law $\GG_E$ over $E^0$,
\index{2GE@$\GG_E$}
represented by the complete Hopf algebra $E^0(\mb{CP}^\infty) \cong
E^0[[x]]$.  Assume that the height of the mod $p$ reduction of $\GG_E$ is
constant and finite.
Then Example~\ref{exam:formal} implies that
taking $\Spec$ of the successive quotients
$R_i = E^0[[x]]/([p^i]_{\GG_E}(x))$ 
gives a $p$-divisible group over $E^0$.
The homotopy equivalence
\[
\holim (\Sigma^\infty_+ B\mu_{p^k})^{\wedge}_p \to
{(\Sigma^\infty_+ BS^1)}^{\wedge}_p
\]
allows us to express the Hopf algebra
$E^0(\mb{CP}^\infty)$ as a limit of the Hopf algebras
$E^0(B\mu_{p^i}) = R_i$ (see, for example, \cite[2.3.1]{sadofsky}).  
\end{exam}

\begin{defn}
A \emph{homomorphism} 
\index{p-divisible group@$p$-divisible group!homomorphism}
of $p$-divisible groups $f: \GG \rightarrow \GG'$ is a
compatible sequence of homomorphisms
$$ f_i: \GG_i \rightarrow \GG'_i. $$
The homomorphism $f$ is an {\em isogeny\/} 
\index{p-divisible group@$p$-divisible group!isogeny}
if there exists a homomorphism
$f': \GG' \rightarrow \GG$ and an integer $k$ such that $ff' = [p^k]$ and
$f'f = [p^k]$. 
\end{defn}

If $\mb G$ is a $p$-divisible group over $S$, it has a Cartier dual
\index{p-divisible group@$p$-divisible group!Cartier dual} \index{2G@$\GG^\vee$}
\index{dual!$p$-divisible group}
\[
\mb G^\vee = \Hom_{\mathit{grp.schm}/S}(\mb G, \mb G_m).
\]
If $S = \Spec(R)$ and $\mb G_i = \Spec(R_i)$, $({\mb G}^\vee)_i$
corresponds to the dual Hopf algebroid $\Hom_R(R_i, R)$.  The
factorization
\[
\mb G_i \xrightarrow{[p]} \mb G_{i-1} \into \mb G_i
\]
of the $p$th power map dualizes to a sequence
\[
\mb G_i^\vee \hookleftarrow \mb G_{i-1}^\vee \xleftarrow{[p]}
\mb G_i^\vee.
\]
The left-hand maps make the family $\{\mb G_i^\vee\}$ into a
$p$-divisible group.  There is a natural isomorphism $\mb G \to (\mb
G^\vee)^\vee$, so $\vee$ is a anti-equivalence of the category of
$p$-divisible groups.

\section{Classification}

A $p$-divisible group is \emph{simple} 
\index{p-divisible group@$p$-divisible group!simple}
if it is not isogenous to a
non-trivial product of $p$-divisible groups.
If $\mb G$ is a simple 
$p$-divisible group over a field $k$ of characteristic $p$, 
there is a short
exact sequence
\[
0 \to \mb G^0 \to \mb G \to \mb G^{et} \to 0,
\]
where $(\mb G^0)_i$ is connected and $(\mb G^{et})_i$ is \'etale
over $k$ \cite{tatepdiv}.  The {\em dimension\/} 
\index{p-divisible group@$p$-divisible group!dimension}
of $\mb G$ is the dimension of the
formal group $\mb G^0$.  
\index{2G0@$\mb G^0$} \index{2Get@$\mb G^{et}$}
Height and dimension are both additive
in short exact sequences and preserved by isogeny.
A $p$-divisible group is classified up to
isogeny by this data.

\begin{thm}[\cite{demazure}, \cite{milne}]\label{thm:pdivclass}
Let $k$ be an algebraically closed field of characteristic $p$.
\begin{enumerate}
\item 
Any $p$-divisible group over $k$ is
isogenous to a product $\prod \mb G_i$ of simple $p$-divisible
groups, unique up to permutation.  

\item Simple $p$-divisible groups $\GG$ 
over $k$ are determined up to isogeny by a pair of relatively
prime integers $(d,h)$, $0 \leq d \leq h$, where $d$ is the
dimension of $\mb G$ and $h \ge 1$ is the height.  The fraction $0
\leq \frac{d}{h} \leq 1$ is called the {\em slope\/} of the $p$-divisible
group.

\item
The dual of such a $p$-divisible group $\mb G$ has height $h$
and dimension $h-d$.  Thus the slope of $\GG^\vee$ is $1 - \frac{d}{h}$.

\item
The endomorphism ring $\End(\mb G)$ is the unique maximal order
in the central division algebra $D$ over 
$\mb Q_p$ of invariant $\frac{d}{h}$.
\end{enumerate}
\end{thm}
  
The isogeny class of a $p$-divisible group is often represented by its
{\em Newton polygon\/}.  
\index{Newton polygon}
One draws a graph in which the horizontal
axis represents total height, and the vertical axis represents total
dimension.  Each simple summand of dimension $d$ and height $h$ is
represented by a line segment of slope $\frac{d}{h}$, 
\index{slope}
with the line
segments arranged in order of increasing slope.

The following is a Newton polygon of a $p$-divisible group of height
$6$ and dimension $3$, with simple summands of slope 
$\frac{1}{2}$,
$\frac{1}{3}$, and $1$.

%%%% ugly ugly ugly but I'm only drawing a couple at most
\vskip 0.7pc
\centerline{
\xy
0;/r1.5pc/:
(0,0);(6,0)**@{.},
(0,1);(6,1)**@{.},
(0,2);(6,2)**@{.},
(0,3);(6,3)**@{.},
(0,0);(0,3)**@{.},
(1,0);(1,3)**@{.},
(2,0);(2,3)**@{.},
(3,0);(3,3)**@{.},
(4,0);(4,3)**@{.},
(5,0);(5,3)**@{.},
(6,0);(6,3)**@{.},
(0,0)*{\bullet};
(3,1)**@{-}*{\bullet};
(5,2)**@{-}*{\bullet};
(6,3)**@{-}*{\bullet};
\endxy
}
\vskip 0.7pc

If $A$ is an abelian variety, the $p$-divisible group of the 
dual abelian variety $A^\vee$ is given by
 $A^\vee(p) = (A(p))^\vee$.  
Every abelian variety $A$ over a field admits a 
{\em polarization\/}, which is a certain type of isogeny 
$\lambda\co A \to A^\vee$ (see Chapter~\ref{chap:pol}, in particular
Remark~\ref{rmk:polexistence}).  
This forces the $p$-divisible group $A(p)$ to have a symmetry condition: 
if $\lambda$ is a slope appearing $m$ times in $A(p)$, so is 
$1 - \lambda$.

\begin{exam}
There are two possible Newton polygons associated to elliptic curves,
corresponding to the supersingular and ordinary types.
\vskip 0.7pc
\centerline{
\xy
0;/r1.5pc/:
(0,0);(2,0)**@{.},
(0,1);(2,1)**@{.},
(0,0);(0,1)**@{.},
(1,0);(1,1)**@{.},
(2,0);(2,1)**@{.},
(3,0);(5,0)**@{.},
(3,1);(5,1)**@{.},
(3,0);(3,1)**@{.},
(4,0);(4,1)**@{.},
(5,0);(5,1)**@{.},
(3,0)*{\bullet};
(4,0)**@{-}*{\bullet};
(5,1)**@{-}*{\bullet};
(0,0)*{\bullet};
(2,1)**@{-}*{\bullet};
\endxy
}
$$ \left( \substack{\text{super-} \\ \text{singular}} \right) 
\qquad \left( \substack{\text{ordinary} \\ \quad} \right) $$
\vskip 0.7pc
\end{exam}

Suppose that $\mb G/S$ is a $p$-divisible group, $F$ is a finite
extension field of $\mb Q$ with ring of integers ${\calO}_F$, and
$i\co {\calO}_F \to \End(\mb G)$ is a ring homomorphism.  Then
the map $i$ factors through the $p$-completion of ${\calO}_F$, and
gives rise to an associated map
\[
i_p \co {\calO}_{F,p} \to \End(\mb G).
\]
Let $\{u_i\}$ be the primes dividing $p$; then ${\calO}_{F,p}
\cong \prod {\calO}_{F,u_i}$.  In particular, the completion
contains idempotents $e_i$ whose images in $\End(\mb G)$ give a direct
sum decomposition
\[
\mb G \cong \bigoplus_i \mb G(u_i).
\]

If $A$ is an abelian variety with a map $\mathcal{O}_F \to \End(A)_{(p)}$
and $u$ is a prime of $F$ dividing $p$, we write $A(u)$ 
\index{2Au@$A(u)$}
for the
corresponding summand of its $p$-divisible group.  

\chapter{The Honda-Tate classification}
\label{chap:tatehonda}

In this chapter we will state the classification, due to Honda and Tate
\cite{honda}, \cite{tate}, of abelian varieties over $\FF_{p^r}$, and over
$\br{\FF}_p$.

\section{Abelian varieties over finite fields}

\index{abelian variety!quasi-homomorphism}
Define the set of \emph{quasi-homomorphisms} between abelian varieties 
$A$ and $A'$ to be the rational vector space
$$ \Hom^0(A,A') = \Hom(A, A') \otimes \QQ. $$
(The \index{2Hom0AA@$\Hom^0(A,A')$}
set $\Hom(A,A')$ is a finitely generated free abelian group.)
Let $\AbVar^0_k$ \index{2Abvar0k@$\AbVar^0_k$} denote the category whose objects are 
abelian varieties over a field $k$, and whose morphisms are the
quasi-homomorphisms.
A \emph{quasi-isogeny} \index{abelian variety!quasi-isogeny}
is an isomorphism in $\AbVar^0_k$.  An
\emph{isogeny} 
\index{abelian variety!isogeny}
is a quasi-isogeny which is an actual homomorphism.

A proof of the following theorem may be found in \cite[IV.19]{mumford}.

\begin{thm}
The category $\AbVar^0_k$ is semisimple.
\end{thm}
\index{abelian variety!simple}

\begin{cor}
For an abelian variety $A \in \AbVar^0_k$, the endomorphism ring
$\End^0(A)$ is semisimple.
\end{cor}

\index{2End0A@$\End^0(A)$}
Suppose that $A$ is a simple abelian variety over $\mb F_q$, where $q
= p^r$.  $A$ admits a {\em polarization\/} $\lambda \co 
A \rightarrow A^\vee$ (Chapter~\ref{chap:pol}).
The quasi-endomorphism ring $E = \End^0(A)$ is of finite dimension
over $\QQ$.  Because $A$ is simple, this is also a division algebra.
The polarization $\lambda$ induces a Rosati involution $\ast$ on $E$ (see
Section~\ref{sec:Rosati}).
Let $M$ \index{2M@$M$} be the center of $E$, and let $m$ be given by 
$$[E:M] = m^2.$$
The field $M$ is a finite extension of $\QQ$.  Let $d$ be the
degree $[M:\QQ]$.

The involution $\ast$ restricts to an involution $c$ on $M$.
Because the Rosati involution is given by conjugation with the
polarization, $c$ \index{2C@$c$}
is independent of the choice of polarization $\lambda$.
Let $M^+$ \index{2M+@$M^+$}
be the subfield of $M$ fixed by $c$.  The involution $\ast$
satisfies a positivity condition that implies that $M^+$ must be
totally real, and if $M \ne M^+$, then $M$ is a
totally imaginary quadratic extension of $M^+$.  In other words, $M$
is a \emph{CM field}.
\index{CM field}

Because $A$ is defined over $\mb F_q$, it has a canonical {\em
Frobenius endomorphism\/} 
\index{Frobenius} \index{1pi@$\pi$}
$\pi \in \End(A) \cap M$.  The ring
$\End(A)$ is finitely generated over $\mb Z$, and so $\pi$ must be in
the ring of integers ${\calO}_M$.

\begin{thm}[Weil, Honda-Tate]
\label{thm:htsimple}
$\quad$
\begin{enumerate}
\item The field $M$ is generated by $\pi$ over $\mb Q$.
\item The algebraic integer $\pi$ is a {\em Weil $q$-integer\/},
  \index{Weil $q$-integer}
  i.e. has absolute value $q^{1/2}$ in any complex embedding $M \into
  \mb C$.
\item The local invariants of $E$ at primes $x$ of $M$ are given by
\begin{alignat*}{2}
\inv_x E & = 1/2, & \qquad & \text{$x$ real}, \\
\inv_x E & = 0, & \qquad & x \not\vert\ p, \\
\inv_x E & = \frac{x(\pi)}{x(q)} [M_x:\mb Q_p], & \qquad & x \vert p.
\end{alignat*}
where 
\index{2Invx@$\inv_x$}
$x(-)$ denotes the valuation associated to a prime $x$.
\item The dimension of $A$ satisfies $2 \dim(A) = d \cdot m$, and $m$
  is the least common denominator of the $\inv_x E$.
\item The natural map
\begin{gather*}
\{ \text{quasi-isogeny classes of simple abelian varieties over $\FF_q$} \}
\\
\downarrow
\\
\{  \text{Weil $q$-integers}\} \Big/ \{\pi \sim i(\pi') \text{ for any } i\co M
\into M'\}
\end{gather*}
is a bijection.
\item For abelian varieties $A$ and $B$ and any prime $\ell$, the map
\[
\Hom_{\FF_q}(A,B)_\ell \to \Hom_{Gal(\br{\FF}_q/\FF_q)}(A(\ell),B(\ell))
\]
is an isomorphism.
\end{enumerate}
\end{thm}

The inclusion map $\mb F_q \to \mb F_{q^r}$ induces an
extension-of-scalars map
\[
A \mapsto A \mathop{\times}_{\Spec(\mb F_q)} \Spec(\mb F_{q^r}),
\]
which preserves quasi-isogenies.  Under the bijection of
Theorem~\ref{thm:htsimple}, the quasi-isogeny class corresponding to
the Weil $q$-integer $\pi$ is taken to the quasi-isogeny class
corresponding to $\pi^r$.

\section{Abelian varieties over $\br{\mb F}_p$}\label{sec:HTFpbar}

In this section we state the classification of abelian varieties over
$\br{\FF}_p$, which is a consequence of the classification over finite
fields.  We follow the treatment given in \cite[Sec.~V.2]{harristaylor}.

Suppose that $A$ is a simple abelian variety over $\br{\FF}_p$, $E =
\End^0(A)$, and $M$ the center of $E$.

Suppose that the prime $p$ splits in $M$ as
$$ (p) = x_1^{e_{x_1}} \cdots x_k^{e_{x_k}} c(x_1)^{e_{x_1}} \cdots
c(x_k)^{e_{x_k}} {x_1'}^{e_{x'_1}} \cdots {x_l'}^{e_{x'_l}} $$
where the primes $x_i'$ are $c$-invariant.  For each prime $x$ of $M$ dividing 
$p$, let $f_x$ denote the
residue degree, and define
\begin{align*}
d_x & = [M_{x}: \QQ_p]
\end{align*}
to be the local degree, so that $d_x = e_xf_x$.

The decomposition 
$M_p = \prod_{x|p} M_x$
induces a decomposition of the $p$-divisible group $A(p)$ (in the
quasi-isogeny category)
$$ A(p) \simeq \bigoplus_{x|p} A(x). 
$$

There is a corresponding 
decomposition of the $p$-complete algebra $E_p$ into simple
local algebras
$$ E_p = \prod_{x|p} E_{x}. $$

The following theorem appears in \cite{waterhousemilne}.

\begin{thm}[Tate]\label{thm:tate}
Let $A$ and $A'$ be abelian varieties over $\br{\FF}_p$.
The canonical map
$$ \Hom(A,A')_p \xrightarrow{\cong} \Hom(A(p),A'(p)) $$
is an isomorphism.
\end{thm}

Therefore, there are
isomorphisms
$$ E_{x} \cong \End^0(A(x)). $$
Because each of the algebras $E_x$ is
simple, each of the $p$-divisible groups $A(x)$ must have pure slope $s_x$
(or equivalently, be quasi-isogenous to a sum of isomorphic 
simple $p$-divisible groups of slope $s_x$).  By relabeling the primes, we
may assume that $s_{x_i} \le s_{c(x_i)}$.
Since $[E_{x}: M_x] = m^2$, $A(x)$ must be
of height $m$ as a $p$-divisible $M_x$-module, and hence it must have
height $d_xm$ as as a $p$-divisible group.  Express the dimension of
the $p$-divisible group $A(x)$ as
$$ \dim A(x) = \eta_x f_x m $$
for $\eta_x \in \QQ$.

We analyze the Newton polygons of these $p$-divisible groups.  For the
primes $x_i$, the polarization $\lambda$ induces a quasi-isogeny
$$ \lambda_*\co A(x_i) \rightarrow A(c(x_i))^\vee. $$
Thus we must have
$$ s_{x_i} = 1 - s_{c(x_i)}. $$
The dimension of $A(x_i)\oplus A(c(x_i))$ may be computed:
\begin{align*}
\dim A(x_i) \oplus A(c(x_i)) 
& = s_{x_i} d_{x_i} m + s_{c(x_i)} d_{c(x_i)} m \\
& = s_{x_i} d_{x_i} m + (1 - s_{x_i}) d_{x_i} m \\
& = d_{x_i} m.
\end{align*}
The Newton polygon for $A(x_i) \oplus A(c(x_i))$ therefore must take the
following form.

\begin{center}
\includegraphics{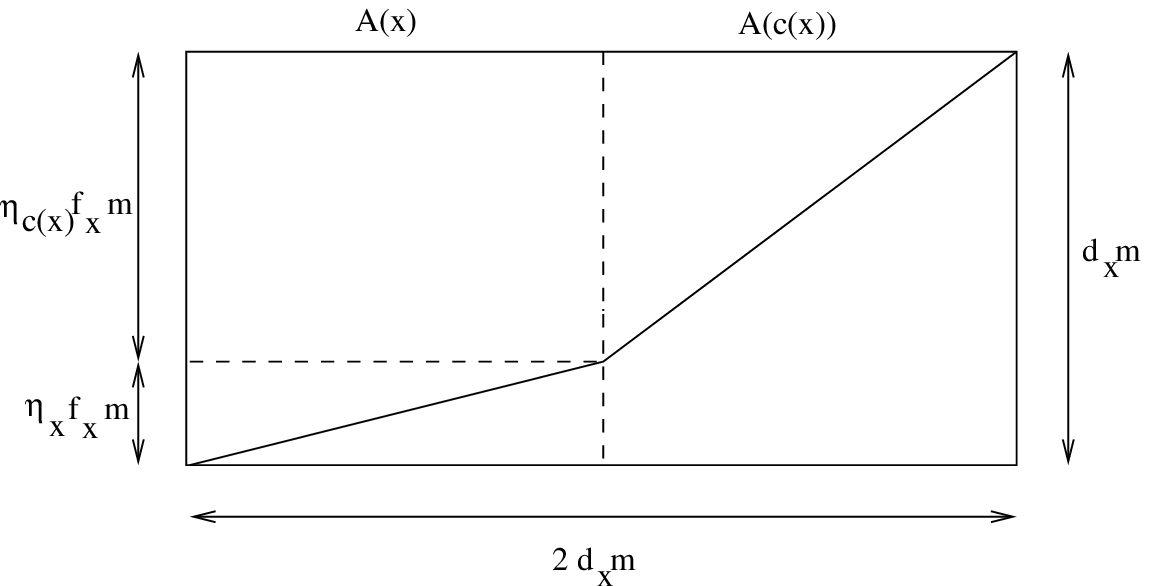}
\end{center}

For the primes $x'_j$, the polarization induces a quasi-isogeny
$$ \lambda\co A(x'_j) \rightarrow A(x'_j)^\vee. $$
This implies that $s_{x'_j} = 1 - s_{x'_j}$, which forces $s_{x'_j} = 1/2$.
The Newton polygon for $A(x'_j)$ therefore takes the following form.

\begin{center}
\includegraphics{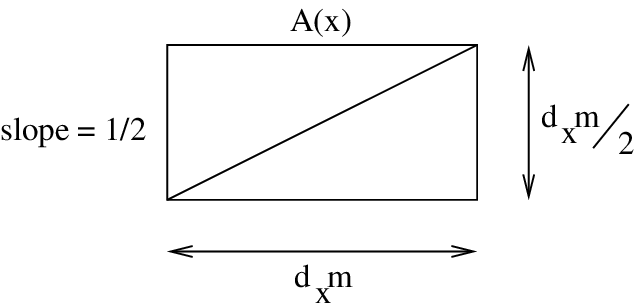}
\end{center}

Since the dimension of the $p$-divisible group $A(p)$ must be the
dimension of $A$, we deduce that
$$ \dim A = dm/2. $$
We also see that for all $x$, the values $\eta_x$ must satisfy
\begin{equation}\label{eq:padictype}
\eta_x / e_x + \eta_{c(x)} / e_x = 1.
\end{equation}
The slopes are recovered by the formula
$$ s_x = \eta_x /e_x. $$

The data $(M,(\eta_x))$ consisting of a CM field $M$ and $(\eta_x)$
satisfying condition~(\ref{eq:padictype}) is called a \emph{$p$-adic type}.
\index{p-adic type@$p$-adic type}
The $p$-adic type associated to $A$ is independent of the polarization
$\lambda$, and invariant under quasi-isogeny.  If $A$ is an abelian
variety defined over $\mb F_{p^r}$ with associated CM field $M$ and
Frobenius morphism $\pi \in M$, the associated $p$-adic type is
$\left(M,\left(x(\pi)/r\right)\right)$.  (Here the valuation $x$ is
normalized so that $x(p) = e_x$.)

We may form a category of $p$-adic types.  A morphism
$$ (M, (\eta_x)) \rightarrow (M', (\eta'_{x'})) $$
consists of an embedding $M \hookrightarrow M'$ such that if $x$ is a prime
of $M$, and $x'$ is a
prime of $M'$ dividing $x$, then we have
$$ \eta'_{x'} = e_{x'/x} \eta_x. $$
A $p$-adic type $(M, (\eta_x))$ is said to be \emph{minimal} 
\index{p-adic type@$p$-adic type!minimal}
if it is
not the target of a non-trivial morphism 
in the category of all $p$-adic types.  The $p$-adic type
associated to a simple abelian variety $A$ is minimal.

\begin{thm}[Honda-Tate]\label{thm:hondatate}
$\quad$
\begin{enumerate}

\item The natural map
\begin{gather*}
\{ \text{quasi-isogeny classes of simple abelian varieties over $\br{\FF}_p$} \}
\\
\downarrow
\\
\{ \text{minimal $p$-adic types} \}
\end{gather*}
is a bijection.  

\item If $A$ is the simple abelian variety associated to a
minimal $p$-adic type $(M,(\eta_x))$, then the local invariants of $E =
\End^0(A)$ are given by
\begin{alignat*}{2}
\inv_x E & = 1/2, & \qquad & \text{$x$ real}, \\
\inv_x E & = 0, & \qquad & x \not\vert\ p, \\
\inv_x E & = \eta_x f_x, & \qquad & x \vert p.
\end{alignat*}
\end{enumerate}
\end{thm}

For our homotopy theoretical applications, we are interested in the most basic
abelian varieties $A$ over $\br{\FF}_p$ for which the $p$-divisible group 
$A(p)$ contains a $1$-dimensional formal group of height $n$.  The
Honda-Tate classification tells us that for $n > 2$, 
the smallest such example
will have $p$-adic type $(M, (\eta_x))$ where
\begin{enumerate}
\item $M$ is a quadratic imaginary extension of $\QQ$ in which $p$ splits
as $x\, c(x)$.
\item $\eta_x = 1/n$ and $\eta_{c(x)} = (n-1)/n$. 
\end{enumerate}
Furthermore, given any quadratic extension $M$ as above with a chosen
prime $x$ over $p$, there exists a \emph{unique} quasi-isogeny class
of abelian variety of dimension $n$ with complex multiplication by $F$
such that $A(x)$ has slope $1/n$.

Let $B$ \index{2B@$B$} 
be a simple algebra whose center $F$ 
\index{2F@$F$}
is a CM field, so that $[B:F]
= s^2$.
We shall want to consider abelian varieties $A$ with the data of a fixed
embedding
$$ i\co B \hookrightarrow \End^0(A). $$
We consider the semisimple category $\AbVar^0_{\br{\FF}_p}(B)$ consisting
of such pairs $(A,i)$, where the morphisms are $B$-linear
quasi-homomorphisms.  We shall refer to such objects as 
\emph{$B$-linear abelian
varieties}.
\index{abelian variety!$B$-linear}

Kottwitz \cite[Sec.~3]{kottwitz} makes some very
general observations concerning this set-up, which we specialize to our
case below, following \cite[V.2]{harristaylor}.  A \emph{$p$-adic type over
$F$}
\index{p-adic type@$p$-adic type!over $F$}
is a $p$-adic type $(L, (\eta_x))$ where $L$ is an extension of $F$.
A $p$-adic type over $F$ is minimal if it is minimal in the category of
$p$-adic types over $F$.
Every simple $B$-linear abelian variety $(A,i)$ is isotypical 
\index{isotypical}
(isogenous to
$A_0^j$ for $A_0$ simple) when viewed
as an object of $\AbVar^0_{\br{\FF}_p}$ \cite[Lem.~3.2]{kottwitz}.

Given a $B$-linear abelian variety $(A,i)$ which is isogenous to
$A_0^j$ with $A_0$ simple, define the following:
\begin{align*}
E & = \End^0(A_0), \\
M & = \text{center of $E$}, \\
D & = \End_B^0(A), \\
L & = \text{center of $D$}, \\
t & = [D:L]^{1/2}.
\end{align*}
Observe that $M$ is naturally contained in $L$, since it is central in
$\End^0(A)$.  Observe also that $F$, when viewed as lying in $\End^0(A)$,
also must be contained in $L$.
Kottwitz \cite[Lem.~3.3]{kottwitz} observes that $L$ is a factor of $F
\otimes_\QQ M$.  He shows that this gives a correspondence
\begin{gather*}
\{ \text{simple $B$-linear abelian varieties $(A,i)$, $A$ isotypical of
type $A_0$} \} \\
\updownarrow \\
\{ \text{fields $L$ occurring as factors of $F\otimes_\QQ M$} \} 
\end{gather*}
Since $L$ is a factor of $F \otimes_\QQ M$, and both $F$ and $M$ are CM
fields, $L$ must be a CM field.

Let $(M, (\nu_y))$ be the $p$-adic type of $A_0$.
We may associate to $A$ the $p$-adic type $(L, (\eta_{x}))$ over
$F$, where, if $x$ is a prime of $L$ lying over a prime $y$ of $M$, 
the value of $\eta_x$ is given by
$$ \eta_x = e_{x/y} \nu_y. $$
If $A$ is simple in $\AbVar_{\br{\FF}_p}^0(B)$, then the associated
$p$-adic type $(L, \eta_x)$ is minimal over $F$. 

If $x$ is a place of $L$ which 
is not invariant under the action of the conjugation $c$, then the
Newton polygon for the $p$-divisible $B$-module $A(x)\oplus A(c(x))$ is
displayed below.

\begin{center}
\includegraphics{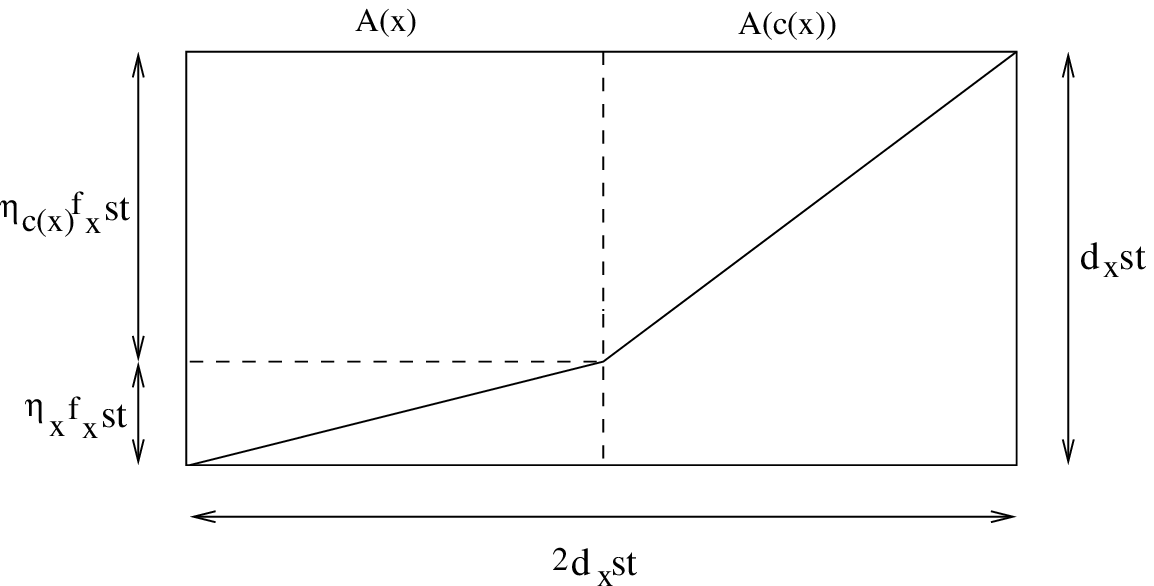}
\end{center}

Otherwise, if $x$ is a place of $L$ which is invariant under the
conjugation action, then the 
Newton polygon of the $p$-divisible group $A(x)$ takes
the following form.

\begin{center}
\includegraphics{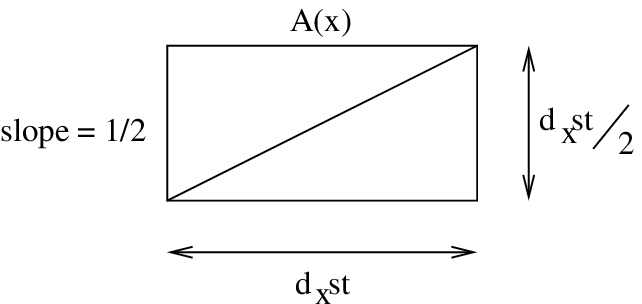}
\end{center}

The dimension of the abelian variety $A$ is related to 
this data by the formula
$$ \dim A = dst/2. $$

Kottwitz proves the following variation of the Honda-Tate classification
for $B$-linear abelian varieties \cite[Lem.~10.13]{kottwitz}.

\begin{thm}[Kottwitz]\label{thm:Bhondatate}
The simple objects of $\AbVar_{\br{\FF}_p}^0(B)$ are classified as follows.
\begin{enumerate}
\item The natural map
\begin{gather*}
\{ \text{quasi-isogeny classes of simple $B$-linear 
abelian varieties over $\br{\FF}_p$} \}
\\
\downarrow
\\
\{ \text{minimal $p$-adic types over $F$} \}
\end{gather*}
is a bijection.

\item If $A$ is the simple $B$-linear abelian variety associated to a
minimal $p$-adic type $(L,(\eta_x))$ over $F$, 
then the local invariants of the central division algebra $D =
\End^0_B(A)$ over $L$ are given by:
\begin{alignat*}{2}
\inv_x D & = 1/2 - \inv_x(B \otimes_F L), & \qquad & \text{$x$ real}, \\
\inv_x D & = -\inv_x(B \otimes_F L), & \qquad & x \not \vert p, \\
\inv_x D & = \eta_x f_x - \inv_x(B \otimes_F L), & \qquad & x \vert p.
\end{alignat*}
\end{enumerate}
\end{thm}

\chapter{Tate modules and level structures}
\label{chap:tate}

\section{Tate modules of abelian varieties}

Let $A$ be an abelian variety over an algebraically closed field $k$.
Assume
$\ell$ is a prime distinct from the characteristic of $k$.  Let
$$ A[\ell^\infty] \subseteq A(k) $$
denote \index{2Al@$A[\ell^\infty]$}
the $\ell$-torsion subgroup of $k$-points of $A$.  
If $A$ is an
abelian variety of dimension $g$, then there is an abstract isomorphism
$$ A[\ell^\infty] = \varinjlim_i A[\ell^i]\cong (\ZZ/\ell^\infty)^{2g}. $$

The \emph{(covariant) $\ell$-adic 
Tate module} 
\index{Tate module}
of $A$ is defined as the inverse limit
$$ T_\ell (A) = \varprojlim_i A[\ell^i] \cong \ZZ_\ell^{2g}. $$
The \index{2TlA@$T_\ell(A)$}
inverse limit is taken over the $\ell$th power maps.  We define
$$ V_\ell (A) = T_\ell(A) \otimes \QQ. $$
Let \index{2VlA@$V_\ell(A)$}
$\mb Z_\ell(1)$ denote the $\ell$-adic Tate module of the
\index{2Zl1@$\mb Z_\ell(1)$}
multiplicative group,
\[
\mb Z_\ell(1) = \varprojlim_i \mu_{\ell^i} \cong \mb Z_\ell.
\]

There are natural isomorphisms
$$ A^\vee[\ell^i] \cong \Hom(A[\ell^i], \mu_{\ell^i}),$$
and taking limits gives rise to a natural isomorphism
\[
T_\ell(A^\vee) \cong \Hom_{\mb Z_\ell}(T_\ell(A), \mb Z_\ell(1)).
\]

It is useful to consider the collection of all Tate modules at once.  Let
$p$ be the characteristic of $k$.  We will denote
$$ T^p(A) = \prod_{\ell \ne p} T_\ell(A). $$
Tensoring \index{2TpA@$T^p(A)$}
with $\QQ$, we get a module over $\AF^{p,\infty}$:
$$ V^p(A) = T^p(A) \otimes \QQ. $$
\index{2VpA@$V^p(A)$}

\section{Virtual subgroups and quasi-isogenies}
\label{sec:isogenysub}

Let $A$ be an abelian variety over $k$, an algebraically closed field
of characteristic $p$.

For a subring $R$ of $\mb Q$, an $R$-isogeny of abelian varieties $A
\to A'$ is a quasi-isogeny that lies in
\[
R \otimes_{\mb Z} \Hom(A,A').
\]
An invertible $R$-isogeny is one whose inverse is also an $R$-isogeny.

Let $A[\mathit{tor}^p]$ 
\index{2Atorp@$A[\mathit{tor}^p]$}
be the subgroup of $A(k)$ consisting of
torsion elements of order prime to $p$.  There is a natural
surjection 
$$ \pi_A\co V^p(A) \to A[\mathit{tor}^p] $$ 
with kernel $T^p(A)$.
For each finite subgroup $H \subset A[tor^p]$ there is a corresponding 
lattice
$\pi_A^{-1}(H) \subset V^p(A)$ such that $T^p(A) \subseteq \pi^{-1}(H)$.

More generally, we can define the set of prime-to-$p$ {\em virtual
subgroups\/} 
\index{virtual subgroup}
of $A$, $\VSub^p(A)$, 
\index{2VsubpA@$\VSub^p(A)$}
to be the set of lattices $K
\subset V^p(A)$.  For any such $K$, there exists some integer $N$ with
$N\cdot T^p(A) \subseteq K$ and $[K:N\cdot T^p(A)] < \infty$.

A $\mb Z_{(p)}$-isogeny $\phi\co A \to A'$ induces an isomorphism $V^p(A) \to
V^p(A')$.  We define the {\em kernel\/} of the quasi-isogeny $\phi$ to
be
\[
\Ker(\phi) = \phi^{-1}(T^p(A')) \in \VSub^p(A).
\]
\index{2Kerp@$\Ker(\phi)$}
\index{abelian variety!quasi-isogeny!kernel of}

Given a virtual subgroup $H \subset V^p(A)$ there is a way to
associate a quasi-isogeny with kernel $H$.  Let $N$ be an integer such
that $T^p(A) \subseteq N^{-1} H$, so the image of $N^{-1} H$ in $A$ is a
finite subgroup $K$.  The quasi-isogeny $\phi_H$ associated to
$H$ is the isogeny
$$ \phi_{H} : A \xrightarrow{\phi_K} A/K \xrightarrow{[N^{-1}]} A/K $$
where $\phi_{K}$ is the quotient isogeny associated to $K$, and
$[N^{-1}]$ is the quasi-isogeny which is inverse to the $N$th power map.
The quasi-isogeny $\phi_H$ actually depends on the choice of
$N$, but a different choice $N'$ yields an isogeny which is
\emph{canonically} isomorphic.  Indeed, consider the
following diagram.
$$
\xymatrix@C+4em@R+1em{
& A \ar[dr]^{\phi_K} \ar[d]^{N} 
\\
A \ar[ru]^{N'} \ar[r]_{NN'} \ar@{=}[d] &
A \ar@{=}[d] \ar[r]_{N^{-1}\phi_{K}} &
A/K \ar@{.>}[d]^{\exists!}_{\cong} 
\\
A \ar[r]^{NN'} \ar[dr]_{N} &
A \ar[r]^{(N')^{-1}\phi_{K'}} &
A/K'
\\
& A \ar[u]_{N'} \ar[ru]_{\phi_{K'}}
}
$$
The composites given by the two rows in the above diagram are 
actual isogenies, and 
the canonical isomorphism in the diagram above exists because the kernels
of these isogenies are equal.  For now we will commit the  
of abuse of defining $A/H$ to be $A/K$.

\section{Level structures}
\label{sec:levelsub}

Let $V$ be a rational vector space of dimension $2n$, 
and fix a $\ZZ$-lattice $L \subset V$.
Define
\begin{align*}
L^p & = \prod_{\ell \ne p} L \otimes \ZZ_\ell, \\
V^p & = V \otimes_\QQ \AF^{p,\infty} = L^p \otimes \QQ.
\end{align*}
\index{2Vp@$V^p$}
\index{2Lp@$L^p$}

Assume that $A$ is an abelian variety of dimension $n$ over an 
algebraically closed field $k$ of characteristic $p$.
By an \emph{integral uniformization} we shall mean an isomorphism
$$ \eta\co L^p \xrightarrow{\cong} T^p(A). $$
By \index{uniformization!integral}
a \emph{rational uniformization} 
\index{uniformization!rational}
we shall mean an isomorphism
$$ \eta\co V^p \xrightarrow{\cong} V^p(A), $$

Let $K^p_0 = \Aut(L^p)$ be  subgroup of $\Aut(V^p)$ given by the lattice
automorphisms:
$$ K^p_0 = \{ g \in \Aut(V^p) \: : \: g(L^p) = L^p \}.  $$
There \index{2Kp0@$K^p_0$}
is a natural action of $\Aut(V^p)$ on the set of rational
uniformizations of $A$.
Let $K^p$ be a subgroup of $\Aut(V^p)$.  A 
\emph{rational $K^p$-level structure} 
\index{level structure!rational}
$[\eta]_{K^p}$ is the $K^p$-orbit of a 
rational uniformization $\eta$.
We say that $[\eta]_{K^p}$ is an \emph{integral $K^p$-level structure}
if $\eta$ is an integral uniformization.
\index{level structure!integral}

Let $\mathcal{L}_{K^p_0}(A)$ 
\index{2LKp0A@$\mathcal{L}_{K^p_0}(A)$}
be the set of rational
$K_0^p$-level structures on $A$.  There
is a natural map
$$ \kappa\co \mathcal{L}_{K^p_0}(A) \rightarrow \VSub^p(A). $$
Given \index{1kappa@$\kappa$}
$[\eta] \in \mathcal{L}_{K^p_0}(A)$, define
$$ \kappa([\eta]) = \eta(L^p). $$
Indeed, there exists a positive integer $N$ such that $\eta(L^p)$ 
contains $N\cdot T^p(A)$.  Since the image $\eta(L^p)$ is
compact, the index $[\eta(L^p):N\cdot T^p(A)]$ must be finite.  Therefore
$\eta(L^p)$ is a virtual subgroup $\kappa(\eta) \in
\VSub^p(A)$.
Automorphisms of $L^p$ do not alter the image $\eta(L^p)$, so the
virtual subgroup $\eta(L^p)$ depends only on the $K_0^p$-orbit of $\eta$.

\begin{lem}\label{lem:integralkappa}
A rational level structure $\eta \co V^p \xrightarrow{\cong} V^p(A)$
restricts to an integral level structure if and only if the virtual
subgroup $\kappa(\eta)$ is $T^p(A)$.
\end{lem}

\section{The Tate representation}

Let $k$ be a field, and let $\ell$ be distinct from
the characteristic of $k$.  If $A$ is an abelian variety over $k$, the
$\ell$-adic Tate module of $A$ is defined by
$$ T_\ell(A) = T_\ell(A \otimes_k \br{k}). $$
The \index{2TlA@$T_\ell(A)$}
Tate module gives a faithful
representation of the category of abelian varieties over $k$ (See, for
instance, \cite[Lem.~9.6]{milneabvar}).

\begin{prop}
Let $A$ and $A'$ be abelian varieties over $k$.
The natural map
$$ \Hom(A,A') \otimes \ZZ_\ell \rightarrow \Hom_{\ZZ_\ell}(T_\ell A, T_\ell
A') $$
is a monomorphism.
\end{prop}

The integral Tate module can be used to isolate the lattice of homomorphisms
using the following lemma.

\begin{lem}\label{lem:tatepullback}
Let $A$ and $A'$ be abelian varieties over $k$.  The diagram
$$
\xymatrix{
\Hom(A,A') \otimes \ZZ_\ell \ar[r] \ar[d]
& \Hom_{\ZZ_\ell}(T_\ell A, T_\ell A') \ar[d] 
\\
\Hom(A,A') \otimes \QQ_\ell \ar[r] 
& \Hom_{\QQ_\ell}(V_\ell A, V_\ell A')
}
$$
is a pullback square.
\end{lem}

\begin{proof}
Suppose that $f\in \Hom(A,A') \otimes \QQ_\ell$ 
has the property that
$f_*(T_\ell A) \subseteq T_\ell A'$.  There is some $i$ so that $g = \ell^i
f$
is contained in $\Hom(A,A') \otimes \ZZ_\ell$.  Then $g_*(T_\ell A)$ is
contained in $\ell^i T_\ell A'$.
Express $g$ in the form $g = g_1 + \ell^i g_2$ for $g_1 \in \Hom(A,A')$ and 
$g_2
\in \Hom(A,A') \otimes \ZZ_\ell$.
We deduce that the map
$$ (g_1)_*\co A[\ell^i] \rightarrow A'[\ell^i] $$
is null.
We therefore have the following factorization.
$$
\xymatrix{
A \ar[r]^{g_1} \ar[d]_{[\ell^i]} 
& A'
\\
A \ar@{.>}[ur]_{f_1}
}
$$
We deduce that $f = f_1 + g_2$ lies in $\Hom(A,A') \otimes \ZZ_\ell$.
\end{proof}

\section{Homomorphisms of abelian schemes}

In this section we explain how the results on abelian varieties given in
the last section extend to abelian schemes.  Fix a connected
base scheme $S$, let $x$ be a point of $S$, 
and let $A$ and $A'$ be abelian schemes over $S$.  As observed in
\cite{fisher}, the following
proposition follows immediately from \cite[Cor 6.2]{mumford2}.

\begin{prop}
The restriction map
$$ \Hom_S(A,A') \rightarrow \Hom_x(A_x, A'_x) $$
is a monomorphism.
\end{prop}

\begin{lem}\label{lem:schemepullback}
Suppose that $S$ is a scheme over $\ZZ[1/\ell]$.
Then the diagram
$$
\xymatrix{
\Hom_S(A,A')_{(\ell)} \ar[r] \ar[d] 
& \Hom_x(A_x, A'_x)_{(\ell)}  \ar[d]
\\
\Hom^0_S(A,A') \ar[r] 
& \Hom^0_x(A_x, A'_x)
}
$$
is a pullback.
\end{lem}

\begin{proof}
Suppose that $f \in \Hom^0_S(A,A')$ has the property that $f_x$ actually
lies in $\Hom_x(A_x,A'_x)_{(\ell)}$.  Let $N$ coprime to $\ell$ and $i$ be
chosen so that $g = \ell^i N f$ lies in $\Hom_S(A,A')$.  Then the
homomorphism of
\'etale group schemes
$$ g_x\co A_x[\ell^i] \rightarrow A'_x[\ell^i] $$
over $x$ is null.  
Consider the map
$$ g\co A[\ell^i] \rightarrow A'[\ell^i] $$
of \'etale group schemes over $S$.  Now every connected component of
$A[\ell^i]$ contains a point of $A_x[\ell^i]$, and 
the zero-section $0_{A'}$ is a connected component of 
$A'[\ell^i]$.  Since $g$ maps every point of $A_x[\ell^i]$ to
the component $0_{A'}$, we deduce that $g$ must map all of $A[\ell^i]$ to
zero.
Therefore, we get a factorization
$$
\xymatrix{
A \ar[r]^g \ar[d]_{[\ell^i]} & A'
\\
A \ar@{.>}[ur]_{g_1} 
}
$$
and we see that $f = N^{-1}g_1$ lies in $\Hom_S(A,A')_{(\ell)}$.
\end{proof}

Lemma~\ref{lem:schemepullback} combines with Lemma~\ref{lem:tatepullback}
to give the following corollary.

\begin{cor}\label{cor:tatepullback}
Let $s$ be a geometric point of $S$.
The square
$$
\xymatrix{
\Hom_S(A,A')_\ell \ar[r] \ar[d] 
&\Hom_{\ZZ_\ell}(T_\ell (A_s), T_\ell (A'_s)) \ar[d] \\
\Hom_S(A,A')_\ell \otimes \QQ \ar[r] 
& \Hom_{\QQ_\ell}(V_\ell(A_s), V_\ell(A'_s))
}
$$
is a pullback.
\end{cor}

\chapter{Polarizations}
\label{chap:pol}

For us, polarizations serve two purposes.  Firstly, the automorphism group of an
abelian variety is often infinite.  However, the subgroup of automorphisms
which preserves a polarization is finite.
Therefore, the moduli of polarized abelian
varieties is more tractable.  Mumford \cite{mumford2}
made use of geometric invariant theory to prove representability of
moduli spaces of polarized abelian varieties, which serves as a
starting point for the representability of the moduli problem associated to 
PEL Shimura
varieties.  Secondly, in Chapter~\ref{chap:deformations} we shall see
that the structure of a polarization on an abelian variety $A$ 
can be used to control the
deformation theory of a
compatible splitting of the $p$-divisible group $A(p)$.  

In this chapter we outline some of the theory of polarizations
of an abelian variety.  We recall the definition of a polarization 
and the associated Rosati involution.  These give
rise to Weil pairings on the Tate module.  We then describe the 
classification of polarized abelian varieties over finite fields.
Our treatment is based
on those of \cite{hida}, \cite{harristaylor}, and \cite{kottwitz}.

\section{Polarizations}

Let $A$ be an abelian variety defined over a field $k$.  The functor
$$ \Pic_{A/k} \co k-\mathrm{Schemes} \rightarrow \text{Abelian groups}, $$
which \index{2PicAk@$\Pic_{A/k}$}
associates to a $k$-scheme
$$ f\co S \rightarrow \Spec (k) $$
the set
$$ \Pic_{A/k}(S) = \left\{ 
\begin{array}{l}
\text{Isomorphism classes of invertible sheaves over
$f^*A$ which} \\
\text{restrict to the trivial line bundle on the subscheme $0_{f^*A}$}
\end{array}
\right\}, $$
is representable by a group scheme $\Pic_{A/k}$ \cite[Sec~13]{mumford}.
The dual abelian variety is defined to be the identity component of this
scheme:
$$ A^\vee = \Pic^0_{A/k}. $$
If \index{2Pic0Ak@$\Pic^0_{A/k}$}
\index{2A@$A^\vee$}
\index{abelian variety!dual}
\index{dual!abelian variety}
$B$ is an abelian variety over $k$, then specifying a $k$-morphism
$$ \alpha\co A \rightarrow B^\vee $$
is the same as giving an isomorphism class of line bundles $L_\alpha$ over 
$A \times_k B$ which
restricts to the trivial line bundle on the subscheme $A \times_k 0_B$.
The morphism $\alpha$ satisfies $\alpha(0_A) = 0_{B^\vee}$ 
if and only if the line
bundle $L_\alpha$ is also trivial when restricted to the subscheme $0_A
\times_k B$.

\begin{prop}
The map 
$$ \Hom(A,B^\vee) \rightarrow 
\left\{ 
\begin{array}{l}
\text{Isomorphism classes of line bundles on} \\ 
\text{$A \times_k B$, trivial on $0_A \times_k B$ and $A \times_k 0_A$} 
\end{array}
\right\}  $$
given by sending a homomorphism $\alpha$ to the corresponding line bundle
$L_\alpha$ is an isomorphism.
\end{prop}

\begin{proof}
We only need to check that every morphism $\alpha\co A \rightarrow B^\vee$
such that $\alpha(0_A) = 0_{B^\vee}$ is a homomorphism.  That is, we need the
following diagram to commute.
$$
\xymatrix@C+1em{
A \times_k A \ar[r]^{\alpha \times \alpha} \ar[d]_{\mu_A} & 
B^\vee \times_k B^\vee \ar[d]^{\mu_{B^\vee}} 
\\
A \ar[r]_{\alpha} & 
B^\vee
}
$$
This is equivalent to verifying that the line bundles
$$ 
L_{\alpha \circ \mu_A} = (\mu_A \times 1)^* L_\alpha \qquad \text{and} \qquad
L_{\mu_{B^\vee} \circ \alpha \times \alpha} = p_1^*L_\alpha \otimes
p_2^*L_\alpha
$$
over $A \times_k A \times_k B$ are isomorphic.  Here
$$ p_i\co A \times_k A \times_k B \rightarrow A \times_k B $$
are the two projections.  However, the two line bundles are easily seen to
agree on the subschemes $0_A \times_k A \times_k B$, $A \times_k 0_A
\times_k B$, and $A \times_k A \times_k 0_B$.  The theorem of the cube
\cite{milneabvar} therefore implies that these line bundles must be
isomorphic over $A \times_k A \times_k B$.
\end{proof}

The identity map $1 \co A^\vee \rightarrow A^\vee$ gives rise to the
\emph{Poincar\'e bundle} $L_{1}$ 
\index{Poincar\'e bundle}
\index{2L1@$L_{1}$}
on $A^\vee \times_k A$.  The line 
bundle $L_1$, when viewed as a bundle over $A \times_k A^\vee$, corresponds
to a canonical isomorphism
$$ A \xrightarrow{\cong} (A^\vee)^\vee. $$ 

Given a homomorphism $\alpha \co A \rightarrow B$, there is a 
dual homomorphism
$$ \alpha^\vee\co B^\vee \rightarrow A^\vee. $$
\index{1alpha@$\alpha^\vee$}
\index{abelian variety!homomorphism!dual}
\begin{enumerate}
\item Pullback along $\alpha$ gives a natural transformation of functors
$$ \alpha^*\co \Pic_{B/k} \rightarrow \Pic_{A/k}, $$
hence a map on representing objects, which restricts to $\alpha^\vee$ on
connected components.

\item Use the canonical isomorphism $B \cong (B^\vee)^\vee$, and take the
corresponding line bundle $L_\alpha$ over $A \times_k B^\vee$.  The pullback
$\tau^*L_\alpha$ under the twist map
$$ \tau\co B^\vee \times_k A \xrightarrow{\cong} A \times_k B^\vee $$
corresponds to the homomorphism
$$ \alpha^\vee\co B^\vee \rightarrow A^\vee. $$
\end{enumerate}

A homomorphism $\lambda\co A \rightarrow A^\vee$ is \emph{symmetric} 
\index{abelian variety!homomorphism!symmetric}
if the
dual homomorphism 
$$ \lambda^\vee\co A \cong (A^\vee)^\vee \rightarrow A^\vee $$
equals $\lambda$.  By the second description of the dual
homomorphism given above, this condition is equivalent to requiring the line
bundle $L_\lambda$ over $A \times_k A$ to satisfy $L_\lambda \cong \tau^*
L_\lambda$, where $\tau\co A \times_k A \rightarrow A \times_k A$ is the
interchange map.

\begin{exam}\label{exam:lambda_L}
Let $L$ be a line bundle over $A$.  Then the line bundle $\mu_A^*L \otimes
p_1^*L^{-1} \otimes p_2^*L^{-1}$ is a degree $0$ 
line bundle over $A \times A$ which
is trivial on $0_A \times_k A$ and $A \times_k 0_A$.  Here $\mu_A\co A
\times_k A \rightarrow A$ is the multiplication map, and $p_i\co A \times_k A
\rightarrow A$ are the two projections.  The corresponding homomorphism
$$ \lambda_L \co A \rightarrow A^\vee $$
is symmetric.
\end{exam}
 
In fact, Theorem~2 in \cite[Sec 20]{mumford}, and the remarks which follow
it, give the following converse.

\begin{thm}\label{thm:symmetric}
If $k$ is algebraically closed, and $\alpha\co A \rightarrow A^\vee$ is a
symmetric homomorphism, then there exists a line bundle $L$ over $A$ so
that $\alpha = \lambda_L$.  The line bundle $L$ satisfies
$$ L^{\otimes 2} = \Delta^* L_\alpha $$
where $\Delta\co A \rightarrow A \times_k A$ is the diagonal.
\end{thm}

The following theorem \cite[Sec~16]{mumford} gives a criterion for
determining if a symmetric homomorphism $A \rightarrow A^\vee$ 
is an isogeny (has finite kernel).

\begin{thm}
Suppose that $L$ is a line bundle on $A$.  Then the following conditions
are equivalent.
\begin{enumerate}
\item The symmetric homomorphism $\lambda_L$ is an isogeny.
\item There is a unique $0 \le i \le \dim A$ 
such that $H_{zar}^i(A,L)$ is non-trivial, and $H_{zar}^j(A,L) = 0$
for $j \ne i$.
\end{enumerate}
\end{thm}

The integer $i = i(L)$ in the previous theorem 
is called the \emph{index} of $L$.  The index 
\index{line bundle!index}
is $0$ if $L$ is
ample \cite[III.5.3]{hartshorne}.

\begin{defn}
Let $A$ be an abelian variety defined over an algebraically closed field
$k$.  A \emph{polarization} 
\index{polarization}
of $A$ is a symmetric isogeny
$$ \lambda = \lambda_L\co A \rightarrow A^\vee $$
determined by a line bundle $L$ which is ample.
\end{defn}

\begin{rmk}\label{rmk:polexistence}
Abelian varieties are in particular projective, and hence admit ample line
bundles.  Therefore, every abelian variety admits a polarization.
\end{rmk}

\begin{rmk}
The theory of symmetric homomorphisms for abelian varieties is analogous to
the theory of symmetric bilinear forms for real vector spaces.  Let $B(-,-)$ be
a symmetric bilinear form on a real vector space $V$.  The symmetric
homomorphism $\lambda_L$ associated to a line bundle $L$ is analogous to
the symmetric bilinear form
$$ B(x,y) = Q(x + y) - Q(x) - Q(y) $$
associated to a quadratic form $Q$.  Theorem~\ref{thm:symmetric} is an
analog to the fact that every form $B$ arises this way, with 
$$ Q(x) = \frac{1}{2}B(x,x). $$
Representing a form $B$ by a hermitian matrix $A$, the index of a
symmetric isogeny is analogous to the number of negative eigenvalues of
$A$.  A polarization is therefore analogous to a positive definite
quadratic form.
\end{rmk}

If
$\lambda = \lambda_L$ is a polarization, then $n \lambda =
\lambda_{L^{\otimes n}}$ is a polarization for $n$ a positive integer.  Two
polarizations $\lambda$ and $\lambda'$ are \emph{equivalent} 
\index{polarization!equivalence of}
if $n \lambda = m \lambda'$ for some positive integers $n$ and $m$.  A
\emph{weak polarization} 
\index{polarization!weak}
is an equivalence class of polarizations.

We defined polarizations only for abelian varieties defined 
over an algebraically
closed field.
If $A$ is an abelian scheme over a scheme $S$, then a polarization of $A$
is an isogeny
$$ \lambda\co A \rightarrow A^\vee $$
which restricts to a polarization on geometric fibers.

\section{The Rosati involution}\label{sec:Rosati}

Let $A$ be an abelian variety over an algebraically closed field $k$, and
let $\lambda = \lambda_L$ be the polarization of $A$ associated to an ample
line bundle $L$.

The polarization $\lambda$ induces an involution $\dag$ 
\index{0dag@$\dag$}
on the semi-simple
$\QQ$-algebra $\End^0(A)$; for $f \in \End^0(A)$, the quasi-endomorphism
$f^\dag$ is given by the
composite
$$ f^\dag \co A \xrightarrow{\lambda} A^\vee \xrightarrow{f^\vee}
A^\vee \xrightarrow{\lambda^{-1}} A. $$ 
This involution is the {\em Rosati involution\/}
\index{Rosati involution}
\index{involution!Rosati}
associated to the polarization $\lambda$.  If the polarization $\lambda$ is
replaced by $n\lambda$ for some $n > 0$, the associated Rosati
involution remains the same.  Therefore, the Rosati involution depends only
on the weak polarization determined by $\lambda$.

The positivity of the line bundle ${L}$ implies that the Rosati
involution is \emph{positive} in the following sense.

\begin{thm}[{\cite[Sec 21]{mumford}}]\label{thm:Rosatipos}
Let $E = \End^0(A)$, and suppose $f \in E$ is nonzero.  Then
\[
\Tr_{E/\mb Q} (f f^\dag) > 0.
\]
Here $\Tr_{E/\QQ}$ is the trace map.
\end{thm}

\section{The Weil pairing}

Suppose that $k$ is algebraically closed, and that 
$\lambda\co A \rightarrow A^\vee$ is a polarization of $A$.
Then $\lambda$ induces a map
$$ \lambda_*\co A[\ell^i] \rightarrow A^\vee[\ell^i] \cong \Hom(A[\ell^i],
\mu_{\ell^i}) $$
whose adjoint gives the \emph{$\lambda$-Weil pairing}
\index{Weil pairing}
$$ \lambda\bra{-,-}_{\ell^i}\co A[\ell^i] \times A[\ell^i] \rightarrow
\mu_{\ell^i}. $$
Taking \index{1lambda@$\lambda\bra{-,-}$}
the inverse limit over $i$, we recover a natural bilinear
pairing 
\[
\lambda\bra{-,-}\co T_\ell(A) \times T_\ell(A) \to \mb Z_\ell(1)
\]
called the {\em $\lambda$-Weil pairing\/}.
The induced pairing 
$$ \lambda\bra{-,-}\co V_\ell(A) \otimes_{\QQ_\ell} V_\ell(A) \rightarrow
\QQ_\ell(1) $$
is non-degenerate, because $\lambda$ is an isogeny.  It can be shown that
the pairing $\lambda\bra{-,-}$ is alternating \cite[Sec~20]{mumford}.

\section{Polarizations of $B$-linear abelian varieties}

Let $A$ be an abelian variety over a field $k$.  Let
$B$ be a simple $\QQ$-algebra.  Recall from Section~\ref{sec:HTFpbar} that
a $B$-linear abelian variety $(A,i)$ is an abelian variety $A/k$ together
with an embedding of rings $i\co B \hookrightarrow \End^0(A)$.  Suppose that
$*$ is an involution on $B$. 

\begin{defn}
A polarization $\lambda$ on $(A,i)$ is
\emph{compatible} 
\index{polarization!compatible}
if the $\lambda$-Rosati involution restricts to
the involution $*$ on $B$.
\end{defn}

Theorem~\ref{thm:Rosatipos} implies that for a
compatible polarization to exist, the involution $*$ must be positive.
Conversely, we have the following lemma.

\begin{lem}[{\cite[Lemma~9.2]{kottwitz}}]
If $*$ is a positive involution on $B$, and $(A,i)$ is a $B$-linear abelian
variety, then there exists a compatible polarization $\lambda$ on $A$.
\end{lem}

The $B$-linear structure on $A$ gives 
the Tate module $V_\ell(A)$ the structure of a $B$-module.  The
compatibility condition on the polarization implies that the $\lambda$-Weil
pairing is \emph{$*$-hermitian}:  
\index{hermitian@$*$-hermitian}
\index{bilinear form!$*$-hermitian}
for all $x,y \in V_\ell(A)$ and $b \in
B$, we have
$$ \lambda\bra{bx,y} = \lambda\bra{x,b^*y}. $$

\section{Induced polarizations}

Let $\lambda = \lambda_L$ be a polarization of an abelian variety $A$, with
associated ample line bundle $L$.  Let 
$$ \alpha\co A' \rightarrow A $$
be an isogeny of abelian varieties.  Then $\alpha$ induces an isogeny
$\alpha^* \lambda$ by the composite
\index{1alpha@$\alpha^* \lambda$}
\index{polarization!induced}
$$ \alpha^* \lambda\co A' \xrightarrow{\alpha} A \xrightarrow{\lambda} A^\vee
\xrightarrow{\alpha^\vee}  (A')^\vee. $$
It is easily verified that $\alpha^* \lambda$ is the symmetric homomorphism
associated to the line bundle $\alpha^*L$.  Since this line bundle is
ample, we see that $\alpha^*\lambda$ is a polarization of $A'$.

Suppose that $(A,\lambda)$ and $(A',\lambda')$ are polarized abelian
varieties.  An isogeny
$$ \alpha\co A' \rightarrow A $$
is an \emph{isometry} 
\index{abelian variety!isogeny!isometry}
if $\alpha^*\lambda = \lambda'$.  If
$\alpha^*\lambda$ is only equivalent to $\lambda'$, we say that $\alpha$ is
a \emph{similitude}.  
\index{similitude}
\index{abelian variety!isogeny!similitude}
The following lemma follows immediately from the
definition of the Rosati involution.

\begin{lem}\label{lem:isometrysimilitude}
Let $A$ be an abelian variety with polarization $\lambda$ and let $\dag$ be
the Rosati involution on $\End^0(A)$.
Let $\alpha\co A \rightarrow A$ be an isogeny.
Then:
\begin{align*}
\text{$\alpha$ is an isometry} & \Leftrightarrow \alpha^\dag\alpha = 1, \\
\text{$\alpha$ is a similitude} & \Leftrightarrow \alpha^\dag\alpha \in
\QQ^\times.
\end{align*}
\end{lem}

An isometry $\alpha\co A \rightarrow A$ of a polarized abelian variety $(A,\lambda)$
induces an isometry 
$$ \alpha_*\co (V_\ell(A), \lambda\bra{-,-}) \rightarrow (V_\ell(A),
\lambda\bra{-,-}) $$
of the $\lambda$-Weil pairing 
on each $\ell$-adic Tate module.  If $\alpha$ is merely a similitude, then
$\alpha_*$ is a similitude of the $\lambda$-Weil pairing: we have
$$ \lambda\bra{\alpha_*(x),\alpha_*(y)} = \nu(\alpha) \lambda \bra{x,y} $$
for some $\nu(\alpha) \in \QQ^\times$, and all $x,y \in V_\ell(A)$.
\index{1nu@$\nu$}

\section{Classification of weak polarizations}\label{sec:polclass}

Assume that $k$ has finite characteristic $p$.
We briefly recall the classification up to isogeny 
of compatible weak polarizations on a
$B$-linear abelian variety $(A,i)$, following \cite[V.3]{harristaylor}.

Consider the following:
\begin{align*}
B & = \text{simple $\QQ$-algebra}, \\
* & = \text{positive involution on $B$}, \\
(A,i) & = \text{$B$-linear abelian variety 
over an algebraically closed field $k$}, \\
\lambda & = \text{compatible polarization of $(A,i)$}.
\end{align*}

Define an algebraic group $H_{(A,i,\lambda)}/\QQ$ 
\index{2HAil@$H_{(A,i,\lambda)}$}
by
$$ H_{(A,i,\lambda)}(R) = \{ h \in (\End^0_B(A) \otimes R)^\times \: : \: 
h^\dag h \in R^\times \}. $$
Here, $\dag$ 
\index{0dag@$\dag$}
is the $\lambda$-Rosati involution.  Note that by
Lemma~\ref{lem:isometrysimilitude}, $H_{(A,i,\lambda)}(\QQ)$ is the group of
quasi-endomorphisms of $A$ which are similitudes of the polarization
$\lambda$.

We shall say that a compatible polarization $\lambda'$ is lies in the same
\emph{similitude class} 
\index{similitude class}
\index{polarization!similitude class}
as
$\lambda$ if there exists a quasi-isometry $\alpha\co A \rightarrow A$ so
that $\alpha^*\lambda$ is equivalent to $\lambda$.

\begin{lem}[{\cite[Lem V.3.1]{harristaylor}}]\label{lem:polclassglobal}
There is a bijective correspondence between similitude classes of compatible weak
polarizations on $A$ and elements in the Galois cohomology kernel:
$$ \ker(H^1(\QQ,H_{(A,i,\lambda)}) \rightarrow H^1(\RR,H_{(A,i,\lambda)})) $$
\end{lem}

\begin{rmk}
Non-abelian Galois $H^1$ is a pointed set. 
The kernel is the collection of elements which
restrict to the distinguished 
element at the completion $\infty$ of $\QQ$.  This
condition reflects the fact that the index of the symmetric homomorphism
$\lambda'$ must agree with the index of $\lambda$ in order to be a
polarization. 
\end{rmk}

Assume now that $\ell$ is a prime different from $p$.
Let $GU_{V_\ell(A)}$ 
\index{2GUVlA@$GU_{V_\ell(A)}$}
be the algebraic group (over $\QQ_\ell$) of $B$-linear 
similitudes of
$(V_\ell(A), \lambda\bra{-,-})$:
$$ GU_{V_\ell(A)}(R) = \left\{ g \in (\End_B(V_\ell(A)) \otimes_{\QQ_\ell}
R)^\times \, : \, 
\begin{array}{l}
\text{there exists $\nu(g) \in R^\times$ such} \\ 
\text{that} \: \lambda\bra{g(x),g(y)} = \nu(g)\lambda\bra{x,y} 
\end{array}
\right\}. 
$$
Similitude classes of non-degenerate $*$-hermitian alternating forms
on $V_\ell(A)$ are classified by elements of $H^1(\QQ_\ell, GU_{V_\ell(A)})$.

Since similitudes of $(A,\lambda)$ induce similitudes of
$(V_\ell(A),\lambda\bra{-,-})$, there is a homomorphism of algebraic groups
$$ H_{(A,i,\lambda)} \times_{\Spec(\QQ)} \Spec(\QQ_\ell) \rightarrow GU_{V_\ell(A)} $$
which induces a map on Galois cohomology
$$ H^1(\QQ,H_{(A,i,\lambda)}) \rightarrow H^1(\QQ_\ell, GU_{V_\ell}). $$

\begin{lem}[{\cite[Lem V.3.1]{harristaylor}}]\label{lem:polclasslocal}
Let $\lambda'$ be a compatible polarization of $(A,i)$.  The image of the
cohomology class $[\lambda'] \in H^1(\QQ, H_{(A,i,\lambda)})$ in
$H^1(\QQ_\ell, GU_{V_\ell}(A))$ quantifies the difference between the
similitude classes of alternating forms represented by $\lambda'\bra{-,-}$
and $\lambda\bra{-,-}$ on $V_\ell(A)$.
\end{lem}

In practice, the cohomology class $[\lambda']$ is often computed by means of
its local invariants in this manner.

\chapter{Forms and involutions}
\label{chap:galois}

In Section~\ref{sec:polclass}, we saw that polarized
$B$-linear abelian varieties up to isogeny were classified by the Galois
cohomology of the similitude group of an alternating form.
In this chapter we will explicitly outline the Galois cohomology
computations relevant for this book.  Most of the results in this chapter
may be found elsewhere (see, for instance, \cite{scharlau}), 
but for convenience we enumerate these
results in one place.

\section{Hermitian forms}

We shall use the following notation.
\begin{align*}
F = & \: \text{quadratic imaginary extension of $\QQ$}. \\
c = & \: \text{conjugation on $F$}. \\
B = & \: \text{central simple algebra over $F$ of dimension $n^2$}. \\
* = & \: \text{positive involution on $B$ of the second kind}. \\
V = & \: \text{free left $B$-module of rank $1$}.  \\
C = & \: \End_B(V) \qquad (\text{noncanonically isomorphic to $B^{op}$}).
\end{align*}
\index{2F@$F$}
\index{2C@$c$}
\index{2B@$B$}
\index{0ast@$*$}
\index{2V@$V$}
\index{2C@$C$}

(An involution of $B$ is said to be \emph{of the second kind} 
\index{involution!of the second kind}
if it
restricts to the conjugation $c$ on $F$.)

\begin{rmk}\label{rmk:involutioncondition}
For $B$ to admit an involution $*$ of the second kind it is necessary
and sufficient \cite{scharlau} that
$$ \inv_x B = 0 $$
for all finite primes $x$ which are not split in $F$ and
$$ \inv_y B + \inv_{y^c} B = 0 $$
for all finite primes $x$ which split as $yy^c$ in $F$.
\end{rmk}

In Section~\ref{sec:polclass} we explained how classification of
compatible weak polarizations on a $B$-linear abelian variety was
equivalent to the classification of certain alternating forms.
For this reason, we are interested in non-degenerate alternating forms
$$ \bra{-,-}\co V \otimes_\QQ V \rightarrow \QQ $$
which are $*$-hermitian, meaning that
$$ \bra{bv,w} = \bra{v, b^*w}. $$
for all $b \in B$ and $v,w \in V$.
Two $*$-hermitian alternating forms $\bra{-,-}$ and $\bra{-,-}'$ on 
$V$ are
said to be \emph{similar} 
\index{bilinear form!similar}
if there exists an endomorphism $\alpha \in C$ and
a unit $\nu(\alpha) \in \QQ^\times$ so that
$$ \bra{\alpha v, \alpha w}' = \nu(\alpha)\bra{ v, w} $$
for all $v,w \in V$.  If $\alpha$ may be chosen so that $\nu(\alpha) = 1$,
then we shall say that $\bra{-,-}$ and $\bra{-,-}'$ are \emph{isometric}.
\index{bilinear form!isometric}

Such an alternating form induces an involution $\iota$ on $C$ via the
\index{1iota@$\iota$}
formula
\[
\bra{c v, w} = \bra{v, c^\iota w}.
\]
Two similar forms induce the same involution $\iota$.

We shall say that a bilinear form
$$ (-,-)\co V \otimes_\QQ V \rightarrow F $$
is \emph{$*$-symmetric} 
\index{bilinear form!$*$-symmetric}
if it is $F$-linear in the first variable, and
$$ (v,w) = (w,v)^c. $$

For definiteness, we shall let
\begin{align*}
V & = B, \\
F & = \QQ(\delta), \\
\end{align*}
where \index{1delta@$\delta$}
$B$ is regarded as a left $B$-module, and $\delta^2 = d$ 
for a negative square-free integer $d$.  There is a natural identification $C =
B^{op}$.

\begin{lem}\label{lem:alternatinghermitian}
There is an one to one
correspondence:
\begin{gather*}
\{ \text{non-degenerate $*$-hermitian alternating forms $\bra{-,-}$ on $V$} \}
\\
\updownarrow
\\
\{ \text{non-degenerate $*$-hermitian $*$-symmetric forms $(-,-)$ on $V$} \}
\\
\end{gather*}
\end{lem}

\begin{proof}
Given an alternating form $\bra{-,-}$, the corresponding $*$-symmetric
form is given by 
$$ (v,w) = \bra{\delta v, w } + \delta \bra{v,w}. $$
Given a $*$-symmetric form $(-,-)$, the corresponding alternating form is
given by
$$ \bra{v,w} = \frac{1}{2d}\Tr_{F/\QQ} \delta (v,w). $$
\end{proof}

The non-degeneracy of the involution $*$ 
allows one to deduce the following lemma, which makes everything more
explicit. 

\begin{lem}\label{lem:betadef}
Let $\bra{-,-}$ be a non-degenerate $*$-hermitian alternating form.
There exists a unique element $\beta \in B = V$ satisfying $\beta^* = -\beta$ 
\index{1beta@$\beta$}
which encodes $\bra{-,-}$:
$$
\bra{x,y} = \Tr_{F/\QQ} \Tr_{B/F} (x \beta y^*).
$$
Regarding $\beta$ as an element $\gamma$ of $C = B^{op}$, we have that
the involution on $C$ is given by 
$$ c^\iota = \gamma^{-1} c^* \gamma. $$
\index{1gamma@$\gamma$}
\end{lem}

\begin{cor}\label{cor:betadef}
There is a $\xi \in B$ 
\index{1xi@$\xi$}
satisfying $\xi^* = \xi$ so that the 
$*$-symmetric pairing $(-,-)$ 
associated to the alternating form $\bra{-,-}$ of
Lemma~\ref{lem:betadef} is given by
$$
(x,y) = \Tr_{B/F} (x \xi y^*).
$$
\end{cor}

\begin{proof}
The element $\xi$ is computed to be $2\delta \beta$.
\end{proof}

For a non-degenerate $*$-hermitian $*$-symmetric form $(-,-)$, we define
the \emph{discriminant}
\index{discriminant}
\index{bilinear form!discriminant}
$$ \mathrm{disc} = N_{B/F}(\xi) \in \QQ^\times/N_{F/\QQ}(F^\times) $$
for $\xi$ the element of Corollary~\ref{cor:betadef}.
The discriminant is an invariant of the isometry
class of $(-,-)$.

Two involutions $\iota$ and $\iota'$ are \emph{equivalent} 
\index{involution!equivalence of}
if there exists
a $c \in C^\times$ so that
$$ c x^\iota c^{-1} = (cxc^{-1})^{\iota'}. $$
We finish this section by explaining how the classification of involutions
up to equivalence is related to the classification of hermitian forms.

\begin{prop}\label{prop:involutionssimilitudeclasses}
The association of an involution on $C$ with a non-degenerate $*$-hermitian 
alternating form on $V$ establishes
a bijective correspondence
\begin{gather*}
\{ \text{Similitude classes of non-degenerate $*$-hermitian 
alternating forms on $V$} \}
\\
\updownarrow
\\
\{ \text{Equivalence classes of involutions of the second kind on $C$} \} 
\end{gather*}
\end{prop}

We first need the following lemma.

\begin{lem}\label{lem:involutionsconjugation}
Let $\iota$ be an involution of the second kind on $C$.  Then there exists
an element $\gamma \in C^\times$ 
\index{1gamma@$\gamma$}
satisfying $\gamma^* = -\gamma$ which
gives $\iota$ by the formula
$$ x^\iota = \gamma^{-1} x^* \gamma. $$
The element $\gamma$ is unique up to multiplication by 
an element in $\QQ^\times$.
\end{lem}

\begin{proof}
The Noether-Skolem theorem implies
that there exists an element $\alpha \in C^\times$, unique up to
$F^\times$-multiple, such that
$$ (x^*)^\iota = \alpha^{-1} x \alpha $$
for all $x \in C$.  
Therefore, we have $x^\iota = \alpha^{-1} x^* \alpha$.
Since $\iota$ is an involution, we determine that 
$\alpha^{-1}\alpha^* = c \in F^\times$. By Hilbert's Theorem 90 there
exists an $a \in F^\times$ so that $\gamma = a\alpha$ satisfies $\gamma^* =
-\gamma$.  Such an element $\gamma$ is determined up to a multiple in
$\QQ^\times$.
\end{proof}

\begin{proof}[Proof of Proposition~\ref{prop:involutionssimilitudeclasses}]
By Lemma~\ref{lem:betadef}, given a non-degenerate $*$-hermitian alternating
form $\bra{-,-}$, there exists an element $\beta$ such that
$$ \bra{x,y} = \bra{x,y}_\beta = \Tr_{F/\QQ} \Tr_{B/F} (x \beta y^*). $$
Letting $\gamma$ be the element $\beta$ regarded as an element of $C =
B^{op}$, the associated involution 
$\iota_\gamma$ 
\index{1iotag@$\iota_\gamma$}
is given by $x^{\iota_\gamma} = \gamma^{-1} x^*
\gamma$.  Changing $\bra{-,-}_\beta$ by a $\QQ^\times$-multiple changes
$\beta$ and $\gamma$ by a $\QQ^\times$-multiple, leaving the involution
$\iota_\gamma$ unchanged.  In the other direction, Lemma~\ref{lem:involutionsconjugation}
associates to an involution $\iota$ an element $\gamma \in C^\times$, 
unique up to $\QQ^\times$-multiple, so that
$$
x^\iota = x^{\iota_{\gamma}} = \gamma^{-1} x^* \gamma
$$
and $\gamma^* = -\gamma$.
Letting
$\beta$ be the element $\gamma$ regarded as an element of $B$, we can then
associate the similitude class of form $\bra{-,-}_\beta$. 
Forms $\bra{-,-}_{\beta}$ 
\index{0brab@$\bra{-,-}_{\beta}$}
and $\bra{-,-}_{\beta'}$ lie in the same
similitude class if and only if the associated involutions $\iota_\gamma$
and $\iota_{\gamma'}$ are equivalent.
\end{proof}

\section{Unitary and similitude groups}

We now fix the following: 
\begin{align*}
\bra{-,-} & = \text{non-degenerate $*$-hermitian alternating pairing on
$V$}, \\ 
(-,-) & = \text{corresponding $*$-symmetric pairing}, \\
(-)^\iota & = \text{involution on $C$ defined by $\bra{cv,w} = \bra{v,
c^\iota w}$}.
\end{align*}
\index{0bra@$\bra{-,-}$}
\index{0bro@$(-,-)$}
\index{1iota@$\iota$}

Associated to this pairing are some group-schemes defined over $\QQ$: 
the \emph{unitary group} $U$ and
the \emph{similitude group} $GU$.  
\index{2U@$U$}
\index{2GU@$GU$}
\index{unitary group}
\index{similitude group}
For a $\QQ$-algebra $R$, the $R$-points
of these groups are given by
\begin{align*}
U(R) & = \{ g \in (C \otimes R)^\times \: : \: g^\iota g = 1 \}, \\
GU(R) & = \{ g \in (C \otimes R)^\times \: : \: g^\iota g \in R^\times \}. 
\end{align*}
Let $V_R$ denote the $C \otimes R$-module $V \otimes R$.  
The $R$-points of $U$ (respectively, $GU$) consist of the isometries
(respectively similitudes) of both pairings 
$(-,-)$ and $\bra{-,-}$ on $V_R$.  The group $GU$ has a \emph{similitude
norm}
\index{similitude norm}
\index{1nu@$\nu$}
$$ \nu : GU \rightarrow \GG_m $$
which, on $R$-points, takes an element $g$ to the quantity $g^\iota g \in
R^\times$.  For $v,w \in V_R$, we have
$$ \bra{gv, gw} = \nu(g)\bra{v,w}. $$

Let $GL_C$ 
\index{2GLC@$GL_C$}
denote the form of $GL_n$ over $F$ whose $R$ points are given by
$$ GL_C(R) = (C \otimes_F R)^\times. $$

\begin{lem}\label{lem:splitF}
There are natural isomorphisms of group schemes over $F$:
\begin{align*}
\Spec(F) \times_{\Spec(\QQ)} U & \cong GL_{C}, \\
\Spec(F) \times_{\Spec(\QQ)} GU & \cong GL_{C} \times \GG_m.
\end{align*}
\end{lem}

\begin{proof}
It suffices to provide a natural isomorphism on $R$-points for 
$F$-algebras $R$.  Since $R$ is an $F$-algebra, there is a decomposition
\begin{align*}
C \otimes_\QQ R 
& \cong C \otimes_F (F \otimes_\QQ R) \\
& \cong (C \otimes_{F} R) \times (C \otimes_{F,c} R).
\end{align*}
The induced involution $\iota$ 
on $(C \otimes_{F} R) \times (C \otimes_{F,c} R)$ is
given by
$$ (c_1 \otimes r_1, c_2 \otimes r_2)^\iota = (c_2^\iota \otimes r_2,
c_1^\iota \otimes r_1). $$
Therefore we have a natural isomorphism
\begin{align*}
U(R)
& = \{ g \in (C \otimes_{\QQ} R)^\times \: : \: g^\iota g = 1 \} \\
& \cong \{ (g_1,g_2) \in (C \otimes_{F} R)^\times \times 
(C \otimes_{F,c} R)^\times \: : \: g_2^\iota g_1 = 1 \} \\
& \cong \{ (g_1, g_1^{-\iota}) \: : \: g_1 \in (C \otimes_{F} R)^\times \} \\
& \cong GL_{C}(R). 
\end{align*}
Similarly, there are natural isomorphisms
\begin{align*}
GU(R)
& = \{ g \in (C \otimes_{\QQ} R)^\times \: : \: g^\iota g = \nu \in
R^\times \} \\
& \cong \{ (g_1, g_1^{-\iota} \nu) \: 
: \: g_1 \in (C \otimes_{F} R)^\times, \: v \in R^\times \} \\
& \cong GL_{C}(R) \times \GG_m(R). 
\end{align*}
\end{proof}

\section{Classification of forms}\label{sec:formclassification}

Let $K$ be either $\QQ$ or $\QQ_x$ for some place $x$ of $\QQ$.
The classification of our various structures are parameterized by
certain Galois cohomology groups:
\begin{align*}
H^1(K, U) & = \left\{ 
\begin{array}{l}
\text{Isometry classes of non-degenerate} \\
\text{$*$-hermitian alternating forms on $V_K$}  \\
\end{array}
\right\},
\\
H^1(K, GU) & = 
\left\{ 
\begin{array}{l}
\text{Similitude classes of non-degenerate} \\ 
\text{$*$-hermitian alternating forms on $V_K$} 
\end{array}
\right\}
\\
& =
\left\{ 
\begin{array}{l}
\text{Equivalence classes of involutions} \\ 
\text{of the second kind on $C \otimes K$} 
\end{array}
\right\}.
\end{align*}

\subsection*{Local case where $x$ is finite and split}

Let $x$ split as $yy^c$.  Then $\mb Q_x$ is an $F$-algebra, and hence
by Lemma~\ref{lem:splitF}, there are
%the isomorphism $C_x \cong C_y \times
%C_{y^c}$ yields 
isomorphisms of group schemes
\begin{align*}
\Spec(\QQ_x) \times_{\Spec(\QQ)} U & \cong GL_{C_y} \\
\Spec(\QQ_x) \times_{\Spec(\QQ)} GU & \cong GL_{C_y} \times \GG_{m}
\end{align*}
If $C$ is a central simple algebra over $\QQ_x$, it is a well-known
generalization of Hilbert's Theorem~90 
that $H^1(\QQ_x, GL_C) = 0$ (see, for instance, the proof of
\cite[Prop.~26.6]{milnealggp}).  
We therefore have:

\begin{lem}\label{lem:localsplitH^1GU}
If $x$ is split in $F$, then we have $H^1(\QQ_x, U) = H^1(\QQ_x, GU) = 0$.
\end{lem}

\subsection*{Local case where $x$ is infinite}

Because $F$ is imaginary, $C$ is necessarily split at $x = \infty$.
The computation of $H^1(\RR, U)$ and $H^1(\RR, GU)$ is then given by the 
classification of hermitian forms
$(-,-)'$ on $W \cong \CC^n$.  These forms are classified up to isometry 
by their signature, 
\index{signature}
and up to similitude by the absolute value of their
signature.  We therefore have the following lemma.

\begin{lem}\label{lem:localinfiniteH^1GU}
The signature induces isomorphisms
\begin{align*}
H^1(\RR, U) & \cong \{ (p,q) \: : \: p+q = n \}, \\
H^1(\RR, GU) & \cong \{ \{ p,q \} \: : \: p+q = n \}. 
\end{align*}
\end{lem}

\subsection*{Local case where $x$ is finite and not split}

Analogously to the global case, relative to the involution $*$ 
we may associate to a form $(-,-)$ on $V_x$
its discriminant $\mathrm{disc} \in \QQ_x^\times/N(F_x^\times) \cong \ZZ/2$.
\index{2Disc@$\mathrm{disc}$}
\index{discriminant}
\index{bilinear form!discriminant}
A proof of the following lemma may be found in \cite{scharlau}.

\begin{lem}
If $x$ does not split in $F$, 
the discriminant gives an isomorphism
$$ H^1(\QQ_x, U) \cong \ZZ/2. $$
\end{lem}

\begin{cor}\label{cor:localnonsplitH^1GU}
If $x$ does not split in $F$,
there are isomorphisms
$$
H^1(\QQ_x, GU) = 
\begin{cases}
\ZZ/2, & \text{$n$ even}, \\
0, & \text{$n$ odd}.
\end{cases}
$$
\end{cor}

\begin{proof}
Let $a \in \QQ^\times$ represent a generator of $\QQ^\times/N(F_x^\times)$.
If $n$ is even, the discriminant of $a(-,-)$ is equal to that of $(-,-)$.
If $n$ is odd, than the discriminant of $a(-,-)$ is not equal to that of
$(-,-)$.
\end{proof}

\subsection*{Global case}

For places $x$ of $\QQ$, we define maps
\begin{align*}
\xi_x\co & H^1(\QQ_x, U) \rightarrow \ZZ/2, \\
\xi'_x\co & H^1(\QQ_x, GU) \rightarrow \ZZ/2, \qquad \text{$n$ even}.
\end{align*}
If \index{1xix@$\xi_x$} \index{1xi'x@$\xi'_x$}
$x$ is finite and split, $H^1(\QQ_x, U) = H^1(\QQ_x, GU) = 0$.
If $x$ is finite and 
not split in $F$, 
let $\xi_x$, $\xi_x'$ be the unique isomorphisms 
\begin{align*}
\xi_x\co & H^1(\QQ_x, U) \xrightarrow[\cong]{\mathrm{disc}}
\QQ_x^\times/N(F_x^\times) \cong \ZZ/2, \\
\xi'_x\co & H^1(\QQ_x, GU) \xrightarrow[\cong]{\mathrm{disc}} 
\QQ_x^\times/N(F_x^\times)
\cong \ZZ/2, \qquad \text{$n$ even}.
\end{align*}
If $n$ is odd we have $H^1(\QQ_x, GU) = 0$.
At the infinite place, we define:
\begin{align*}
\xi_\infty((p,q)) & \equiv q \pmod 2, \\
\xi'_\infty(\{p,q\}) & \equiv p \equiv q \pmod 2, \qquad \text{$n$ even}.
\end{align*}

The following result may be translated from the classification of global 
hermitian forms given in \cite{scharlau} (see also \cite{clozel}).

\begin{thm}
There is a short exact sequence 
$$
0 \rightarrow H^1(\QQ, U) \rightarrow \bigoplus_x H^1(\QQ_x, U)
\xrightarrow{\sum \xi_x}
\ZZ/2 \rightarrow 0.
$$
\end{thm}

\begin{proof}
Scharlau \cite{scharlau} shows that non-degenerate 
$*$-symmetric $*$-hermitian forms on $V$ are classified by their
discriminant $\mathrm{disc} \in \QQ^\times/N(F^\times)$ and 
signature $(p,q)$.
Global properties of the Hilbert symbol give rise to a short exact
sequence
\begin{equation}\label{eq:normSES}
0 \rightarrow \QQ^\times/N(F^\times) \rightarrow \bigoplus_x
\QQ_x^\times/N(F_x^\times) \xrightarrow{\sum_x} \ZZ/2 \rightarrow 0.
\end{equation}
The result follows from the local computations, together with the obvious
relation in $\RR^\times/N(\CC^\times) \cong \ZZ/2$:
$$ \mathrm{disc}_\infty \equiv q \pmod 2. $$
\end{proof}

\begin{cor}\label{cor:globalH^1GU}
If $n$ is even, then there is a short exact sequence 
$$
0 \rightarrow H^1(\QQ, GU) \rightarrow \bigoplus_x H^1(\QQ_x, GU) 
\xrightarrow{\sum \xi'_x}
\ZZ/2 \rightarrow 0.
$$
If $n$ is odd then the absolute value of the signature gives an isomorphism
$$
H^1(\QQ,GU) \xrightarrow{\cong} H^1(\RR, GU).
$$
\end{cor}

\begin{proof}
If a form $(-,-)$ has discriminant
$\mathrm{disc}$, the form $r(-,-)$ has discriminant $r^n \cdot
\mathrm{disc}$,
for $r \in \QQ^\times$.
Since $(\QQ^\times)^2$ is contained in $N(F^\times)$, if $n$ is even the
discriminant is an invariant of the similitude class of $(-,-)$.  However,
if $n$ is
odd, then the short exact sequence (\ref{eq:normSES}) shows that there
exists a global multiple $r$ which kills the local invariants
$\mathrm{disc}_x$ for finite places $x$, with the effect of possibly
changing the sign of the signature.
\end{proof}

\chapter{Shimura varieties of type $U(1,n-1)$}
\label{chap:shimura}

In this chapter we describe the Shimura stacks we wish to study.  The
integral version described here is essentially due to Kottwitz, and the
exposition closely follows \cite{kottwitz}, \cite{harristaylor}, and
\cite{hida}.

\section{Motivation}

In Section~\ref{sec:HTFpbar}, we saw that the simplest examples of 
abelian varieties $A$ 
over $\br{\FF}_p$, whose 
$p$-divisible group $A(p)$ contains a $1$-dimensional
summand of slope $1/n$, had dimension $n$ and 
complex multiplication by $F$, 
a quadratic imaginary
extension of $\QQ$ in which $p$ splits.  
For such abelian varieties, the $1$-dimensional
summand of $A(p)$ is given by $A(u)$, where $u$ is a prime of $F$ dividing
$p$.

Let $(A,i)$ be an $n^2$-dimensional $M_n(F)$-linear abelian variety.
The action of $M_n(F)$
makes $A$ isogenous to $A_0^n$, where $A_0$ is an $n$-dimensional abelian
variety with complex multiplication by $F$.
Thus there  is an equivalence between 
the quasi-isogeny category of $n$-dimensional abelian varieties $A_0$ with 
complex multiplication by $F$ and the quasi-isogeny 
category of $n^2$ dimensional 
$M_n(F)$-linear abelian varieties $(A,i)$. 

More generally, we saw in Section~\ref{sec:HTFpbar} that 
if $B$ is any central simple algebra over $F$ of dimension
$n^2$, then there exist simple $B$-linear abelian varieties with
$A(u)$ of dimension $n$ and slope $1/n$.  The $B$-linear structure on $A$
induces a $B_u$-linear structure on $A(u)$.  If we assume that $B$ is split
over $u$, then the $p$-divisible group $A(u)$ is isogenous to a product
$(\epsilon A(u))^n$, where $\epsilon$ is an idempotent of $B_u = M_n(F_u)$
and $\epsilon A(u)$ is a $1$-dimensional $p$-divisible group of slope
$1/n$.

Let $*$ be a positive involution on $B$.
Introducing the structure of a compatible polarization $\lambda$ on our 
$n^2$-dimensional $B$-linear abelian varieties, 
one could form a moduli stack of tuples $(A,i,\lambda)$ for which the
$p$-divisible group $A(u)$ is $n$-dimensional.

These moduli stacks will, in general, have
infinitely many components.
However, Lemma~\ref{lem:polclasslocal} indicates that the classification of
such polarizations is controlled by the local similitude classes of 
$\lambda$-Weil pairings on the Tate modules $V_\ell(A)$.
We may pick out finitely many components of our moduli space by fixing a
global pairing $\bra{-,-}$, and restricting our attention to only those
tuples $(A,i,\lambda)$ for which the Weil pairings $\lambda\bra{-,-}$ on
$V_\ell(A)$ are similar to the local pairing $\bra{-,-}_\ell$ for each
$\ell \ne p$.  Such moduli stacks are instances of \emph{Shimura stacks}.
\index{Shimura stack}

\section{Initial data}\label{sec:data}

Begin with the following data.
\begin{align*}
F & = 
\text{quadratic imaginary extension of $\QQ$, such that $p$ splits as $u
u^c$.}
\index{2F@$F$}
\index{2U@$u$}
\index{2Uc@$u^c$}
\\
\calO_F & = 
\text{ring of integers of $F$.}
\\
B & = 
\text{central simple algebra over $F$, $\dim_F B = n^2$, split over $u$ and
$u^c$.} 
\\
(-)^* & = 
\text{positive involution on $B$ of the second kind, i.e.:} 
\\
& \qquad 1. \: 
\Tr_{F/\QQ} \Tr_{B/F} (xx^*) > 0 \: \text{for $x \ne 0$.} 
\\
& \qquad 2. \:
\text{$\ast$ restricts to conjugation on $F$.} 
\\
\calO_B & = 
\text{maximal order in $B$ such that $\calO_{B,(p)}$ is preserved under $*$.} 
\\
V & = 
\text{free left $B$-module.}
\\
\bra{-,-} & = 
\text{$\QQ$-valued non-degenerate alternating form on $V$} \\
& \qquad \text{which is
$*$-hermitian.  (This means $\bra{\alpha x, y} = \bra{x, \alpha^* y}$.} \\
L & = 
\text{$\calO_B$-lattice in $V$,
% } \\
% & \qquad \text{
$\bra{-,-}$ restricts to give integer values on $L$,} \\
& \qquad \text{and makes $L_{(p)}$ self-dual.}
\end{align*}
\index{2OF@$\mc{O}_F$}
\index{2B@$B$}
\index{0ast@$*$}
\index{2OB@$\mc{O}_B$}
\index{2V@$V$}
\index{0bra@$\bra{-,-}$}
\index{2L@$L$}
\index{involution!positive}
\index{involution!of the second kind}

From this data we define:
\begin{align*}
L^p & = 
\prod_{\ell \ne p} L \otimes \ZZ_\ell,
\\
V^p & = V \otimes_\QQ \AF^{p,\infty},
\\
C & = 
\End_B(V),
\\
(-)^\iota & =
\text{involution on $C$ defined by $\bra{av, w} = \bra{v, a^\iota w}$,}
\\
\calO_C & = 
\text{order of elements $x \in C$ such that $x(L) \subseteq L$},
\\
GU(R) & = 
\{ g \in (C \otimes_\QQ R)^\times \: : \: g^\iota g \in R^\times \},
\\
U(R) & = 
\{ g \in (C \otimes_\QQ R)^\times \: : \: g^\iota g = 1 \},
\\
K^{p}_0 & =
\{ g \in GU(\AF^{p,\infty}) \: : \: g(L^p) = L^p \}. 
\end{align*}
\index{2Lp@$L^p$}
\index{2Vp@$V^p$}
\index{2C@$C$}
\index{1iota@$\iota$}
\index{2OC@$\mc{O}_C$}
\index{2GU@$GU$}
\index{2U@$U$}
\index{2Kp0@$K^p_0$}

We are interested in the case where we have:
\begin{align*}
V & = B, \\
L & = \calO_B, \\
U(\RR) & \cong U(1,n-1).
\end{align*}
It then follows that we have
\begin{align*}
C & \cong B^{op}, \\
\calO_C & \cong \calO_B^{op} \qquad \text{(by maximality of $\calO_B$)}.
\end{align*}
In our case $*$ may equally well be regarded as an involution $*$ on $C$,
and there exists (Lemma~\ref{lem:betadef}) 
an element $\beta \in B$ 
\index{1beta@$\beta$}
satisfying $\beta^* = -\beta$ which
encodes $\bra{-,-}$:
$$
\bra{x,y} = \Tr_{F/\QQ} \Tr_{B/F}(x \beta y^*).
$$
Let $\gamma$ be the element 
\index{1gamma@$\gamma$}
$\beta$ regarded as an element of $C$.  Then
$\iota$ is given by
$$
z^\iota = \gamma^{-1} z^* \gamma.
$$

Tensoring with $\RR$, and taking the complex embedding of $F$ which sends
$\delta$ to a negative multiple of $i$,
\index{1delta@$\delta$}
we may identify the completions of our simple algebras with matrix
algebras over $\CC$:
\begin{align*}
B_\infty & = M_n(\CC), \\ 
* & = \text{conjugate transpose}, \\
C_\infty & = M_n(\CC) \quad \text{(identified with $B$ through the
transpose)}, \\
\beta & = 
\begin{bmatrix}
e_1 i &   &   & 
\\
  & -e_2 i&   & 
\\
  &   & \ddots &
\\
  &   &   & -e_n i
\end{bmatrix}.  \\
\end{align*}
Here the $e_i$'s are positive real numbers.
\index{2Ei@$e_i$}

Because $B$ was assumed to be split over $u$, we may fix an isomorphism
$$ \mc{O}_{B,u} \cong M_n(\mc{O}_{F,u}) = M_n(\ZZ_p). $$
Let $\epsilon \in {\calO}_{B,u}$ 
\index{1epsilon@$\epsilon$}
be the projection associated by this
isomorphism to the
matrix which has a $1$ in the $(1,1)$ entry, and zeros elsewhere.

\begin{rmk}
Giving the data $(B, (-)^*, \bra{-,-})$ is essentially the same as
specifying a form of the similitude group $GU$.  The forms of $GU$
are classified by the Galois cohomology group $H^1(\QQ, PGU)$, where the
algebraic group $PGU$
\index{2PGU@$PGU$}
is the quotient of $GU$ by the subgroup $T_F$, 
\index{2TF@$T_F$}
where
$$ T_F(R) = (F \otimes_\QQ R)^\times. $$
Clozel describes this computation  in \cite{clozel}: an element in
$H^1(\QQ,PGU)$ corresponds uniquely 
to the local invariants of the division algebra
$B$, as well as the difference between the classes of $H^1(\QQ, GU)$
determined  by the involution $*$ and the pairing $\bra{-,-}$ (see
Section~\ref{sec:formclassification}).
\end{rmk}

\section{Statement of the moduli problem}\label{sec:moduli}

Assume that $S$ is a scheme, and that $A$ is an abelian
scheme over $S$.
A \emph{polarization} $\lambda\co A \rightarrow A^\vee$ of $A$ 
\index{polarization!of an abelian scheme}
is an isogeny
which restricts to a polarization on each of the geometric fibers of $A/S$.
For $R$ a ring contained in $\RR$, an $R$-polarization 
\index{polarization!$R$-polarization}
is an $R$-isogeny
$$ \lambda \in \Hom(A,A^\vee) \otimes R $$
which is a positive linear combination of polarizations.
A $\QQ$-polarization defines a $\lambda$-\emph{Rosati 
involution} 
\index{involution!Rosati}
\index{Rosati involution}
$\dag$ on
\index{0dag@$\dag$}
$\End^0(A)$ by $f^\dag = \lambda^{-1} f^\vee \lambda$.

We shall define two functors
$$ \mathcal{X}, \mathcal{X}'\co \{\text{locally noetherian 
formal schemes}/\Spf(\ZZ_p)\}
\rightarrow \text{groupoids}. $$
We \index{2X@$\mc{X}$} \index{2X'@$\mc{X}'$}
shall then show that these functors are equivalent.

\begin{rmk}
Every formal scheme $S$ over $\Spf(\ZZ_p)$ is a formal colimit $S =
\varinjlim S_i$, where $S_i$ is a scheme on which $p$ is locally nilpotent.  It
therefore suffices to define the functors $\mc{X}$ and $\mc{X}'$ for locally
noetherian schemes on which $p$ is locally nilpotent.
\end{rmk}

\subsection*{The functor $\mathcal{X}$}\label{sec:functorX}

Assume that $S$ is a locally noetherian scheme on which $p$ is locally
nilpotent.
The objects of the groupoid $\mathcal{X}(S)$ consist of
\index{2X@$\mc{X}$}
tuples of data
$(A,i,\lambda)$
as follows.

\begin{tabular}{lp{19pc}}
$A$ & is an abelian scheme over $S$ of dimension $n^2$.\\
$\lambda\co A \rightarrow A^\vee$& is a $\ZZ_{(p)}$-polarization,
 with Rosati involution $\dag$ on $\End(A)_{(p)}$.\\
$i\co \calO_{B,(p)} \hookrightarrow \End(A)_{(p)}$ & is an
  inclusion of rings, satisfying $i(b^*) = i(b)^\dag$ in $\End(A)_{(p)}$,
  such that $\epsilon A(u)$ is $1$-dimensional.
\end{tabular}
\index{2A@$A$}
\index{2I@$i$}
\index{1lambda@$\lambda$}
\index{0dag@$\dag$}

We impose the following additional restrictions on $(A,i,\lambda)$. 
Choose a geometric point $s$ in each 
component of $S$. We require that for each of these points
there \emph{exists} an ${\calO}_B$-linear integral uniformization
$$
\eta \co L^p \xrightarrow{\cong} T^p(A_s)
$$
so \index{1eta@$\eta$}
that when tensored with $\QQ$, $\eta$ sends $\bra{-,-}$ to an 
$(\AF^{p,\infty})^\times$-multiple of
the $\lambda$-Weil pairing. 
%\begin{enumerate}
%\item $\eta$ is $\calO_B$-linear.
%\item $\eta$ is $\pi_1(S,s)$-invariant.
%\item when tensored with $\QQ$, $\eta$ sends $(-,-)$ to an 
%$(\AF^{p,\infty})^\times$-multiple of
%the $\lambda$-Weil pairing. 
%\end{enumerate}
We do not fix this uniformization as part of the data.

We pause to explain this last restriction.  By
Lemma~\ref{lem:polclasslocal}, the isogeny classes of
compatible polarizations $\lambda$ on a fixed $B$-linear abelian variety
$A$ over an algebraically closed field $k$
are determined by the elements of $H^1(\QQ_\ell, GU)$ given at each $\ell
\ne p$ by the $\lambda$-Weil pairing on $V_\ell(A)$.
%The
%$\pi_1(S,s)$-invariance implies that
The existence of the isomorphism
$\eta$ is independent of the point $s$.

A morphism
$$ (A, i, \lambda) \rightarrow (A', i', \lambda') $$
in $\mathcal{X}(S)$ consists of an isomorphism of abelian schemes over $S$
$$ \alpha\co A \xrightarrow{\cong} A' $$
such that
\begin{alignat*}{2}
\lambda & = r\alpha^\vee \lambda' \alpha, 
\quad && 
r \in \ZZ_{(p)}^\times, \\
i'(b)\alpha & = \alpha i(b), 
\quad &&
b \in \calO_B. \\
\end{alignat*}
In particular, the isomorphism class of $(A,i,\lambda)$ depends only on
the weak polarization determined by $\lambda$.

\subsection*{The functor $\mathcal{X}'$}

\index{2X'@$\mc{X}'$}
The functor $\mathcal{X}$ classifies $\calO_B$-linear polarized 
abelian varieties up
to isomorphism (with certain restrictions on the slopes of the
$p$-divisible group and the Weil pairings).  We shall now introduce a
different functor $\mathcal{X}'$ which classifies 
$\calO_{B,(p)}$-linear polarized
abelian varieties with rational level structure up to isogeny.  The functor
$\mathcal{X}'$ will be shown to be equivalent to $\mathcal{X}$.

Assume that $S$ is a locally noetherian 
connected scheme on which $p$ is locally nilpotent. 
Fix a geometric point $s$ of $S$.
The objects of the groupoid $\mathcal{X}'(S)$ consist of
tuples of data
$(A,i, \lambda,[\eta])$
as follows.

\begin{center}
\begin{tabular}{lp{19pc}}
$A$& is an abelian scheme over $S$ of dimension $n^2$.\\
$\lambda\co A \rightarrow A^\vee$ &is a $\ZZ_{(p)}$-polarization,
  with Rosati involution $\dag$ on $\End(A)_{(p)}$.\\
$i\co \calO_{B,(p)} \hookrightarrow \End(A)_{(p)}$ &is an inclusion of rings,
  satisfying $i(b^*) = i(b)^\dag$ in $\End(A)_{(p)}$,
  such that $\epsilon A(u)$ is $1$-dimensional. \\
$[\eta]$ &is a rational $K^p_0$ level structure, i.e. the $K^{p}_0$-orbit
of a rational uniformization $\eta\co V^p \xrightarrow{\cong}
V^p(A_s)$,
such that $\eta$ is ${\calO}_{B,(p)}$-linear, 
$[\eta]$ is $\pi_1(S,s)$-invariant, and 
$\eta$ sends $\bra{-,-}$ to an $(\AF^{p,\infty})^\times$-multiple
of the $\lambda$-Weil pairing.
\end{tabular}
\end{center}
\index{2A@$A$}
\index{2I@$i$}
\index{1lambda@$\lambda$}
\index{0dag@$\dag$}

A morphism
$$ (A, i, \lambda, \eta) \rightarrow (A', i', \lambda',\eta') $$
in $\mathcal{X}'(S)$ consists of a $\mb Z_{(p)}$-isogeny of
abelian schemes over $S$ 
$$ \alpha\co A \xrightarrow{\simeq} A' $$
such that
\begin{alignat*}{2}
\lambda & = r\alpha^\vee \lambda' \alpha, 
\quad && 
r \in \ZZ_{(p)}^\times, \\
i'(b)\alpha & = \alpha i(b), 
\quad &&
b \in \calO_{B,(p)}, \\
[\eta'] & = \alpha_* [\eta].
\quad &&
\end{alignat*}

The $\pi_1(S,s)$-invariance of the level structure $[\eta]$ implies
that the objects are independent of the choice of the point $s$ in
$S$.  The functor $\mathcal{X}'$ extends to schemes which are not
connected by taking products over the values on the
components.

\section{Equivalence of the moduli problems}\label{sec:equivmoduli}

Fix a locally noetherian 
connected scheme $S$ on which $p$ is locally nilpotent, 
with geometric 
point $s$ as before.  For each object $(A,i, \lambda)$ of
$\mathcal{X}(S)$, choose an $\calO_B$-linear similitude
$$ \eta\co L^p \xrightarrow{\cong} T^p(A_s). $$
Then $(A, i_{(p)}, \lambda, [\eta_\QQ])$ is an object of $\mathcal{X}'(S)$.
Different choices of $\eta$ yield canonically equivalent objects 
$(A,i_{(p)},\lambda, [\eta_\QQ])$, because any two choices of $\eta$ will
necessarily differ by an element of $K^{p}_0$.  Thus we have produced a
functor
$$ F_S\co \mathcal{X}(S) \rightarrow \mathcal{X}'(S) $$
which is natural in $S$.  The rest of this section will be devoted to
proving the following theorem.

\begin{thm}\label{thm:equivmoduli}
The functor $F_S$ is an equivalence of categories.
\end{thm}

\subsection*{$F_S$ is essentially surjective}

Suppose that $(A, i, \lambda, [\eta])$ is an object of
$\mathcal{X}'(S)$.  In particular, choose a rational uniformization
$\eta$. We need to show that this data is isogenous to a new set of
data $(A', i', \lambda', [\eta'])$ where $\eta'$ lifts to an integral
level structure.

Define $L_s$ to be the virtual subgroup $\kappa(\eta) \in \VSub^p(A_s)$ 
of Section~\ref{sec:levelsub}.  The $\pi_1(S,s)$ invariance of $[\eta]$
implies that this virtual subgroup is invariant under $\pi_1(S,s)$.
Therefore, the virtual subgroup $L_s$ extends to a local system of
virtual subgroups $L$ of $A$.  Explicitly, 
there exists an integer $N$ prime to $p$ such that $N^{-1}L_s$ corresponds to
an actual subgroup scheme $H_s$ of $A_s$, and
a finite subgroup scheme $H$ of $A$ over $S$ extending $H_s$.
Define $A'$ to be the quotient abelian variety, with quotient isogeny
$$ \alpha' \co A \rightarrow A/H = A'. $$

The $\mb Z_{(p)}$-isogeny $\alpha = N^{-1}\alpha'$ induces a 
$\mb Z_{(p)}$-isogeny $\alpha_s\co
A_s \to A'_s$.  By construction, $\alpha_s(L_s) = T^p(A'_s)$.

The $\ZZ_{(p)}$-polarization $\lambda'$ on $A'$ is 
defined to be the $\ZZ_{(p)}$-isogeny that makes the following diagram commute.
$$
\xymatrix{
A \ar[d]_{\alpha} \ar[r]^{\lambda} &
A^\vee 
\\
A' \ar[r]_{\lambda'} &
(A')^\vee \ar[u]_{\alpha^\vee}
} $$

Define $\calO_{B,(p)}$-multiplication on $A'$
$$ i' \co \calO_{B,(p)} \hookrightarrow \End(A')_{(p)} $$
by the formula $i'(b) = \alpha i(b) \alpha^{-1}$.
Compatibility of $i'$ with the $\lambda'$-Rosati involution
$\dag'$ is easily checked:
\begin{align*}
i'(b)^{\dag'} 
& = (\lambda')^{-1} i'(b)^\vee \lambda' \\ 
& = ((\alpha^\vee)^{-1} \lambda \alpha^{-1})^{-1}(\alpha i(b) \alpha^{-1})^\vee
((\alpha^\vee)^{-1} \lambda \alpha^{-1}) \\
& = \alpha \lambda^{-1} i(b)^\vee \lambda \alpha^{-1} \\
& = \alpha i(b)^\dag \alpha^{-1} \\
& = \alpha i(b^*) \alpha^{-1} \\
& = i'(b^*).
\end{align*}

The rational level structure $[\eta']$ is defined to be the the orbit of
the composite
$$ \eta'\co V^p \xrightarrow{\eta} V^p(A_s) \xrightarrow{\alpha_s}
V^p(A'_s). $$
By construction, the $\mb Z_{(p)}$-isogeny $\alpha$ gives an
isomorphism from $(A,i,\lambda,[\eta])$ to $(A',i',\lambda',[\eta'])$.
We wish to show that $\eta'$ lifts to an integral uniformization.  By
Lemma~\ref{lem:integralkappa}, it suffices to show that the virtual
subgroup $\kappa(\eta') \in \VSub^p(A'_s)$ is $T^p(A_s)$.

However, by definition $\kappa(\eta') = \alpha_s(\eta(L^p)) =
\alpha_s(L_s)$, and we have already shown that $\alpha_s(L_s) =
T^p(A_s)$.

We will have shown that $(A',i',\lambda')$ may be regarded as an object of
$\mathcal{X}(S)$, provided that the $\calO_{B,(p)}$-multiplication on
$\End(A')_{(p)}$ given by
$i'$ lifts to $\calO_B$-multiplication on $\End(A')$.  
The existence of an $\calO_B$-action on $A'$ follows from the following
pair of pullback squares (the righthand square is a pullback by
Corollary~\ref{cor:tatepullback}).
$$
\xymatrix@C+1em{
\End(A') \ar[r] \ar[d] &
\: \prod_{\ell \ne p} \End(A')_\ell \: \ar@{^{(}->}[r] \ar[d] &
\End(T^p A'_s) \ar[d] 
\\
\End(A')_{(p)} \ar[r] &
\: \left( \prod_{\ell \ne p} \End(A')_\ell \right) \otimes \QQ \: 
\ar@{^{(}->}[r] &
\End (V^p A'_s) }
$$

\subsection*{$F_S$ is fully faithful}

The faithfulness of $F_S$ is obvious.  We just need to show that it is
full.
Suppose that $(A,\lambda,i,\eta)$ and $(A',\lambda',i',\eta')$ are two
objects of $\mathcal{X}'(S)$ which are in the image of $F_S$.  Without
loss of generality assume that the rational level structures $\eta$
and $\eta'$ lift to integral level structures.  Suppose that
$$ f\co (A,\lambda, i, \eta) \rightarrow (A',\lambda', i', \eta') $$
is a morphism in $\mathcal{X}'(S)$.  In particular, $f\co A \rightarrow A'$
is a $\ZZ_{(p)}$-isogeny.  
Because $\eta$ and $\eta'$ are integral level structures, we may conclude
that the induced map
$$ f_*\co V^p(A) \xrightarrow{\cong} V^p(A') $$
lifts to an isomorphism
$$ f_*\co T^p(A) \xrightarrow{\cong} T^p(A'). $$
The pair of pullback squares (the righthand square is a pullback by 
Corollary~\ref{cor:tatepullback})
$$
\xymatrix@C+.5em{
\Hom(A,A') \ar[r] \ar[d] &
\: \prod_{\ell \ne p} \Hom(A,A')_\ell \: \ar@{^{(}->}[r] \ar[d] &
\Hom(T^p A_s,T^p A'_s) \ar[d] 
\\
\Hom(A,A')_{(p)} \ar[r] &
\: \left( \prod_{\ell \ne p} \Hom(A,A')_\ell \right) \otimes \QQ \: 
\ar@{^{(}->}[r] &
\Hom (V^p A_s,V^p A'_s) }
$$
allows us to conclude that $f$ is actually an isomorphism of abelian
varieties.

\section{Moduli problems with level structure}

Let ${K}^p$ 
\index{2Kp@$K^p$}
be an open subgroup of $K^p_0$.  Associated to $K^p$ are
variants of the functors $\mc{X}$ and $\mc{X'}$.

For a locally noetherian connected scheme $S$ on which $p$ is
locally nilpotent,
with geometric point $s \rightarrow S$, define a groupoid
$\mc{X}_{K^p}(S)$ 
\index{2XKp@$\mc{X}_{K^p}$}
to have objects $(A,i,\lambda, [\eta]_{K^p})$ where:
\begin{center}
\begin{tabular}{lp{21pc}}
$(A,i,\lambda)$ & is an object of $\mc{X}(S)$, \\
$[\eta]_{K^p}$ & is an integral $K^p$ level structure, i.e. the $K^{p}$-orbit
of an integral uniformization $\eta\co L^p \xrightarrow{\cong}
T^p(A_s)$,
such that $\eta$ is ${\calO}_{B}$-linear, 
$[\eta]_{K^p}$ is $\pi_1(S,s)$-invariant, and 
$\eta$ sends $\bra{-,-}$ to an $(\AF^{p,\infty})^\times$-multiple
of the $\lambda$-Weil pairing.
\end{tabular}
\end{center}
The \index{1etaKp@$[\eta]_{K^p}$}
morphisms of $\mc{X}_{K^p}(S)$ are those isomorphisms in $\mc{X}(S)$
which preserve the level structure.

Define a groupoid
$\mc{X'}_{K^p}(S)$ 
\index{2X'Kp@$\mc{X}'_{K^p}$}
to have objects $(A,i,\lambda, [\eta]_{K^p})$ where:
\begin{center}
\begin{tabular}{lp{21pc}}
$(A,i,\lambda, [\eta]_{K^p_0})$ & is an object of $\mc{X'}(S)$, \\
$[\eta]_{K^p}$ & is a rational $K^p$ level structure, i.e. the 
$K^{p}$-orbit
of a rational uniformization $\eta\co V^p \xrightarrow{\cong}
V^p(A_s)$,
such that $\eta$ is $B$-linear, 
$[\eta]_{K^p}$ is $\pi_1(S,s)$-invariant, and 
$\eta$ sends $\bra{-,-}$ to an $(\AF^{p,\infty})^\times$-multiple
of the $\lambda$-Weil pairing.
\end{tabular}
\end{center}
The \index{1etaKp@$[\eta]_{K^p}$}
morphisms of $\mc{X'}_{K^p}(S)$ are those prime-to-$p$ isogenies in 
$\mc{X'}(S)$
which preserve the $K^p$ level structures.

Clearly, the functors $\mc{X}$ and $\mc{X}'$ are recovered by taking $K^p
= K^p_0$.

Since any integral level structure determines a rational level structure,
the functors $F_S$ lift to give functors
$$ F_{K^p,S}\co \mc{X}_{K^p}(S) \rightarrow \mc{X'}_{K^p}(S) $$
which are natural in $S$.
Theorem~\ref{thm:equivmoduli} generalizes to the following
theorem.

\begin{thm}\label{thm:levelequivalence}
The functors $F_{K^p,S}$ are equivalences of categories.
\end{thm}

\begin{rmk}
The functor $\mc{X}'_{K^p}$ has the advantage that it may be defined for
\emph{any} compact open subgroup $K^p$ of $GU(\AF^{p,\infty})$, not just
those contained in $K^p_0$.
\end{rmk}

\section{Shimura stacks}

If $K'^p < K^p$ are open compact subgroups of $GU(\AF^{p,\infty})$, 
then there are natural
transformations
\begin{equation}\label{eq:f_def}
f_{K'^p, K^p}\co \mc{X'}_{K'^p}(-) \rightarrow \mc{X'}_{K^p}(-)
\end{equation}
given by sending a tuple $(A,i,\lambda, [\eta]_{K'^p})$ to the tuple
$(A,i,\lambda, [\eta]_{K^p})$.

The functors $\mc{X'}_{K^p}$ are representable.  For $K^p$ sufficiently
small, the representability of $\mc{X}'_{K^p}$ by a scheme is discussed in 
\cite[Sec.~5,6]{kottwitz}.  For general subgroups $K^p$, the
representability of $\mc{X}'_{K^p}$ by an algebraic stack is discussed in
\cite[Sec.~7.1.2]{hida}.  The following theorem represents a compilation of
these results.

\begin{thm}\label{thm:shimurastack}
$\quad$
\begin{enumerate}
\item 
The functor $\mc{X'}_{K^p}$ is representable by 
a Deligne-Mumford stack $\Sh(K^p)$ over $\ZZ_p$.  
\index{Shimura stack}
\index{2ShKp@$\Sh(K^p)$}

\item
The map 
$$ f_{K'^p, K^p}\co \Sh(K'^p) \rightarrow \Sh(K^p) $$
induced by the natural transformation (\ref{eq:f_def}) is \'etale, of degree
equal to $[K^p : K'^p]$.  If $K'^p$ is normal in $K^p$, then the covering
is Galois with Galois group $K^p/K'^p$.

\item
For $K^p$ sufficiently small, $\Sh(K^p)$ is a quasi-projective scheme
over $\ZZ_p$.  It is projective if $B$ is a division algebra.

\end{enumerate}
\end{thm}

\begin{rmk}
The moduli problem considered in either of these
sources has a slightly different formulation than presented here, due to
the fact that the authors of \cite{kottwitz}, \cite{hida} 
are working over $\mc{O}_{F,(u)}$.
Specializing to $\mc{O}_{F,u} = \ZZ_p$ produces an equivalent moduli
problem (see Remark~\ref{rmk:det}).
\end{rmk}

\chapter{Deformation theory}
\label{chap:deformations}

\section{Deformations of $p$-divisible groups}

We briefly summarize the deformation theory of
$1$-dimensional $p$-divisible groups, 
restricting to the case when the dimension of the
formal group is $1$, following \cite[II.1]{harristaylor}, 
\cite[Prop.~4.5]{drinfeld}.
\index{p-divisible group@$p$-divisible group!deformation theory}

Let $k$ be a field of characteristic $p$.  
Fix a $p$-divisible group $\br{\GG}$ 
\index{2G@$\br{\GG}$}
over $k$.  
Consider the functor
$$ \Def_{\br{\GG}}\co 
\{ \text{noetherian complete local rings} \} 
\rightarrow \{ \mathrm{groupoids} \} $$
whose \index{2DefG@$\Def_{\br{\GG}}$}
$R$-objects are given by
$$ \Def_{\br{\GG}}(R) = \left\{ (\GG, i, \alpha) \: : \: 
\begin{array}{l}
i\co k \rightarrow R/\mf{m}_R, \\
\GG \: \text{a $p$-divisible over $R$}, \\
\alpha\co \GG \big\vert_{\Spec(R/\mf{m}_R)} \xrightarrow{\cong} i^* \br{\GG}
\end{array}
\right\},
$$
and whose $R$-morphisms are those isomorphisms of $p$-divisible groups that
restrict to the identity on $\br{\GG}$.

Lubin and Tate showed that deformations of 
$1$-dimensional formal groups of finite height 
have no automorphisms, and that this functor is represented by a formal scheme.

\begin{thm}[Lubin-Tate \cite{lubintate}]\label{thm:lubintate}
Suppose that $\br{\GG}$ is a $1$-dimensional formal group of finite height $h$. 
The functor $\Def_{\br{\GG}}(-)$ is representable by the formal scheme 
$\Def_{\br{\GG}} = \Spf(B)$ where 
$B = W(k)[[u_1, \ldots, u_{h-1}]]$.
\end{thm}

Given a $p$-divisible group $\GG = \varinjlim \GG_i$ 
over $R$, it is convenient to consider its associated 
sheaf of $\ZZ_p$-modules on $\Spf(R)$ in the fppf topology.  Given a 
local $R$-algebra $T$, the sections are given by
$$ \GG(T) = \varprojlim_j \varinjlim_i \GG_i(T/\mf{m}_T^j). $$
The \index{2GT@$\GG(T)$} \index{p-divisible group@$p$-divisible group!fppf sheaf}
resulting functor
$$ \{ \text{$p$-divisible groups} \} \rightarrow \{ \text{fppf sheaves of
$\ZZ_p$-modules} \} $$
is fully faithful.

\begin{lem}\label{lem:GG^0(T)}
Suppose $\GG$ is a $1$-dimensional formal group over $R$ of finite height, 
with coordinate $x$, so that with respect to 
this coordinate the sum formula is given by 
$$ x_1 +_\GG x_2 = \sum_{i,j} a_{i,j} x_1^i x_2^j $$
and that $T$ is a local $R$-algebra.  Then
there is an isomorphism of groups
$$ \GG(T) \cong (\mf{m}_T, +_{\GG}). $$
\end{lem}

\begin{proof}
We have an isomorphism of topological
rings (see Example~\ref{exam:formal}):
$$ T[[x]] \cong \varprojlim_i T[[x]]/[p^i]_\GG(x). $$
The $T$-sections of $\GG$ are then given by the formula
\begin{align*}
\GG(T) & = \varprojlim_j \varinjlim_i \Alg_T(T[[x]]/[p^i]_\GG(x),
T/\mf{m}_T^j) \\ 
& = \Alg_T^c(T[[x]], T) \\
& \cong \mf{m}_T.
\end{align*}
\end{proof}

We are now in a position to give a description of the deformation theory of
$1$-dimensional $p$-divisible groups.

\begin{thm}\label{thm:pdivdef}
Suppose that the $\br{\GG}$ is a $1$-dimensional $p$-divisible group over
$k$ of finite height $h$.  Then the functor $\Def_{\br{\GG}}(-)$ is
representable by the formal scheme $\Def_{\br{\GG}} = \Spf(B)$, where
$$ B = W(k)[[u_1, \ldots, u_{h-1}]]. $$
\end{thm}
\index{2Ui@$u_i$}

\begin{proof}[Sketch proof of Theorem~\ref{thm:pdivdef}]
Let $(\GG, i, \alpha)$ be an object of $\Def_{\br{\GG}}(T)$.  
Then there is a short exact sequence (\cite[Cor.~II.1.2]{harristaylor})
$$ 0 \rightarrow \GG^0 \rightarrow \GG \rightarrow \GG^{et} \rightarrow 0
$$
where $\GG^0$ is formal of height $k$ and $\GG^{et}$ is ind-\'etale of height
$h-k$.  
By Theorem~\ref{thm:lubintate}, the deformation $\GG^0$ is
classified by a map of local rings
$$ W(k)[[u_1, \ldots, u_{k-1}]] \rightarrow T. $$
The deformation $\GG^{et}$ is unique up to isomorphism.  The extension is
classified by an element of 
$$ \Ext^1_{\ZZ_p}(\GG^{et}, \GG^0), $$
where the $\Ext$ group is taken in the category of sheaves of
$\ZZ_p$-modules.
Let $T\GG^{et}$ denote the sheaf given by $\varprojlim \GG^{et}[p^i]$.
The short exact sequence of sheaves
$$ 0 \rightarrow T\GG^{et} \rightarrow T\GG^{et} \otimes_{\ZZ_p} \QQ_p
\rightarrow \GG^{et} \rightarrow 0 $$
gives rise to a long exact sequence of $\Ext$ groups.  By descent, we may
assume that $T$ is separably closed, so that $\GG^{et}$ is abstractly
isomorphic to $(\QQ_p/\ZZ_p)^{h-k}_T$.
The boundary homomorphism gives an isomorphism
$$ \Hom_{\ZZ_p}((\ZZ_p)^{h-k}_T, \GG^0) \xrightarrow{\cong}
\Ext^1_{\ZZ_p}(\GG^{et}, \GG^0). $$
Therefore we have,  
\begin{align*}
\Def_{\br{\GG}}(T) & \cong \Def_{\br{\GG}^0}(T) \times \Hom_{\ZZ_p}(\GG^{et},
\GG^0) \\
& \cong \Def_{\br{\GG}^0}(T) \times \Hom_{\ZZ_p}((\ZZ_p)^{h-k}_T, \GG^0) \\
& \cong \Def_{\br{\GG}^0}(T) \times \mf{m}_T^{h-k}
\end{align*}
The last isomorphism is given by Lemma~\ref{lem:GG^0(T)}.
We deduce that there is an isomorphism
$$ \Def_{\br{\GG}}(T) \cong \mathrm{Ring}^c(W(k)[[u_1, \ldots,
u_{k-1}]][[u_{k}, \ldots, u_{h-1}]], T). $$
\end{proof}

\section{Serre-Tate theory}

We recall the Serre-Tate theorem \cite{katz}. 

\begin{thm}[Serre-Tate]
Suppose $j\co S \into R$ is a closed embedding of schemes with nilpotent
ideal sheaf ${\mc{I}}$ such that $p$ is locally nilpotent on $R$.
Then the diagram of categories
\xym{
\{\text{Abelian varieties over $R$}\} \ar[r]^{j^*} 
\ar[d]^{(-)(p)} &
\{\text{Abelian varieties over $S$}\} \ar[d]^{(-)(p)} \\
\{\text{$p$-divisible groups over $R$}\} \ar[r]^{j^*} &
\{\text{$p$-divisible groups over $S$}\}
}
is a pullback.  In other words, there is an equivalence of categories
between the category of abelian varieties over $R$ and the category of
tuples
\[
\left\{(A,\mb G, \phi)\ \Bigg|
\text{  \begin{tabular}{l}
    $A$ an abelian variety over $S$, \\
    $\mb G$ a $p$-divisible group over $R$, \\
    $\phi\co A[p^\infty] \to j^* \mb G$ an isomorphism
  \end{tabular}}
\right\}.
\]
\end{thm}
\index{abelian variety!deformation theory}

\section{Deformation theory of points of $\Sh$}

Let $K^p$ be a compact open subgroup of $GU(\AF^{p,\infty})$.
Consider the functor
\begin{align*}
\Phi\co \mc{X}_{K^p}(S) & \rightarrow \text{$p$-divisible groups}/S, \\
(A,i,\lambda, [\eta]) & \mapsto \epsilon A(u) 
\end{align*}
where $\epsilon A(u)$ is the $1$-dimensional summand of the $p$-divisible
group $A(p)$.
In this section we explain how the deformation theory of points of
$\Sh(K^p)$ is
controlled by the deformation theory of the associated 
$1$-dimensional $p$-divisible groups.  The material in this section is
well known (see, for instance, \cite[Lem.~III.4.1]{harristaylor}).

\begin{thm}\label{thm:defsh}
Suppose $j\co S \into R$ is a closed embedding of schemes with nilpotent
ideal sheaf ${\mc{I}}$ such that $p$ is locally nilpotent on $R$.
Then the diagram of categories
\xym{
\mc{X}'_{K^p}(R) \ar[r]^{j^*} 
\ar[d]^{\Phi} &
\mc{X}'_{K^p}(S) \ar[d]^{\Phi} \\
\{\text{$p$-divisible groups over $R$}\} \ar[r]^{j^*} &
\{\text{$p$-divisible groups over $S$}\}
}
is a pullback.  Equivalently, there is an equivalence of categories
between the category $\mc{X}'_{K^p}(R)$ and the category of
tuples
\[
\left\{(A, i, \lambda, [\eta],\mb G, \alpha)\ \: : \:
\text{  \begin{tabular}{l}
    $(A,i,\lambda, [\eta])$ an object of $\mc{X}_{K^p}(S)$, \\
    $\mb G$ a $p$-divisible group over $R$, \\
    $\alpha \co \epsilon A(u) \to j^* \mb G$ an isomorphism.
  \end{tabular}}
\right\}
\]
\end{thm}

\begin{proof}
The map $i$ takes conjugation on ${\calO}_F$ to the Rosati
involution.  Therefore, the isomorphism
\[
\lambda\co A(p) \to A^\vee(p)
\]
breaks up as a direct sum of the following two isomorphisms:
\begin{align*}
A(u^c) &\to (A^\vee)(u) \cong A(u)^\vee,\\
A(u) &\to (A^\vee)(u^c) \cong A(u^c)^\vee.
\end{align*}
In particular, $A(u)$ and $A(u^c)$ have the same
height with opposite slope decompositions.

The splitting ${\calO}_{B,u} \cong M_n(\mb Z_p)$
gives rise to a further decomposition
\[
A(u) \cong (\epsilon A(u))^n,
\]
viewed as column vectors acted on by the matrix ring.

The summand $\epsilon A(u)$ of $A(p)$ is $1$-dimensional
by assumption.  The associated formal group $A(p)^0$ has a
corresponding canonical $1$-dimensional summand 
$\epsilon A(u)^0$ of height less than or equal to $n$.

A deformation of
$(A,\lambda,i,[\eta])$ to $R$ consists of element $(\tilde A, \tilde
\lambda, \tilde i, [\tilde \eta]) \in \mathcal{X}'_{K^p}(R)$ and an
isomorphism $\phi\co A \to j^* \tilde A$ respecting the structure.
This determines a deformation $\epsilon \tilde A(u)$ of the
$p$-divisible group $\epsilon A(u)$.

The Serre-Tate theorem provides a functorial equivalence between
deformations of $(A,\lambda,i,[\eta])$ and deformations of the
associated $1$-dimensional summand of the $p$-divisible group, as
follows.

Let $\mb G$ be a deformation of $\epsilon A(u)$ to $R$ with
isomorphism $\psi\co \epsilon A(u) \to j^* \mb G$.  The isomorphism
\[
\psi^n\co A(u) \cong (\epsilon A(u))^n \to j^*(\mb G^n)
\]
makes $\mb G^n$ into a deformation of $A(u)$.  The action
of $({\calO}_B)_u \cong M_n(\mb Z_p)$ extends uniquely to an action
on $\mb G^n$.  

The dual $p$-divisible group $(\mb G^\vee)^n$ is a deformation of
$A(u^c) \cong A(u)^\vee$.  It carries an action of
$\mathcal{O}_{B,u^c}$.  We obtain an isomorphism
\[
\phi\co A(p) \cong A(u) \times A(u^c) \to
j^*(\mb G^n \times (\mb G^\vee)^n)
\]
respecting the action of ${\calO}_{B,(p)}$.

The Serre-Tate theorem then implies that for each deformation, there
is an abelian variety $\tilde A$ together with isomorphisms
\begin{eqnarray*}
f\co j^*(\tilde A) &\to& A,\\
g\co \tilde A(p) &\to& \mb G^n \times (\mb G^\vee)^n
\end{eqnarray*}
such that $\phi \circ (f)(p) = j^* g$, and the action of 
${\calO}_{B,(p)}$ extends to a map
\[
\tilde i\co {\calO}_{B,(p)} \to \End(\tilde A)_{(p)}.
\]

The twist isomorphism
\[
\tau\co \mb G^n \times (\mb G^\vee)^n \to (\mb G^\vee)^n \times \mb G^n
\]
is symmetric, and extends the isomorphism
\[
A(u) \times A(u^c) \to A^\vee(u) \times
A^\vee(u^c).
\]

The abelian variety $\tilde A^\vee$ provides a lifting of
the data $(A^\vee, (\mb G^\vee)^n \times \mb G^n, (\phi^\vee)^{-1})$.
The polarization $\lambda$ on $A$ and the twist morphism $\tau$ lift
uniquely to a map $\tilde \lambda\co \tilde A \to \tilde A^\vee$.  As
both $\lambda$ and $\tau$ are symmetric, so is $\tilde \lambda$. 
The maps $\lambda$ and $\tau$ both $*$-commute with the action of
${\calO}_{B,(p)}$, and hence so does $\tilde \lambda$.

The positivity of $\tilde \lambda$ is characterized by positivity at the
geometric points $s$ of $R$, and there is a bijective correspondence
between geometric points of $S$ and $R$; therefore, $\tilde \lambda$
is positive.  Similarly, any level structure $[\eta]$ is determined by
its values on the geometric points $s \in S$, and so $[\eta]$ extends
uniquely.

Therefore, the extension $\mb G$ determines a unique deformation
$(\tilde A, \tilde \lambda, \tilde i, [\eta])$ of the tuple
$(A,\lambda,i,[\eta])$ to $R$.

This extension is functorial, as follows.  Let $f\co
(A,\lambda,i,[\eta]) \to (A',\lambda',i',[\eta'])$ be a morphism in
$\mathcal{X}'_{K^p}(S)$, and $(\mb G,\psi)$, $(\mb G',\psi')$ 
corresponding deformations with a choice of extension $f\co \mb G \to
\mb G'$.

The deformation $\mb G^n \times (\mb G^\vee)^n$ of $A(p)$,
together with its action of ${\calO}_{B,(p)}$ and twist morphism, is
functorial in $\mb G$.  In other words, if $f\co A \to A'$ is an honest
isogeny of abelian varieties, we have a morphism of Serre-Tate data 
\[
(A,\mb G^n \times (\mb G^\vee)^n, \phi) \to 
(A',(\mb G')^n \times ((\mb G')^\vee)^n, \phi').
\]
There is a corresponding map of liftings $\tilde A \to \tilde A'$.  If
$f\co A \to A'$ is instead a quasi-isogeny, there exists a diagram of
isogenies
\[
A \overfrom^{[n]} A \xrightarrow{n f} A',
\]
where $n$ is relatively prime to $p$, and hence a diagram of lifts
\[
\tilde A \xleftarrow{[n]} \tilde A \xrightarrow{n \tilde f} \tilde A'.
\]
Therefore, the category of elements of ${\mc{X}}'_{K^p}(S)$ together with
deformations of their $p$-divisible group to $R$ is equivalent to the
category ${\mc{X}}'(R)$.
\end{proof}

\begin{cor}\label{cor:defsh}
Assume that $K^p$ is sufficiently small so that $\Sh(K^p)$ is a scheme.
Let $s$ be a point of $\Sh(K^p)^\wedge_p$, classifying
the tuple $(A,i,\lambda, [\eta])$.  Then the map
$$ \Phi\co \Sh(K^p)^\wedge_s \rightarrow \Def_{\epsilon A(u)} $$
is an isomorphism of formal schemes.
\end{cor}

The infinitesimal criterion of smoothness \cite[Prop.~17.14.2]{EGA4.4} allows
us to deduce the following.

\begin{cor}
$\Sh(K^p)$ is smooth of relative dimension $n-1$ over $\Spec(\mb
Z_p)$.
\end{cor}

\chapter{Topological automorphic forms}\label{chap:taf}

\section{The generalized Hopkins-Miller theorem}

The following terminology is introduced in \cite{lurie}.

\begin{defn}\label{defn:weaklyevenperiodic}
A homotopy commutative ring spectrum $E$ is \emph{weakly even periodic} if
\index{weakly even periodic}
\begin{enumerate}
\item The homotopy groups $\pi_* E$ are concentrated in even degrees.

\item
The maps
$$ E^2 \otimes_{E^0} E^{n} \rightarrow E^{n+2} $$
are isomorphisms.
\end{enumerate}
\end{defn}

Let $E$ be a weakly even periodic ring spectrum.  The usual 
$H$-space structure on $\CC P^\infty$ gives $\GG_E = 
\Spf(E^0(\CC P^\infty))$ the
structure of a group object in the category of formal schemes.
Here, $\Spf$ is taken by regarding $E^0(\CC P^\infty)$ as $\varprojlim_n
E^0(\CC P^n)$ (condition~(1) above implies that there is no
$\varprojlim^1$-term).

\begin{lem}
The formal scheme $\GG_E$ 
\index{2GE@$\GG_E$}
is a formal group over
$\Spec(E^0)$, i.e. it admits
the structure of a formal group law Zariski locally over $\Spec(E^0)$.
\end{lem}

\begin{proof}
Condition $(2)$ of Definition~\ref{defn:weaklyevenperiodic} implies that
$E^2$ is an invertible $E^0$-module with inverse $E^{-2}$.
Therefore, $E^2$ is locally free (see \cite[Thm~11.6(a)]{eisenbud}, where
the noetherian hypothesis is not used).  
Let $\Spec(R) \subseteq \Spec (E^0)$ be
an affine Zariski open so that $R \otimes_{E^0} E^{-2} \cong R\{u\}$ 
is a free $R$-module
of rank $1$.  Then there are isomorphisms
$$ R \otimes_{E^0} E^*(\CC P^\infty) \cong R \otimes_{E^0} E^*[[x]] \cong
E^*[[xu]] $$
where $x$ lies in degree $2$.  The coordinate $xu$ gives $R \otimes_{E^0}
E^0(\CC P^\infty)$ the structure of a formal group law.
\end{proof}

Let
$$ \omega_{\GG_E} = (\Lie(\GG_E))^* $$
be \index{1omegaG@$\omega_\GG$}
the line bundle over $\Spec(E^0)$ of invariant $1$-forms on $\GG_E$.
There is
a canonical isomorphism 
$$ \pi_{2t} E \cong \Gamma \omega_{\GG_E}^{\otimes t}. $$

Goerss and Hopkins
\cite{goersshopkins} extended work of Hopkins and Miller \cite{rezk} to
provide a partial converse to this construction.

\begin{thm}[Goerss-Hopkins-Miller]\label{thm:goersshopkinsmiller}
Let $k$ be a perfect field of characteristic $p$ 
and let $\br{\GG}$ be a formal group of finite height
over $k$.  Then there is an even periodic $E_\infty$-ring spectrum
$E_{\br{\GG}}$ such that
\index{2EG@$E_{\br{\GG}}$}
\begin{enumerate}
\item $E_{\br{\GG}}^0 = \mc{O}_{\Def_{\br{\GG}}}$, the coordinate ring of
the formal scheme $\Def_{\br{\GG}}$.
\item $\GG_E$ is a universal deformation of $\br{\GG}$.
\end{enumerate}
The construction is functorial in pairs $(k, \br{\GG})$: given another pair
$(k',\br{\GG}')$, a map of 
fields 
$i\co k' \rightarrow k$, and an isomorphism $\alpha\co 
\br{\GG} \xrightarrow{\cong} i^*\br{\GG}'$, there is an induced map of 
$E_\infty$-ring
spectra
$$ (i,\alpha)^*\co E_{\br{\GG}'} \rightarrow E_{\br{\GG}}. $$
\end{thm}

Our goal in this section is to state a generalization, recently
announced by Lurie, which gives functorial liftings of rings with
$1$-dimensional $p$-divisible groups.

\begin{thm}[Lurie]\label{thm:lurie}
Let $A$ be a local ring with maximal ideal $\mf{m}_A$ and residue field $k$ 
perfect of characteristic $p$, and suppose $X$
is a locally noetherian separated Deligne-Mumford stack over $\Spec(A)$.  
Suppose that $\mb G \to X$ is a
$p$-divisible group over $X$ of constant height $h$ and dimension $1$.
Suppose that for some \'etale cover $\pi: \td{X} \rightarrow X$ where $\td{X}$
is a scheme, the following condition is satisfied:
\begin{description}
\item[$(*)$]
for every point $x \in \td{X}^\wedge_{\mf{m}_A}$ the induced map
$$ \td{X}^\wedge_x \rightarrow \Def_{\pi^*\GG_x} $$
classifying the deformation $(\pi^* \GG) \big\vert_{\td{X}^\wedge_x}$ of 
$(\pi^*\GG)_x$ is an isomorphism of formal schemes.
\end{description}
Then there
exists a presheaf of $E_\infty$-ring spectra $\mc{E}_{\GG}$ 
\index{2EG@$\mc{E}_{\GG}$}
on the \'etale site of $X^\wedge_{\mf{m}_A}$, such that
\begin{enumerate}
\item $\mc{E}_{\GG}$ satisfies homotopy descent: it is (locally) fibrant in the
Jardine model structure \cite{jardine}, \cite{duggerhollanderisaksen}.
\item For every formal affine \'etale open $f\co \Spf(R) \rightarrow
X^\wedge_{\mf{m}_A}$,
the $E_\infty$-ring spectrum of sections $\mc{E}_{\GG}(R)$ is weakly even
periodic, with $\mc{E}_{\GG}(R)^0 = R$.
There is an isomorphism
$$ \gamma_f\co f^*\GG^0 \xrightarrow{\cong} \GG_{\mc{E}_{\GG}(R)} $$
natural in $f$.
\end{enumerate}
The construction is functorial in $(X,\GG)$: given another pair
$(X',\GG')$, a morphism $g\co X \rightarrow X'$ of Deligne-Mumford stacks 
over $\Spec(A)$, 
and an isomorphism of $p$-divisible groups
$\alpha\co \GG
\xrightarrow{\cong} g^* \GG'$, there is an induced map of presheaves of 
$E_\infty$-ring spectra
$$ (g,\alpha)^*\co g^*\mc{E}_{\GG'} \rightarrow \mc{E}_{\GG}. $$
\end{thm}

\begin{rmk}\label{rmk:lurie}
Suppose that $\GG$ and $X$ satisfy the hypotheses of
Theorem~\ref{thm:lurie}.
The functoriality of the sections of the presheaf $\mc{E}_{\GG}$ is
subsumed by the functoriality of the presheaf in $\GG$. 
Indeed, suppose that $g\co U \rightarrow X^\wedge_{\mf{M}_A}$ is 
an \'etale open.
\begin{enumerate}
\item The pullback $g^*\GG$ satisfies condition $(*)$ of 
Theorem~\ref{thm:lurie}, and hence gives rise to a presheaf
$\mc{E}_{g^*\GG}$ on $U$.

\item The spectrum of sections $\mc{E}_{\GG}(U)$ is given by 
the global sections $\mc{E}_{g^*\GG}(U)$.

\item The map 
$$ (g, \mathrm{Id})^* \co g^* \mc{E}_{\GG} \rightarrow \mc{E}_{g^*\GG} $$
induces a map
$$ g^*\co \mc{E}_{\GG}(X) \rightarrow \mc{E}_{g^*\GG}(U) = \mc{E}_{\GG}(U) $$
which agrees with the restriction map for the presheaf $\mc{E}_{\GG}$.
\end{enumerate}
\end{rmk}

Suppose that
$$ f: \Spf(R) \rightarrow X^\wedge_{\mf{m}_A} $$
is a formal affine \'etale open, and suppose that the invertible sheaf $f^*
\omega_{\GG^0}$ is trivial.  By Condition~(2) above, this is equivalent
to asserting that the spectrum of sections $\mc{E}_{\GG}(R)$ is even
periodic.  Then Condition~(2) actually determines the homotopy type of the
spectrum $\mc{E}_{\GG}(R)$.  Indeed, by choosing a coordinate $T$ of
the formal group $f^*\GG^0$ over  $\Spf(R)$, we get a map of rings
$$ \rho_{f^*\GG^0, T} \co MUP^0 \rightarrow R = \mc{E}_{\GG}(R)^0 $$
classifying the pair $(f^*\GG^0, T)$.  Here, $MUP$ 
\index{2MUP@$MUP$}
is the periodic
complex cobordism spectrum \cite[Ex.~8.5]{strickland}.  The ring $MUP^0$ is
the Lazard ring, and there is a canonical map
$$ \Spec(MUP^0) \rightarrow \mc{M}_{FG} $$
where $\mc{M}_{FG}$ is the moduli stack of formal groups.
\index{2MFG@$\mc{M}_{FG}$}

\begin{lem}\label{lem:flat}
The composite
$$ \Spec(R) \xrightarrow{\rho_{f^*\GG^0, T}} \Spec(MUP^0) \rightarrow
\mc{M}_{FG} $$ 
is flat.
\end{lem}

\begin{proof}
We check Landweber's criterion \cite{landweber}, \cite{naumannFG}.  
Let $I_{p,n}$ be the
ideal $(p, v_1, \ldots, v_{n-1})$ of $MUP^0$.  We must show that for each
$n$, the map
$$ \cdot v_n\co R /I_{p,n} \rightarrow R / I_{p,n} $$ is injective.
Suppose that $r$ is an element of $R/ I_{p,n}$.
Assuming that $r$ is non-zero, we need to verify that $v_n r$ is nonzero.
Since $f$ is \'etale and $X$ is locally noetherian, $R$ must be noetherian.
Therefore there is a 
maximal ideal $\mf{m}$ of $R$ so that the image of 
$r$ in $R^\wedge_{\mf{m}}/I_{p,n}$ is non-zero.  Let $x$ be the closed
point of $\Spf(R^\wedge_\mf{m})$, and consider the map
$$ \phi \co \Spf(R^\wedge_{\mf{m}}) \rightarrow \Def_{f^*\GG_x}
$$
classifying the deformation $f^*\GG \big\vert_{\Spf(R^\wedge_{\mf{m}})}$.
Since condition $(*)$ is an \'etale local condition, the map $\phi$ is an
isomorphism of formal schemes, and gives an 
isomorphism of $MUP^0$-algebras
$$ \phi^* \co 
W(k(x))[[u_1, \ldots, u_{h-1}]] \xrightarrow{\cong} R^\wedge_{\mf{m}}. $$
Let $h' \le h$ be the height of the formal group $\GG_{f(x)}^0$.
Lemma~6.10 of \cite{rezk} implies that multiplication by $v_n$ is injective
on the subalgebra $W(k(x))[[u_1, \ldots, u_{h'-1}]]/I_{p,n}$ classifying
deformations of the formal group $\GG_{f(x)}^0$, hence multiplication by
$v_n$ is injective on 
$$ W(k(x))[[u_1, \ldots, u_{h-1}]]/I_{p,n} = (W(k(x))[[u_1, \ldots,
u_{h'-1}]]/I_{p,n})[[u_{h'}, \ldots, u_{h-1}]]. $$
Thus the image of $v_n r$ in this ring is non-zero, implying
that $v_n r$ must be non-zero.
\end{proof}

We conclude that the theories $\mc{E}_{\GG}(R)$ are Landweber exact
\cite{landweber}, \cite{naumannFG}.

\begin{cor}\label{cor:LEFT}
There is an isomorphism of cohomology theories
$$ R \otimes_{MUP^0} MUP^*(X) \cong \mc{E}_{\GG}(R)^*(X). $$
\end{cor}

\begin{rmk}
Hovey and Strickland \cite[Prop.~2.20]{hoveystrickland} showed that for 
Landweber exact
cohomology theories such as 
$\mc{E}_{\GG}(R)$ which have no odd dimensional homotopy
groups,
Corollary~\ref{cor:LEFT} actually determines the homotopy type of these
spectra.
\end{rmk}

\begin{cor}\label{cor:GHML}
Suppose that $k$ is a perfect field of characteristic $p$ and that $\br{\GG}$
is a formal group of finite height over $k$.  Let $\GG/ \Def_{\br{\GG}}$ be the universal
deformation of $\br{\GG}$.  Then there is an
equivalence
$$ E_{\br{\GG}} \simeq \mc{E}_{\GG}(\Def_{\br{\GG}}). $$
\end{cor}

\begin{proof}
Both spectra are Landweber exact cohomology theories whose associated
formal groups are isomorphic.
\end{proof}

\section{The descent spectral sequence}

Let $X$ and $\GG$ be as in Theorem~\ref{thm:lurie}.  
Let $\omega_{\GG^0}$ be the line bundle of invariant $1$-forms on the formal
subgroup $\GG^0$.
In the previous section we have identified the homotopy groups of the 
spectrum of sections of $\mc{E}_{\GG}(U)$ for formal affine \'etale opens
$$ f\co U \rightarrow X^\wedge_{\mf{m}_A}. $$
Namely, for $k$ odd,
$\pi_k(\mc{E}_{\GG}(R)) = 0$ whereas
$$ \pi_{2t}(\mc{E}_{\GG}(U)) = H^0(U, \omega_{f^*\GG^0}^{\otimes t}). $$
For an arbitrary \'etale open
$$ f\co U \rightarrow X^\wedge_{\mf{m}_A} $$
we have the following.

\begin{thm}\label{thm:dss}
There is a conditionally convergent spectral sequence 
$$ E_2^{s,2t} = H_{zar}^s(U, \omega_{f^* \GG^0}^{\otimes t}) 
\Rightarrow \pi_{2t-s}(\mc{E}_{\GG}(U)). $$
\end{thm}

\begin{proof}
Let $U' \rightarrow U$ be a formal affine \'etale open.
Consider the cosimplicial object given by the Cech nerve:
$$ \mc{E}_{\GG}(U'^{\bullet+1}) = \left\{ \mc{E}_{\GG}(U') \Rightarrow 
\mc{E}_{\GG}(U' \times_U U') \Rrightarrow \cdots \right\}
$$
Because the sheaf $\mc{E}_{\GG}$ satisfies homotopy descent, the map
$$ \mc{E}_{\GG}(U) \rightarrow \holim_\Delta \mc{E}_{\GG}(U'^{\bullet+1}) $$
is an equivalence.
The spectral sequence of the theorem is the Bousfield-Kan spectral sequence
for computing the homotopy groups of the homotopy totalization.  
Since $U'$ is an affine formal scheme and $X$ is separated, 
the pullbacks $U' \times_U \cdots \times_U U'$ are affine
formal schemes. 
Therefore, 
the $E_2^{*,2t}$-term is the cohomology of the cosimplicial abelian group
$$ \omega_{f^*\GG}^{\otimes t}(U') \Rightarrow 
\omega_{f^*\GG}^{\otimes t}(U' \times_U U')
\Rrightarrow \cdots. $$
Since the coherent sheaf $\omega^{\otimes t}_{f^*\GG}$ has no higher cohomology when
restricted to $U' \times_U \cdots \times_U U'$, 
we deduce that the cohomology of this cosimplicial
abelian group computes the sheaf cohomology of $\omega_{f^*\GG}^{\otimes t}$ 
over $U$.
\end{proof}

We are able to deduce the following generalization of
Corollary~\ref{cor:LEFT}.

\begin{prop}\label{prop:LEFT}
Suppose that 
$$ f\co U \rightarrow X^\wedge_{\mf{m}_A} $$
is an \'etale open,
and that for every quasi-coherent sheaf $\mc{F}$ on $U$
$$ H^s_{zar}(U, \mc{F}) = 0 $$
for $s > 0$.  Then the spectrum $E = \mc{E}_\GG(U)$ is Landweber exact.
\end{prop}

\begin{proof}
Let $g \co U' \rightarrow U$ be a formal affine \'etale cover.  
By Lemma~\ref{lem:flat} the composite
$$ U' \xrightarrow{g} U \xrightarrow{f} X^\wedge_{\mf{m}_A} \rightarrow 
\mc{M}_{FG} $$
is flat.  Since $g$ is faithfully flat, we deduce that the composite
$$ h \co U \xrightarrow{f} X^\wedge_{\mf{m}_A} \rightarrow \mc{M}_{FG} $$
is flat.

By hypothesis the descent spectral sequence of Theorem~\ref{thm:dss}
collapses to give
$$ \pi_k(E) \cong
\begin{cases}
H^0(U, \omega^{\otimes t}), & k = 2t, \\
0, & \text{$k$ odd}.
\end{cases}
$$
Since $\pi_*E$ is concentrated in even degrees, $E$ is complex orientable.
Let $X$ be a spectrum, and let $\mc{MU}^{ev}(X)$ (respectively 
$\mc{MU}^{odd}(X)$) be the sheaf over
$\mc{M}_{FG}$ determined by the $MU_*MU$-comodule $MU_{2*}(X)$
(respectively $MU_{2*+1}(X)$).  We have
$$  
(MU_*(X) \otimes_{MU_*} E_*)_{k} =
\begin{cases}
H^0(U, h_{coh}^*(\mc{MU}^{ev} \otimes \omega^{t})), & k = 2t, \\
H^0(U, h_{coh}^*(\mc{MU}^{odd} \otimes \omega^{t})), & k = 2t+1. \\
\end{cases}
$$
Since $h_{coh}^*$ and $H^0_{coh}(U,-)$ are both exact functors, we deduce
that $MU_*(-) \otimes_{MU_*} E_*$ is a homology theory, and hence that
$E$ is Landweber exact.
\end{proof}

\section{Application to Shimura stacks}\label{sec:TAFdef}

Let $\Sh(K^p)$ be the Shimura stack over $\ZZ_p$ of Chapter~\ref{chap:shimura}.
Define a line bundle $\omega$ over $\Sh(K^p)$ as follows: for an $S$-point
of $\Sh(K^p)$ classifying a tuple $(A,i,\lambda,[\eta])$ over $S$, let 
$$ \omega = \omega_{\epsilon A(u)^0}, $$
the \index{1omega@$\omega$}
module of invariant $1$-forms on the formal part of the $1$-dimensional
summand $\epsilon A(u)$ of the $p$-divisible group $A(p)$.

Corollary~\ref{cor:defsh} implies that 
the $p$-divisible group $\epsilon A(u)$ on the Shimura stack $\Sh(K^p)$
satisfies condition $(*)$ of Theorem~\ref{thm:lurie}.  Thus, there is a
presheaf $\mc{E}(K^p)$
\index{2EKp@$\mc{E}(K^p)$}
of $E_\infty$-ring spectra on the \'etale site of
$\Sh(K^p)^\wedge_p$, satisfying the following.
\begin{enumerate}
\item For a formal affine \'etale open
$$ f\co \Spf(R) \rightarrow \Sh(K^p)^\wedge_p $$
classifying the tuple $(A,i,\lambda, [\eta])$ over $\Spf(R)$, the spectrum
of sections $\mc{E}(K^p)(R)$ is a weakly even periodic $E_\infty$-ring spectrum,
equipped with an isomorphism
$$ \GG_{\mc{E}(K^p)(R)} \cong \epsilon A(u)^0 $$
of its formal group with the formal subgroup $\epsilon A(u)^0$ of $A(p)$.
\item For an \'etale open $f\co U \rightarrow \Sh(K^p)^\wedge_p$, there is a
spectral sequence
$$ H^s(U, \omega^{\otimes t}) \Rightarrow
\pi_{2t-s}(\mc{E}(K^p)(U)). $$
\end{enumerate}

In particular, we define a spectrum of \emph{topological automorphic forms} 
\index{topological automorphic forms}
\index{2TAFKp@$\TAF(K^p)$}
by
$$ \TAF(K^p) = \mc{E}(K^p)(\Sh(K^p)^\wedge_p), $$
with descent spectral sequence
\begin{equation}\label{eq:TAFdescent}
H^s(\Sh(K^p)^\wedge_p, \omega^{\otimes t}) \Rightarrow \pi_{2t-s}(\TAF(K^p)).
\end{equation}

The following lemma follows immediately from Remark~\ref{rmk:lurie}.

\begin{lem}\label{lem:TAFsection}
Let $K'^p$ be an open subgroup of $K^p$.  
The spectrum $\TAF(K'^p)$ is given as the sections of the sheaf
$\mc{E}(K^p)$ on the \'etale cover
$$
f_{K'^p, K^p}\co \Sh(K'^p) \rightarrow \Sh(K^p) 
$$
of Theorem~\ref{thm:shimurastack}.
\end{lem}

\chapter{Relationship to automorphic forms}\label{chap:autform}

Assume that $K^p$ is sufficiently small, so that $\Sh(K^p)$ is a
quasi-projective scheme.
We briefly explain the connection with the classical theory of holomorphic
automorphic forms.  Namely, we will explain that there is a variant
$\Sh(K^p)_F$
\index{2ShKpF@$\Sh(K^p)_F$}
of the
Shimura variety defined over $F$, so that there is an isomorphism
$$ \Sh(K^p) \otimes_{\mc{O}_{F,u}} F_u \cong \Sh(K^p)_F \otimes_F F_u. $$
There is therefore a zig-zag of maps
$$
\xymatrix{
H^0(\Sh(K^p), \omega^{\otimes k}) \ar[r] & 
H^0(\Sh(K^p)_{F_u}, \omega^{\otimes k}) 
\\
& H^0(\Sh(K^p)_{F}, \omega^{\otimes k}) \ar[r] \ar[u] &
H^0(\Sh(K^p)_\CC, \omega^{\otimes k})
}
$$
where
\begin{align*}
\Sh(K^p)_{F_u} & = \Sh(K^p)_F \otimes_F F_u, \\
\Sh(K^p)_\CC & = \Sh(K^p)_F \otimes_F \CC. 
\end{align*}
We \index{2ShKpFu@$\Sh(K^p)_{F_u}$} \index{2ShKpC@$\Sh(K^p)_\CC$}
will explain that the sections
$$ H^0(\Sh(K^p)_\CC, \omega^{\otimes k}) $$
give instances of classical holomorphic automorphic forms for $GU$.

\begin{rmk}
There is a further variant $\Sh(K^p)_{(u)}$ 
\index{2ShKpu@$\Sh(K^p)_{(u)}$}
of the Shimura variety $\Sh(K^p)$
which is defined over the local ring $\mc{O}_{F,(u)}$ studied in
\cite{kottwitz} and \cite{hida}.  The varieties $\Sh(K^p)$ and 
$\Sh(K^p)_{F}$ are obtained from $\Sh(K^p)_{(u)}$ by base change to
$\mc{O}_{F,u}$ and $\mc{F}$, respectively  (see Remarks~\ref{rmk:det} and
\ref{rmk:tr}).
\end{rmk}

\section{Alternate description of $\Sh(K^p)$}\label{sec:ulocalSh}

For an abelian scheme $A$, let $\Lie A$ 
\index{2LieA@$\Lie A$}
denote the tangent sheaf at the
identity section.  For a $p$-divisible group $\GG$ over a locally $p$-nilpotent
base, let $\Lie{\GG}$ 
\index{2LieG@$\Lie \GG$}
be the
sheaf of invariant vector fields in the formal group $\GG^0$.

In Section~\ref{sec:moduli} we defined a functor
$$
\mc{X}_{K^p}'\co \left\{ 
\genfrac{}{}{0pt}{}{\text{locally noetherian schemes}}{\text{on which $p$ is 
locally nilpotent}} \right\} \rightarrow \{ \text{groupoids} \}.
$$
The functor of points of the Deligne-Mumford stack $\Sh(K^p)$ 
restricts to give the functor
$\mc{X}'_{K^p}$.  
In order to treat the characteristic zero points of the stack, we
must describe an extension 
$$ \mc{X}''_{K^p} \co \{ \text{locally noetherian schemes}/\Spec(\ZZ_p) \}
\rightarrow \{ \text{groupoids} \} $$
of the functor $\mc{X}'_{K^p}$.
\index{2X''Kp@$\mc{X}''_{K^p}$}

For a locally noetherian connected scheme $S$, the objects of the 
groupoid $\mc{X}''_{K^p}$ are tuples $(A,i,\lambda, [\eta])$ just as in the
definition of $\mc{X}'_{K^p}$ with one modification: the condition that
$\epsilon A(u)$ is $1$-dimensional must be replaced with the following
condition on $\Lie A$: 
\begin{description}
\item[$(*)$] The coherent sheaf $\Lie A \otimes_{\mc{O}_{F,p}} 
\mc{O}_{F, u}$ is locally free of rank $n$.
\end{description}
The tensor product in condition $(*)$ above is taken with respect to the
action of $\mc{O}_{F,p}$ on $\Lie A$ induced from the $B$-linear structure
$i$ on $A$.

\begin{lem}
There is a natural isomorphism of functors
$$ \mc{X}'_{K^p} \cong \mc{X}''_{K^p} $$
(where $\mc{X}''_{K^p}$ has been restricted to the domain of $\mc{X}'_{K^p}$).
\end{lem}

\begin{proof}
If $S$ is a locally noetherian connected scheme 
for which $p$ is locally nilpotent, 
we must show that a tuple $(A,i,\lambda, [\eta])$ is an object of
$\mc{X}'_{K^p}(S)$ if and only if it is an object of $\mc{X}''_{K^p}(S)$.
However, since $p$ is locally nilpotent on $S$, there is a natural
isomorphism
$$ \Lie A(p) \cong \Lie A $$
which induces a natural isomorphism
$$ \Lie A(u) \cong \Lie A \otimes_{\mc{O}_{F,p}} \mc{O}_{F,u}. $$
The coherent sheaf $\Lie A(u)$ is locally free of rank $n$ if and only if
the summand $\epsilon \Lie A(u)$ is locally free of rank $1$.  The latter
is equivalent to the formal part of the 
$p$-divisible group $\epsilon A(u)$ having dimension equal to $1$.
\end{proof}

\begin{rmk}\label{rmk:det}
In \cite{kottwitz} and \cite{hida}, a slightly different functor
$$ \mc{X}'''_{K^p}: \{ \text{locally noetherian schemes}/\Spec(\mc{O}_{F,(u)})
\}
\rightarrow \{ \text{groupoids} \} $$
is studied.  \index{2X'''Kp@$\mc{X}'''_{K^p}$}
The $S$-points are still given by tuples $(A,i, \lambda,
[\eta])$, but
condition $(*)$ is replaced with a ``determinant condition''
\cite[p390]{kottwitz}, \cite[p308]{hida}.  
\index{determinant condition}
When $S$ is
a scheme over $\mc{O}_{F,u}$, the determinant condition implies that there
is a local isomorphism of $\mc{O}_{B,p}$-linear
$\mc{O}_S$-modules
$$ \Lie A \cong \mc{O}_{F,u}^n \otimes_{\mc{O}_{F,u}} \mc{O}_S \oplus
(\mc{O}_{F,u^c}^n)^{n-1} \otimes_{c, \mc{O}_{F,u}} \mc{O}_S. $$
Here we are implicitly using our fixed isomorphism $\mc{O}_{B,p} 
\cong M_n(\mc{O}_{F,p})$, and
$\mc{O}_{F,u}^n$ (respectively $\mc{O}_{F,u^c}^n$) is the 
standard representation of the algebra $M_n(\mc{O}_{F,u})$ (respectively 
$M_n(\mc{O}_{F,u^c})$).  
The notation $\otimes_{c,\mc{O}_{F,u}}$ above indicates that we are
regarding $(\mc{O}^n_{F,u^c})^{n-1}$ as an $\mc{O}_{F,u}$-module through the
conjugation homomorphism.
Asserting the existence of such a local isomorphism
is equivalent to asserting that $\Lie A \otimes_{\mc{O}_{F,p}} \mc{O}_{F,u}$
is locally free of rank $n$.  The functor $\mc{X}'''_{K^p}$ is proven in
\cite[Sec.~5,6]{kottwitz} to be representable by a quasi-projective variety
$\Sh(K^p)_{(u)}$ 
\index{2ShKpu@$\Sh(K^p)_{(u)}$}
when $K^p$ is sufficiently small.  For general $K^p$, the
representing object $\Sh(K^p)_{(u)}$ 
still exists as a Deligne-Mumford stack \cite[Sec~7.1.2]{hida}.
The Shimura stack $\Sh(K^p)$ which represents the 
functors $\mc{X}'_{K^p}$ and $\mc{X}''_{K^p}$ is obtained from $\Sh(K^p)_{(u)}$
by defining
$$ \Sh(K^p) := \Sh(K^p)_{(u)} \otimes_{\mc{O}_{F,(u)}} \mc{O}_{F,u}. $$
\end{rmk}

\section{Description of $\Sh(K^p)_F$}

We now describe the stack $\Sh(K^p)_F$
\index{2ShKpF@$\Sh(K^p)_F$}
over $F$.
Let $S$ be a locally noetherian connected scheme over $F$.  
An object of the groupoid of $S$-points $\Sh(K^p)_{F}(S)$
consists of a tuple $(A,i,\lambda, [\eta])$ just as in the definition of
the functor $\mc{X}'_{K^p}$ except that the condition that
$\epsilon A(u)$ is $1$-dimensional must be replaced with a suitable
condition on the locally free sheaf $\Lie A$ on $S$.  
The field $F$ acts on $\Lie A$ 
in two different ways:
both through the complex multiplication $i$ of $F$ on $A$, and through
structure of $S$ as a scheme over $F$.
This induces a splitting of the sheaf $\Lie A$:
$$ \Lie A \cong \Lie A^+ \oplus \Lie A^- $$
where \index{2LieA+@$\Lie A^+$} \index{2LieA-@$\Lie A^-$}
$\Lie A^+$ is the subsheaf where the $F$-linear structures agree, and
$\Lie A^-$ is the subsheaf where the $F$-linear structures are conjugate.  
We require:
\begin{description}
\item[$(*)$] the sheaf $\Lie A^+$ must be locally free of rank $n$.
\end{description}

\begin{lem}
There is an equivalence of stacks
$$ \Sh(K^p) \otimes_{\mc{O}_{F,u}} F_u \simeq \Sh(K^p)_{F} \otimes_F F_u. $$
\end{lem}

\begin{proof}
For an $S$-object $(A,i, \lambda, [\eta])$ of either stack,
the sheaf $\Lie A^+$ is naturally isomorphic to $\Lie A
\otimes_{\mc{O}_{F,p}} \mc{O}_{F,u}$.
\end{proof}

\begin{rmk}\label{rmk:tr}
When working in characteristic $0$, the determinant condition of
Remark~\ref{rmk:det} can be replaced by a ``trace condition''.  The points
of the resulting stack $\Sh(K^p)_F = \Sh(K^p)_{(u)} \otimes_{\mc{O}_{F,(u)}}
F$ are precisely those tuples $(A,i,\lambda, [\eta])$ which satisfy
condition $(*)$ above (see \cite[Lemma~III.1.2]{harristaylor}).
\end{rmk}

\section{Description of $\Sh(K^p)_\CC$}\label{sec:complexpoints}

Let $K^p$ be sufficiently small so that $\Sh(K^p)_F$ is a scheme.
Fix an embedding
$$ j: F \hookrightarrow \CC. $$
We \index{2J@$j$}
shall explain how the pullback 
$$ \Sh(K^p)_{\CC} = \Sh(K^p)_{F} \times_{\Spec(F)} \Spec(\CC) $$
admits \index{2ShKpC@$\Sh(K^p)_\CC$}
a classical description (Theorem~\ref{thm:shimuraadelic}).
The data of Section~\ref{sec:data} allows us to construct explicit 
points of $\Sh(K^p)_\CC$.
Let $W$ 
\index{2W@$W$}
be the $2n^2$ dimensional real vector space given by 
$V_\infty = V \otimes_\QQ \RR$ with lattice $L$.  A complex structure $J$ 
\index{2J@$J$}
\index{complex structure}
is an
element of $GL(W)$ satisfying $J^2 = -1$.  For each such $J$, 
the real vector space
$W$ admits the structure of a complex vector space, which we denote
$W_J$.  
\index{2WJ@$W_J$}
Associated to $J$ is a polarized abelian variety $A_J$
\index{2AJ@$A_J$}
given by the quotient
$$
A_J = W_J/L.
$$
The dual abelian variety is given by
$$ A^\vee_J = W_J^*/L^* $$
where 
\begin{align*}
W_J^* & = \{ \alpha\co W_J \rightarrow \CC \: : \: 
\text{$\alpha$ conjugate linear} \},
\\
L^* & = \{ \alpha \in W_J^* \: : \: \im \alpha(L) \subseteq \ZZ
\}.
\end{align*}
\index{2W*@$W^*$}
\index{2L*@$L^*$}

Associated to the non-degenerate alternating form $\bra{-,-}$ on $V_\infty$ 
is a non-degenerate hermitian form
$(-,-)_J$
\index{0broJ@$(-,-)_J$}
on $W_J$, given by the correspondence of
Lemma~\ref{lem:alternatinghermitian}:
$$ (v,w)_J = \bra{Jv, w} + i\bra{v,w}. $$
We get an induced complex-linear isomorphism:
\begin{align*}
\td{\lambda}_J\co W_J & \rightarrow W_J^*, \\
v & \mapsto (v,-)_J. 
\end{align*}
The integrality conditions imposed on $\bra{-,-}$ with respect to the lattice 
$L$ (Section~\ref{sec:data}) 
imply that
$\td{\lambda}_J$ 
\index{1lambdaJ@$\td{\lambda}_J$}
descends to a prime-to-$p$ isogeny of abelian varieties:
$$ \lambda_J\co A_J \rightarrow A_J^\vee. $$
In \index{1lambdaJ@$\lambda_J$}
order for the isogeny $\lambda_J$ to be a polarization, the alternating
form $\bra{-,-}$ must be a \emph{Riemann form}: 
\index{Riemann form}
\index{bilinear form!Riemann form}
we require that the form
$\bra{-,J-}$ be symmetric and positive.

The $B$-action on $V$ gives rise to an embedding
$$ B \hookrightarrow \Aut_\RR (W). $$
In order for $B$ to act by complex linear maps, we must have
$$ J \in C_\infty \subset \Aut_\RR(W). $$
Here, $C_\infty$ 
\index{2Cinfty@$C_\infty$}
is the algebra $C \otimes_\QQ \RR$, where $C$ is the
algebra of Section~\ref{sec:data}.
Since $\mc{O}_B$ acts on $L$, we get an inclusion
$$ i\co \mc{O}_B \hookrightarrow \End(A_J). $$
Let \index{2I@$i$}
$(-)^\iota$ be the involution on $C$ of Section~\ref{sec:data}.  Then
in order for the involution $*$ on $B$ to agree with the $\lambda_J$-Rosati
involution, we must have $J^\iota = -J$.

The $B$-linear polarized abelian variety $(A_J, i, \lambda_J)$ admits a
canonical integral uniformization $\eta_1$ by the composite
\begin{align*}
\eta_1\co L^p & = \varprojlim_{(N,p) = 1} L/NL \\
& \cong \varprojlim_{(N,p) = 1} N^{-1}L/L \\
& = \varprojlim_{(N,p) = 1} A_J[N] \\
& = T^p(A_J).
\end{align*}
For \index{1eta1@$\eta_1$}
$g \in GU(\AF^{p,\infty})$, we let $\eta_g$ 
\index{1etag@$\eta_g$}
denote the associated
rational uniformization
$$
\eta_g\co V^p \xrightarrow{g} V^p \xrightarrow{\eta_1} V^p(A_J).
$$

Let $W^+_J \subseteq W_J$ 
\index{2WJ+@$W_J^+$}
denote the subspace where the complex structure $J$ agrees with
the complex structure $W_J$ inherits from the $F$-linear structure of $W$
through our chosen isomorphism $j: F \otimes_\QQ \RR \xrightarrow{\cong}
\CC$.  
Then requiring that $\Lie A^+$ is $n$-dimensional amounts to
requiring that $W_J^+$ is $n$-dimensional.

Let $\mc{H}$
\index{2H@$\mc{H}$}
be the set of \emph{compatible complex structures}:
\index{complex structure!compatible}
$$ \mc{H} = \left\{ J \in C_\infty \: : \: 
\begin{array}{l}
J^2 = -1, \\
\bra{-,J-} \: \text{is symmetric and positive}. \\
\end{array}
\right\}
$$

\begin{lem}
Suppose that $J$ is a compatible complex structure.  Then we have:
\begin{gather*}
J^\iota = -J, \\
\dim_\CC W_J^+ = n.
\end{gather*}
\end{lem}

\begin{proof}
Using the fact that $\bra{-,J-}$ is symmetric, we have
\begin{align*}
\bra{x, Jy} & = \bra{y, Jx} \\
& = \bra{J^\iota y, x} \\
& = -\bra{x, J^\iota y}.
\end{align*}
Since $\bra{-,-}$ is non-degenerate, it follows that $J^\iota = -J$.

Our assumptions on the alternating form $\bra{-,-}$ imply that 
the associated $*$-symmetric form
$$ (x,y) = \bra{\delta x, y} + \delta\bra{x,y} $$
on $W$ has signature $(n,(n-1)n)$.  Here, $F = \QQ(\delta)$, where, 
under our fixed complex embedding, $\delta$ is equal to $-ai$ for $a \in
\RR_{>0}$.
Suppose that $x$ is a non-zero element of $W_J^+$.  Then we have
\begin{align*}
(x,x) & = \bra{\delta x,x} + \delta\bra{x,x} \\
& = \bra{-ai x, x} \\
& = a\bra{x, ix } \\
& = a\bra{x, Jx} > 0,
\end{align*}
since we have assumed that $\bra{-,J-}$ is positive definite.  We deduce,
using the signature of $(-,-)$, that $\dim W_J^+ \le n$.  However, the
same argument shows that $(-,-)$ is negative definite when restricted
to $W_J^-$, which implies that $\dim W_J^- \le (n-1)n$.  Since
$$ \dim W_J^+ + \dim W_J^- = n^2, $$
we deduce that $\dim W_J^+ = n$.
\end{proof}

We deduce the following lemma.

\begin{lem}
For $J \in \mc{H}$, the tuple $(A_J, i, \lambda_J, [\eta_g]_{K^p})$ 
gives a $\CC$-point
of $\Sh(K^p)_\CC$.
\end{lem}

The group
$$ GU(\RR)^+ = \{ g \in GU(\RR) \: : \: g^\iota g > 0 \} $$
acts \index{2GUR+@$GU(\RR)^+$}
on $\mc{H}$ by conjugation.  Under the identification $U(\RR) \cong
U(1,n-1)$ of Section~\ref{sec:data}, we have the following
lemma.

\begin{lem}
The action of $U(\RR)$ on $\mc{H}$ is transitive, and the stabilizer of
some fixed $J_0
\in \mc{H}$ is given by a maximal compact subgroup 
$$ K_\infty = U(W_J^+) \times U(W_J^-) \cong U(1) \times U(n-1). $$
Thus \index{2Kinfty@$K_\infty$}
there is an isomorphism $\mc{H} \cong U(\RR)/K_\infty$.
\end{lem}

\begin{rmk}
The symmetric hermitian domain $\mc{H} = U(\RR)/K_\infty$ 
admits a canonical complex
structure \cite[VIII]{helgason}.
\end{rmk}

Define a group
$$ GU(\ZZ_{(p)})^+ = \{ g \in GU(\QQ)^+ \: : \: g(L_{(p)}) = L_{(p)} \}. $$
Every \index{2GUZp+@$GU(\ZZ_{(p)})^+$}
element $\gamma \in GU(\ZZ_{(p)})^+$ gives rise to a $\ZZ_{(p)}$-isogeny
$$ \gamma\co (A_{J}, i, \lambda_{J}, [\eta_g]_{K^p}) 
\xrightarrow{\cong} (A_{\gamma J \gamma^{-1}}, i, \lambda_{\gamma J
\gamma^{-1}}, [\eta_{\gamma g}]_{K^p}). $$
We have the following theorem \cite[Sec.~8]{kottwitz}.

\begin{thm}\label{thm:shimuraadelic}
The map
$$ GU(\ZZ_{(p)})^+ \backslash (GU(\AF^{p,\infty})/K^p \times \mc{H})
\rightarrow \Sh(K^p)_\CC $$
sending a point $[g,J]$ to the point that classifies $(A_J, i, \lambda_J,
[\eta_g]_{K^p})$
is an isomorphism of complex analytic manifolds.
\end{thm}

\begin{rmk}
The statement of Theorem~\ref{thm:shimuraadelic} differs mildly that of
\cite{kottwitz}, where it is proven that $\Sh(K^p)_\CC$ is isomorphic to
$$ \coprod_{\ker^1 (\QQ, GU)} GU(\QQ)^+ 
\backslash (GU(\AF^\infty)/K^pK_p \times \mc{H}) $$
where 
\begin{align*}
\ker^1(\QQ, GU) & = \ker(H^1(\QQ, GU) \rightarrow \bigoplus_v H^1(\QQ_v,
GU)), \\
K_p & = \{ g \in GU(\QQ_p) \: : \: g(L_p) = L_p\}.
\end{align*}
In our case, the $\ker^1$-term is trivial
(Corollary~\ref{cor:globalH^1GU}).  Moreover, we may remove the
$p$-component
from the adelic quotient by combining Lemma~7.2 of \cite{kottwitz} with the
fact that $H^1(\QQ_p, GU)$ is trivial (Lemma~\ref{lem:localsplitH^1GU}).
\end{rmk}

\section{Automorphic forms}

Let $T \cong (S^1)^n$ 
\index{2T@$T$}
be the maximal
torus of $K_\infty = U(1) \times U(n-1)$, and let $X \cong \ZZ^n$ be 
\index{2X@$X$}
the lattice of characters
of $T$.  
Let $V_\kappa$ 
\index{2Vkappa@$V_\kappa$}
be the finite dimensional irreducible complex representation 
of the group
$K_\infty$ of highest weight $\kappa \in X$.
\index{1kappa@$\kappa$}
The collection of all highest weights of
irreducible representations forms a cone 
$$ X^+ = \ZZ \times \ZZ \times \mb{N}^{n-2} \subseteq X$$ 
(provided $n > 1$).  
Here, the first factor of $\ZZ$ corresponds to the powers of the
standard representation on $U(1)$ while the second
factor of $\ZZ$ corresponds to the powers of the determinant representation
on $U(n-1)$.  The factor $\mb{N}^{n-2}$ is the space of
weights of irreducible representations of $SU(n-1)$.

Associated to $V_\kappa$ is a holomorphic 
vector bundle $\mc{V}_\kappa$ on $\Sh(K^p)_\CC$, 
given by
the Borel construction.
\begin{gather*}
\mc{V}_\kappa = GU(\ZZ_{(p)})^+ \backslash (GU(\AF^{p,\infty})/K^p \times U(\RR)
\times_{K_\infty} V_\kappa) \\
\qquad \downarrow \\
\Sh(K^p)_\CC = GU(\ZZ_{(p)})^+ \backslash (GU(\AF^{p,\infty})/K^p \times U(\RR)
/K_\infty)
\end{gather*}

By a \emph{weakly holomorphic automorphic form} 
\index{automorphic form!weakly holomorphic}
of weight 
$\kappa$
\index{automorphic form!weight}
we shall mean a holomorphic
section of $\mc{V}_\kappa$.  A \emph{holomorphic automorphic form} 
\index{automorphic form!holomorphic}
of
weight $\kappa$ is a weakly holomorphic automorphic form
which satisfies certain growth conditions \cite{borel}.  If
$\Sh(K^p)_\CC$ is compact (i.e. if $B$ is a division algebra --- see
Theorem~\ref{thm:shimurastack}) then these growth conditions are always
satisfied, and every weakly holomorphic automorphic form is holomorphic.

Let $\omega_{F}$ 
\index{1omegaF@$\omega_F$}
denote the line bundle $(\epsilon \Lie A^+)^*$
on $\Sh(K^p)_{F}$.  Then
$\omega_{F}$ restricts to $\omega$ on $\Sh(K^p)_{F_u}$.  The pullback
$\omega_\CC$ 
\index{1omegaC@$\omega_\CC$}
of $\omega_{F}$ to $\Sh(K^p)_\CC$ is the line bundle whose fiber
over 
$$ [g,J] \in GU(\ZZ_{(p)})^+ \backslash 
(GU(\AF^{p,\infty})/K^p \times \mc{H}) \cong \Sh(K^p)_\CC $$
is given by
$$ (\omega_\CC)_{[g,J]} = (\Lie A_J^+)^* = (W_J^+)^*. $$
We therefore have the following lemma.

\begin{lem}
There is an isomorphism of vector bundles $\omega_\CC^{\otimes k} \cong
\mc{V}_{-\ul{k}}$, where $-\ul{k}$ 
\index{2K@$\ul{k}$}
is the weight
$$ (-k, \underbrace{0, \ldots, 0}_{n-1}) \in X^+. $$
\end{lem}

We therefore have a zig-zag
$$
\xymatrix{
H^0(\Sh(K^p), \omega^{\otimes k}) \ar[r] & 
H^0(\Sh(K^p)_{F_u}, \omega^{\otimes k})
\\
& H^0(\Sh(K^p)_{F}, \omega_F^{\otimes k}) \ar[u] \ar[r] & 
H^0(\Sh(K^p)_\CC, \omega_\CC^{\otimes k})  
}
$$
relating the $E_2$-term of the descent spectral sequence for $\TAF_{K^p}$
(\ref{eq:TAFdescent}) with the space of weakly holomorphic automorphic
forms of weight $-\ul{k}$.

\chapter{Smooth $G$-spectra}\label{chap:smGspt}

In this chapter we introduce the notion of a \emph{smooth $G$-spectrum} for
a totally disconnected locally compact group $G$.  
\index{2G@$G$}
We shall see in the next
chapter that the collection of spectra $\{ \TAF(K^p) \}$, as 
$K^p$ varies over the compact open subgroups of $G = GU(\AF^{p,\infty})$, 
gives rise to a smooth $G$-spectrum $V_{GU}$.  

\section{Smooth $G$-sets}

Recall that a $G$-set $X$ is \emph{smooth} 
\index{smooth!$G$-set}
if the stabilizer of every $x \in
X$ is open.  Such a $G$-set $X$ then satisfies
$$ X = \varinjlim_{H \le_{o} G} X^H $$
where \index{0leo@$\le_o$}
$H$ runs over the open subgroups of $G$.
Let $\Set^{sm}_{G}$ be the category whose objects are the smooth
$G$-sets and
\index{2SetsmG@$\Set^{sm}_{G}$}
whose morphisms are $G$-equivariant maps.

Fix an infinite cardinality $\alpha$ 
\index{1alpha@$\alpha$}
larger than the cardinality of $G$.
Let $\Set^{sm,\alpha}_G$ 
\index{2SetsmGalpha@$\Set^{sm,\alpha}_G$}
be the (essentially) \emph{small} 
subcategory of $\Set_G^{sm}$ consisting of
those smooth $G$-sets which have cardinality less than $\alpha$.
We may endow $\Set^{sm,\alpha}_G$ with
the structure of a Grothendieck topology by declaring that surjections are
covers.  The only nontrivial thing to verify is that $\Set^{sm,\alpha}_G$
contains pullbacks, but this is accomplished by the following lemma.

\begin{lem}\label{lem:Gsetpullback}
Suppose that $H_1$ and $H_2$ are open subgroups of $H$, an open 
subgroup of $G$.  Then the following diagram is a pullback diagram of
smooth $G$-sets.
$$
\xymatrix{
\coprod\limits_{h \in H_1 \backslash H / H_2} G/(H_1 \cap hH_2h^{-1}) 
\ar[r]^-r \ar[d]
& G/H_2 \ar[d] \\
G/H_1 \ar[r] 
& G/H
}
$$
Here all of the maps are the evident surjections except for the map
$r$, which is given by the following formula:
$$ r(g(H_1 \cap hH_2h^{-1})) = ghH_2. $$
\end{lem}

The site $\Set^{sm,\alpha}_G$ has enough points, with one unique
point given by the filtering system $\{ G/H \}_{H \le_{o} G}$.

Every smooth $G$-set is a
disjoint union of smooth $G$-orbits, and (since $G$-orbits are small) the
natural map
\begin{equation}\label{eq:mapsmGset}
\prod_i \coprod_j \Map_G(G/H_i, G/H_j) \rightarrow \Map_G\left(\coprod_i
G/H_i, \coprod_j G/H_j \right)
\end{equation}
is an isomorphism.
Therefore, to describe the morphisms in $\Set^{sm}_G$ 
it suffices to describe the
morphisms between smooth $G$-orbits. 
Let $H$ and $H'$ be open subgroups of $G$, and consider the subset
$$ {}_{H}N_{H'} = \{ g \in G \: : \: g^{-1}Hg \le H' \}. $$
We \index{2NHH'@${}_{H}N_{H'}$}
have an isomorphism
\begin{equation}\label{eq:mapGorbit}
\Map_G(G/H, G/H') \cong (H \backslash {}_{H}N_{H'} / H')^{op} 
\end{equation}
where for $g \in {}_HN_{H'}$, the double coset $HgH'$ corresponds to the map
$$ R_g\co xH \mapsto xgH'. $$

Suppose that $\mc{X}$ 
\index{2X@$\mc{X}$}
is a sheaf on the site $\Set^{sm,\alpha}_{G}$.  Then the
colimit
$$ X = \varinjlim_{H \le_o G} \mc{X}(G/H) $$
taken \index{2X@$X$}
over open subgroups $H$ is a smooth $G$-set with $H$-fixed
points
$$ X^H = \mc{X}(G/H). $$
Conversely, any smooth $G$-set $X$ gives rise to a sheaf $\mc{X}$ on
$\Set^{sm,\alpha}_{G}$ whose sections on a $G$-set 
$S = \amalg_i G/H_i$ are given
by
$$ \mc{X}(S) = \prod_i X^{H_i}. $$
These constructions give the following lemma.

\begin{lem}\label{lem:Gset}
The category of smooth $G$-sets is equivalent to the category of sheaves of
sets on
the site $\Set^{sm,\alpha}_G$.
\end{lem}

There is an adjoint pair $(U, (-)^{sm})$ of functors
$$
U : \Set^{sm}_G \leftrightarrows \Set_G : (-)^{sm}
$$
between the category of smooth $G$-sets and the category of all $G$-sets,
where
\begin{align*}
U(X) & = X, \\
X^{sm} & = \varinjlim_{H \le_o G} X^H.
\end{align*}
The category of smooth $G$-sets is complete and cocomplete, with colimits
formed in the underlying category of sets, and limits formed by applying
the functor $(-)^{sm}$ to the underlying limit of sets.  Finite limits are
formed in the underlying category of sets.

\begin{rmk}
The properties listed above make the adjoint pair $(U, (-)^{sm})$ a
geometric morphism of topoi.
\end{rmk}

\section{The category of simplicial smooth $G$-sets}\label{sec:smGspace}

It is convenient to construct model categories of discrete $G$-spaces and 
discrete $G$-spectra for
a profinite group $G$ by considering the category of
simplicial (pre)-sheaves of sets  
on the site of finite discrete $G$-spaces (see, for instance
\cite{jardineetale}).  An alternative, and equivalent approach, was
proposed by Goerss \cite{goerssgalois}, who considered simplicial objects in the
category of discrete $G$-sets.  In this section we modify these
frameworks to the case of $G$ a totally disconnected locally compact group.
We investigate these model structures in some detail in preparation for
some technical applications in
Chapter~\ref{chap:Qspectra}.

\begin{defn}
A \emph{simplicial smooth $G$-set} is a simplicial object in the category
$\Set^{sm}_G$.
\end{defn}
\index{smooth!$G$-set!simplicial}

Consider the following categories:
\begin{align*}
s\Set_G^{sm} & = \text{simplicial smooth $G$-sets,} \\
s\Pre(\Set_G^{sm,\alpha}) & = \text{simplicial presheaves of 
sets on $\Set_G^{sm, \alpha}$,} \\
s\Shv(\Set_G^{sm,\alpha}) & = \text{simplicial sheaves of sets on
$\Set^{sm, \alpha}_G$.}
\end{align*}
\index{2SSetGsm@$s\Set_G^{sm}$}
\index{2SPreSetGsmalpha@$s\Pre(\Set_G^{sm,\alpha})$}
\index{2SShvSetGsmalpha@$s\Shv(\Set_G^{sm,\alpha})$}

Consider the diagram of functors
\begin{equation}
\xymatrix{
s\Pre(\Set_G^{sm, \alpha}) 
\ar@<.5ex>[r]^-{\mc{L}^2} \ar@<-.5ex>@{<-}[r]_-{U} &
s\Shv(\Set_G^{sm, \alpha}) \ar@<.5ex>[r]^-{L} \ar@<-.5ex>@{<-}[r]_-{R} & 
s\Set_G^{sm}
}
\end{equation}
where
\begin{align*}
\mc{L}^2 & = \text{sheafification}, \\
U & = \text{forgetful functor}, \\
L(\mc{X}) & = \varinjlim_{H \le_{o} G} \mc{X}(G/H), \\
R(X)(\amalg_i G/K_i) & = \prod_i X^{K_i}.
\end{align*}
The \index{2L2@$\mc{L}^2$} \index{2U@$U$} \index{2L@$L$} \index{2R@$R$}
pairs $(\mc{L}^2, U)$ and $(L,R)$ are adjoint pairs of functors.
Lemma~\ref{lem:Gset} has the following corollary.

\begin{cor}\label{cor:Gspace}
The adjoint pair $(L,R)$ gives an adjoint equivalence of categories.
\end{cor}

The categories of 
simplicial presheaves and sheaves of sets 
admit injective local model structures 
such that the adjoint pair 
$(\mc{L}^2, U)$ forms a Quillen equivalence \cite{jardinepresheaf}.
The category $s\Set_G^{sm}$ inherits a Quillen equivalent model structure
under the adjoint equivalence $(L,R)$. We shall call this induced model
structure on $s\Set_G^{sm}$ the \emph{injective local model structure}.
\index{injective local model structure}

\begin{lem}\label{lem:smmodelspace}
A map in $s\Set_G^{sm}$ is a weak equivalence (cofibration) in the
injective local model structure 
if and only if it induces a weak
equivalence (monomorphism) on the underlying simplicial set.
\end{lem}

\begin{proof}
The local weak equivalences in $s\Pre(\Set^{sm,\alpha}_G)$ and 
$s\Shv(\Set^{sm, \alpha}_G)$ are the stalkwise
weak equivalences, and the functor $L$ gives the stalk.    
A cofibration of simplicial presheaves is a sectionwise monomorphism. The 
result for cofibrations follows immediately from the definition of
the functor $L$, and the fact that for a simplicial sheaf $\mc{X}$, and 
open subgroups $U \le V \le G$, the
induced map
$$ \mc{X}(G/V) \rightarrow \mc{X}(G/U) $$
is a monomorphism.
\end{proof}

The category of simplicial sheaves on $\Set_G^{sm,\alpha}$ also admits a
\emph{projective local model structure} \cite{blander}.  
\index{projective local model structure}
In this model
structure, the weak equivalences are still the local (stalkwise) weak
equivalences, but the projective cofibrations are generated by the 
inclusions of the
representable simplicial sheaves: the generating cofibrations are morphisms
of the form
$$ i_{n,Z}: 
\partial \Delta^n \times \Map_G(-, Z) \hookrightarrow \Delta^n \times
\Map_G(-, Z) $$
for \index{2InZ@$i_{n,Z}$}
$Z \in \Set_G^{sm,\alpha}$.  The projective local fibrations are
determined.  The category of
simplicial smooth $G$-sets inherits a projective local model structure from
the adjoint equivalence $(L,R)$.

\begin{lem}
A map in $s\Set_G^{sm}$ is a projective cofibration if and only if it is a
monomorphism.
\end{lem}

\begin{proof}
The class of projective cofibrations is generated by the monomorphisms
$$ Li_{n,Z}: Z \times \partial \Delta^n \hookrightarrow Z \times \Delta^n
$$
for $Z \in \Set_G^{sm,\alpha}$.  It clearly suffices to consider the
generating morphisms $Li_{n,G/H}$ for open subgroups $H \le G$.
Because the maps $Li_{n,G/H}$ are monomorphisms, every projective
cofibration is a monomorphism.
Conversely, given an inclusion $j : X \hookrightarrow Y$ of simplicial
smooth $G$-sets, we have
$$ Y = \varinjlim_k Y^{[k]} $$
where $Y^{[-1]} = X$ and $Y^{[k]}$ is given by the pushout
$$
\xymatrix{
\coprod_{NY_k - NX_k} \partial \Delta^n \ar[r] \ar[d] 
& Y^{[k-1]} \ar[d] \\
\coprod_{NY_k - NX_k} \Delta^n \ar[r] 
& Y^{[k]}
}
$$
where $NX_k$ and $NY_k$ are the (smooth $G$-sets of) nondegenerate
$k$-simplices of $X$ and $Y$, respectively.
This pushout may be rewritten in terms of $G$-orbits of non-degenerate
$k$-simplices as follows.
$$
\xymatrix{
\coprod_{[\sigma] \in (NY_k - NX_k)/G} G/H_\sigma \times 
\partial \Delta^n \ar[r] \ar[d]_{\amalg Li_{n,G/H_\sigma}} 
& Y^{[k-1]} \ar[d] \\
\coprod_{[\sigma] \in (NY_k - NX_k)/G} G/H_\sigma \times \Delta^n \ar[r] 
& Y^{[k]}
}
$$
Here, for a representative $\sigma$ of an $G$-orbit $[\sigma]$, the group
$H_\sigma$ is the (open) stabilizer of $\sigma$.  We have shown that the
class of monomorphisms is contained in the category of cofibrations generated
by the morphisms $Li_{n,Z}$.
\end{proof}

We have shown that the cofibrations and weak equivalences are the same in the 
injective and projective local model category structures on $s\Set_G^{sm}$,
giving the following surprising corollary.

\begin{cor}
The injective and projective local model category structures on
$s\Set_G^{sm}$ agree.
\end{cor}

It is shown in \cite{blander} 
that a map between simplicial sheaves is a
projective local fibration if and only if it is a projective local
fibration on the level of underlying presheaves.  In
\cite[Thm~1.3]{duggerhollanderisaksen}, an explicit characterization of the 
projective fibrant presheaves is given.  This fibrancy condition, when
translated into $s\Set_G^{sm}$ using the functor $L$, gives the following
characterization of the fibrant objects of $s\Set_G^{sm}$.

\begin{lem}\label{lem:fibrantsmGspace}
A simplicial smooth $G$-set $X$ is fibrant if and only if:
\begin{enumerate}
\item the $H$-fixed points $X^H$ is a Kan complex for every open subgroup
$H \le G$,
\item for every open subgroup $H \le G$ and every hypercover 
$$ \cdots \coprod_{\alpha \in I_2} \td{G/U}_{\alpha,2} \Rrightarrow 
\coprod_{\alpha \in I_1} \td{G/U}_{\alpha,1} \Rightarrow \coprod_{\alpha
\in I_0} \td{G/U}_{\alpha,0} \rightarrow \td{G/H} $$
(where $\td{G/U}$ denotes the representable sheaf associated to a smooth
orbit $G/U$) the induced map
$$ X^H \rightarrow \holim_{\Delta} 
\prod_{\alpha \in I_\bullet} X^{U_{\alpha,\bullet}} $$
is a weak equivalence.
\end{enumerate}
\end{lem}

\begin{rmk}
Since filtered colimits of Kan complexes are Kan complexes, the underlying
simplicial set of a fibrant simplicial smooth $G$-set is a Kan complex.
\end{rmk}

\section{The category of smooth $G$-spectra}\label{sec:smGsp}

The methods of Section~\ref{sec:smGspace} extend in a straightforward manner to give
model categories of smooth $G$-spectra, as investigated in the case of $G$
profinite by \cite{jardineetale} and \cite{davis}.  These smooth
$G$-spectra are intended to generalize 
``naive equivariant stable homotopy theory'': the fibrant objects will not
distinguish between fixed points and homotopy fixed points.

\begin{defn}
A \emph{smooth $G$-spectrum} 
\index{smooth!$G$-spectrum}
is a Bousfield-Friedlander spectrum of smooth
$G$-sets.  It consists of a sequence $\{ X_i \}_{i \ge 0}$ of pointed 
simplicial smooth $G$-sets together with $G$-equivariant structure maps
$$ \Sigma X_i \rightarrow X_{i+1}. $$
\end{defn}

\begin{rmk}
If $G$ is a profinite group, then the notion of a smooth $G$-spectrum
coincides with that of a discrete $G$-spectrum studied in \cite{davis}.
\end{rmk}

Consider the following categories (where \emph{spectrum} means 
Bousfield-Friedlander spectrum of simplicial sets):
\begin{align*}
\Sp_G^{sm} & = \text{category of smooth $G$-spectra,} \\
\Pre\Sp(\Set_G^{sm,\alpha}) & = \text{category of presheaves of 
spectra on $\Set_G^{sm, \alpha}$,} \\
\Shv\Sp(\Set_G^{sm,\alpha}) & = \text{category of sheaves of spectra on
$\Set^{sm, \alpha}_G$.}
\end{align*}
\index{2SpGsm@$\Sp_G^{sm}$}
\index{2PreSpSetGsmalpha@$\Pre\Sp(\Set_G^{sm,\alpha})$}
\index{2ShvSpSetGsmalpha@$\Shv\Sp(\Set_G^{sm,\alpha})$}

Just as in Section~\ref{sec:smGspace}, we have a diagram of functors
\begin{equation}
\xymatrix{
\Pre\Sp(\Set_G^{sm, \alpha}) 
\ar@<.5ex>[r]^-{\mc{L}^2} \ar@<-.5ex>@{<-}[r]_-{U} &
\Shv\Sp(\Set_G^{sm, \alpha}) \ar@<.5ex>[r]^-{L} \ar@<-.5ex>@{<-}[r]_-{R} & 
\Sp_G^{sm}
}
\end{equation}
where
\begin{align*}
\mc{L}^2 & = \text{sheafification}, \\
U & = \text{forgetful functor}, \\
L(\mc{X}) & = \varinjlim_{H \le_{o} G} \mc{X}(G/H), \\
R(X)(\amalg_i G/K_i) & = \prod_i X^{K_i}.
\end{align*}
The \index{2L2@$\mc{L}^2$} \index{2U@$U$} \index{2L@$L$} \index{2R@$R$}
pairs $(\mc{L}^2, U)$ and $(L,R)$ are adjoint pairs of functors,
and the adjoint pair $(L,R)$ gives an adjoint equivalence of categories.

The categories of 
presheaves and sheaves of spectra admit injective local model structures
\index{injective local model structure}
such that the adjoint pair 
$(\mc{L}^2, U)$ forms a Quillen equivalence \cite{goerssjardine}.
The category $\Sp_G^{sm}$ inherits a Quillen equivalent model structure
under the adjoint equivalence $(L,R)$.

\begin{lem}\label{lem:smmodel}
A map in $\Sp_G^{sm}$ is a weak equivalence (cofibration) 
if and only if it induces a stable
equivalence (cofibration) on the underlying spectrum.
\end{lem}

\begin{proof}
The weak equivalences in $\Pre\Sp(\Set^{sm,\alpha}_G)$ and 
$\Shv\Sp(\Set^{sm, \alpha}_G)$ are the stalkwise
stable equivalences.  
For a sheaf of spectra $\mc{E}$, the underlying
spectrum of $L\mc{E}$ is the stalk of $\mc{E}$ at the (unique) point of the
site $\Set_G^{sm,\alpha}$.
An argument analogous to that given in \cite[Thm~3.6]{davis} proves the
statement concerning cofibrations.
\end{proof}

\section{Smooth homotopy fixed points}

For a smooth $G$-spectrum $X$, we define the fixed point spectrum by
taking the fixed points level-wise
$$ (X^G)_i = (X_i)^G. $$
The \index{2XG@$X^G$}
$G$-fixed point functor is right adjoint to the functor ${\it triv}$
\index{2Triv@$\mathit{triv}$}
which
associates to a spectrum the associated smooth $G$-spectrum with trivial
$G$-action.
$$ {\it triv}\co \Sp \leftrightarrows \Sp_G \co (-)^G. $$
Since the functor ${\it triv}$ preserves cofibrations and weak
equivalences, we have the following lemma.

\begin{lem}\label{lem:Quillenfixedpoint}
The adjoint functors $({\it triv}, (-)^G)$ form a Quillen pair.
\end{lem}

Let $\beta_{G,X}\co X \rightarrow X_{fG}$ 
\index{2XfG@$X_{fG}$}
denote a functorial 
fibrant replacement functor
for the model category $\Sp^{sm}_G$, where $\beta_{G,X}$ is a trivial
cofibration of smooth $G$-spectra.
The homotopy fixed point functor $(-)^{hG}$
\index{smooth!homotopy fixed points}
is the Quillen right derived functor of $(-)^G$, and is thus given by
$$ X^{hG} = (X_{fG})^G. $$
\index{2XhG@$X^{hG}$}

\section{Restriction, induction, and coinduction}

Fix an open subgroup $H$ of $G$.
For a smooth $H$-set $Z$, we define the coinduced smooth $G$-set by
$$ \CoInd_H^G Z = \Map_H(G,Z)^{sm} := \varinjlim_{U \le_o G} \Map_H(G/U, Z)
$$
where \index{2CoInd@$\CoInd$} \index{2Mapsm@$\Map_G(-,-)^{sm}$} \index{coinduction}
the colimit is taken over open subgroups.  An element $g \in G$ sends
an element $\alpha \in \Map_H(G/U,Z)$ to an element $g\cdot \alpha \in
\Map_H(G/gUg^{-1}, Z)$ by the formula
$$ (g \cdot \alpha)(g'(gUg^{-1})) = \alpha(g'gU). $$
This construction extends to simplicial smooth $G$-sets and smooth
$G$-spectra in the obvious manner to give a functor
$$ \CoInd_H^G \co \Sp^{sm}_H \rightarrow \Sp^{sm}_G. $$
Let $\Res_H^G$ be the restriction functor:
$$ \Res_H^G \co \Sp_G^{sm} \rightarrow \Sp_H^{sm}. $$
The \index{2Res@$\Res$} \index{restriction}
functors $(\Res_H^G, \CoInd_H^G)$ form an adjoint pair.
Since $\Res_H^G$ preserves cofibrations and weak equivalences, we have the
following lemma.

\begin{lem}
The adjoint functors $(\Res_H^G, \CoInd_H^G)$ form a Quillen pair.
\end{lem}

The Quillen pair $(\Res_H^G, \CoInd_H^G)$ gives rise to a derived 
pair $(L\Res_H^G, R\CoInd_H^G)$.  Since the functor $\Res_H^G$ preserves
all weak equivalences, there are equivalences $L \Res_H^G X \simeq \Res_H^G
X$ for all smooth $G$-spectra $X$.

\begin{lem}[Shapiro's Lemma]\label{lem:Shapiro}
Let $X$ be a smooth $H$-spectrum, and suppose that $H$ is an open
subgroup of $G$.  Then there is an equivalence
$$ (R\CoInd_H^G X)^{hG} \xrightarrow{\simeq} X^{hH}. $$
\end{lem}

\begin{proof}
The lemma follows immediately from the following commutative diagram of
functors.
$$
\xymatrix{
\Sp_H^{sm} \ar[rr]^{\CoInd_H^G} \ar[dr]_{(-)^H} && 
\Sp_G^{sm} \ar[dl]^{(-)^G} 
\\
& \Sp
}
$$
\end{proof}

Define the 
induction
functor 
on a smooth $H$-spectrum $Y$ to be the Borel construction
$$ \Ind_H^G Y = G_+ \wedge_H Y. $$
Here, \index{induction} \index{2Ind@$\Ind$}
the Borel construction is taken regarding $G$ and $H$ as being
discrete groups, but this is easily seen to produce a smooth $G$-spectrum
since $H$ is an open subgroup of $G$.
The induction and restriction functors form an adjoint pair 
$(\Ind_H^G, \Res_H^G)$
$$ \Ind_H^G\co \Sp^{sm}_{H} \leftrightarrows \Sp^{sm}_{G} \co \Res_H^G. $$
Since non-equivariantly we have an isomorphism
$$ \Ind_H^G Y \cong G/H_+ \wedge Y, $$
we have the following lemma.

\begin{lem}\label{lem:resfibrant}
The functor $\Ind_H^G$ preserves cofibrations and weak equivalences, and
therefore $\Res_H^G$ preserves fibrant objects.
\end{lem}

\begin{defn}
Let $X$ be a smooth $G$-spectrum.
Define the smooth homotopy $H$-fixed points of $X$ by 
$$ X^{hH} := (\Res_H^G X)^{hH} $$
\end{defn}

Lemma~\ref{lem:resfibrant} has the following corollary.

\begin{cor}\label{cor:resfibrant}
For a smooth $G$-spectrum $X$, the spectrum $X^{hH}$ is equivalent to the
$H$-fixed points $(X_{fG})^H$.
\end{cor}

\section{Descent from compact open subgroups}\label{sec:descent}

In this section we will explain how the homotopy type of a smooth
$G$-spectrum can be reconstructed using only its homotopy fixed points for
\emph{compact} open subgroups of $G$ (Construction~\ref{const:V}).  This
construction will be used in Section~\ref{sec:GUaction} to construct a
smooth $GU(\AF^{p,\infty})$-spectrum $V_{GU}$ from the collection of
spectra $\{ \TAF(K^p) \}$, where $K^p$ ranges over the compact open
subgroups of $GU(\AF^{p,\infty})$. 

Let $\Set_G^{co,\alpha}$ 
\index{2SetGcoalpha@$\Set_G^{co,\alpha}$}
be the full subcategory of $\Set_G^{sm, \alpha}$
consisting of objects $X$ whose stabilizers are all compact.  Every object
$Z$ of $\Set_G^{co, \alpha}$ is therefore isomorphic to a disjoint union of
$G$-orbits
$$ Z = \coprod_i G/K_i $$
where the subgroups $K_i$ are compact and open.
Lemma~\ref{lem:Gsetpullback} implies that $\Set^{co, \alpha}$ is a
Grothendieck site, with covers given by surjections.  The site
$\Set_G^{co, \alpha}$ also possesses enough points, with a unique point
given by the filtering system $\{ G/K \}_{K}$ where $K$ ranges over the
compact open subgroups of $G$.

Just as in Section~\ref{sec:smGsp}, there is an 
adjoint equivalence $(L',R')$
$$ L' \co \Shv\Sp(\Set_G^{co,\alpha}) \leftrightarrows \Sp_G^{sm} \co R' $$
where 
\begin{align*}
L'\mc{X} & = \varinjlim_{K \le_{co} G} \mc{X}(G/K), \\
(R' X)(\amalg_i G/K_i) & = \prod_i X^{K_i}.
\end{align*}
The \index{2L'@$L'$} \index{2R'@$R'$}
colimit in the definition of $L'$ is taken over compact open subgroups
of $G$.

The adjoint equivalence $(L',R')$ induces a model structure on $\Sp^{sm}_G$
from the Jardine model structure on $\Shv\Sp(\Set_G^{co,\alpha})$.
Precisely the same arguments proving Lemma~\ref{lem:smmodel} apply to prove
the following lemma.

\begin{lem}
In the model structure on $\Sp_G^{sm}$ induced from the Jardine model
structure on $\Shv\Sp(\Set_G^{co,\alpha})$, a map is a weak equivalence
(cofibration) if and only if it is a weak equivalence (cofibration) on the
underlying spectrum.
\end{lem}

\begin{cor}
The model structure on $\Sp_G^{sm}$ induced from the Jardine model
structure on $\Shv\Sp(\Set_G^{co,\alpha})$ is identical to that induced
from the the model structure on $\Shv\Sp(\Set^{sm,\alpha}_G)$.
\end{cor}

\begin{const}\label{const:V}
Begin with a presheaf of spectra 
$\mc{X}$
on the site $\Set_G^{co,\alpha}$.  We associate to $\mc{X}$ a 
smooth $G$-spectrum 
$$ V_{\mc{X}} = L'(\mc{L}^2\mc{X})_f, $$
the \index{2VX@$V_{\mc{X}}$}
functor $L'$ applied to the fibrant replacement of the sheafification
of $\mc{X}$. 
\end{const}

\begin{lem}
The smooth $G$-spectrum $V_{\mc{X}}$ is fibrant.
\end{lem}

\begin{proof}
The fibrant objects in $\Sp_G^{sm}$ are precisely those objects $X$ 
for which $R'X$ is
fibrant in $\Shv\Sp(\Set_G^{co,\alpha})$.  Since $(L',R')$ form an equivalence of
categories, there is an isomorphism
$$ R'V_{\mc{X}} = R'L'(\mc{L}^2\mc{X})_f \cong (\mc{L}^2\mc{X})_f. $$
\end{proof}

\begin{lem}\label{lem:cofixedpoints}
Suppose that $\mc{X}$ is a fibrant presheaf of spectra on
$\Set^{co,\alpha}_G$, and that $K$ is a compact open subgroup of $G$.  
Then there
is an equivalence of spectra
$$ \mc{X}(G/K) \xrightarrow{\simeq} V_{\mc{X}}^K. $$ 
\end{lem}

\begin{proof}
Consider the Quillen equivalence given by the adjoint pair
$$ \mc{L}^2 \co \Pre\Sp(\Set_G^{co,\alpha}) \leftrightarrows
\Shv\Sp(\Set_G^{co,\alpha})\co U $$
where $\mc{L}^2$ is the sheafification functor and $U$ is the forgetful
functor.
\index{2L2@$\mc{L}^2$}
\index{2U@$U$}
The functor $U$ preserves all weak equivalences.

Define $\phi$ to be the composite
$$ \phi\co \mc{X} \xrightarrow{\eta} U\mc{L}^2 \mc{X} \xrightarrow{U(\alpha)} U
(\mc{L}^2 \mc{X})_f $$
in the category $\Pre\Sp(\Set_G^{co,\alpha})$,
where $\eta$ is the unit of the adjunction and
$\alpha$ is the fibrant replacement morphism for the sheaf $\mc{L}^2
\mc{X}$.  The map $\eta$ is always a weak equivalence, and, since $U$
preserves all weak equivalences, the map $U(\alpha)$ is a weak equivalence.
Thus $\phi$ is a weak equivalence of presheaves of spectra.

The functor $U$, being the right adjoint of a Quillen pair, preserves
fibrant objects, and thus $\phi$ is a weak equivalence between fibrant
objects.  The sections functor
\begin{align*} 
\Gamma_{G/K}\co \Pre\Sp(\Set_G^{co,\alpha}) & \rightarrow \Sp \\
\mc{E} & \mapsto \mc{E}(G/K) 
\end{align*}
is also the right adjoint of a Quillen pair, and hence preserves weak
equivalences between fibrant objects.  We therefore have a string of
equivalences
\begin{align*}
\mc{X}(G/K) 
\xrightarrow[\simeq]{\Gamma_{G/K}\phi} & \: U(\mc{L}^2 \mc{X})_f(G/K) \\
= & \: (\mc{L}^2 \mc{X})_f(G/K) \\
\cong & \: (R'L'\mc{L}^2 \mc{X})_f)(G/K) \\
= & \: V_{\mc{X}}^K.
\end{align*}
\end{proof}

\section{Transfer maps and the Burnside category}\label{sec:hecke}

In this section we will construct transfer maps between the homotopy fixed
points of compact open subgroups.  We will then explain how 
this extra structure assembles to give this collection of homotopy fixed
point spectra
the structure of a Mackey functor.

\subsection*{Construction of transfers}

We will now explain how to use Shapiro's Lemma to construct transfer
maps for smooth $G$-spectra.  

Suppose $K$ is a compact open subgroup of $G$ with open subgroup
$K'$.  There is a continuous based map $K_+ \to K'_+$, equivariant with
respect to the left and right $K'$-actions, given by
\[
k \mapsto 
\begin{cases}
k & \text{if } k \in K', \\
* & \text{otherwise.}
\end{cases}
\]

Given a smooth $K'$-spectrum $X$, this gives rise to a map of $K'$-spectra
\[
X \cong \Map_{K',*}(K'_+, X)^{sm} \to \Map_{K',*}(K_+,X)^{sm} \cong \Res^K_{K'}
\CoInd^K_{K'} X.
\]
The adjoint of this map is a natural map of $K$-spectra
\[
i^K_{K'}\co \Ind^K_{K'} X \to \CoInd^K_{K'} X.
\]
On the level of underlying spectra, this map is an inclusion
\[
\bigvee_{K/K'} X \to \prod_{K/K'} X,
\]
and hence is a natural weak equivalence of smooth $K$-spectra.

Consider the counit map
\[
\varepsilon^K_{K'}\co \Ind^K_{K'} \Res^K_{K'} X \to X.
\]
Application of the right derived functor of smooth coinduction gives a
diagram of maps
\[
R \CoInd^G_{K'} X \overfrom^{\sim} R \CoInd^G_K \Ind^K_{K'} X
\to R \CoInd^G_K X
\]
of smooth $G$-spectra.  Upon taking homotopy fixed points,
Lemma~\ref{lem:Shapiro} gives rise to natural maps
\[
X^{hK'} \overfrom^{\sim} (\Ind^K_{K'} \Res^G_{K'} X)^{hK} \to X^{hK}.
\]

\begin{defn}
The \emph{transfer} 
\index{transfer}
associated to the inclusion of $K'$ in $K$ is
the corresponding map 
\index{2Tr@$\Tr$}
$\Tr^K_{K'}\co X^{hK'} \to X^{hK}$ in the homotopy
category.
\end{defn}

\subsection*{The Burnside category}

The treatment in this section follows \cite{minami}.
Define the Burnside category
\index{Burnside category}
${\mc{M}}^{co}_G$
\index{2McoG@${\mc{M}}^{co}_G$}
of $G$-sets to be the a category with objects given by
$$ 
\mathrm{Ob}(\mc{M}_G^{co}) 
= \{ Z \in
\Set^{co,\alpha}_G \: : \: |Z/G| < \infty \}. 
$$
Thus the objects of $\mc{M}^{co}_G$ are $G$-sets with finitely many orbits and
compact open stabilizers.

We now define the morphisms.
Let $S,T$ be objects of $\mc{M}^{co}_G$.  Let $\Cor(S,T)$ 
\index{2CorST@$\Cor(S,T)$}
denote the set of
correspondences given by
isomorphism classes of diagrams in $\Set^{co}_G$ of the form
$$
\xymatrix{
& X \ar[dl]_p \ar[dr]^q \\
S && T.
}
$$
The set $\Cor(S,T)$ admits the structure of a commutative monoid through
disjoint unions:
$$
\xymatrix{
& X \ar[dl]_p \ar[dr]^q \\
S && T
}
\begin{array}{c}
\\ \\
+ 
\end{array}
\xymatrix{
& X' \ar[dl]_{p'} \ar[dr]^{q'} \\
S && T
}
\begin{array}{c}
\\ \\
=
\end{array}
\xymatrix{
& X \amalg X' \ar[dl]_{p \amalg p'} \ar[dr]^{q \amalg q'} \\
S && T
}
$$
The morphisms of the category $\mc{M}_G^{co}$ from $S$ to $T$ 
are defined to be the Grothendieck group
of the monoid $\Cor(S,T)$:
$$ \Hom_{\mc{M}^{co}_G}(S,T) = \mathrm{Groth}(\Cor(S,T)). $$

Composition in $\mc{M}_G^{co}$ is given by:
$$
\xymatrix{
& Y \ar[dl] \ar[dr] \\
T && U
}
\begin{array}{c}
\\ \\
\circ
\end{array}
\xymatrix{
& X \ar[dl] \ar[dr] \\
S && T
}
\begin{array}{c}
\\ \\
=
\end{array}
\xymatrix{
&& Z \ar[dl] \ar[dr] \\
& X \ar[dl] \ar[dr] && Y \ar[dl] \ar[dr] \\
S && T && U
}
$$
where $Z$ is the pullback $X \times_T Y$.

The morphisms of $\mc{M}^{co}_G$ 
are generated by the right multiplication
by $g$ maps 
$$ \rho_g\co G/gKg^{-1} \to G/K, $$
the \index{1rhog@$\rho_g$}
projections
$$ r^K_{K'}\co G/K' \to G/K, $$
and \index{2RKK'@$r^K_{K'}$}
the transfers 
$$ i^K_{K'}\co G/K \to
G/K' $$ 
for \index{2IKK'@$i^K_{K'}$}
$g \in G$, $K' < K$, satisfying the following relations.
\begin{enumerate}
\item $\rho_g \circ \rho_h = \rho_{hg}$
\item $r^{K}_{K'} \circ r^{K'}_{K''} = r^{K}_{K''}$
\item $i^{K'}_{K''} \circ i^{K}_{K'} = i^{K}_{K''}$
\item $\rho_g \circ r^{gKg^{-1}}_{gK'g^{-1}} = r^K_{K'} \circ \rho_g$
\item $\rho_g \circ i^{gKg^{-1}}_{gK'g^{-1}} = i^K_{K'} \circ \rho_g$
\item $i^H_L \circ r^H_K = \sum_{KxL} r^L_{x^{-1} K x \cap L} \rho_x
    \circ i^K_{K \cap xLx^{-1}}$
\end{enumerate}

\begin{defn}
A \emph{Mackey functor} (for the group $G$) is an additive functor
$$ \mc{M}^{co}_G \rightarrow \text{Abelian groups}$$
which takes coproducts to direct sums.
\index{Mackey functor}
\end{defn}

\subsection*{The Mackey functor associated to a smooth $G$-spectrum}

For a smooth $G$-spectrum $X$ we get an induced functor
$$ \mc{X}: (\mc{M}^{co}_G)^{op} \rightarrow \mathrm{Ho}(\Sp) $$
from the Burnside category to the homotopy category of spectra as follows.
On objects, we define
$$ \mc{X}(\amalg_i G/K_i) = \prod_i X^{hK_i}. $$
To a correspondence of the form
\begin{equation}\label{eq:specialcorr}
\xymatrix{
& G/K'' \ar[dl]_{r_{K''}^{K}} \ar[dr]^{r^{K'}_{g^{-1}K''g} \circ \rho_g} \\
G/K && G/K'
}
\end{equation}
we assign the composite
$$ X^{hK'} \xrightarrow{\Res^{K'}_{g^{-1}K''g}}
X^{hg^{-1}K''g} \xrightarrow{g_*} X^{hK''} \xrightarrow{\Tr^K_{K''}}
X^{hK}.
$$
Every correspondence between $G/K$ and $G/K'$ is a disjoint union of
correspondences of the form (\ref{eq:specialcorr}), and we assign to it the
corresponding sum of morphisms in the homotopy category of spectra.  Thus
we have defined a map of commutative monoids
$$ \Cor(G/K, G/K') \rightarrow \Hom_{\mathrm{Ho}(\Sp)}(X^{hK'}, X^{hK}) $$
which extends to a homomorphism
$$ \Hom_{\mc{M}^{co}_G}(G/K, G/K') 
\rightarrow \Hom_{\mathrm{Ho}(\Sp)}(X^{hK'}, X^{hK}). $$
Extending additively we have homomorphisms
$$ \Hom_{\mc{M}^{co}_G}(S, T) 
\rightarrow \Hom_{\mathrm{Ho}(\Sp)}(\mc{X}(T), \mc{X}(S)) $$
for every pair of objects $S,T \in \mc{M}^{co}_G$.

In view of the definition of composition in $\mc{M}^{co}_G$, the following
lemma implies the functoriality of $\mc{X}$.  The lemma is easily verified
using the definition of transfer.

\begin{lem}
Suppose that $H_1$ and $H_2$ are open subgroups of $H$, a compact open 
subgroup of $G$.  Consider the following pullback diagram of
Lemma~\ref{lem:Gsetpullback}.
$$
\xymatrix@C+6em{
\coprod\limits_{h \in H_1 \backslash H / H_2} G/(H_1 \cap hH_2h^{-1}) 
\ar[r]^-{\coprod r^{H_2}_{h^{-1}H_1h \cap H_2} \circ\rho_h} \ar[d]
& G/H_2 \ar[d] \\
G/H_1 \ar[r] 
& G/H
}
$$
Then the following diagram commutes in $\mathrm{Ho}(\Sp)$.
$$
\xymatrix@C+6em{
\prod\limits_{h \in H_1 \backslash H / H_2} X^{hH_1 \cap hH_2h^{-1}} 
\ar@{<-}[r]^-{\sum h_* \circ \Res^{H_2}_{h^{-1}H_1h \cap H_2}} \ar[d]_{\sum
\Tr^{H_1}_{H_1 \cap hH_2h^{-1}}}
& G/H_2 \ar[d]^{\Tr^H_{H_2}} \\
G/H_1 \ar@{<-}[r]_{\Res^{H}_{H_1}} 
& G/H
}
$$
\end{lem}

The previous
relations imply the following.

\begin{prop}
The conjugation, restriction, and transfer maps satisfy the following
relations.
\begin{enumerate}
\item $h_* \circ g_* = (hg)_*$
\item $\Res^{K'}_{K''} \circ \Res^K_{K'} = \Res^{K}_{K''}$
\item $\Tr^K_{K'} \circ \Tr^{K'}_{K''} = \Tr^{K}_{K''}$
\item $\Res^{gKg^{-1}}_{gK'g^{-1}} \circ g_* = g_* \circ \Res^K_{K'}$
\item $\Tr^{gKg^{-1}}_{gK'g^{-1}} \circ g_* = g_* \circ \Tr^K_{K'}$
\item $\Res^H_K \circ \Tr^H_L = \sum_{KxL} \Tr^K_{K \cap xLx^{-1}}
  \circ x_* \circ \Res^L_{x^{-1} K x \cap L}$
\end{enumerate}
\end{prop}

\begin{cor}
The functor
\begin{align*}
\pi_* \mc{X}(-) \co  \mc{M}^{co}_G & \rightarrow \text{Abelian groups} \\
\amalg_i G/K_i & \mapsto \bigoplus_i \pi_*(X^{hK_i}) 
\end{align*}
is a Mackey functor.
\end{cor}

For any $G$-module $M$, the fixed points form a Mackey functor.
\begin{align*}
M^{(-)} \co (\mc{M}_G^{co})^{op} & \rightarrow \text{Abelian Groups} \\
\amalg_i G/K_i & \mapsto \bigoplus_i M^{K_i}
\end{align*}
The transfers
$$ \Tr_{K'}^K : M^{K'} \rightarrow M^K $$
are given by
$$ \Tr_{K'}^K(m) = \sum_{[g] \in K/K'} gm. $$

Consider the canonical ``edge homomorphisms''
$$ \epsilon_K \co \pi_*(X^{hK}) \rightarrow \pi_*(X)^K. $$
The \index{1epsilonkappa@$\epsilon_\kappa$}
following lemma is easily verified by comparing the definitions of
transfers.

\begin{lem}\label{lem:epsilon}
The homomorphisms $\epsilon_{(-)}$ assemble to give a natural
transformation of Mackey functors
$$ \epsilon_{(-)} \co \pi_*(X^{h(-)}) \rightarrow \pi_*(X)^{(-)}. $$
\end{lem}

\chapter{Operations on $\TAF$}\label{chap:ops}

In this chapter we will describe three classes of cohomology operations,
which extend classical operations in the theory of automorphic
forms.
\begin{description}
\item[(1) Restriction] For compact open subgroups $K'^p \le K^p$, a map of
$E_\infty$-ring spectra
$$ \Res^{K^p}_{K'^p} \co \TAF(K^p) \rightarrow \TAF(K'^p). $$

\item[(2) Action of $GU(\AF^{p,\infty})$]
For $g \in GU(\AF^{p,\infty})$ and $K^p$ as above, 
a map of $E_\infty$-ring spectra
$$ g_*\co \TAF(K^p) \rightarrow \TAF(gK^pg^{-1}). $$

\item[(3) Transfer] For $K^p$ and $K'^p$ as above, a map in the homotopy
category
spectra
$$ \Tr_{K'^p}^{K^p}\co \TAF(K'^p) \rightarrow \TAF(K^p). $$
\end{description}

Operations of type $(1)$ and $(2)$ may be most elegantly formulated as a
functor on an orbit category, as we describe in
Section~\ref{sec:GUaction}.  This structure gives rise to a smooth
$GU(\AF^{p,\infty})$-spectrum $V_{GU}$.  
The operations of type $(3)$ then come automatically (see
Section~\ref{sec:hecke}).  All three types of operations encode the fact
that $\TAF(-)$ induces a Mackey functor from the Burnside category into the 
homotopy category of spectra.  

The relationship of this Mackey functor structure to the
theory of Hecke algebras is discussed in Section~\ref{sec:heckeTAF}.  
While we are able to lift Hecke operators on automorphic forms to
cohomology operations on $\TAF(K^p)$, we are unable to determine in general
if this induces an action of the corresponding Hecke algebra (we do treat a
very restrictive case in Propositions~\ref{prop:cohomologicalmackey} and 
\ref{prop:heckealgebra}).  We are \emph{not} advocating that having the
Hecke algebra act on $\TAF(K)$ is necessarily 
the ``right'' point of view.  The most
important structure is the Mackey functor structure.

\section{The $E_\infty$-action of
$GU(\AF^{p,\infty})$}\label{sec:GUaction}

Let $\Set^{co,\alpha}_{GU(\AF^{p,\infty})}$ be the Grothendieck site of
Section~\ref{sec:descent}.

\begin{prop}\label{prop:Einftyops}
The assignment 
$$ \coprod_i GU(\AF^{p,\infty})/K^p_i \mapsto \prod_i \TAF(K^p_i) $$
gives a presheaf 
$$ \TAF(-) \co (\Set^{co,\alpha}_{GU(\AF^{p,\infty})})^{op} \rightarrow 
\text{$E_\infty$-ring spectra}. $$
\end{prop}

\begin{proof}
In light of the isomorphism~(\ref{eq:mapsmGset}), 
it suffices to establish functoriality on orbits in
$\Set^{co,\alpha}_{GU(\AF^{p,\infty})}$.
Let $\Orb^{co}_{GU(\AF^{p,\infty})}$ 
\index{2Orbco@$\Orb^{co}$}
be the subcategory of orbits. 
Consider the functor
\begin{align*}
\Orb^{co}_{GU(\AF^{p,\infty})} & \rightarrow \{\text{$p$-divisible
groups}\}, \\
K^p & \mapsto (\Sh(K^p), \epsilon \mathbf{A}(u)).
\end{align*}
Here, $(\mathbf{A}, \mathbf{i},
\pmb{\lambda}, [\pmb{\eta}]_{K^p})$ is the universal tuple on $\Sh(K^p)$.

Using the isomorphism (\ref{eq:mapGorbit}),
functoriality of this correspondence is established as follows: for open
compact subgroups $g^{-1}K'^pg \le K^p$, the map
$$ g_{K'^p, K^p}\co \Sh(K'^p) \rightarrow \Sh(K^p) $$
classifies the tuple $(\mathbf{A}', \mathbf{i}',
\pmb{\lambda}', [\pmb{\eta}'g]_{K^p})$, where 
$(\mathbf{A}', \mathbf{i}',
\pmb{\lambda}', [\pmb{\eta}']_{K'^p})$ is the universal tuple on
$\Sh(K'^p)$.  In particular, there is a canonical isomorphism of
$p$-divisible groups
$$ \alpha_{can}\co \epsilon\mathbf{A}'(u) \xrightarrow{\cong} (g_{K'^p,
K^p})^* \epsilon\mathbf{A}(u)  $$
giving a map
$$ (g_{K'^p,K^p}, \alpha_{can})\co 
(\Sh(K'^p), \epsilon \mathbf{A}'(u)) \rightarrow
(\Sh(K^p), \epsilon \mathbf{A}(u)) $$ 
in the category of $p$-divisible groups.  The functoriality
Theorem~\ref{thm:lurie} gives an induced map of presheaves of
$E_\infty$-ring spectra
$$ (g_{K'^p, K^p}, \alpha_{can})^* \co (g_{K'^p, K^p})^*
\mc{E}(K^p) \rightarrow \mc{E}(K'^p) $$
and hence a map 
$$ g_*\co \TAF(K^p) \rightarrow \TAF(K'^p). $$
\end{proof}

\begin{lem}
The presheaf $\TAF(-)$ on the site $\Set_{GU(\AF^{p,\infty})}^{co, \alpha}$
is fibrant.
\end{lem}

\begin{proof}
Fibrancy is shown to be a local condition in \cite{duggerhollanderisaksen},
and hence follows from Lemma~\ref{lem:TAFsection} and the fibrancy of the
presheaves $\mc{E}(K^p)$ of Section~\ref{sec:TAFdef}.
\end{proof}

Construction~\ref{const:V} associates to the presheaf $\TAF(-)$ a fibrant smooth
$GU(\AF^{p,\infty})$-spectrum
$$ V_{GU} = V_{\TAF(-)}. $$
By \index{2VGU@$V_{GU}$}
Lemma~\ref{lem:cofixedpoints}, for any compact open 
subgroup $K^p < GU(\AF^{p,\infty})$ we have an equivalence
$$ \TAF(K^p) \xrightarrow{\simeq} V_{GU}^{hK^p}. $$

\section{Hecke operators}\label{sec:heckeTAF}

Let $G = GU(\AF^{p,\infty})$.
\index{2G@$G$}
Let $\mc{M} = \mc{M}^{co}_{G}$ 
\index{2M@$\mc{M}$}
denote the Burnside
category.  In Section~\ref{sec:hecke} we observed that the homotopy fixed
points of a smooth $G$-spectrum assembled to yield a Mackey functor in the
stable homotopy category.  In
particular, since the spectra $\TAF(-)$ are given as homotopy fixed point
spectra, our Mackey functor takes the form:
\begin{align*}
\mc{M}^{op} & \rightarrow \mathrm{Ho}(\Sp) \\
\amalg_i G/K_i & \mapsto \prod_i \TAF(K_i)
\end{align*}

Let $\mc{H}$ denote the Hecke category of $G$.  
\index{2H@$\mc{H}$}
The objects of $\mc{H}$ are the same as the objects of $\mc{M}$.  The
morphisms are additively determined by
$$ \Hom_{\mc{H}}(G/K, G/K') =
\Hom_{\ZZ[G]}(\ZZ[G/K], \ZZ[G/K']) \cong (\ZZ[G/K'])^K. $$
This morphism space is a free abelian group with basis
$$ \left\{ T_{[g]} = 
\sum_{h \in K/K\cap gK'g^{-1}} hgK' \: : \: [g] \in K \backslash 
G/K'\right\} \subset \ZZ[G/K']^K. $$
The \index{2Tg@$T_{[g]}$}
endomorphisms in $\mc{H}$ of $G/K$ give the Hecke
algebra 
\index{Hecke algebra}
of the pair $(G,K)$.  

There is a natural functor
$$ F : \mc{M} \rightarrow \mc{H} $$
(see \cite{minami})
\index{2F@$F$}
which associates to the correspondence
$$
\xymatrix{
& X \ar[dl]_f \ar[dr]^g \\
S && T
}
$$
the morphism given by
\begin{align*}
\ZZ[S] & \rightarrow \ZZ[T] \\
s & \mapsto \sum_{x \in f^{-1}(s)} g(x)
\end{align*}
for elements $s \in S$.  The functor $F$ is full, but not  
faithful.  The
generators $T_{[g]} \in \Hom_\mc{H}(G/K, G/K')$ 
may be lifted to correspondences
$$
\xymatrix{
& G/K \cap gK'g^{-1} \ar[dl]_{r_{K \cap gK'g^{-1}}^{K}} \ar[dr]^{\rho_g
\circ
r_{K \cap gK'g^{-1}}^{gK'g^{-1}}} \\
G/K && G/K'
}
$$
but this does \emph{not} induce a splitting of $F$.

\begin{defn}[Hecke Operators on $\TAF$]\label{defn:hecke}
For $[g] \in K\backslash G/K'$ define the \emph{Hecke operator}
$$ T_{[g]} : \TAF(K') \rightarrow \TAF(K) $$
to be the composite:
$$ \TAF(K') \xrightarrow{g_*} \TAF(gK'g^{-1}) \xrightarrow{\Res_{K \cap
gK'g^{-1}}^{gK'g^{-1}}} \TAF(K \cap gK'g^{-1}) \xrightarrow{\Tr_{K \cap
gK'g^{-1}}^K} \TAF(K) $$
\end{defn}
\index{Hecke operator}

\begin{ques}\label{ques:hecke}
Do these Hecke operators 
induce a factorization of the following form?
$$
\xymatrix@C+2em{
\mc{M}^{op} \ar[r]^{\TAF(-)} \ar[d]_{F} & \mathrm{Ho}(\Sp) \\
\mc{H}^{op} \ar@{.>}[ur] 
}
$$
In particular, does this induce a map of \emph{rings} 
from the Hecke algebra for
$(G,K)$ into the ring of stable cohomology operations on $\TAF(K)$?
\end{ques}

We are unable to resolve this question in general. We can give
some partial results.

\begin{prop}\label{prop:cohomologicalmackey}
Suppose that
for every compact open subgroup $K$ of $G$, and every quasi-coherent sheaf
$\mc{F}$ on $\Sh(K)^\wedge_p$,
the sheaf cohomology groups
$$ H^s(\Sh(K)^\wedge_p, \mc{F}) $$
vanish for $s > 0$.
Then the Mackey functor $\TAF(-)$ factors through the Hecke
category.
\end{prop}

\begin{rmk}
The hypotheses of Proposition~\ref{prop:cohomologicalmackey} are extremely
restrictive.  There are morally two ways in which the hypotheses can fail to be
satisfied.
\begin{enumerate}
\item The stack $\Sh(K)$ can fail to have an \'etale covering space
which is an affine scheme.
\item The stack $\Sh(K)$ can have stacky points whose automorphism group
contains $p$-torsion.
\end{enumerate}
This immediately rules out the case where the algebra $B$ is a division
algebra (Theorem~\ref{thm:shimurastack}(3)).  Note that the
(non-compactified) moduli stack of
elliptic curves $\mc{M}_{ell}$ \emph{does} satisfy the hypotheses of
Proposition~\ref{prop:cohomologicalmackey} for $p \ge 5$.  (See
Remark~\ref{rmk:baker} for more discussion on the modular case.)  
\end{rmk}

\begin{proof}[Proof of Proposition~\ref{prop:cohomologicalmackey}]
Fix a compact open subgroup $K$.
The descent spectral sequence (\ref{eq:TAFdescent})
$$
H^s(\Sh(K)^\wedge_p, \omega^{\otimes t}) \Rightarrow \pi_{2t-s}(\TAF(K)) 
$$
collapses to give an edge isomorphism
$$ \pi_{2t}(\TAF(K)) \xrightarrow{\cong} 
H^0(\Sh(K)^\wedge_p, \omega^{\otimes t}). $$
By Proposition~\ref{prop:LEFT},
our hypotheses guarantee that the spectrum $\TAF(K)$ is Landweber exact.
Picking a complex orientation for $\TAF(K)$, the ring $\pi_*\TAF(K)$
becomes an $MU_*$-algebra, and the natural map
$$ MU_*MU \otimes_{MU_*} \pi_*\TAF(K) \rightarrow MU_*\TAF(K) $$
is an isomorphism, giving (see \cite{adams})
\begin{equation}\label{eq:MUTAF}
MU_*\TAF(K) \cong \pi_* \TAF(K) [b_1, b_2, ...].
\end{equation}

By Landweber exactness, for any pair of compact open subgroups $K$,$K'$, 
the natural homomorphism
$$ [\TAF(K), \TAF(K')] \rightarrow \Hom_{MU_*MU}(MU_*\TAF(K), MU_*\TAF(K')) $$
is an isomorphism \cite[Cor.~2.17]{hoveystrickland}.
Yoshida (\cite{yoshida}, \cite[Thm.~3.7]{minami}) proved that a Mackey
functor
$$ \mc{M}^{op} \rightarrow \text{Abelian Groups} $$
factors through the Hecke category if and only if for every $K' \le K$ we
have $\Tr_{K'}^K \Res_{K'}^K = \abs{K:K'}$.  It therefore suffices to prove
that the composite
$$ MU_*\TAF(K) \xrightarrow{MU_*\Res_{K'}^{K}} MU_*\TAF(K')
\xrightarrow{MU_*\Tr_{K'}^K} MU_*\TAF(K) $$
is multiplication by $\abs{K:K'}$.

Recall that for $K'$, an open normal subgroup of $K$, the morphism
$$ \Sh(K') \rightarrow \Sh(K) $$
is a Galois covering space with Galois group $K/K'$
(Theorem~\ref{thm:shimurastack}).  We therefore have
\begin{align*}
\pi_*\TAF(K) & \cong (\colim_{K' \le_o K} \pi_*\TAF(K'))^{K} \\
& \cong (\pi_*V_{GU})^{K}.
\end{align*}
where the colimit is taken over open subgroups of $K$.
We deduce from (\ref{eq:MUTAF}) that the map 
$$ MU_*(\TAF(K)) = MU_*(V_{GU}^{hK}) \rightarrow MU_*(V_{GU}) $$
is a monomorphism. Consider the following diagram (where $V = V_{GU}$).

\begin{landscape}
$$
\xymatrix{
MU_*V^{hK} \ar[r]^{\Res_{K'}^K} \ar@{^{(}->}[ddd] &
MU_*V^{hK'} \ar[r]_-{(1)}^-{\cong} \ar@{^{(}->}[ddd] \ar@/^3pc/[rrr]^{\Tr_{K'}^{K}} &
MU_*(\CoInd_{K'}^{K} V)^{hK} \ar[d] &
MU_*(\Ind_{K'}^K V)^{hK} \ar[l]_{\cong} \ar[d] \ar[r]_-{(2)} &
MU_*V^{hK} \ar@{^{(}->}[ddd]
\\
&& 
MU_*\CoInd_{K'}^{K} V \ar[d]_{(3)}^{\cong} &
MU_*\Ind_{K'}^K V \ar[l]_{\cong} \ar[d]_{(4)}^{\cong} &
\\
&&
MU_*\Map(K/K', V) \ar[d]_\cong &
MU_*K/K'_+ \wedge V \ar[l]^{\cong} \ar[d]^\cong
\\
MU_*V \ar@{=}[r] &
MU_*V \ar[r]_-{\mathrm{diag}} &
\prod\limits_{K/K'} MU_*V &
\bigoplus\limits_{K/K'} MU_*V \ar[l]^{\cong} \ar[r]_-{\mathrm{add}} &
MU_*V
}
$$
\end{landscape}

Here $(1)$ is the isomorphism given by Lemma~\ref{lem:Shapiro} and $(2)$ is
the counit of the adjunction.  The isomorphisms $(3)$ and $(4)$ are induced
by the canonical $K$-equivariant isomorphisms:
\begin{align*}
\CoInd_{K'}^K V \cong \Map(K/K', V), \\
\Ind_{K'}^K V \cong K/K'_+ \wedge V,
\end{align*}
where $K$ acts on $\Map(K/K', V)$ by conjugation and on $K/K'_+ \wedge V$
diagonally.  The bottom horizontal 
composite is multiplication by $\abs{K: K'}$.
Since the leftmost and rightmost vertical arrows are monomorphisms, we
deduce that the topmost horizontal composite is multiplication by $\abs{K:
K'}$.
\end{proof}

Precisely the same methods prove the following variant.

\begin{prop}\label{prop:heckealgebra}
Suppose that $K$ is a fixed compact open subgroup, and that
for every compact open subgroup $K'$ of $K$, and every quasi-coherent sheaf
$\mc{F}$ on $\Sh(K')^\wedge_p$,
the sheaf cohomology groups
$$ H^s(\Sh(K')^\wedge_p, \mc{F}) $$
vanish for $s > 0$.
Then the Hecke operators induce a homomorphism of \emph{rings}
$$ \ZZ[G/K]^K \hookrightarrow [\TAF(K), \TAF(K)]. $$
Here $\ZZ[G/K]^K \cong \Hom_{\ZZ[G]}(\ZZ[G/K], \ZZ[G/K])$ is the Hecke
algebra for $(G,K)$.
\end{prop}

\begin{rmk}\label{rmk:baker}
This is an analog of the work of Baker \cite{baker}.  Baker proved that for
$p \ge 5$, the spectra
$\TMF_{(p)}$ admit an action of the Hecke algebra for the pair
$(GL_2(\AF^{p,\infty}), GL_2(\widehat{\ZZ}^p))$.  The authors do not know
if his results extend to the cases where $p = 2$ or $3$.  These are
precisely the cases where the corresponding moduli stack has nontrivial
higher sheaf cohomology.
\end{rmk}

We note that regardless of the outcome of
Question~\ref{ques:hecke}, Lemma~\ref{lem:epsilon} implies that the
edge homomorphism
$$ \pi_{2t}\TAF(K) \rightarrow (\pi_{2t}V_{GU})^K = 
H^0(\Sh(K)_p^\wedge, \omega^{\otimes t}) $$
commutes with the action of the Hecke operators (as given in
Definition~\ref{defn:hecke}).  
This relates the action of the
Hecke operators on $\TAF(K)$ to the action of the classical Hecke operators
of the corresponding space of automorphic forms.

\chapter{Buildings}
\label{chap:buildings}

In this chapter we give explicit descriptions for the Bruhat-Tits 
buildings of the
local forms of the group $GU$ and $U$.  
The buildings are certain 
finite dimensional contractible simplicial complexes on
which the groups $GU(\QQ_\ell)$ and $U(\QQ_\ell)$ act with compact
stabilizers.
Our treatment follows closely that
of \cite{abramenkonebe}.

\section{Terminology}

A \emph{semi-simplicial set} 
\index{semi-simplicial set}
is a simplicial set without
degeneracies.  Thus a semi-simplicial set
$\{ X_\bullet, d_\bullet \}$ consists of a sequence of sets $X_k$ of
$k$-simplices, with face maps $d_\bullet: X_k \rightarrow X_{k-1}$
satisfying the simplicial face relations.

A \emph{simplicial complex} 
\index{simplicial complex}
$\{X_\bullet, \le \}$ consists of a sequence of
sets $X_k$ of $k$-simplices together with a poset structure on the set
$\coprod_k X_k$ which encodes face containment.  Thus a semi-simplicial set
is a simplicial complex with a compatible ordering on the vertices of its
simplices.

For a finite dimensional simplicial complex $X_\bullet$, let $d$ denote the
maximal dimension of its simplices.  A $d$-simplex in $X_\bullet$ is
called a \emph{chamber}.  
\index{chamber}
The simplicial complex $X_\bullet$ is called a
\emph{chamber complex} 
\index{chamber complex}
if every simplex is contained in a chamber, and
given any two chambers $a$ and $b$ there is a sequence of chambers
$$ a = a_0, a_1, \ldots, a_k = b $$
such that for each $i$, the chambers $a_i$ and $a_{i-1}$ intersect in a
$(d-1)$-dimensional simplex (a \emph{panel}).  
\index{panel}
A chamber complex is \emph{thin}
\index{chamber complex!thin}
if each panel is contained in exactly $2$ chambers.  A chamber complex  is 
\emph{thick}
\index{chamber complex!thick}
\index{building!thick}
if each panel is contained in at least $3$ chambers.  

For a thin $d$-dimensional 
chamber complex $X_\bullet$, a \emph{folding} 
\index{chamber complex!folding}
is a morphism of
simplicial complexes
$$ f : X_\bullet \rightarrow X_\bullet $$
satisfying $f^2 = f$ for which every chamber in the image of $f$ is in the
image of exactly $2$ chambers.  The complex $X_\bullet$ is called a 
\emph{Coxeter
complex} 
\index{Coxeter complex}
if, for every pair of chambers $a$ and $b$, there exists a folding
$f$ for which $f(a) = b$.  

A \emph{building} 
\index{building}
$\mc{B}$ is a chamber complex which is a union of Coxeter
complexes.  These subcomplexes are the
\emph{apartments} 
\index{apartments}
of $\mc{B}$, and are required to satisfy the following
two conditions: 
\begin{enumerate}
\item
For any pair of chambers, there must exist an apartment containing both of
them.

\item
For any pair of apartments $\mc{A}_1$ and $\mc{A}_2$, there must
be an isomorphism of simplicial complexes
$$ \phi: \mc{A}_1 \rightarrow \mc{A}_2 $$
such that $\phi$ restricts to the identity on $\mc{A}_1 \cap \mc{A}_2$.
\end{enumerate}
A building $\mc{B}$ is said to be \emph{affine} 
\index{building!affine}
if its apartments are
tessellations of a Euclidean space by reflection hyperplanes of an
irreducible affine reflection group.

\begin{rmk}
All of the buildings we shall consider will be affine.  The buildings for
$SL$ and $U$ will be thick, but the buildings for $GL$ and $GU$ will
\emph{not} be thick.
\end{rmk}

\section{The buildings for $GL$ and $SL$}

Let $K$ 
\index{2K@$K$}
be a discretely valued local field with ring of integers $\calO$
\index{2O@$\mc{O}$}
and uniformizer $\pi$.
\index{1pi@$\pi$}
Let $W$ be a $K$-vector space of dimension $n$.
\index{2W@$W$}

\subsection*{The building for $GL(W)$}

The building $\mathcal{B}(GL(W))$ 
\index{2BGL@$\mc{B}(GL)$}
\index{building!for $GL$}
is 
the geometric realization of a
semi-simplicial set $\mathcal{B}(GL(W))_\bullet$.  The
$k$-simplices are given by sets of lattice chains in $W$:
$$ \mathcal{B}(GL(W))_k = 
\{ L_0 < L_1 < \cdots < L_k \le \pi^{-1} L_0
\}.
$$
The $i$th face map is given by deleting the $i$th lattice from the chain.
The group $GL(W)$ acts on the building simplicially by permuting the
lattices.

A basis $\mathbf{v} = (v_1, \ldots, v_n)$ of 
$W$ gives rise to a maximal simplex
(chamber)
$$
\mathcal{C}(\mathbf{v}) = (
L_0(\mathbf{v}) < L_1(\mathbf{v}) < \cdots < L_n(\mathbf{v}) = \pi^{-1}
L_0(\mathbf{v}) )
$$
where
$$ L_i(\mathbf{v}) = 
\pi^{-1} \calO v_1 \oplus \cdots \oplus \pi^{-1} \calO v_i \oplus \calO v_{i+1}
\oplus \cdots \oplus \calO v_n.
$$
In particular,
the building $\mathcal{B}(GL(W))$ is $n$-dimensional.  

The following lemma is clear.

\begin{lem}
The action of $GL(W)$ on $\mathcal{B}(GL(W))$ 
is transitive on vertices and chambers.
\end{lem}

\subsection*{The building for $SL(W)$}

Two lattice chains $\{L_i\}$, $\{L'_i\}$ 
are \emph{homothetic} 
\index{homothetic}
if there exists an integer $m$ so
that for all $i$, $L'_i = \pi^m L_i$.

The building for $SL(W)$ is the quotient of $\mathcal{B}(GL(W))$ given
by taking homothety classes of lattice chains.  
The building $\mathcal{B}(SL(W))$ 
\index{2BSL@$\mc{B}(SL)$}
\index{building!for $SL$}
may be described as the simplicial complex
whose $k$-simplices 
are given by homothety classes of \emph{certain}
lattice chains:
$$
\mathcal{B}(SL(W))_k = 
\{ [L_0 < L_1 < \cdots < L_k < \pi^{-1} L_0] \}.
$$
The face containment relations of $\mathcal{B}(SL(W))$ 
are given by chain containment up to
homothety.

The building for
$SL(W)$ is $n-1$ dimensional.
There is a natural projection of simplicial complexes
$$ 
\mathcal{B}(GL(W)) \rightarrow \mathcal{B}(SL(W))
$$
and the building 
$\mathcal{B}(GL(W))$ is homeomorphic to $\mathcal{B}(SL(W)) \times \RR$.  
The $k$-simplices $\{L_i\}$ of $\mathcal{B}(GL(W))$ for which $L_k =
\pi^{-1} L_0$ become degenerate in $\mathcal{B}(SL(W))$.  

\section{The buildings for $U$ and $GU$}

Suppose now that $K$ 
\index{2K@$K$}
is a quadratic extension of $\QQ_\ell$.  Let $(-,-)$
\index{0bro@$(-,-)$}
be a non-degenerate hermitian form on $W$.  

\subsection*{Totally isotropic subspaces}

A vector $v \in W$ is \emph{isotropic}
\index{isotropic!vector}
if $(v,v) = 0$.  A subspace of $W$ is \emph{anisotropic}
\index{anisotropic}
if it contains no non-zero
isotropic vectors. A subspace $W'$ of $W$
is \emph{totally isotropic} 
\index{isotropic!subspace}
\index{totally isotropic}
if $(v,w) = 0$ for all $v,w \in W'$.
A subspace $W'$ of $W$ is \emph{hyperbolic} 
\index{hyperbolic!subspace}
if it admits a basis $(v_1, \ldots,
v_{2k})$ such that with respect to this basis the restricted form
$(-,-)\vert_{W'}$ is given by a matrix
$$
\begin{bmatrix}
&& * \\
& \reflectbox{$\ddots$} & \\
* & & 
\end{bmatrix}
$$
where an entry is nonzero if and only if it lies on the reverse diagonal.
Such a basis is a $\emph{hyperbolic basis}$.
\index{hyperbolic!basis}

The \emph{Witt index} 
\index{Witt index}
$r$ of $(-,-)$ is the dimension of the maximally totally
isotropic subspace of $W$.  The dimension of a maximal hyperbolic subspace
$V$ of $W$ is $2r$, and the subspace $V^\perp$ is anisotropic.
The dimensions of the buildings for $GU(W)$ and $U(W)$ are functions of the Witt
index.

Recall from Chapter~\ref{chap:galois} that
isometry classes of quadratic hermitian forms are classified by their
discriminant $\mathrm{disc}$ in $\QQ_\ell^\times/N(K^\times)$.  Let $a \in
\QQ_\ell^\times$ be a generator of $\QQ_\ell^\times/N(K^\times)$.
The following lemma is easily deduced.

\begin{lem}\label{lem:isotropy}
There 
exists a basis so that the form $(-,-)$ is given by the 
matrices in the following table.

\begin{tabular}{c|c|c|}
&
$\mathrm{disc} = (-1)^{\lfloor n/2 \rfloor} \in \QQ^\times/N(K^\times)$ &
$\mathrm{disc} \ne (-1)^{\lfloor n/2 \rfloor} \in \QQ^\times/N(K^\times)$
\\
\hline
&& \\
$n$ even &
$\left[
\begin{array}{ccc}
& & 1 \\
& \reflectbox{$\ddots$} & \\
1 & & 
\end{array}
\right]$ &
$\left[
\begin{array}{ccc|cc}
& & 1 & & \\
& \reflectbox{$\ddots$} & & & \\
1 & & & & \\
\hline
&&& 1 & \\
& & & & -a \\
\end{array}
\right]$
\\
& $r = n/2$ & $r = (n-2)/2$ \\
\hline
&& \\
$n$ odd &
$\left[
\begin{array}{ccc|c}
& & 1 & \\
& \reflectbox{$\ddots$} & & \\
1 & & & \\
\hline
&&& 1 \\
\end{array}
\right]$ &
$\left[
\begin{array}{ccc|c}
& & 1 & \\
& \reflectbox{$\ddots$} & & \\
1 & & & \\
\hline
&&& a \\
\end{array}
\right]$
\\
& $r = (n-1)/2$ & $r = (n-1)/2$ \\
\hline
\end{tabular}
\end{lem}

\subsection*{The building for $U$}

Given an $\calO$-lattice $L$ of $W$, the dual lattice 
\index{2L@$L^\#$}
\index{dual!lattice}
\index{lattice!dual}
$L^\# \subset W$ is the
$\calO$-lattice given by 
$$ L^\# = \{ w \in W \: : \: (w,L) \subseteq \calO \}. $$
The dual lattice construction has the following properties:
\begin{enumerate}
\item $L^{\# \#} = L$, 
\item if $L_0 < L_1$, then $L_1^\# < L_0^\#$,
\item $(\pi L)^\# = \pi^{-1}L^\#$.
\end{enumerate}
Thus, there is an induced involution
$$ \iota \co \mathcal{B}(SL(W)) \rightarrow \mathcal{B}(SL(W)). $$
The building $\mathcal{B}(U(W)) \subseteq \mathcal{B}(SL(W))$ 
\index{building!for $U$}
\index{2BU@$\mc{B}(U)$}
consists of the fixed points of this
involution.  Since the involution $\iota$ does not necessarily restrict to
the identity on invariant chambers, the subspace $\mathcal{B}(U(W))$ is not
necessarily a simplicial subcomplex.  

However, the building
$\mathcal{B}(U)$ does admit the structure of a semi-simplicial set whose
simplices are self-dual lattice chains.
The $k$-simplices are
given by
$$
\mathcal{B}(U(W))_k = \{ L_0 < L_1 < \cdots < L_k \le L_k^\# < \cdots <
L_0^\# \le \pi^{-1} L_0 \}.
$$
The face maps $d_i$ are obtained by deletion of the lattices $L_i$ and
$L_i^\#$ from the
chain.
In particular, the vertices are given by lattices $L$ satisfying $L \le
L^\# \le \pi^{-1} L$.  
We shall say that such lattices are \emph{preferred}.
\index{preferred lattice}
\index{lattice!preferred}
Thus the $k$-simplices of $\mathcal{B}(U(W))$ are given by chains of
preferred lattices:
$$
\mathcal{B}(U(W))_k = \{ L_0 < L_1 < \cdots < L_k < \pi^{-1}L_0 \: : \:
\text{$L_i$ are preferred} \}.
$$

Preferred lattices satisfy the following properties.
\begin{enumerate}
\item There is at most one preferred representative in a homothety lattice 
class.
\item If $L$ is preferred, then either $L = L^\#$, or $L^\#$ is not
preferred.
\end{enumerate}
These properties allow one to relate the simplicial description of
$\mathcal{B}(U(W))$ with the geometric fixed points of the involution
$\iota$ on $\mathcal{B}(SL(W))$.  A combinatorial vertex corresponding to a
preferred lattice $L$ corresponds to the midpoint of the edge joining $[L]$
and $[L^\#]$.

It is easily checked that for $g \in GL(W)$, one has $(g(L))^\# = g^{-\iota}
(L^\#)$.  Therefore the building $\mathcal{B}(U(W))$ admits a simplicial
action by $U(W)$.  Let $GU^1(W)$ be the subgroup
$$ GU^1(W) = 
\ker \left( GU(W) \xrightarrow{\nu} 
\QQ_\ell^\times \xrightarrow{\nu_\ell} \ZZ \right) $$
of \index{2GU1@$GU^1$}
similitudes of $W$ with similitude norm of valuation $0$.  Then the
following lemma is immediate.

\begin{lem}
The action of $U(W)$ on $\mc{B}(U(W))$ extends to an action of $GU^1(W)$.
\end{lem}

In order to get an explicit description of the chambers of
$\mathcal{B}(U(W))$ we recall
the following lemma of \cite[Sec.~5]{abramenkonebe}.

\begin{lem}\label{lem:abramenkonebe}
Suppose that the residue characteristic $K$ does not equal $2$.  Let $W'$
be an anisotropic subspace of $W$.  Then there is a unique preferred 
lattice $X$ contained in $W'$.
This lattice is given by 
$$ X = \{ w \in W' \: : \: (w,w) \in \calO\}. $$
\end{lem}

Let $\mathbf{v} = (v_1, \ldots, v_{2r})$ be a hyperbolic basis of a 
maximal hyperbolic 
subspace $V$ of
$W$ such that $(v_i, v_{2r+1-i}) = 1$.  We shall say that a hyperbolic
basis with this property is \emph{normalized}, 
\index{hyperbolic!basis!normalized}
and clearly every 
hyperbolic basis can be
normalized.  Let $X$ be the
unique preferred lattice in the anisotropic subspace $V^\perp$. 
Then we may associate to $\mathbf{v}$ a chamber
$$ \mathcal{C}(\mathbf{v}) = 
( L_0(\mathbf{v}) < \cdots < L_r(\mathbf{v}) \le 
L_r^\#(\mathbf{v}) < \cdots <
L_0^\#(\mathbf{v}) \le \pi^{-1} L_0(\mathbf{v}) )
$$
where
$$
L_i(\mathbf{v}) = \pi \calO v_1 + \cdots \pi \calO v_{r-i} + \calO v_{r-i+1} +
\cdots \calO v_{2r} + X.
$$
In particular the building $\mathcal{B}(U(W))$ has dimension $r$.

Unlike the building for $SL(W)$, the group $U(W)$ does not act transitively
on vertices.  However, every chamber of $\mathcal{B}(U(W))$ can be shown to
take the form $\mathcal{C}(\mathbf{v})$ for some normalized hyperbolic
basis $\mathbf{v}$.  Witt's theorem and Lemma~\ref{lem:abramenkonebe} may
then be used to derive the following lemma.

\begin{lem}
The group $U(W)$ acts transitively on the set of chambers of the building
$\mathcal{B}(U(W))$.
\end{lem}

\begin{lem}
The stabilizers of the action of $U(W)$ and $GU^1(W)$ on $\mc{B}(U(W))$ are
all compact and open.
\end{lem}

\begin{proof}
We give the proof for $U(W)$; the argument for $GU^1(W)$ is identical.
Let $x = (L)$ be a vertex of $\mathcal{B}(U(W))$ given by a preferred
lattice $L$, 
and let $K_x \subset U(W)$ be its stabilizer.  Then $K_x$ is given by the
intersection $GL(L) \cap U(W)$ in $GL(W)$.  Since $U(W)$ is a closed
subgroup of $GL(W)$ and $GL(L)$ is compact and open, the subgroup $K_x$ is
compact and open in $U(W)$.  Since the action is simplicial, 
it follows that all of the stabilizers of
simplices in $\mathcal{B}(U(W))$ are compact and open.
\end{proof}

\subsection*{The building for $GU(W)$}

The group $GU(W)$ does not preserve the class of preferred lattices, but we
do have the following lemma.

\begin{lem}\label{lem:preferredGU}
Let $g$ be an element of $GU(W)$ with similitude norm $\nu(g)$ of valuation
$k$.
Suppose that $L$ is homothetic to a preferred lattice.  Then $(g(L))^{k\#}$
is homothetic to a preferred lattice.
\end{lem}

\begin{proof}
Suppose that $L = \pi^{b}\br{L}$ where $\br{L}$ is a
preferred lattice.
Observe that we have
$$ g(\br{L}^\#) = (g^{-\iota}(\br{L}))^\# = (\nu^{-1}g(\br{L}))^\# =
\pi^{k} (g(\br{L}))^\#. $$

Suppose that $k = 2a$ is even.  Then the element $g_1 = \pi^{-a}g$ has
similitude norm of valuation zero, and $g_1(\br{L})$ is preferred.
The lattice $g(L)$ is homothetic to $g_1(\br{L})$.

Suppose that $k = 2a-1$ is odd.  
The element $g_1 = \pi^{-a}g$ has similitude norm
of valuation $-1$.  Applying
$g_1$ to the chain $\pi L^\# \le L \le L^\#$, and using the fact that 
$(g_1(\br{L}))^\# = \pi g_1(\br{L}^\#)$, we see that $(g_1(\br{L}))^\#$ is
preferred.  The lattice $(g(L))^\#$ is therefore homothetic to a preferred
lattice.
\end{proof}

The action of $U(W)$ on the building $\mathcal{B}(U(W))$ can be extended to
an action of the group $GU(W)$.  The action is determined by its
effect on
vertices, where we regard the vertices as preferred homothety classes of
lattices.  The action is given by
$$ g \cdot [L] = 
\begin{cases}
[g(L)] & \text{if $\nu(g)$ has even valuation,} \\
[(g(L))^\#] & \text{if $\nu(g)$ has odd valuation}.
\end{cases}
$$
Lemma~\ref{lem:preferredGU} implies that $g \cdot [L]$ is a preferred
homothety class.
This action is not associative on lattices, but is associative on homothety
classes.  This is an action of simplicial complexes, but does not preserve
the ordering on the vertices.

The building $\mathcal{B}(SL(W))$ carries a natural action by $GU(W)$
through the natural action of $GL(W)$.  We have the following lemma.

\begin{lem}
The inclusion of the fixed point space $j \co \mathcal{B}(U(W)) \hookrightarrow
\mathcal{B}(SL(W))$ is $GU(W)$-equivariant.
\end{lem}

\begin{proof}
It suffices to check equivariance on vertices.
Recall that for a preferred lattice $L$, the inclusion
$j$ maps $[L]$ to the midpoint $\mathrm{mid}([L],[L^\#])$.  
Suppose that $g$ is an element of $GU(W)$, with similitude norm of
valuation $k$.  Then we have
\begin{align*}
j(g\cdot [L]) 
& = j([(g(L))^{k\#}]) \\
& = \mathrm{mid}([(g(L))^{k\#}], [(g(L))^{(k+1)\#}]) \\
& = \mathrm{mid}([g(L)], [(g(L))^\#]) \\
& = g \cdot j([L]).
\end{align*}
\end{proof}

We define $\mathcal{B}(GU(W))$ by the following $GU(W)$-equivariant
\index{building!for $GU$}
\index{2BGU@$\mc{B}(GU)$}
pullback.
$$
\xymatrix{
\mathcal{B}(GU(W)) \ar@{^{(}->}[r] \ar[d] &
\mathcal{B}(GL(W)) \ar[d] \\
\mathcal{B}(U(W)) \ar@{^{(}->}[r] &
\mathcal{B}(SL(W))
}
$$
We see that $\mathcal{B}(GU(W))$ is homeomorphic to $\mathcal{B}(U(W))
\times \RR$.
In particular, the dimension of $\mathcal{B}(GU(W))$ is $r+1$.

Since $GU(W)$ is a closed subgroup of $GL(W)$, and the stabilizer in
$GL(W)$ of every
point of $\mathcal{B}(GL(W))$ is open and compact, we have the following
lemma.

\begin{lem}
The stabilizer in $GU(W)$ of every point of $\mathcal{B}(GU(W))$ is open and
compact.
\end{lem}

\begin{rmk}
The buildings $\mathcal{B}(GU(W))$ and $\mathcal{B}(GL(W))$ are easily seen
(using a simplicial prism decomposition)
to be \emph{equivariantly} homeomorphic to the products
\begin{align*}
\mathcal{B}(GU(W)) & \approx \mathcal{B}(U(W)) \times \RR_{\nu}, \\
\mathcal{B}(GL(W)) & \approx \mathcal{B}(SL(W)) \times \RR_{\det}.
\end{align*}
Here, $\RR_\nu$ is the similitude line, with $GU(W)$ action
$g(x) = v_\pi(\nu(g)) \cdot x$, and $\RR_{\det}$ is the determinant line,
with $GL(W)$ action $g(x) = v_\pi(\det(g)) \cdot x$.  
The embedding of $\mathcal{B}(GU(W))$ into $\mathcal{B}(GL(W))$ is given by
the embedding of $\mathcal{B}(U(W))$ into $\mathcal{B}(SL(W))$, together
with the $GU(W)$-equivariant homeomorphism 
\begin{align*}
\RR_\nu & \xrightarrow{\approx} \RR_{\det}, \\
x & \mapsto n/f \cdot x.
\end{align*}
Here, $f$ is the residue field degree of $K$ over $\QQ_\ell$.
\end{rmk}

\chapter{Hypercohomology of adele groups}\label{chap:Qspectra}

\section{Definition of $Q_{GU}$ and $Q_U$}

Lemma~\ref{lem:cofixedpoints} has the following corollary.

\begin{cor}
For compact open $K^p \le GU(\AF^{p,\infty})$ there are equivalences
$$ \TAF(K^p) \simeq V_{GU}^{hK^p}. $$
\end{cor}

\begin{defn} 
Let $K^p = \prod_{\ell \ne p} K_\ell$ 
be a compact open subgroup of $GU(\AF^{p,\infty})$, and 
let $S$ be a set of primes not containing $p$.
Let $K^{p,S}$ be the product 
$$ K^{p,S} = \prod_{\ell \not\in S \cup \{ p \} } K_\ell. $$
\begin{enumerate}
\item  
Define $Q_{GU}(K^{p,S})$ to be the
homotopy fixed point spectrum
$$ Q_{GU}(K^{p,S}) = V_{GU}^{hK^{p,S}_+} $$
where \index{2QGU@$Q_{GU}$}
$$ K^{p,S}_+ := GU(\AF_S)K^{p,S}. $$
If \index{2KpS+@$K^{p,S}_+$}
$S$ consists of all primes different from $\ell$, we
denote this spectrum $Q_{GU}$. 

\item 
Let $GU^1(\AF_S)$ be the open subgroup
\index{2GU1@$GU^1$}
$$ GU^1(\AF_S) = \ker \left( GU(\AF_S) \xrightarrow{\nu} {\prod_{\ell \in S}}'
\QQ_\ell^\times \xrightarrow{\nu_\ell} \prod_{\ell \in S} \ZZ \right) $$
of similitudes $g$ whose similitude norm $\nu(g)$ 
has valuation $0$ at every place in
$S$.
Define $Q_U(K^{p,S})$ 
\index{2QU@$Q_U$}
to be the homotopy fixed point spectrum
$$ Q_U(K^{p,S}) = V_{GU}^{hK^{p,S}_{1,+}}, $$
where
$$ K^{p,S}_{1,+} := GU^1(\AF_S)K^{p,S}. $$
If \index{2KpS1+@$K^{p,S}_{1,+}$}
$S$ consists of all primes different from $\ell$, we
denote this spectrum $Q_{U}$. 
\end{enumerate}
\end{defn}

\section{The semi-cosimplicial resolution}

Fix a prime $\ell \ne p$ and a compact open subgroup $K^{p,\ell} \subset
GU(\AF^{p,\ell,\infty})$.  
\index{2Kpl@$K^{p,\ell}$}
We shall consider two cases simultaneously:

\begin{center}
\begin{tabular}{r|c|c}
& Case I & Case II \\
\hline
&& \\
$G_\ell =$ & $GU(\QQ_\ell)$ & $GU^1(\QQ_\ell)$ \\
&& \\
$\mc{B} =$ & $\mc{B}(GU(\QQ_\ell))$ & $\mc{B}(U(\QQ_\ell))$ \\
&& \\
$Q =$ & $Q_{GU}(K^{p,\ell})$ & $Q_U(K^{p,\ell})$
\end{tabular}
\end{center}
In \index{2Gl@$G_\ell$} \index{2B@$\mc{B}$} \index{2Q@$Q$}
either case, $\mc{B} = \abs{\mc{B}_\bullet}$ 
\index{2B@$\mc{B}_\bullet$}
is the realization of a 
finite dimensional contractible semi-simplicial
set $\mc{B}_\bullet$, and the group $G_\ell$ acts on $\mc{B}$ 
simplicially with compact open
stabilizers.  Throughout this section, we will let $\mc{B}'_\bullet$ be
the simplicial set generated by the semi-simplicial set $\mc{B}_\bullet$ by
including degeneracies.  
\index{2B'@$\mc{B}'_\bullet$}

Let $G$ denote the open subgroup 
$$ G = G_\ell K^{p,\ell} \subset GU(\AF^{p,\infty}), $$
so \index{2G@$G$}
that $Q = V_{GU}^{hG}$.
The action of the group $G_\ell$ on $\mc{B}$ naturally extends to an action
of the group $G$, where we simply let the factors $K_{\ell'}$ act trivially
for $\ell' \ne p, \ell$.  The stabilizers of this action are still compact
and open.

\subsection*{The canonical contracting homotopy}

We recall some of the material from \cite{garrett}, in particular
Section~13.7.  
These apartments of $\mc{B}$
possess \emph{canonical metrics}, 
\index{building!canonical metric}
and simplicial automorphisms of these are
isometries.  The building $\mc{B}$ inherits a metric as follows: for $x, y
\in \mc{B}$, let $\mc{A}$ be an apartment containing $x$ and $y$.  The
distance between $x$ and $y$ is equal to the distance taken in the
apartment $\mc{A}$ with respect to its canonical metric.  Since the group $G$ acts simplicially on $\mc{B}$, and
sends apartments to apartments, $G$ acts by isometries on $\mc{B}$.

Any two points $x, y \in \mc{B}$, are joined by a unique geodesic
$$ \gamma_{x,y}: I \rightarrow \mc{B}. $$
Let \index{1gammaxy@$\gamma_{x,y}$}
$\mc{A}$ be an apartment of $\mc{B}$ which contains both $x$ and $y$.
Then the geodesic $\gamma_{x,y}$ is given by the affine combination
$$ \gamma_{x,y}(t) = (1-t)x + ty $$
in the affine space $\mc{A}$.

A contracting homotopy for $\mc{B}$ is very simple to describe
\cite[Sec.~14.4]{garrett}.  Fix a point $x_0$ of $\mc{B}$.  A contracting
homotopy
$$ H: \mc{B} \times I \rightarrow \mc{B} $$
is given by
\begin{equation}\label{eq:contracting}
H_t(x) = \gamma_{x,x_0}(t).
\end{equation}
\index{2H@$H$}
\index{building!contracting homotopy}

\begin{rmk}
The argument \cite[Sec.~14.4]{garrett} that the contracting homotopy $H$ is
continuous
applies to \emph{thick} buildings.  While
the building $\mc{B}(U(\QQ_\ell))$ is thick if $\ell$ does not split in $F$, 
the building $\mc{B}(U(\QQ_\ell))$ in the split case, as well as the building
$\mc{B}(GU(\QQ_\ell))$ in either case, is not thick.  However, because
these other buildings are
products of thick buildings with an affine space, the continuity of $H$ is
easily seen to extend to these cases.
\end{rmk}

\subsection*{A technical lemma}

Let $E$ be a smooth $G$-spectrum, and suppose that $X$ is a simplicial
$G$-set.  For an open subgroup $U \le G$, define $\ul{\Map}_U(X, E)$ 
\index{2MapU@$\ul{\Map}_U(-,-)$}
to be
the spectrum whose $i$th space is the simplicial mapping 
space of $U$-equivariant maps: 
$$ \ul{\Map}_U(X, E)_i = \ul{\Map}_U(X,E_i). $$
Define a smooth $G$-spectrum
$$ \ul{\Map}(X, E)^{sm} := \varinjlim_{U \le_o G} 
\ul{\Map}_{U}(X, E) $$
(where \index{2Mapsm@$\ul{\Map}(-,-)^{sm}$}
the colimit is taken over open subgroups of $G$).
\index{0leo@$\le_o$}
The group $G$ acts on the maps in the colimit by conjugation.
In this section we will prove the following technical lemma.

\begin{lem}\label{lem:technicallemma}
The map of smooth $G$-spectra
$$ E \rightarrow \ul{\Map}(\mc{B}'_\bullet, E)^{sm} $$
induced by the map $\mc{B}'_\bullet \rightarrow \ast$
is an equivalence.
\end{lem}

The key point to proving Lemma~\ref{lem:technicallemma} 
is the following topological lemma.

\begin{lem}\label{lem:technicallemmaspace}
Suppose that $X$ is a topological space with a smooth $G$-action.  
Then the inclusion of the constant maps induces a natural inclusion
$$ X \rightarrow \Map(\mc{B}, X)^{sm} := \varinjlim_{U \le_o G}
\Map_U(\mc{B}, X) $$
which is the inclusion of a deformation retract.  Here, $\Map_U(\mc{B}, X)$
\index{2MapU@$\Map_U(-,-)$}
is given the subspace topology with respect to the mapping
space $\Map(\mc{B}, X)$, 
and $\Map(\mc{B}, X)^{sm}$ 
\index{2Mapsm@$\Map(-,-)^{sm}$}
is given the topology of
the union of the spaces $\Map_U(\mc{B}, X)$.
\end{lem}

\begin{proof}
Fix a point $x_0$ in $\mc{B}$ with compact open stabilizer $K_{x_0} \le
G$.
Let $H$ be the contracting homotopy defined by (\ref{eq:contracting}).
Then $H$ induces a deformation retract
$$ H' \co \Map(\mc{B}, X) \times I \rightarrow \Map(\mc{B}, X) $$
given by
$$ H'_t(f)(x) = f(H_t(x)) = f(\gamma_{x,x_0}(t)) $$
where $\gamma_{x,x_0}$ is the unique geodesic from $x$ to $x_0$.  
We claim that for each open subgroup $U$ of $G$,
$H'$ restricts to give a map 
\begin{equation}\label{eq:H'onUequis} 
H' \co \Map_U(\mc{B}, X) \times I \rightarrow \Map_{U \cap
K_{x_0}}(\mc{B}, X).
\end{equation}
Indeed, let $f$ be an
element of $\Map_U(\mc{B}, X)$.  Let $g$ be an element of $U \cap K_{x_0}$.  
Since $G$ acts
by isometries, the induced map
$$ g \co \mc{B} \rightarrow \mc{B} $$
sends geodesics to geodesics.  Hence we have
\begin{align*}
g(H_t'(f)(x)) & = g(f(\gamma_{x,x_0}(t))) \\
& = f(g\gamma_{x,x_0}(t)) \\
& = f(\gamma_{gx,gx_0}(t)) \\
& = f(\gamma_{gx, x_0}(t)) \\
& = H_t'(f)(gx).
\end{align*}
We have shown that $H_t'(f)$ is $U \cap K_{x_0}$-equivariant.

By taking the colimit of the maps $H'_t$ of (\ref{eq:H'onUequis}), we
obtain a deformation retract
$$ H' \co \Map(\mc{B}, X)^{sm} \times I \rightarrow \Map(\mc{B}, X)^{sm}. $$
\end{proof}

We now explain how, using the geometric realization and singular functors, 
the topological result of Lemma~\ref{lem:technicallemmaspace} gives an
analogous simplicial result.
The category $s\Set$ of simplicial sets and the category $\Top$ of
topological spaces are Quillen equivalent by adjoint functors $(\abs{-}, S)$
$$ \abs{-} \co s\Set \leftrightarrows \Top \co S $$
where $S(-)$ 
\index{2S@$S(-)$}
\index{0abs@$\abs{-}$}
is the singular complex functor and $\abs{-}$ is geometric
realization.  The geometric realization of a simplicial smooth
$G$-set is a $G$-space with smooth $G$-action.  Since the
realization of the $k$-simplex $\abs{\Delta^k}$ is compact, it is easily
verified that the singular chains on a $G$-space with smooth
$G$-action is a simplicial smooth $G$-set.

\begin{lem}\label{lem:realizationfibrant}
Suppose that $X$ is a fibrant simplicial smooth $G$-set.  Then $S\abs{X}$,
the singular complex of the geometric realization of $X$, is also fibrant.
\end{lem}

\begin{proof}
By Lemma~\ref{lem:fibrantsmGspace}, we must verify that for every open
subgroup $H \le G$, the $H$-fixed points
of $S\abs{X}$ is a Kan complex, and that $S\abs{X}$ satisfies homotopy descent
with respect to hypercovers of $G/H$.  Since both the functors $S(-)$ and
$\abs{-}$ both preserve fixed points, we see that $(S\abs{X})^H = S
\abs{X^H}$ is a Kan complex.  For a hypercover
$\{G/U_{\alpha,\bullet}\}_{\alpha \in I_\bullet}$, we must verify that the
map
$$ (S\abs{X})^H \rightarrow \holim_\Delta \prod_{\alpha \in
I_\bullet}(S\abs{X})^{U_{\alpha,\bullet}} $$
is a weak equivalence.  The functors $S(-)$ and
$\abs{-}$ commute with fixed points, $S(-)$ commutes with arbitrary products and
totalization, and $\abs{-}$ commutes (up to homotopy equivalence) 
with arbitrary products of Kan
complexes, so we need the map
$$ S\abs{X^H} \rightarrow S \holim_\Delta \abs{\prod_{\alpha \in
I_\bullet}X^{U_{\alpha,\bullet}}} $$
to be a weak equivalence.  It suffices to show that the map
$$ \abs{X^H} \rightarrow \holim_\Delta \abs{\prod_{\alpha \in
I_\bullet}X^{U_{\alpha,\bullet}}} $$
is a weak equivalence, or equivalently, since $(\abs{-},S(-))$ form Quillen
equivalence, that the map
$$ X^H \rightarrow S\holim_\Delta \abs{\prod_{\alpha \in
I_\bullet}X^{U_{\alpha,\bullet}}} \xrightarrow{\simeq} \holim_\Delta S\abs{\prod_{\alpha \in
I_\bullet}X^{U_{\alpha,\bullet}}} $$
is a weak equivalence.  This follows from the fact that the map
$$ 
\prod_{\alpha \in
I_\bullet}X^{U_{\alpha,\bullet}} \rightarrow 
S\abs{\prod_{\alpha \in
I_\bullet}X^{U_{\alpha,\bullet}}}
$$
is a level-wise weak equivalence of cosimplicial Kan complexes, and that
the map
$$ X^H \rightarrow \holim_\Delta \prod_{\alpha \in
I_\bullet}X^{U_{\alpha,\bullet}} $$
is a weak equivalence, since $X$ is fibrant.
\end{proof}

Lemma~\ref{lem:technicallemma} follows immediately from the following
lemma.

\begin{lem}
Let $X$ be a simplicial smooth $G$-set.
The natural map
$$ X \rightarrow \ul{\Map}(\mc{B}'_\bullet, X)^{sm} $$
is a weak equivalence of simplicial sets.
\end{lem}

\begin{proof}
By Lemma~\ref{lem:technicallemmaspace}, the map
$$ \abs{X} \xrightarrow{\simeq} \Map(\mc{B}, \abs{X})^{sm} $$
is an equivalence.  Since geometric realization is the left adjoint of a
Quillen equivalence, we see that the adjoint
$$ X \xrightarrow{\simeq} S\Map(\mc{B}, \abs{X})^{sm} $$
is an equivalence.  We have the following sequence of isomorphisms:
\begin{align*}
S\Map(\abs{\mc{B}'_\bullet}, \abs{X})^{sm}
& = S \varinjlim_{U \le_o G} \Map_U(\abs{\mc{B}'_\bullet}, \abs{X}) \\
& \cong \Map^{\Top}(\abs{\Delta^\bullet}, 
\varinjlim_{U \le_o G} \Map_U(\abs{\mc{B}'_\bullet}, \abs{X})) \\
& \cong  
\varinjlim_{U \le_o G} 
\Map^{\Top}(\abs{\Delta^\bullet},\Map_U(\abs{\mc{B}'_\bullet}, \abs{X})) \\
& \cong  
\varinjlim_{U \le_o G} 
\Map^{\Top}_U(\abs{\Delta^\bullet} \times \abs{\mc{B}'_\bullet}, \abs{X}) \\
& \cong  
\varinjlim_{U \le_o G} 
\Map^{\Top}_U(\abs{\Delta^\bullet \times \mc{B}'_\bullet}, \abs{X}) \\
& \cong  
\varinjlim_{U \le_o G} 
\Map^{s\Set}_U(\Delta^\bullet \times \mc{B}'_\bullet, S\abs{X}) \\
& =   
\ul{\Map}(\mc{B}'_\bullet, S\abs{X})^{sm}.
\end{align*}
This isomorphism fits into the following commutative diagram.
$$
\xymatrix{
X \ar[r] \ar[d]_\simeq &
\ul{\Map}(\mc{B}'_\bullet, X)^{sm} \ar[d]^{\eta_*} 
\\
S\Map(\mc{B}, \abs{X})^{sm} \ar[r]_\cong &
\ul{\Map}(\mc{B}'_\bullet, S\abs{X})^{sm}
}
$$
We just need to verify that the map $\eta_*$ is a weak equivalence.
For any simplicial smooth $G$-set $Z$,
the functor $- \times Z$ preserves weak equivalences and cofibrations.
Therefore,
the functors $(- \times Z, \ul{\Map}(Z, -)^{sm})$ form a Quillen adjoint
pair on $s\Set_G^{sm}$.  In particular, the functor
$\ul{\Map}(\mc{B}'_\bullet, -)^{sm}$ preserves weak equivalences
between fibrant simplicial smooth $G$-sets.  By Lemma~\ref{lem:realizationfibrant}, 
the map
$$ \eta: X \rightarrow S\abs{X} $$
is a weak equivalence between fibrant simplicial smooth $G$-sets.
\end{proof}

\begin{rmk}
The authors do not know a purely simplicial argument to prove
Lemma~\ref{lem:technicallemma}.  Clearly, much of the work in this section
could be eliminated by working with spectra of topological
spaces with smooth $G$-action as opposed to 
spectra of simplicial smooth $G$-sets.  Our reason for choosing to work
simplicially is that we are simply 
unaware of a treatment of local model structures on categories of sheaves
of topological spaces in the literature.
\end{rmk}

\subsection*{The semi-cosimplicial resolution}

Recall from the beginning of this section that $Q$ is the hypercohomology
spectrum $Q_{GU}(K^{p,\ell})$ or $Q_U(K^{p,\ell})$.
We state our main theorem describing a finite semi-cosimplicial resolution
of $Q$.

\begin{thm}\label{thm:mainthm}
There is a
semi-cosimplicial spectrum $Q^\bullet$ 
\index{2Q@$Q^\bullet$}
of length $d = \dim \mc{B}$
\index{2D@$d$}
whose $s$th term ($0 \le s \le d$)
is given by
$$ Q^s = \prod_{[\sigma]} \TAF(K(\sigma)). $$
The product ranges over $G_\ell$ orbits of 
$s$-simplices $[\sigma]$ in the building $\mathcal{B}$.  
The groups $K(\sigma)$ 
\index{2Ks@$K(\sigma)$}
are given by
$$ K(\sigma) = K^{p,\ell} K_\ell(\sigma) $$
where $K_\ell(\sigma)$ 
\index{2Kls@$K_\ell(\sigma)$}
is the subgroup of $G_\ell$ which stabilizes
$\sigma$.
There is an equivalence
$$ Q \simeq \holim Q^\bullet. $$
\end{thm} 

\begin{proof}
By Lemma~\ref{lem:technicallemma}, the map
$$ r : V_{GU} \rightarrow \ul{\Map}(\mc{B}'_\bullet, V_{GU})^{sm} $$
is an equivalence of smooth $G$-spectra.  
The functor
$$ \ul{\Map}(\mc{B}'_\bullet, -)^{sm} : \Sp_G^{sm} \rightarrow \Sp_G^{sm} $$
is right Quillen adjoint to the functor $- \wedge (\mc{B}'_\bullet)_+$. 
Therefore, since $V_{GU}$ is a fibrant smooth $G$-spectrum, the spectrum
$\ul{\Map}(\mc{B}'_\bullet, V_{GU})^{sm}$ is a fibrant smooth $G$-spectrum.
We have the following sequence of $G$-equivariant isomorphisms.
\begin{align*}
\ul{\Map}(\mc{B}'_\bullet, V_{GU})^{sm}  
& = \varinjlim_{U \le_o G} \ul{\Map}_U(\mc{B}'_\bullet, V_{GU}) \\
& \cong \varinjlim_{U \le_o G} \holim \Map_U(\mc{B}'_\bullet, V_{GU}) \\
& \cong \varinjlim_{U \le_o G} \holim \Map_U(\mc{B}_\bullet, V_{GU}) \\
& \cong \holim \varinjlim_{U \le_o G} \Map_U(\mc{B}_\bullet, V_{GU}) \\
& \cong \holim \Map(\mc{B}_\bullet, V_{GU})^{sm}.
\end{align*}
Here, we were able to commute the homotopy limit past the colimit because 
it is the homotopy limit of a \emph{finite length} semi-cosimplicial
spectrum.  
Taking $G$-fixed points of both sides, we have a map
$$ r^G: V_{GU}^{G} \rightarrow (\holim \Map(\mc{B}_\bullet,
V_{GU})^{sm})^G. $$
The map $r^G$ is a weak equivalence because $r$ is a weak equivalence
between fibrant smooth $G$-spectra.  
There are isomorphisms
\begin{align*}
(\holim \Map(\mc{B}_\bullet, V_{GU})^{sm})^G
& \cong \holim \Map_G(\mc{B}_\bullet, V_{GU}) \\
& \cong \holim \prod_{[\sigma] \in \mc{B}_\bullet/G} \Map_G(G/K(\sigma), V_{GU}) 
\\
& \cong \holim \prod_{[\sigma] \in \mc{B}_\bullet/G} V_{GU}^{K(\sigma)}.
\end{align*}
By Lemma~\ref{lem:cofixedpoints}, there are equivalences
$$ V_{GU}^{K(\sigma)} \simeq \TAF(K(\sigma)). $$
\end{proof}

\begin{rmk}
The coface maps of the semi-cosimplicial spectrum $Q^\bullet$ are all
instances of the $E_\infty$ operations arising from
Proposition~\ref{prop:Einftyops}.  Thus, the spectrum $Q^\bullet$ is
actually a semi-cosimplicial $E_\infty$-ring spectrum, and the
totalization $Q$ therefore inherits the structure of an $E_\infty$-ring
spectrum.
\end{rmk}

\chapter{$K(n)$-local theory}\label{chap:K(n)local}

Fix a compact open subgroup $K^p$ of $GU(\AF^{p,\infty})$ so that
$\Sh(K^p)$ is a scheme.  
Let 
$$ \mbf{A}_{univ} = (A_{univ}, i_{univ}, \lambda_{univ}, [\eta_{univ}]) $$ 
be the universal
tuple over $\Sh(K^p)$.
Let
$\Sh(K^p)^{[n]}$ be the reduced closed subscheme of $\Sh(K^p) \otimes_{\ZZ_p}
\FF_p$ 
\index{2ShKpn@$\Sh(K^p)^{[n]}$}
where the the formal group
$\epsilon A_{univ}(u)^0$ has height $n$ 
(see \cite[Lem.~II.1.1]{harristaylor}).  

\section{Endomorphisms of mod $p$ points}

Suppose that 
\index{2A@$\mbf{A}$}
$\mbf{A} = (A,i,\lambda, [\eta])$ is an element of 
$\Sh(K^p)^{[n]}(\br{\FF}_p)$.

We make the following definitions.
\begin{align*}
D & = \End^0_B(A), \\
\mc{O}_D & = \End_B(A), \\
\dag & = \text{$\lambda$-Rosati involution on $D$}.
\end{align*}
By \index{2D@$D$} \index{2OD@$\mc{O}_D$} \index{0dag@$\dag$}
Theorem~\ref{thm:Bhondatate}, $D$ is a central simple algebra over $F$
of dimension $n^2$,
with invariants given by:
\begin{alignat*}{2}
\inv_x D & = -\inv_x B, & \qquad & x \not | p, \\
\inv_u D & = 1/n, & & \\
\inv_{\br{u}} D & = (n-1)/n.
\end{alignat*}
The ring $\mc{O}_D$ is an order in $D$.  Since the $p$-divisible
$\mc{O}_B$-module of $A$ takes the form
$$ A(p) \cong (\epsilon A(u))^n \times (\epsilon A(u^c)^n), $$
Theorems~\ref{thm:pdivclass}(4) and \ref{thm:tate} combine to show that the
order $\mc{O}_{D,(p)}$ is maximal at $p$.  Because the polarization
$\lambda$ is prime-to-$p$, the order $\mc{O}_{D,(p)}$ is preserved by the
Rosati involution $\dag$.  

Define algebraic groups $GU_\mbf{A}$ 
\index{2GUA@$GU_\mbf{A}$}
and $U_\mbf{A}$ 
\index{2UA@$U_\mbf{A}$}
over $\ZZ_{(p)}$ by
\begin{align*}
GU_\mbf{A}(R) & = \{ g \in (\mc{O}_{D,(p)} \otimes_{\ZZ_{(p)}} R)^\times \: : \:
g^\dag g \in R^\times \}, \\
U_\mbf{A}(R) & = \{ g \in (\mc{O}_{D,(p)} \otimes_{\ZZ_{(p)}} R)^\times \: : \:
g^\dag g = 1 \}.
\end{align*}
There is a short exact sequence
$$ 1 \rightarrow U_\mbf{A} \rightarrow GU_\mbf{A} \xrightarrow{\nu} \GG_m
\rightarrow 1. $$
Observe that we have (Lemma~\ref{lem:isometrysimilitude})
\begin{align*}
GU_\mbf{A}(\ZZ_{(p)}) = & \text{group of prime-to-$p$ 
quasi-similitudes of $(A,i,\lambda)$}, \\
U_\mbf{A}(\ZZ_{(p)}) = & \text{group of prime-to-$p$ 
quasi-isometries of $(A,i,\lambda)$}.
\end{align*}
The rational uniformization induces isomorphisms
\begin{align*}
\eta_*: GU(\QQ_\ell) \xrightarrow{\cong} GU_\mbf{A}(\QQ_\ell), \\
\eta_*: U(\QQ_\ell) \xrightarrow{\cong} U_\mbf{A}(\QQ_\ell),
\end{align*}
for \index{1etaast@$\eta_*$}
primes $\ell \ne p$.  
The maximality of the $\mc{O}_{F,(p)}$-order
$\mc{O}_{D,(p)}$, together with Lemma~\ref{lem:splitF}, gives the following
lemma.

\begin{lem}\label{lem:Z_ppoints}
The induced action of $U_\mbf{A}(\ZZ_{(p)})$ on 
the summand $\epsilon A(u)$ of the
$p$-divisible group $A(p)$ induces an isomorphism
$$
U_\mbf{A}(\ZZ_p) \xrightarrow{\cong} \MS_n 
$$
where $\MS_n = \Aut(\epsilon A(u))$ 
\index{2Sn@$\MS_n$}
is the $n$th Morava stabilizer group.
The similitude norm 
gives a split short exact sequence
$$ 1 \rightarrow U_{\mbf{A}}(\ZZ_p) \rightarrow GU_{\mbf{A}}(\ZZ_p) 
\xrightarrow{\nu}
\ZZ_p^\times \rightarrow 1. $$
\end{lem}

Define $\Gamma$ 
\index{1gamma@$\Gamma$}
to be the quasi-isometry group $U_\mbf{A}(\ZZ_{(p)})$.
The action of $\Gamma$ on the Tate module $V^{p,S}(A)$ induces an inclusion
$$ i_\eta: \Gamma \hookrightarrow U(\AF^{p,S,\infty}) \subset
GU(\AF^{p,\infty}). $$
For \index{2Ieta@$i_\eta$}
any subgroup $K$ of $GU(\AF^{p,S,\infty})$ let $\Gamma(K)$ 
\index{1gammaK@$\Gamma(K)$}
be the
subgroup of $\Gamma$ given by the intersection $\Gamma \cap K$.

\begin{prop}\label{prop:aut}
Suppose that $K^p$ is an arbitrary open compact subgroup of
$GU(\AF^{p,\infty})$, so that $\Sh(K^p)$ is not necessarily a scheme.  Then
the automorphisms of the $\br{\FF}_p$-point $\mbf{A} = 
(A,i,\lambda, [\eta]_{K^p})$
are given by $\Gamma(K^p)$.
\end{prop}

\begin{proof}
By definition, we have
$$ \Aut(\mbf{A}) = GU_\mbf{A}(\ZZ_{(p)}) \cap K^p
\subset GU(\AF^{p,\infty}) $$
(using the isomorphism $\eta_*: GU(\AF^{p,\infty}) \xrightarrow{\cong}
GU_\mbf{A}(\AF^{p,\infty})$).  The similitude norm restricts to give
a homomorphism
$$ \nu: \Aut(\mbf{A}) \rightarrow \{\pm 1\} = (\widehat{\ZZ}^p)^\times \cap
\ZZ_{(p)}^\times \subset (\AF^{p,\infty})^\times. $$
However, by the positivity of the Rosati involution
(Theorem~\ref{thm:Rosatipos}), the similitude norm $\nu$ cannot be
negative.  We therefore deduce
$$ \Aut(\mbf{A}) = U_\mbf{A}(\ZZ_{(p)}) \cap K^p = \Gamma(K^p). $$
\end{proof}

\section{Approximation results}

In this section we compile various approximation results that we shall
appeal to in later sections.  These results both will allow us to compare the
spectrum $Q_U(K^{p,S})_{K(n)}$ with the $K(n)$-local sphere, 
as well as to manipulate certain adelic
quotients.

For each prime $x \ne p$, let $K_{1,x} < GU(\QQ_x)$ be the image of the subgroup
$$ \{ g \in \mc{O}_{D,x}^\times \: : \: g^\dag g \in \ZZ_\ell^\times \} <
GU_\mbf{A}(\QQ_\ell) $$
under the isomorphism $\eta_*^{-1}: GU_\mbf{A}(\QQ_x) \xrightarrow{\cong}
GU(\QQ_x)$.
For a set of primes $S$ not containing $p$, define the group $K_1^{p,S}$ by 
$$ K_1^{p,S} = \prod_{x \not\in \{p\} \cup S} K_{1,x}. $$
Naumann \cite[Cor.~21, Rmk.~22]{naumann} proves the following theorem, 
quantifying, at least in certain situations, the degree to which the 
group $\Gamma(K_1^{p,S})$ approximates the Morava stabilizer group $\MS_n$.

\begin{thm}[Naumann]\label{thm:naumann}
Suppose that
\begin{enumerate}
\item The polarization $\lambda$ is principal,
\item The order $\mc{O}_D$ is maximal.
\end{enumerate}
Then we have the following.
\begin{description}
\item[Case $\mathbf{n}$ odd]
There exists a prime $\ell \ne p$ which splits in $F$ so that the group
$\Gamma(K_1^{p,\ell})$ is dense in $\MS_n$.

\item[Case $\mathbf{n}$ even]
Suppose that $p \ne 2$.
Then there exists a prime $\ell$ which splits in 
$F$ such that the closure of the group $\Gamma(K_1^{p,2,\ell})$ in $\MS_n$
is of index less than or equal to the order of the unit group 
$\mc{O}_F^\times$.  
\end{description}
\end{thm}

\begin{rmk}
Since $F$ is a quadratic imaginary extension, $\mc{O}_F^\times$
is of order $2$, $4$, or $6$.
\end{rmk}

A key observation of Naumann is the following proposition.

\begin{prop}\label{prop:SESdense}
Suppose that we have a map of short exact sequences groups
$$
\xymatrix{
1 \ar[r] &
H' \ar[r] \ar@{^{(}->}[d] &
H \ar[r] \ar@{^{(}->}[d] &
H'' \ar[r] \ar@{^{(}->}[d] &
1 \\
1 \ar[r] &
G' \ar[r] &
G \ar[r]_\pi &
G'' \ar[r] &
1
}
$$
where the bottom row is a short exact sequence of first countable
topological groups.  Assume that $H'$ is dense in $G'$ and that 
there exists an open subgroup $U$ of $G''$
so that $\pi^{-1}(U)$ is compact in $G$.  Then $H$ is dense in
$\pi^{-1}(\br{H}'')$, where $\br{H}''$ is the closure of $H''$ in $G''$.
\end{prop}

Naumann's methods may be used to prove the following easier proposition.

\begin{prop}\label{prop:plocaldense}$\quad$
\begin{enumerate}
\item The group $U_\mbf{A}(\QQ)$ is dense in $U_\mbf{A}(\QQ_p)$.
\item The group $\Gamma = U_\mbf{A}(\ZZ_{(p)})$ is dense in $\MS_n =
U_\mbf{A}(\ZZ_p)$.
\item The group $GU_\mbf{A}(\QQ)$ is dense in $GU_\mbf{A}(\QQ_p)$.
\end{enumerate}
\end{prop}

\begin{proof}  We prove statements (1) and (2) simultaneously.
Let $SU_\mbf{A}$ be the kernel of the reduced norm:
$$ 1 \rightarrow SU_\mbf{A} \rightarrow U_\mbf{A} \xrightarrow{N_{D/F}} T
\rightarrow 1 $$
where the algebraic group $T/\ZZ_{(p)}$ is given by
$$ T(R) = \{ t \in (\mc{O}_{F,(p)} \otimes_{\ZZ_{(p)}} R)^\times \: : \:
N_{F/\QQ}(t) = 1 \}. $$
Using the fact that there is a pullback \cite[1.4.2]{platonov}
\begin{equation}\label{diag:pullbackO_Du}
\xymatrix{
\mc{O}_{D,u} \ar[r]^{N_{D/F}} \ar[d] &
\mc{O}_{F,u} \ar[d] \\
D_u \ar[r]_{N_{D/F}} &
F_u 
}
\end{equation}
we deduce that the following diagram is a pullback.
\begin{equation}\label{diag:pullbackU}
\xymatrix{
U_\mbf{A}(\ZZ_{(p)}) \ar[r]^{N_{D/F}} \ar[d] &
T(\ZZ_{(p)}) \ar[d] \\
U_\mbf{A} (\QQ) \ar[r]_{N_{D/F}} &
T(\QQ) 
}
\end{equation}

The weak approximation theorem \cite[Lemma~7.2]{platonov} 
implies that the embedding
$$ SU_\mbf{A}(\QQ) \hookrightarrow SU_\mbf{A}(\QQ_p) $$
is dense.  However, the pullback in diagram~(\ref{diag:pullbackO_Du}) implies
that $SU_\mbf{A}(\ZZ_p) = SU_\mbf{A}(\QQ_p)$ and the pullback in
diagram~(\ref{diag:pullbackU}) implies that $SU_\mbf{A}(\ZZ_{(p)}) = 
SU_\mbf{A}(\QQ)$.  So
we actually have determined that the embedding
$$ SU_\mbf{A}(\ZZ_{(p)}) \hookrightarrow SU_\mbf{A}(\ZZ_p) $$
is dense.
In \cite{naumann}, it is established that the
following maps are surjections:
\begin{align*}
N_{D/F}: & U_\mbf{A}(\QQ) \rightarrow T(\QQ), \\
N_{D/F}: & U_\mbf{A}(\ZZ_p) \rightarrow T(\ZZ_p).
\end{align*}
The pullback of Diagram~(\ref{diag:pullbackU}) implies that the  map
$$
N_{D/F}: U_\mbf{A}(\ZZ_{(p)}) \rightarrow T(\ZZ_{(p)}) 
$$
is a surjection.  We therefore have the following 
diagrams of short exact sequences.
$$
\xymatrix{
1 \ar[r] & 
SU_\mbf{A}(\ZZ_{(p)}) \ar[r] \ar[d] &
U_\mbf{A}(\ZZ_{(p)}) \ar[r] \ar[d] &
T(\ZZ_{(p)}) \ar[r] \ar[d] &
1 \\
1 \ar[r] & 
SU_\mbf{A}(\ZZ_{p}) \ar[r] &
U_\mbf{A}(\ZZ_{p}) \ar[r] &
T(\ZZ_{p}) \ar[r] &
1
\\
1 \ar[r] & 
SU_\mbf{A}(\QQ) \ar[r] \ar[d] &
U_\mbf{A}(\QQ) \ar[r] \ar[d] &
T(\QQ) \ar[r] \ar[d] &
1 \\
1 \ar[r] & 
SU_\mbf{A}(\QQ_{p}) \ar[r] &
U_\mbf{A}(\QQ_{p}) \ar[r] &
T(\QQ_{p}) \ar[r] &
1
}
$$
The groups $T(\ZZ_{(p)})$ and $T(\QQ)$ are dense in $T(\ZZ_p)$ and
$T(\QQ_p)$, respectively (see, for instance,
\cite{naumann}).  Therefore, (1) and (2) follow from
Proposition~\ref{prop:SESdense}.  (To verify that the second diagram above
satisfies the hypotheses of Proposition~\ref{prop:SESdense} we need to
again appeal to the pullback (\ref{diag:pullbackO_Du}).)

By (1) and Proposition~\ref{prop:SESdense}, to prove (3) it suffices to
prove that the image of the similitude norm
$$ \nu: GU_{\mbf{A}}(\QQ) \rightarrow \QQ^\times $$
is dense in $\QQ_p^\times$.  Because the similitude norm restricts to the
norm $N_{F/\QQ}$ on the subgroup $F^\times \le GU(\QQ)$, it suffices to
prove that the image of the norm
$$ N_{F/\QQ}: F^\times \rightarrow \QQ^\times $$
is dense in $\QQ_p^\times$.  Let $D$ be the absolute value of the 
discriminant of $F$, and let 
$$ \chi: (\ZZ/D)^\times \rightarrow \{ \pm 1\} $$
be the corresponding Dirichlet character, so that a prime $q$ splits in $F$
if and only if $\chi(q) = 1$.

Assume that $p$ is odd.  Fix a prime $\ell$ which splits in $F$
and is a generator of $\ZZ_p^\times$.  Such a prime exists because the
former represents a congruence condition modulo $D$, whereas the latter
represents a congruence condition modulo $p^2$, and since $p$ was assumed
to split, it is coprime to $D$.  To prove the image of
$N_{F/\QQ}$ is dense in $\QQ_p^\times$, it suffices to prove that $\ell$ and $p$
are in the image.  By the fundamental short exact sequence
$$ 1 \rightarrow \QQ^\times/N(F^\times) \rightarrow \bigoplus_{x}
\QQ_x^\times/N(F_x^\times) \rightarrow \ZZ/2 \rightarrow 1 $$
it suffices to prove that for all non-split $q$, both $p$ and $\ell$ are
zero in the group $\QQ_q^\times/N(F_q^\times)$.  If $q$ is inert in $F$, then this
follows from the fact that $p$ and $\ell$ are coprime to $q$.  If $q$ is
ramified, this follows from the fact that the kernel of the composite
$$ \ZZ_q^\times \rightarrow (\ZZ/D)^\times \xrightarrow{\chi} \{ \pm 1 \}
$$
is equal to $N(F_q^\times) \cap \ZZ_q^\times$.
The case of $p = 2$ is similar, but because $\ZZ_2^\times$ is not cyclic,
two generating split primes must be used instead.
\end{proof}

\begin{lem}\label{lem:valuationimage}
The images of the composites
\begin{gather*}
GU_\mbf{A}(\QQ) \xrightarrow{\nu} \QQ^\times \xrightarrow{\oplus \nu_\ell}
\bigoplus_{\text{$\ell$ prime}} \ZZ 
\\
GU_\mbf{A}(\AF^\infty) \xrightarrow{\nu} (\AF^\infty)^\times 
\xrightarrow{\oplus \nu_\ell}
\bigoplus_{\text{$\ell$ prime}} \ZZ 
\end{gather*}
are equal.
\end{lem}

\begin{proof}
The result follows from applying Galois cohomology computations of
Section~\ref{sec:formclassification} to the map of exact sequences.
\begin{equation}\label{diag:SESgaloisH^1}
\xymatrix{
GU_\mbf{A}(\QQ) \ar[r]^\nu \ar[d] &
\QQ^\times \ar[r] \ar[d] &
H^1(\QQ, U_\mbf{A}) \ar[r] \ar[d] &
H^1(\QQ, GU_\mbf{A}) \ar[d] 
\\
GU_\mbf{A}(\AF^\infty) \ar[r]_\nu &
(\AF^\infty)^\times \ar[r] &
H^1(\AF^\infty, U_\mbf{A}) \ar[r] &
H^1(\AF^\infty, GU_\mbf{A}) 
}
\end{equation}
If $n$ is even, then the image of $\nu: GU(\QQ) \rightarrow \QQ^\times$ is
seen to be $(\QQ^\times)^+$, and the theorem is clear.
If $n$ is odd, then 
we have the following map of exact sequences.
$$
\xymatrix@C+1em{
GU_\mbf{A}(\QQ) \ar[r]^{\oplus \nu_\ell \nu} \ar[d] &
\bigoplus\limits_{\text{$\ell$ prime}} \ZZ \ar[r] \ar@{=}[d] &
\bigoplus\limits_{\substack{\text{$\ell$ inert and} \\ \text{unramified in
$F$}}} \ZZ/2 \ar@{=}[d] \ar[r] &
0
\\
GU_\mbf{A}(\AF^\infty) \ar[r]_{\oplus \nu_\ell \nu} &
\bigoplus\limits_{\text{$\ell$ prime}} \ZZ \ar[r] &
\bigoplus\limits_{\substack{\text{$\ell$ inert and} \\ \text{unramified in
$F$}}} \ZZ/2 \ar[r] &
0
}
$$
\end{proof}

\section{The height $n$ locus of $\Sh(K^p)$}

Let $K^p$ be sufficiently small so that $\Sh(K^p)$ is a scheme.
Combining Corollary~\ref{cor:defsh} with Corollary~II.1.4 of
\cite{harristaylor}, we see that
the subscheme $\Sh(K^p)^{[n]}$ 
is either empty, or smooth of dimension zero.  
We therefore have the following lemma.

\begin{lem}
The scheme $\Sh(K^p)^{[n]}$ is \'etale over $\Spec(\FF_p)$.
\end{lem}

The scheme $\Sh(K^p)^{[n]}$ is therefore affine.  Its structure is
completely determined by 
\begin{enumerate}
\item
the set of $\br{\FF}_p$-points
$\Sh(K^p)^{[n]}(\br{\FF}_p)$,

\item 
the action of the Galois group $Gal$ of $\FF_p$ on 
$$ \Sh(K^p)^{[n]} \otimes_{\FF_p} \br{\FF}_p =
\coprod_{\Sh(K^p)^{[n]}(\br{\FF}_p)} \Spec(\br{\FF}_p). $$
\end{enumerate}
Item (1) above is described by Corollary~\ref{cor:Fpbarpoints}.  Item
(2) is a serious arithmetic question related to the zeta function of
the Shimura variety --- we do not investigate it here.

\subsection*{Existence of mod $p$ points}

The following proposition appears in \cite[Cor~V.4.5]{harristaylor} in the
case where $B$ is a division algebra.

\begin{prop}\label{prop:Fpbarnonempty}
The set $\Sh(K^p)^{[n]}(\br{\FF}_p)$ is non-empty.  
\end{prop}

\begin{proof}
By Theorem~\ref{thm:Bhondatate}, there exists a $B$-linear abelian variety
$(A,i)$ over $\br{\FF}_p$ of dimension $n^2$ associated to the minimal
$p$-adic type $(F,\zeta)$ where
\begin{align*}
\zeta_u & = 1/n, \\
\zeta_{u^c} & = (n-1)/n.
\end{align*}

The $p$-adic type of $(A, i)$ determines the slopes of the $p$-divisible
group $A(p)$, and there is therefore an isogeny of $B$-linear 
$p$-divisible groups
$$ \phi: (A(p),i_*) \rightarrow (\GG, i') $$
such that
$$ \GG = (\epsilon \GG(u))^n \oplus (\epsilon \GG(u^c))^n. $$
Since the $p$-completion of the order $\mc{O}_B$ is given by
$$ \mc{O}_{B,p} = M_n(\mc{O}_{F,u}) \times M_n(\mc{O}_{F,u^c}), $$
the inclusion of rings
$$ i': B \hookrightarrow \End^0(\GG) $$ 
lifts to an inclusion
$$ i': \mc{O}_{B,(p)} \hookrightarrow \End(\GG). $$
Since $\phi$ is an isogeny, there exists a $k$ such that $\ker \phi$ is
contained in the finite group-scheme $A[p^k]$ of $p^k$-torsion points of
$A$.  We define an isogenous $B$-linear abelian variety $(A',i')$ 
by taking the
quotient
$$ (A,i) \rightarrow (A/\ker \phi,i') = (A', i'). $$
There is a canonical isomorphism of $B$-linear $p$-divisible groups
$$ (A'(p), (i')_*) \cong (\GG, i'). $$
By Theorem~\ref{thm:tate}, the inclusion of rings
$$ i': B \hookrightarrow \End^0(A') $$
lifts to an inclusion
$$ i': \mc{O}_{B,(p)} \hookrightarrow \End(A')_{(p)}. $$
Choose a compatible polarization 
$$ \lambda: A' \rightarrow (A')^\vee, $$
and let $\lambda_*: A(p) \rightarrow A(p)^\vee$ be the induced map of
$p$-divisible groups.
(According to
Lemma~\ref{lem:polclassglobal}, compatible polarizations exist.)
Since $\lambda$ is compatible, there exist finite subgroups $K_u < A(u)$
and $K_{u^c} < A(\br{u})$ such that
$$ \ker \lambda_* = K_u \oplus K_{u^c} < A(u) \oplus A(u^c) = A(p). $$
Define $A''$ to be the quotient $A'/K_u$, with quotient isogeny $q: A'
\rightarrow A''$.  Compatibility of the polarization implies that the
$p$-divisible group $A'(u)/K_u \oplus A(\br{u})$ inherits a 
$\mc{O}_{B,(p)}$-linear
structure, giving an inclusion
$$ i'': \mc{O}_{B,(p)} \hookrightarrow \End(A''). $$
The quasi-isogeny
$$ \lambda' = (q^{-1})^\vee \lambda q^{-1}: A'' \rightarrow (A'')^\vee $$
is easily seen to be prime-to-$p$, and thus gives a prime-to-$p$
compatible polarization on $(A'',i'')$.

In order to produce a point of $\Sh(K^p)^{[n]}(\br{\FF}_p)$, we need to
show that the polarization $\lambda'$ can be altered so that there exists a
similitude between the
Weil pairing on the Tate module $V^p(A'')$ and our fixed pairing
$\bra{-,-}$ on $V^p$.    If $n$ is odd, then for all $\ell$, $H^1(\QQ_\ell, GU) = 0$
(Lemma~\ref{lem:localsplitH^1GU} and
Corollary~\ref{cor:localnonsplitH^1GU}), and so \emph{any} two
non-degenerate $*$-hermitian alternating forms on $V_\ell \cong
V_{\ell}(A'')$ are similar, so there exists a uniformization
$$ \eta: (V^p, \bra{-,-}) \xrightarrow{\simeq} (V^p(A''),
\lambda'\bra{-,-}). $$

We must work harder if $n$ is even. 
Using Theorem~\ref{thm:defsh}, there exists a
deformation $(\td{A}'', \td{i}'', \td{\lambda}')$ of the tuple  $(A'',
i'', \lambda')$ over the Witt ring $W(\br{\FF}_p)$.  Choose an embedding of
the field
$F_u^{nr} \cong W(\br{\FF}_p) \otimes \QQ$ in $\CC$ so that the 
following diagram
of fields commutes.
$$
\xymatrix{
F \ar[r] \ar[d] & F_u \ar[d] \\
\CC & F_u^{nr} \ar[l]
} 
$$
Under this inclusion, we may pull back $(\td{A}'', \td{i}'', \td{\lambda})$
to 
a polarized $\mc{O}_{B,(p)}$-linear abelian variety 
$(A''_\CC, i''_\CC, \lambda'_\CC)$
over $\CC$.  The methods of Section~\ref{sec:complexpoints} imply that
there is a non-degenerate $*$-hermitian alternating form $\bra{-,-}'$ on
$V$, a lattice $L' \subset V$, and a compatible 
complex structure $J$ on $V_\infty = V
\otimes_\QQ \RR$ so that:
\begin{align*}
A''_\CC & = V_\infty/L', \\
i''_\CC & = \text{induced from $B$-module structure of $V$}, \\
\lambda''_\CC & = \text{polarization associated to Riemann form
$\bra{-,-}'$.}
\end{align*}
Thus, there is a canonical uniformization
$$ \eta_1: (V^p, \bra{-,-}') \xrightarrow{\simeq} (V^p(A''_\CC), 
\lambda'_\CC\bra{-,-}).
$$
Furthermore, there are isomorphisms (using proper-smooth base change):
\begin{align*}
(V^p(A''_\CC), \lambda'_\CC\bra{-,-}) 
& \cong (H^1_{et}(A''_\CC, \AF^{p,\infty})^*, \lambda'_\CC\bra{-,-}) \\
& \cong (H^1_{et}(\td{A}'', \AF^{p,\infty})^*, \td{\lambda}'\bra{-,-}) \\
& \cong (H^1_{et}(A'', \AF^{p,\infty})^*, \lambda'\bra{-,-}) \\
& \cong (V^p(A''), \lambda'\bra{-,-}).
\end{align*}

Let $GU_{\mbf{A}''}$ be the group of quasi-similitudes of
$\mbf{A}'' = (A'',i'',\lambda')$.
By Lemma~\ref{lem:polclassglobal},
any other isogeny class of weak polarization $\lambda''$ 
is determined by an element $[\lambda''] \in H^1(\QQ, GU_{\mbf{A}''})$ 
such that
$$ [\lambda']_\infty = [\lambda'']_\infty \in H^1(\RR, GU_{\mbf{A}''}). $$ 
The calculations of Section~\ref{sec:formclassification} give isomorphisms:
\begin{align*}
\mathrm{disc}_{\dag'}: H^1(\QQ_\ell, GU_{\mbf{A}''}) & \xrightarrow{\cong}
\QQ_\ell^\times/N_{F/\QQ}(F^\times), \\
\mathrm{disc}_{*}: H^1(\QQ_\ell, GU) & \xrightarrow{\cong}
\QQ_\ell^\times/N_{F/\QQ}(F^\times).
\end{align*}
In particular, $H^1(\QQ_p, GU) = H^1(\QQ_p, GU_{\mbf{A}''}) = 0$. 
Here, $\mathrm{disc}_{\dag'}$ denotes the discriminant taken relative to
the $\lambda'$-Rosati involution $(-)^{\dag'}$, and $\mathrm{disc}_*$ 
denotes the
discriminant taken relative to the involution $(-)^*$.  We identify the
Galois cohomologies $H^1(\QQ_\ell, GU_{\mbf{A}''})$ and $H^1(\QQ_\ell,
GU)$ in this manner, and use these isomorphisms to endow these isomorphic
sets with the structure of a group.
By Lemma~\ref{lem:polclasslocal}, the image of $[\lambda'']$ in
$H^1(\QQ_\ell, GU)$ is given by
$$ [\lambda'']_\ell = [\lambda'\bra{-,-}_\ell] + [\lambda''\bra{-,-}_\ell], $$
the difference of the classes represented by the Weil pairings on
$V_\ell(A'')$.
Let $[\bra{-,-}']$ be the element of $H^1(\QQ, GU)$ representing the
similitude class of the form $\bra{-,-}'$.   
We wish to demonstrate that there is such a class $[\lambda'']$ so that
$$ [\lambda''\bra{-,-}_\ell] = [\bra{-,-}]_\ell \in H^1(\QQ_\ell, GU) $$
for every $\ell \ne p$.  There are short exact sequences
(Corollary~\ref{cor:globalH^1GU})
\begin{gather*}
0 \rightarrow H^1(\QQ, GU) \rightarrow \bigoplus_x H^1(\QQ_x, GU) 
\xrightarrow{\sum \xi'_x}
\ZZ/2 \rightarrow 0, \\
0 \rightarrow H^1(\QQ, GU_{\mbf{A}''}) \rightarrow \bigoplus_x H^1(\QQ_x,
GU_{\mbf{A}''}) 
\xrightarrow{\sum \xi'_x}
\ZZ/2 \rightarrow 0.
\end{gather*}
By the positivity of the Rosati involution (Theorem~\ref{thm:Rosatipos}),
we have $\xi'_\infty[\lambda'] = 0$.
It therefore suffices to show that
$$ \sum_\ell (\xi'_\ell([\bra{-,-}]_\ell) +
\xi'_\ell([\lambda'\bra{-,-}_\ell)]) = 0. $$
Since the complex structure $J$ on $V_\infty$ given above is compatible
with $\bra{-,-}'$, the form $\bra{-,-}'$ has signature $\{1, n-1\}$, and
hence, by Lemma~\ref{lem:localinfiniteH^1GU}, we have
$$ [\bra{-,-}]_\infty = [\bra{-,-}']_\infty \in H^1(\RR, GU) $$
We compute
\begin{align*}
\sum_\ell \xi'_\ell([\bra{-,-}]_\ell) + 
\sum_\ell \xi'_\ell([\lambda'\bra{-,-}_\ell])
& = \sum_\ell \xi'_\ell([\bra{-,-}]_\ell) + 
\sum_\ell \xi'_\ell([\bra{-,-}']_\ell) \\
& = \xi'_\infty([\bra{-,-}]_\infty) + \xi'_\infty([\bra{-,-}']_\infty) \\
& = 0.
\end{align*}
Thus there exist a class $[\lambda''] \in H^1(\QQ, GU)$ with
$[\lambda''\bra{-,-}_\ell] = [\bra{-,-}]_\ell$ and $[\lambda'']_\infty =
[\lambda']_\infty$.  There exists a corresponding polarization $\lambda''$.
Using the same methods used to
construct $\lambda'$, we may assume without loss of generality that
$\lambda''$ is prime-to-$p$.
\end{proof}

\subsection*{Calculation of $\Sh(K^p)^{[n]}(\br{\FF}_p)$}

We shall make use of the following lemma.

\begin{prop}
Given any two tuples
$$ (A, i, \lambda, [\eta]), (A', i', \lambda', [\eta']) \in 
\Sh(K^p)^{[n]}(\br{\FF}_p)$$  
there exists a
prime-to-$p$ isogeny 
$$ (A, i, \lambda) \rightarrow (A', i', \lambda') $$
of weakly polarized $B$-linear abelian varieties.
\end{prop}

\begin{proof}
The completion $\mc{O}_{D,u}$ is the unique maximal order of
the central $\QQ_p$-division algebra $D_u$ of Hasse invariant $1/n$.
Let $S$ be a uniformizer of the maximal ideal $\mf{m}_{D,u}$ of
$\mc{O}_{D,u}$.
There is an isomorphism $U_\mbf{A}(\QQ_p) \cong D^\times_u$ 
(Lemma~\ref{lem:Z_ppoints}).
By Proposition~\ref{prop:plocaldense}(1), there exists an element $T \in
U_\mbf{A}(\QQ)$ such that the image $T_u \in D_u^\times$ satisfies
$$ T_u = S(1+Sx) $$
for some element $x \in D_u$.  In particular, $T$ is a $B$-linear 
quasi-isometry
$$ T: (A,i,\lambda) \rightarrow (A,i,\lambda) $$
for whose norm has $u$-valuation $\nu_u(N_{D/F}(T)) = 1$.
Lemma~\ref{lem:valuationimage} implies that there exists an element $R \in
GU_\mbf{A}(\QQ)$ such that $p$-adic valuation of the similitude norm
satisfies $\nu_p(\nu(R)) = 1$.

By Theorem~\ref{thm:hondatate} and Lemma~\ref{lem:polclassglobal}, there
exists a $B$-linear quasi-similitude
$$ \alpha: (A, i, \lambda) \rightarrow (A', i', \lambda'). $$
There exists an integer $e_1$ so that 
$$ (\alpha R^{e_1})^* \lambda' = c\lambda $$
for some $c \in \ZZ_{(p)}^\times$.  By altering $\lambda'$ within its
$\ZZ_{(p)}^\times$-weak polarization class, we may assume that $c = 1$.
Then $\alpha R^{e_1}$ is an quasi-isometry between $B$-linear 
abelian varieties with
prime-to-$p$ polarizations.
There exists an integer $e_2$ so that induced quasi-isogeny  
$$ (\alpha R^{e_1})_*S^{e_2}: \epsilon A(u) \rightarrow \epsilon A'(u) $$
is an isomorphism of formal groups.
By Theorem~\ref{thm:tate} the following diagram is a pullback.
$$
\xymatrix{
\mathrm{Isom}(\mbf{A}, \mbf{A'})_{(p)} \ar[d] \ar[r] &
\mathrm{Isog}(\epsilon A(u), \epsilon A'(u)) \ar[d] \\
\mathrm{Isom}(\mbf{A}, \mbf{A'})\otimes \QQ \ar[r] &
\mathrm{Isog}(\epsilon A(u), \epsilon A'(u)) \otimes \QQ
}
$$
Here, $\mathrm{Isom}(-,-)$ denotes the isometries between $B$-linear
polarized abelian varieties, and $\mathrm{Isog}$ denotes isogenies between
$p$-typical groups.
We deduce that
$$ \alpha R^{e_1} T^{e_2} : (A,i,\lambda) \rightarrow (A',i',\lambda') $$
is a prime-to-$p$ isometry of $B$-linear polarized abelian varieties.
\end{proof}

\begin{cor}\label{cor:Fpbarpoints}
The map
$$ GU_\mbf{A}(\ZZ_{(p)}) \backslash GU(\AF^{p,\infty}) /K^p \rightarrow
\Sh(K^p)^{[n]}(\br{\FF}_p) $$
given by sending a double coset $[g]$ to the tuple 
$(A,i,\lambda,[\eta g]_{K^p})$ is an isomorphism.
\end{cor}

We give an alternative characterizations of the adelic quotient of
Corollary~\ref{cor:Fpbarpoints}.

\begin{lem}\label{lem:adelicquotients}
$\quad$
\begin{enumerate}
\item
The natural map
$$ \Gamma \backslash GU^1(\AF^{p,\infty}) \rightarrow GU_\mbf{A}(\ZZ_{(p)}) 
\backslash GU(\AF^{p,\infty}) $$
is an isomorphism.
\item
The natural map
$$ \Gamma \backslash GU^1(\AF^{p,\infty})/K^p \rightarrow GU_\mbf{A}(\ZZ_{(p)}) 
\backslash GU(\AF^{p,\infty})/K^p $$
is an isomorphism.

\item The natural map
$$ GU_\mbf{A}(\ZZ_{(p)}) \backslash GU(\AF^{p,\infty})/K^p \rightarrow 
GU_\mbf{A}(\QQ) \backslash GU_\mbf{A}(\AF^{\infty})/K^pGU_{\mbf{A}}(\ZZ_p)
$$
is an isomorphism.
\end{enumerate}
\end{lem}

\begin{proof}
(1) follows from Lemma~\ref{lem:valuationimage}, and the fact that
$\Gamma = GU_\mbf{A}(\ZZ_{(p)}) \cap GU^1(\AF^{p,\infty})$
(Theorem~\ref{thm:Rosatipos}), and    
(1) implies (2). 
%and (2) implies (3).  
Statement (3) follows from
Proposition~\ref{prop:plocaldense}(3).  
\end{proof}

\begin{rmk}
Lemma~\ref{lem:adelicquotients} relates the number of $\br{\FF}_p$-points
of $\Sh(K^p)^{[n]}$ to a class number of $GU_\mbf{A}$.
\end{rmk}

\section{$K(n)$-local $\TAF$}

In this section, fix an open compact subgroup $K^p$ of
$GU(\AF^{p,\infty})$,
sufficiently small so that $\Sh(K^p)$ is a scheme.

\begin{lem}\label{lem:TAFE(n)local}
Let $U \rightarrow \Sh(K^p)^\wedge_p$ be an \'etale open.
Then the spectrum of sections $\mc{E}(K^p)(U)$ is $E(n)$-local.  In
particular, the spectrum $\TAF(K^p)$ is $E(n)$-local. 
\end{lem}

\begin{proof}
Fix a formal affine \'etale open
$$ f: U' = \Spf(R) \rightarrow U. $$
Let $(A_f,i_f,\lambda_f, [\eta_f])/R$ be the tuple classified by $f$.
The spectrum of sections $\mc{E}(K^p)(U')$ is Landweber exact
(Corollary~\ref{cor:LEFT}), therefore
it is an $MU$-module spectrum (\cite[Thm~2.8]{hoveystrickland}).
Greenlees and May \cite[Thm~6.1]{greenleesmay} express 
the Bousfield localization at $E(n)$ 
as the localization with respect to the regular
ideal $I_{n+1} = (p, v_1, \ldots, v_{n})$ of $MU_*$:
$$
\mc{E}(K^p)(U')_{E(n)} \simeq \mc{E}(K^p)(U')[I_{n+1}^{-1}].
$$
There is a spectral sequence \cite[Thm~5.1]{greenleesmay}:
$$ H^s(\Spec(R) - V, \omega^{\otimes t}) \Rightarrow
\pi_{2t-s}\mc{E}(K^p)(U')[I_{n+1}^{-1}]. $$
Here $V$ is the locus of $\Spec(R/p)$ where the formal group $\epsilon
A_f[u]^0$ is of height greater than $n$.  However, the formal group cannot
have height greater that the height of the height $n$ 
$p$-divisible group $\epsilon A(u)$.  Therefore $V$ is empty, and the
spectral sequence collapses, because $\Spec(R)$ is affine.
We conclude that the map
$$ \pi_* \mc{E}(K^p)(U') \rightarrow \pi_* \mc{E}(K^p)(U')_{E(n)} $$
is an isomorphism, so $\mc{E}(K^p)(U')$ must be $E(n)$-local.

Because $\Sh(K^p)$ is quasi-projective, $U$ is separated, hence all of the
terms of the \v Cech nerve $(U'^{\bullet+1})$
$$ U'^{\bullet+1} = \{ U' \Leftarrow U' \times_{U} U' \Lleftarrow
U' \times_{U} U' \times_{U} U' \cdots \} $$
are affine formal schemes.
Because the sheaf $\mc{E}(K^p)$ satisfies homotopy descent, the map
\begin{equation}\label{eq:cech}
\mc{E}(K^p)(U) \rightarrow \holim_\Delta \mc{E}(K^p)(U'^{\bullet+1})
\end{equation}
is an equivalence.
Since the terms of the homotopy totalization are $E(n)$-local, we determine
that $\mc{E}(K^p)(U)$ is $E(n)$-local.
\end{proof}

Let $Gal$ 
\index{2Gal@$Gal$}
be the Galois group $Gal(\br{\FF}_p/\FF_p) \cong \widehat{\ZZ}$,
generated by the Frobenius automorphism $\Fr$.
\index{2Fr@$\Fr$}
Let $\Sh(K^p)_{\FF_{p^k}}$ 
\index{2ShKpFq@$\Sh(K^p)_{\FF_q}$}
be the Galois cover 
$$ \Sh(K^p)_{\FF_{p^k}} = \Sh(K^p) \otimes_{\ZZ_p} W(\FF_{p^k}) \rightarrow
\Sh(K^p). $$ 
Define spectra
$$ \TAF(K^p)_{\FF_{p^k}} = \mc{E}(K^p)((\Sh(K^p)_{\FF_{p^k}})^\wedge_p). $$
The \index{2TAFKpFq@$\TAF(K^p)_{\FF_q}$}
functoriality of the presheaf $\mc{E}(K^p)$ gives rise to a
contravariant functor $\Phi$ on
the subcategory $\Orb^{sm}_{Gal} \subset \Set^{sm}_{Gal}$ of smooth (i.e.
finite)
$Gal$-orbits by
$$ \Phi: Gal/H \mapsto \TAF(K^p)_{\br{\FF}_p^H}. $$
Let 
$$ \TAF(K^p)_{\br{\FF}_p} = V_\Phi \simeq \varinjlim_{k}
\TAF(K^p)_{\FF_{p^k}} $$ 
be \index{2TAFKpFpbar@$\TAF(K^p)_{\br{\FF}_p}$}
the associated associated smooth
$Gal$-spectrum of Construction~\ref{const:V}.

\begin{lem}
The spectrum $\TAF(K^p)_{\br{\FF}_p}$ is $E(n)$-local.
\end{lem}

\begin{proof}
The underlying spectrum of $\TAF(K^p)_{\br{\FF}_p}$ is weakly equivalent to
a colimit of $E(n)$-local spectra.  Since localization with respect to
$E(n)$ is smashing, the colimit is $E(n)$-local.
\end{proof}

By Lemma~\ref{lem:cofixedpoints}, we
can recover $\TAF(K^p)$ by taking $Gal$-homotopy fixed points:
$$ \TAF(K^p) \simeq (\TAF(K^p)_{\br{\FF}_p})^{hGal}. $$
By Theorem~5.2.5 of \cite{markbldg}, we have the following lemma

\begin{lem}\label{lem:Galfixed}
The sequence
$$ \TAF(K^p) \rightarrow \TAF(K^p)_{\br{\FF}_p} \xrightarrow{\Fr-1}
\TAF(K^p)_{\br{\FF}_p} $$
is a fiber sequence.
\end{lem}

For a multi-index $J = (j_0, j_1, \ldots, j_{n-1})$, let
$M(J)$ 
\index{2MI@$M(I)$}
denote the
corresponding generalized Moore spectrum with $BP$-homology
$$ BP_*M(J) = BP_*/(p^{j_0}, v_1^{j_1}, \ldots, v_{n-1}^{j_{n-1}}). $$
The periodicity theorem of Hopkins and Smith \cite{hopkinssmith} 
guarantees the existence of
generalized Moore spectra $M(J)$ for a cofinal collection of multi-indices
$J$.
The $K(n)$-localization of an $E(n)$-local spectrum $X$ may be calculated by
by the following proposition (\cite[Prop~7.10]{hoveystrickland}).

\begin{prop}\label{prop:hoveystrickland}
Suppose that $X$ is an $E(n)$-local spectrum.  Then there is an equivalence
$$ X_{K(n)} \simeq \holim_J X \wedge M(J). $$
\end{prop}

There is an isomorphism
$$ \Sh(K^p)^{[n]}_{\FF_{p^k}} \cong \coprod_{i \in I(K^p,k)}
\Spec(\FF_{p^{k_i}}) $$
for some finite index set $I(K^p,k)$.  The numbers $k_i$ are all greater
than or equal to $k$.  Let $\br{\GG}_i$ be the restriction of the
formal group $\epsilon A(u)^0$ to the $i$th factor.

\begin{prop}\label{prop:K(n)localTAF}
There is an canonical equivalence
$$ \TAF(K^p)_{\FF_{p^k},K(n)} \xrightarrow{\simeq} 
\prod_{i \in I(K^p,k)} E_{\br{\GG}_i}. $$
Here $E_{\br{\GG}_i}$ is the Morava $E$-theory spectrum associated to the
height $n$ formal group $\br{\GG}_i$ over $\FF_{p^{k_i}}$.
\end{prop}

\begin{proof}
Let $(\Sh(K^p)_{\FF_{p^k}})^\wedge_{I_n}$ denote the completion of 
$\Sh(K^p)_{\FF_{p^k}}$ at the
subscheme $\Sh(K^p)^{[n]}_{\FF_{p^k}}$.  By Corollary~\ref{cor:defsh}, 
there is an isomorphism
$$ (\Sh(K^p)_{\FF_{p^k}})^\wedge_{I_n} \cong \coprod_{i \in I(K^p, k)}
\Def_{\br{\GG}_i}. $$
Let $g$ be the inclusion
$$ g: (\Sh(K^p)_{\FF_{p^k}})^\wedge_{I_n} \rightarrow 
(\Sh(K^p)_{\FF_{p^k}})^\wedge_p. $$
Let $f : U \rightarrow (\Sh(K^p)_{\FF_{p^k}})^\wedge_p$ be a formal
affine \'etale cover, and let $U^{\bullet+1}$ be its \v Cech nerve.  The
cover $U$ pulls back to an \'etale cover $\widehat{f}: \widehat{U}
\rightarrow (\Sh(K^p)_{\FF_{p^k}})^\wedge_{I_n}$, with \v Cech nerve
$\widehat{U}^{\bullet+1}$.  There are isomorphisms
$$ \widehat{U}^{s} \cong \coprod_{i \in I(K^p,k,s)} \Def_{\br{\GG}_{s,i}}
$$ 
for some finite index set $I(K^p, k, s)$.
The unit of the adjunction $(g^*, g_*)$ gives rise to a map of spectra of
sections
$$ \Res_g: \mc{E}(K^p)(U^{s}) \rightarrow (g^*\mc{E}(K^p))(\widehat{U}^s)
$$
The functoriality of Theorem~\ref{thm:lurie}, together with
Corollary~\ref{cor:GHML}, gives a map
$$ (g,Id)^*: (g^*\mc{E}(K^p))(\widehat{U}^s) \rightarrow \prod_{i \in
I(K^p,k,s)} E_{\br{\GG}_{s,i}}. $$
Since $U^s$ is a affine formal scheme, the homotopy groups of 
$\mc{E}(K^p)(U^s)$ is
given by
$$ \pi_k \mc{E}(K^p)(U^s) \cong 
\begin{cases}
\omega^{\otimes k/2} (U^s) & \text{$k$ even}, \\
0 & \text{$k$ odd}.
\end{cases}
$$
We have the following diagram.
$$
\xymatrix{
\pi_{2t} \mc{E}(K^p)(U^{s}) \ar[r]^{\Res_g} \ar@{=}[d] &
\pi_{2t} (g^*\mc{E}(K^p))(\widehat{U}^s) \ar[r]^{(g,Id)^*} \ar@{=}[d] &
\pi_{2t} \prod_{i \in
I(K^p,k,s)} E_{\br{\GG}_{s,i}} \ar@{=}[d]
\\
\omega^{\otimes t}(U^s) \ar[r]_{\Res_g} &
(g^* \omega^{\otimes t})(\widehat{U}^s) \ar[r] &
(g_{coh}^* \omega^{\otimes t})(\widehat{U}^s)
}
$$
Here, $g^*$ denotes the pullback of sheaves of spectra/abelian groups whereas
$g^*_{coh}$ denotes the pullback of coherent sheaves.
The bottom row induces an isomorphism
$$ \omega^{\otimes t}(U^s)^\wedge_{I_n} \xrightarrow{\cong} 
(g_{coh}^* \omega^{\otimes
t})(\widehat{U}^s), $$
where we have completed the $MU_*$-module $\omega^{\otimes t}(U^s)$ at the
ideal 
$$ I_n = (p, v_1, \ldots, v_{n-1}) \subset MU_*. $$  
Therefore, we
deduce that the map $\Res_g \circ (g,Id)^*$ induces an equivalence
$$ \holim_J \mc{E}(K^p)(U^{s})\wedge M(J) \xrightarrow{\simeq} 
\prod_{i \in
I(K^p,k,s)} E_{\br{\GG}_{s,i}}. 
$$
Taking $\holim_\Delta$, we get an equivalence
\begin{equation}\label{eq:completetot}
\holim_{\Delta} \holim_J \mc{E}(K^p)(U^{\bullet+1})\wedge M(J) \xrightarrow{\simeq} 
\holim_\Delta \prod_{i \in
I(K^p,k,\bullet+1)} E_{\br{\GG}_{\bullet+1,i}}. 
\end{equation}
(The cosimplicial structure of the target is induced by the simplicial
structure of the \v Cech nerve $\widehat{U}^{\bullet+1}$ using the
functoriality of Theorem~\ref{thm:goersshopkinsmiller}.)
Using Lemma~\ref{lem:TAFE(n)local}, 
Proposition~\ref{prop:hoveystrickland}, the homotopy descent property
for the presheaf $\mc{E}(K^p)$, and the fact that the
complexes $M(J)$ are finite, we have equivalences
\begin{align*}
\holim_{\Delta} \holim_J \mc{E}(K^p)(U^{\bullet+1})\wedge M(J)
& \simeq \holim_J M(J) \wedge \holim_{\Delta}  \mc{E}(K^p)(U^{\bullet+1}) \\
& \simeq \holim_J M(J) \wedge \TAF(K^p)_{\FF_{p^k}} \\
& \simeq \TAF(K^p)_{\FF_{p^k},K(n)}.
\end{align*}
Since the coherent sheaf $g_{coh}^* \omega^{\otimes t}$ satisfies \'etale
descent, the cohomology of the cosimplicial abelian group 
$(g_{coh}^* \omega^{\otimes t})(\widehat{U}^{\bullet+1})$ is given by 
$$
\pi^s (g_{coh}^* \omega^{\otimes t})(\widehat{U}^{\bullet+1})
= 
\begin{cases}
(g_{coh}^* \omega^{\otimes t})((\Sh(K^p)_{\FF_{p^k}})^\wedge_{I_n}) & 
s = 0, \\
0 & s \ne 0.
\end{cases}
$$
We therefore deduce that the map
$$
\prod_{i \in I(K^p,k)} E_{\br{\GG}_i} \rightarrow
\holim_\Delta \prod_{i \in
I(K^p,k,\bullet+1)} E_{\br{\GG}_{\bullet+1,i}}
$$
is an equivalence.
The equivalence (\ref{eq:completetot}) therefore gives an equivalence
$$ \TAF(K^p)_{\FF_{p^k},K(n)} \xrightarrow{\simeq} 
\prod_{i \in I(K^p,k)} E_{\br{\GG}_i}. $$
\end{proof}

Let $H_n$ be the height $n$ Honda formal group over $\br{\FF}_p$.  Let
$E_n = E_{H_n}$ denote the associated Morava $E$-theory.

\begin{cor}\label{cor:K(n)localTAFFpbar}
Assume that $\Sh(K^p)^{[n]}$ is a scheme.
There is an equivalence
$$ \TAF(K^p)_{\br{\FF}_p, K(n)} \simeq \prod_{\Gamma \backslash
GU^1(\AF^{p,\infty})/K^p} E_n. $$
\end{cor}

\begin{proof}
Let $k$ be sufficiently large so that 
$$ \Sh(K^p)^{[n]}(\FF_{p^k}) = \Sh(K^p)^{[n]}(\br{\FF}_{p})
\cong \Gamma \backslash GU^1(\AF^{p,\infty})/K^p. $$
(The last isomorphism is Corollary~\ref{cor:Fpbarpoints}.)  
Then there is an isomorphism:
$$ \Sh(K^p)^{[n]} \cong \coprod_{\Gamma \backslash GU^1(\AF^{p,\infty})/K^p} 
\Spec(\FF_{p^k}). $$
Then, by
Proposition~\ref{prop:K(n)localTAF}, there is an equivalence
$$ TAF(K^p)_{\FF_{p^k},K(n)} \xrightarrow{\simeq}
\prod_{i \in \Gamma \backslash GU^1(\AF^{p,\infty})/K^p} E_{(\br{\GG}_i,
\FF_{p^k})}. $$
Taking the colimit over $k$ and $K(n)$-localizing gives an equivalence
\begin{equation}\label{eq:colimequiv}
(\varinjlim_k TAF(K^p)_{\FF_{p^k},K(n)})_{K(n)} \xrightarrow{\simeq}
\prod_{i \in \Gamma \backslash GU^1(\AF^{p,\infty})/K^p} (\varinjlim_k 
E_{(\br{\GG}_i,\FF_{p^k})})_{K(n)}.
\end{equation}
Because $K(n)$-equivalences are preserved under colimits, we have
\begin{align*}
(\varinjlim_k TAF(K^p)_{\FF_{p^k},K(n)})_{K(n)}
& \simeq (\varinjlim_k TAF(K^p)_{\FF_{p^k}})_{K(n)} \\
& \simeq \TAF(K^p)_{\br{\FF}_p,K(n)}.
\end{align*}
Each height $n$ formal group $\br{\GG}_i$ is isomorphic to $H_n$ over
$\br{\FF}_p$.  Therefore, using Proposition~\ref{prop:hoveystrickland} and the
functoriality of the Goerss-Hopkins-Miller theorem
(Theorem~\ref{thm:goersshopkinsmiller}), there are
equivalences
\begin{align*}
(\varinjlim_k E_{(\br{\GG}_i, \br{\FF}_p^k)})_{K(n)} 
& \simeq \holim_J \varinjlim_k E_{(\br{\GG}_i, \br{\FF}_p^k)} \wedge M(J) \\
& \simeq E_n.
\end{align*}
The equivalence~(\ref{eq:colimequiv}) therefore gives an equivalence
$$ 
\TAF(K^p)_{\br{\FF}_p, K(n)} \xrightarrow{\simeq} \prod_{\Gamma \backslash
GU^1(\AF^{p,\infty})/K^p} E_n. $$
\end{proof}

Taking $Gal$-homotopy fixed points, we arrive at the following.

\begin{cor}\label{cor:K(n)localTAF}
Assume that $\Sh(K^p)^{[n]}$ is a scheme.
There is an equivalence
$$ TAF(K^p)_{K(n)} \simeq \left( \prod_{\Gamma \backslash
GU^1(\AF^{p,\infty})/K^p} E_n \right)^{hGal} . $$
\end{cor}

\begin{rmk}
In Corollary~\ref{cor:K(n)localTAF},
the action of the Galois group $Gal$ is typically non-trivial on the
index set
$$ \Sh(K^p)^{[n]}(\br{\FF}_{p})
\cong \Gamma \backslash GU^1(\AF^{p,\infty})/K^p $$
and the action of $Gal$ on the components $E_n$ is not the typical one.
This is because the equivalence
$$
(\varinjlim_k E_{(\br{\GG}_i, \br{\FF}_p^k)})_{K(n)} 
\simeq E_n
$$
appearing in the proof of Corollary~\ref{cor:K(n)localTAFFpbar}
is not $Gal$-equivariant; the formal group $\GG_i$ may not even be defined
over $\FF_p$.
\end{rmk}

\section{$K(n)$-local $Q_{U}$}

Let $G$ 
\index{2G@$G$}
and $G^1$ 
\index{2G1@$G^1$}
denote the groups $GU(\AF^{p,\infty})$ and 
$GU^1(\AF^{p,\infty})$, respectively.
In Section~\ref{sec:GUaction}, using a construction described in
Section~\ref{sec:descent}, we produced a fibrant smooth
$G$-spectrum $V = V_{GU}$ such that
\index{2V@$V$}
$$ V \simeq \varinjlim_{K^p} \TAF(K^p). $$
The colimit is taken over compact open subgroups $K^p$ of
$G$.  
There is a cofinal collection of subgroups $K^p$ such
that $\Sh(K^p)^{[n]}$ is a scheme.
Let $V_{\br{\FF}_p}$ be a fibrant smooth 
\index{2VFpbar@$V_{\br{\FF}_p}$}
$G \times Gal$-spectrum such that
$$ V_{\br{\FF}_p} \simeq \varinjlim_{K^p,k} \TAF(K^p)_{\br{\FF}_{p^k}}. $$

Let $\Sp_{\Gamma}$ 
\index{2SpGamma@$\Sp_\Gamma$}
be the category of $\Gamma$-equivariant spectra, with
the injective model structure.  The cofibrations and weak equivalences in
this model structure are detected on the underlying category of spectra.
Consider the adjoint pair 
$(\Res_\Gamma^{G^1}, 
\Map_{\Gamma}(G^1, -)^{sm})$:
$$ \Res_\Gamma^{G^1} : \Sp^{sm}_{G^1}
\leftrightarrows \Sp_\Gamma : \Map_{\Gamma}(G^1, -)^{sm}
$$
defined by
\begin{align*}
\Res_\Gamma^{G^1} X & = X, 
\\
\Map_{\Gamma}(G^1, Y)^{sm} & = 
\varinjlim_{H \le_o G^1} 
\Map_{\Gamma}(G^1/H, Y).
\end{align*}
The \index{2Res@$\Res$}
$G^1$-action on 
$\Map_{\Gamma}(G^1, Y)^{sm}$ is by precomposition with
\index{2MapGammasm@$\Map_{\Gamma}(-,-)^{sm}$}
right multiplication.

The following double coset formula is very useful.

\begin{lem}\label{lem:doublecoset}
Let $H$ and $K$ be subgroups of a group $G$.
Let $Y$ be an $H$-spectrum.  
There is an isomorphism
$$ \Map_{H} (G/K, Y) \cong \prod_{[g] \in 
H \backslash G/K} Y^{H \cap gKg^{-1}}. $$
\end{lem}

\begin{lem}\label{lem:doublecosetspecial}
Let $Y$ be a $\Gamma$-spectrum.
Suppose that $K^p$ be sufficiently small so that $\Sh(K^p)^{[n]}$ is a scheme.
Then there is an isomorphism
$$ \Map_\Gamma(G^1/K^p, Y) \cong \prod_{\Gamma \backslash 
G^1/K^p} Y. $$
\end{lem}

\begin{proof}
By Proposition~\ref{prop:aut}, for every element $g \in
G^1$, the group 
$$ \Gamma(gK^pg^{-1}) = \Gamma \cap
gK^pg^{-1} $$ 
is the automorphism group of a $\br{\FF}_p$-point of the stack
$\Sh(K^p)^{[n]}$.  Because these automorphism groups are trivial, the lemma
follows immediately from Lemma~\ref{lem:doublecoset}.
\end{proof}

\begin{lem}\label{lem:Gammamodel}
$\quad$
\begin{enumerate}
\item
The functors $(\Res_\Gamma^{G^1}, 
\Map_{\Gamma}(G^1, -)^{sm})$ form a Quillen pair.  

\item
The functor $\Map_{\Gamma}(G^1, -)^{sm}$ takes
fibrant $\Gamma$-spectra to fibrant smooth $G^1$-spectra.

\item
The functor $\Map_{\Gamma}(G^1, -)^{sm}$ preserves all
weak equivalences.
\end{enumerate}
\end{lem}

\begin{proof}
Using Lemma~\ref{lem:smmodel}, it is clear that the functor
$\Res_\Gamma^{G^1}$ preserves cofibrations and weak
equivalences, and this proves (1).  Statement (2) is an immediate
consequence of statement (1).  Because there is a natural isomorphism
$$ \Map_{\Gamma}(G^1, -)^{sm} \cong \varinjlim_{K^p} 
\Map_{\Gamma}(G^1/K^p, -), $$
where the colimit is taken over $K^p$ such that $\Sh(K^p)^{[n]}$ is a
scheme,
Statement~(3) follows from Lemma~\ref{lem:doublecosetspecial}.
\end{proof}

Fix an $\br{\FF}_p$-point $(A,i,\lambda, [\eta])$ of $\Sh(K^p)^{[n]}$ 
(where $K^p$ is some compact open subgroup of $G^1$).
Fix an isomorphism $\alpha: \epsilon A(u) \cong H_n$, where $H_n$ is the Honda
height $n$ formal group.  The action of $\Gamma$ on $\epsilon A(u)$ and
the isomorphism $\alpha$ gives an action of $\Gamma$ on the spectrum $E_n$.
The following proposition is immediate from 
Corollary~\ref{cor:K(n)localTAFFpbar} and
Lemma~\ref{lem:doublecosetspecial}.

\begin{prop}\label{prop:K(n)localV}
There is an $K(n)$-equivalence
$$ V_{\br{\FF}_p} \rightarrow
\Map_{\Gamma}(G^1, E_n)^{sm} $$
of smooth $G^1 \times Gal$-spectra.
\end{prop}

We now are able to prove our main $K(n)$-local result.

\begin{thm}
Let $U$ be an open subgroup of $G^1$.  There is an equivalence
$$ (V^{hU})_{K(n)} \xrightarrow{\simeq} \left( \prod_{[g] \in \Gamma \backslash
G^1/U} E_n^{h\Gamma(gUg^{-1})} \right)^{hGal}. $$
\end{thm}

\begin{proof}
It suffices to prove that there is a $Gal$-equivariant equivalence
$$ (V_{\br{\FF}_p}^{hU})_{K(n)} \xrightarrow{\simeq} 
\prod_{[g] \in \Gamma \backslash
G^1/U} E_n^{h\Gamma(gUg^{-1})}. $$
The result is then obtained by taking homotopy fixed points with respect to
$\Fr \in Gal$ (Lemma~\ref{lem:Galfixed}).
By Proposition~\ref{prop:hoveystrickland}, a map between $E(n)$-local
spectra which is an $M(J)$-equivalence for every $J$ is a
$K(n)$-equivalence (it actually suffices to only check this for a single
$J$).  Since localization with respect to $E(n)$ is smashing, colimits of
$E(n)$-local spectra are $E(n)$-local.  Therefore, the spectra
$V_{\br{\FF}_p}$ and $\Map_\Gamma(G^1, E_n)^{sm}$ are $E(n)$-local.  Since
$M(J)$ is finite, Proposition~\ref{prop:K(n)localV} implies that there is a 
$Gal$-equivariant $K(n)$-local equivalence
$$ V_{\br{\FF}_p}^{hU} \rightarrow (\Map_\Gamma(G^1, E_n)^{sm})^{hU}. $$
Let $E_n'$ be a fibrant replacement for $E_n$ in the category of
$\Gamma$-equivariant spectra.
By Lemma~\ref{lem:Gammamodel}, the spectrum $\Map_\Gamma(G^1, E_n)^{sm}$ is
fibrant as a smooth $G^1$-spectrum, and the map
$$ \Map_\Gamma(G^1, E_n)^{sm} \rightarrow \Map_\Gamma(G^1, E_n')^{sm} $$
is a weak equivalence.  

Using Corollary~\ref{cor:resfibrant} and 
Lemma~\ref{lem:doublecoset}, we have the following sequence of
equivalences:
\begin{align*}
(\Map_\Gamma(G^1, E_n)^{sm})^{hU}
& \simeq (\Map_\Gamma(G^1, E'_n)^{sm})^{U} \\
& \cong \Map_\Gamma(G^1/U, E'_n) \\
& \cong \prod_{[g] \in \Gamma \backslash G^1/U} (E'_n)^{\Gamma(gUg^{-1})} \\
& \simeq \prod_{[g] \in \Gamma \backslash G^1/U} (E_n)^{h\Gamma(gUg^{-1})}.
\end{align*}
\end{proof}

Specializing to the cases of $U = K^p$ and $U = K^{p,S}$, we have the
following corollary.

\begin{cor}\label{cor:K(n)localTAFQ}
Let $K^p$ be an open compact subgroup of $GU(\AF^{p,\infty})$.
There are equivalences
\begin{align*}
\TAF(K^p)_{K(n)} & \simeq \left( \prod_{[g] \in \Gamma \backslash
GU^1(\AF^{p,\infty})/K^p} E_n^{h\Gamma(gK^pg^{-1})} \right)^{hGal}, \\
Q_U(K^{p,S})_{K(n)} & \simeq \left( \prod_{[g] \in \Gamma \backslash
GU^1(\AF^{p,S,\infty})/K^{p,S}} E_n^{h\Gamma(gK^{p,S}g^{-1})}
\right)^{hGal}.
\end{align*}
\end{cor}

\chapter{Example: chromatic level $1$}\label{chap:examples}

In this chapter we provide some analysis of the spectrum $\TAF$ and
the associated homotopy fixed point spectrum at chromatic filtration
$1$.  In particular, we find that these spectra are closely related to
the $K(1)$-local sphere.

\section{Unit groups and the $K(1)$-local
sphere}\label{sec:unitgroups}

In this section, we indicate how we can recover a description of the
$K(1)$-local sphere by making use of units in a field extension of
the rationals.

Let $E_1$ be the Lubin-Tate spectrum whose homotopy groups are 
$W(\br{\FF}_p)[u^{\pm 1}]$ with $|u| = 2$. 
We regard $E_1$ as being the Hopkins-Miller spectrum 
associated to the formal multiplicative group $\widehat{\GG}_m$
over $\bar{\FF}_p$. Since the formal group $\widehat{\GG}_m$ is defined
over $\FF_p$, the spectrum $E_1$ possesses an action of
$$ Gal = Gal(\bar{\FF}_p/\FF_p) $$
by $E_\infty$ ring maps.  The work of Goerss and Hopkins
\cite{goersshopkins}
specializes to prove that there is an
isomorphism 
$$ \GG_1 \cong \Aut_{E_\infty}(E_1) $$
where 
$$ \GG_1 = \MS_1 \times Gal \cong \ZZ_p^\times \times \widehat{\ZZ}. $$
Specifically, 
\index{2Gn@$\GG_n$}
for $k \in \mb
Z_p^\times$, there is an Adams operation $\psi^k\co E_1 \to E_1$ such
that $\psi^k_*(u) = k u$.
The $Gal$-homotopy fixed points of the spectrum $E_1$ is the 
spectrum $KU_p$, the $p$-completion of the 
complex $K$-theory spectrum.  The action of
$\ZZ_p^\times$ on $E_1$ descends to an action on $KU_p$, and the Adams
operations $\psi^k_*$ restrict to give the usual Adams operations on
$p$-adic $K$-theory.

The product decomposition
$$ \Aut_{E_\infty}(E_1) \cong \MS_1 \times Gal $$
is not canonical.  Rather, there is a canonical short exact sequence of
profinite groups
\begin{equation}\label{eq:G1SES} 
1 \rightarrow \MS_1 \rightarrow \Aut_{E_\infty}(E_1) \rightarrow Gal
\rightarrow 1 
\end{equation}
and the choice of formal group $\widehat{\GG}_m$ over $\FF_p$ gives rise to
a splitting.  More generally, Morava \cite{morava} 
studied \emph{forms of $K$-theory}:
\index{K-theory@$K$-theory!form of}
$p$-complete ring spectra $K'$ such that there exists an isomorphism of
multiplicative cohomology theories
$$ KU_p^*(-) \otimes_{\ZZ_p} W(\bar{\FF}_p) \xrightarrow{\cong} {K'}^*(-)
\otimes_{\ZZ_p} W(\bar{\FF}_p). $$
Morava showed that there was an isomorphism
$$ \{ \text{forms of $KU_p$} \} \xrightarrow{\cong} H^1_c(Gal;
\ZZ_p^\times) \cong \ZZ_p^\times. $$
Using Goerss-Hopkins-Miller theory, one can strengthen Morava's theorem to
prove that there is an isomorphism
$$ \{ \text{$E_\infty$ forms of $KU_p$} \} \xrightarrow{\cong} H^1_c(Gal;
\ZZ_p^\times). $$
Given a Galois cohomology class 
$$ \alpha \in H^1(Gal; \ZZ_p^\times) $$
we may regard it as giving a splitting of the short exact sequence
(\ref{eq:G1SES}), and thus an inclusion
$$ \iota_\alpha: Gal \hookrightarrow \Aut_{E_\infty}(E_1). $$
The $E_\infty$-form $K_\alpha$ associated to the cohomology class $\alpha$ 
is given as the
homotopy fixed points of the \emph{new} Galois action on $E_1$
induced by the inclusion $\iota_\alpha$:
\begin{equation}\label{eq:Kalpha}
K_\alpha = E_1^{h_\alpha Gal}. \
\end{equation}
\index{2Kalpha@$K_\alpha$}

The $K(1)$-local sphere is known to be homotopy equivalent to the
fiber of the map $\psi^k - 1\co {KO}_p \to KO_p$, where $k$
is any topological generator of the group $\mb Z_p^\times/\{\pm 1\}$.
Using the equivalence of Devinatz-Hopkins \cite{devinatzhopkins}
$$ S_{K(1)} \xrightarrow{\simeq} E_1^{h\GG_1} $$
and (\ref{eq:Kalpha}), we can substitute the $KO_p$-spectrum with the
fixed points of any form
of $K$-theory to give a fiber sequence
$$ S_{K(1)} \rightarrow K_\alpha^{h\{\pm 1\}} \xrightarrow{\psi^k - 1}
K_\alpha^{h\{\pm 1\}}. $$

Let $F$ 
\index{2F@$F$}
be a quadratic imaginary extension field of $\mb Q$, and let
$p$ be a prime of $\mb Q$ that splits in $F$ as $u u^c$.  
\index{2U@$u$}
\index{2Uc@$u^c$}
This
corresponds to the existence of an embedding $u\co F \to \mb Q_p$.  If
${\calO}_F$ is the ring of integers of $F$, there is a corresponding
embedding $u\co {\calO}_F \to \mb Z_p$.
\index{2OF@$\mc{O}_F$}

The Dirichlet unit theorem tells us that the unit group of 
${\calO}_F$ is finite.  As the extension is quadratic, there are only three
possibilities.  If $F = \mb Q(i)$, ${\calO}_F^\times$ is the group of
fourth roots of unity.  If $F = \mb Q(\omega)$, where $\omega$ is a
third root of unity, then ${\calO}_F^\times$ is the group of sixth
roots of unity. In either of these two exceptional cases, the primes
$2$ and $3$ are nonsplit.  In any other case, ${\calO}_F^\times$ is
$\{\pm 1\}$.

Fix a form of $K$-theory $K_\alpha$.
The map ${\calO}_F \to \mb Z_p$ gives an action of the (finite) unit
group of ${\calO}_F$ on $K_\alpha$.  
If $p \neq 2$, this action factors through
the
action of the roots of unity $\mu_{p-1} \subset \mb Z_p^\times$, and
the homotopy fixed point set of the action of ${\calO}_F^\times$ on
$\mb Z_p$ is a wedge of suspensions (forms of) of Adams summands.  
If $p = 2$,
then ${\calO}_F^\times = \{\pm 1\}$, and the homotopy fixed point set
is $K_\alpha^{\{\pm 1\}}$, additively equivalent to $KO_2$.

However, the image of ${\calO}_F$ in $\mb Z_p$ contains more units
than merely this finite subgroup.  Let $S_F$ be a finite set of primes
of ${\calO}_F$ that do not divide $p$.  Then there is an extension
map $S_F^{-1} {\calO}_F \to \mb Z_p$, and an action of $(S_F^{-1} 
{\calO}_F)^\times$ on $K_\alpha$.

Let $\Cl(F)$ denote the ideal class group of $F$.
\index{2ClF@$\Cl(F)$}
There is
an exact sequence
\begin{equation}\label{eq:locseq}
0 \to {\calO}_F^\times \to (S_F^{-1}{\calO}_F)^\times \longoverto^{\oplus
  \nu} \bigoplus_{S_F} \mb Z \to \Cl(F) \to \Cl(F)/\langle S_F \rangle
\to 0.
\end{equation}
The class group of a number field is finite, so 
$(S_F^{-1} {\calO}_F)^\times$ is 
isomorphic to ${\calO}_F^\times \oplus \mb Z^{|S_F|}$.

In particular, if we invert only 1 prime $w$ (so $|S_F| = 1$), there
exists a smallest integer $d$ such that $w^d = (\kappa)$ as ideals for
some $\kappa \in {\calO}_F$.  The homotopy fixed point set of $K_\alpha$
under the action of $(S_F^{-1} {\calO}_F)^\times$ is the spectrum
\[
\left(K_\alpha^{h{\calO}_F^\times}\right)^{h\langle \kappa\rangle},
\]
or equivalently the homotopy fiber of the map 
$[\kappa] - 1\co K_\alpha^{h{\calO}_F^\times} \to K_\alpha^{h{\calO}_F^\times}$.  
If $\kappa$ is a
topological generator of $\mb Z_p^\times/{\calO}_F^\times$, the
resulting homotopy fixed point spectrum is, in fact, the $K(1)$-local
sphere.

Global class field theory and the Chebotarev density theorem show that
prime ideals of $F$ are uniformly distributed in the class group, and
principal prime ideals are uniformly distributed in any congruence
condition.  Therefore, it is always possible to pick a principal prime
ideal $(\kappa)$ of ${\calO}_F$ such that $\kappa$ maps to a
topological generator of the (topologically cyclic) group $\mb
Z_p^\times/\{\pm 1\}$.  

One the other hand, suppose one inverts a set $S_F$ of primes with $|S_F|
> 1$.  Let $H$ be the closure of the image of 
$(S_F^{-1} {\calO}_F)^\times$ in 
$\mb Z_p^\times/{\calO}_F^\times$.  The group $H$ is cyclic,
and hence there is a decomposition
\[
(S_F^{-1} {\calO}_F)^\times \cong {\calO}_F^\times \times \langle x
\rangle \times \mb Z^{|S_F| - 1}
\]
such that the image of $x$ is a topological generator of $H$.  Then
there is a decomposition of the homotopy fixed point set of $K_\alpha$
under $(S_F^{-1} {\calO}_F)^\times$ as
\[
\left( K_\alpha^{h{\calO}_F^\times \times \langle x \rangle}\right)^{h \mb
  Z^{|S_F| - 1}}.
\]
However, $\mb Z^{|S_F|-1}$ acts through $H$, which acts trivially on the
homotopy fixed point spectrum.  Therefore, the homotopy fixed point
spectrum is the function spectrum
\[
F\left((B\mb Z^{|S_F|-1})_+, K_\alpha^{h{\calO}_F^\times \times \langle x
    \rangle}\right),
\]
and therefore decomposes as a wedge of suspensions of 
$K_\alpha^{h {\calO}_F^\times \times \langle x \rangle}$.

\section{Topological automorphic forms in chromatic filtration $1$}

We now analyze the homotopy fixed point spectra associated to height
$1$ Shimura varieties.  In this case, much of the required data
becomes redundant.

Let $F$ be a quadratic imaginary extension of $\mb Q$ in which
$p$ splits.  The central simple algebra over $F$ of degree $n^2$ must
be $F$ itself in this case, and the maximal order must be 
${\calO}_F$.

A $\mb Q$-valued nondegenerate hermitian alternating form on $F$ is of
the form
\[
(x,y) = \Tr_{F/\mb Q} (x \beta y^*)
\]
for \index{0bro@$(-,-)$}
some nonzero $\beta \in F$ 
\index{1beta@$\beta$}
such that $\beta^* = -\beta$,
i.e. $\beta$ is purely imaginary.  However, any two such $\beta$
differ by a scalar multiple, and so there is a unique similitude class
of such pairings on $F$.  The involution induced on $F$ is, of
necessity, complex conjugation.

The group scheme $GU$ satisfies $GU(R) \cong (R \otimes F)^\times$, and
the group scheme $U$ is the associated unitary group.  
There are isomorphisms
\begin{align*}
GU(\AF^{p,\infty}) & \cong (\AF_F^{u,u^c,\infty})^\times \\
& \cong {\prod}'_{w \not \, \vert p} F_w^\times \\
GU^1(\AF^{p,\infty}) & \cong 
\left({\prod_{\text{$\ell$ split}} }' 
\QQ_\ell^\times \times \ZZ_\ell^\times
\right)
\times 
\left( \prod_{\text{$\ell$ nonsplit}} \mc{O}_{F,\ell}^\times \right)
\end{align*}

There is a
unique maximal open compact subgroup of $GU(\mb A^{p,\infty})$, 
namely the group of units
\[
K_0^{p} = \prod_{w \not\,|\, p} {\calO}_{F,w}^\times.
\]
\index{2Kp0@$K^p_0$}

Every elliptic curve $E$ 
\index{2E@$E$}
comes equipped with a canonical weak
polarization which is an isomorphism, and any other weak polarization
differs by multiplication by a rational number.  A map 
${\calO}_{F,(p)} \to \End(E)_{(p)}$ is equivalent to a map $F \to
\End^0(E)$, as follows.  The existence of the map implies that the
elliptic curve is ordinary, as $p$ is split in $F$.  Using the
pullback
$$
\xymatrix{
\End(E)_{(p)} \ar[r] \ar[d] 
& \End(E(u)) \times \End(E(u^c)) \ar[d] 
\\
\End^0(E) \ar[r]
& \End^0(E(u)) \times \End^0(E(u^c))
}
$$
and the isomorphism 
$$ \End(E(u)) \times \End(E(u^c)) \cong \mb Z_p \times \mb Z_p, $$
we see the map ${\calO}_{F,(p)} \to \End^0(E)$ factors through
$\End(E)_{(p)}$.

There are two choices for the map ${\calO}_{F,(p)} \to
\End(E)_{(p)}$ that differ by complex conjugation.  However, there
exists a unique choice such that the corresponding summand
$E(u)$ is the formal summand and $E(u^c)$ is the
\'etale summand.

Therefore, the Shimura variety $\Sh(K_0^p)$ associated to this data
classifies elliptic curves with complex multiplication by $F$.  The
automorphism group of any such object is the unit group 
${\calO}_F^\times$.  The moduli of elliptic curves with complex
multiplication by $F$ breaks up geometrically as a disjoint union
indexed by the class group of $F$, as follows.
\[
\Sh(K_0^p) \times_{\Spec(\mb Z_p)} \Spec(W(\bar{\mb F}_p)) \cong
\coprod_{\Cl(F)} \Spec(\mb W(\bar{\mb F}_p)) \mmod {\calO}_F^\times.
\]
Here $\mmod$ denotes the stack quotient by a trivial group action.
\index{0mmod@$\mmod$}

The main theorem of complex multiplication for elliptic curves has the
following consequence.  Let $K$ 
\index{2K@$K$}
be the Hilbert class field of $F$,
i.e. the unique totally unramified extension of $F$, which has Galois
group $\Cl(F)$, and let ${\calO}_K$ be its ring of integers.  Then we
have the following isomorphism:
\begin{equation}\label{eq:hilbertclassfield}
\Sh(K_0^p) \cong \Spec({\calO}_{K,u}) \mmod {\calO}_F^\times \cong
\coprod_{\mf P | u} \Spec({\calO}_{K,\mf P}) \mmod {\calO}_F^\times.
\end{equation}
Here the coproduct is over primes $\mf P$ of $K$ dividing $u$.
Isomorphism~(\ref{eq:hilbertclassfield}) follows 
from the fact that there is a map
\[
{\calO}_F[j(E)] \to {\calO}_{K}
\]
which an isomorphism at $p$, where $j(E)$ is the $j$-invariant of any
elliptic curve with complex multiplication by $F$.

The spectrum of topological automorphic forms associated to the Shimura
variety is the spectrum
\[
\TAF(K_0^p) \simeq \left( \prod_{\Cl(F)} 
E_1^{h {\calO}_F^\times} \right)^{hGal}.
\]

The spectrum of topological automorphic forms can be expressed
as follows.  For a prime $\mf P$ of $K$ dividing $u$, let $D_{\mf P}$ be
the decomposition group, and let $f$ denote its order.  Under the isomorphism
${\rm Gal}(K/F) \cong \Cl(F)$, $D_{\mf P}$ is the subgroup
generated by $u$.  Because $K$ is unramified over $F$, the integer $f$ is
equal to the residue class degree of $\mf{P}$ over $u$.  
In particular, there are isomorphisms
$$ \mc{O}_{K,u} \cong \prod_{\mf P \vert u} \mc{O}_{K,\mf P} \cong
\prod_{\Cl(F)/D_{\mf P}} W(\FF_{p^j}). $$
We deduce that all of the $\bar{\FF}_p$ points of $\Sh(K_0^p)$ are defined
over $\FF_{p^j}$, and that, for each prime $\mf P$, there exists an
elliptic curve $E_{\mf P}/\FF_{p^j}$ with complex multiplication by $F$,
such that the isomorphism classes of $\Sh(K^p_0)(\bar{\FF}_p)$ are
represented by the set of elliptic curves
$$ \{ E_{\mf P}^{(p^i)} \: : \: \mf P \vert u, \: 0 \le i < f \} $$
(where $E^{(p^i)}_{\mf P}$ denotes the pullback of $E_{\mf P}$ 
over the $i$th power of
the Frobenius).

Let $K(E_{\mf P})$ 
\index{2KEP@$K(E_\mf{P})$}
be the form of $KU_p \otimes_{\ZZ_p} W(\FF_{p^j})$
corresponding to the unique deformation of the height $1$ formal group
$\widehat{E}_{\mf P}$ over $W(\FF_{p^j})$.
The homotopy of $K(E_{\mf P})$ is given by
$$ \pi_* K(E_{\mf P}) \cong \mc{O}_{K,\mf{P}}[u^{\pm 1}] \cong
W(\FF_{p^f})[u^{\pm 1}]. $$
We have the following:
\begin{align*}
\TAF(K_0^p) 
& \simeq \prod_{\mf P | u} K(E_{\mf P})^{h {\calO}_F^\times}. 
\end{align*}

Let $S$ be a finite set of primes of $\QQ$.  Let $S_F$ be the collection of
primes of $F$ dividing the primes in $S$.

\subsection*{Case where the primes of $S$ do not split in $F$}

Suppose that each of the primes in $S$ does not split in $F$.  Then there is
an equivalence
$$ Q_{GU}(K_0^{p,S}) = \left( \prod_{\Cl(F)/\bra{S}}
E_1^{h(S^{-1}\mc{O}_F)^\times} \right)^{hGal} $$

In particular, if $\ell$ is chosen to be a topological generator of
$\ZZ_p^\times$, and is inert in $F$, then there is an equivalence
$$ Q_{GU}(K_0^{p,\ell}) \xrightarrow{\simeq} \prod_{\Cl(F)/D_{\mf P}}
S_{K(1),\FF_{p^f}}
$$
Here, $S_{K(1), \FF_{p^f}}$ is the Galois extension of $S_{K(1)} =
(E_1^{h\MS_1})^{hGal_{\FF_p}}$ given by the fixed point spectrum
$$ S_{K(1),\FF_{p^f}} = (E_1^{h\MS_1})^{hGal_{\FF_{p^f}}}. $$

The building $\mc{B}(GU)$ 
for $GU(\QQ_\ell) = \mc{O}_{F,\ell}^\times$ is homeomorphic to the real
line $\RR$.  The group $\mc{O}_{F,\ell}^\times$ 
acts by translation by $\ell$-adic
valuation.  The resulting building decomposition is a product of $J$-fiber
sequences
$$ Q_{GU}(K_0^{p,\ell}) \rightarrow \prod_{\mf P | u}
K(E_{\mf P})^{h\mc{O}_F^\times} \xrightarrow{\prod \psi^\ell - 1}
\prod_{\mf P | u}
K(E_\mf{P})^{h\mc{O}_F^\times}. $$

\subsection*{Case where the primes in $S$ split in $F$}

Suppose that each of the primes in $S$ is split in $F$.  
Let $S'_F \subset S_F$
be a set containing exactly one prime dividing $\ell$ for each $\ell \in
S$.
The group $\Gamma(K_0^{p,S})$ is given by 
$$ \Gamma(K_0^{p,S}) = 
(S^{-1}\mc{O}_F)^{N=1}, $$ 
the subgroup of $(S^{-1}\mc{O}_F)^\times$
for which the norm $N = N_{F/\QQ}$ is $1$.
There is an exact sequence
$$ 0 \rightarrow (S^{-1}\mc{O}_F)^{N=1} \rightarrow (S^{-1}\mc{O}_F)^\times
\xrightarrow{\oplus \nu_\ell N} \bigoplus_{S} \ZZ. $$
It follows from (\ref{eq:locseq}) that there is an exact sequence
$$ 0 \rightarrow \mc{O}_F^\times \rightarrow (S^{-1} \mc{O}_F)^{N=1}
\xrightarrow{\oplus \nu_w} \bigoplus_{w \in S'_F} \ZZ \xrightarrow{\kappa}
\Cl(F) $$
where
$$ \kappa(\sum_{w \in S'_F} n_w (w)) = \sum_{w \in S'_F} n_w [w] - n_w
[w^c] \in \Cl(F). $$
In particular, there is an isomorphism
$$ \Gamma(K_0^{p,S}) \cong \mc{O}_F^{\times} \times \ZZ^{\abs{S}}. $$
Following the techniques of Section~\ref{sec:unitgroups}, the closure $H$ of
the image of $\Gamma(K_0^{p,S})$ in $\ZZ_p^\times/ \mc{O}_{F}^{\times}$ 
is cyclic.  Choosing a generator $x \in \Gamma(K_0^{p,S})$ gives a
decomposition
$$ \Gamma(K_0^{p,S}) \cong \mc{O}_F^\times \times \bra{x} \times
\ZZ^{\abs{S}-1}. 
$$
We therefore have an equivalence
\begin{align*}
Q_U(K_0^{p,S}) & \simeq \left( \prod_{\Cl(F)/\bra{S'_F}}
F(B\ZZ^{\abs{S}-1}_+,  E_1^{h\mc{O}_F^\times \times \bra{x}})
\right)^{hGal}.
\end{align*}

Assume now that $S$ consists of a single prime $\ell$ which splits as
$ww^c$, such that:
\begin{enumerate}
\item the ideal $w = (t)$ is principal,
\item the element $q = t/t^c \in \Gamma(K_0^{p,\ell}) 
\subset \ZZ_p^\times$ is a
topological generator of $\ZZ_p^\times/\mc{O}_F^\times$.
\end{enumerate}
Infinitely many such primes can be shown to exist, using class field theory
and the Chebotarev density theorem.  Condition $(2)$ implies that
$\Gamma(K_0^{p,\ell})$ is dense in $\ZZ_p^\times$.
Then there is an equivalence
$$ Q_U(K_0^{p,\ell}) \simeq \prod_{\Cl(F)/D_{\mf P}} S_{K(1),\FF_{p^f}}. $$

The building $\mc{B}(U)$ for the group 
$$ U(\QQ_\ell) = (F_\ell)^{N = 1} \cong F_w^\times $$ 
is homeomorphic to $\RR$.  An element $g \in U(\QQ_\ell)$ acts on
$\mc{B}(U)$ by translation by $\nu_w(g)$.  The building gives a
fiber sequence
$$ Q_U(K_0^{p,\ell}) \rightarrow \prod_{\mf P | u} K(E_\mf{P})^{h\mc{O}_F^\times}
\xrightarrow{\prod \psi^q-1}
\prod_{\mf P | u} K(E_{\mf{P}})^{h\mc{O}_F^\times}. $$

\backmatter

\bibliographystyle{amsalpha}
\bibliography{assembly3}
\printindex

\end{document}